\documentclass[a4paper,11pt]{amsart}
\usepackage{graphicx}
\usepackage{subcaption}
\usepackage{float}
\usepackage{overpic}
\usepackage{amsmath,amsthm,amssymb}
\usepackage{slashed}
\usepackage{comment}
\usepackage[all]{xy}
\newcommand{\mychoice}[2]{#1
}
\mychoice{
\newcommand{\plabel}[1]{ \label{#1}}
\newcommand{\gbibitem}[1]{ \bibitem{#1}}
\newcommand{\snewpage}{}
\newcommand{\scleardoublepage}{}

\specialcomment{commentx}{   }{   }
\specialcomment{commenty}{   }{   }
\excludecomment{commentx}
\excludecomment{commenty}
\newcommand{\commleaveout}[1]{}
}{
\usepackage{xcolor}
\newcommand{\plabel}[1]{ \label{#1}\rlap{\smash{${}^{^{[#1]}}$}}}
\newcommand{\gbibitem}[1]{ \bibitem{#1}\rlap{\smash{${}^{^{[#1]}}$}}}
\newcommand{\snewpage}{\newpage}
\newcommand{\scleardoublepage}{\cleardoublepage}

\newenvironment{commentx}{\color{magenta} }{\color{black} }
\newenvironment{commenty}{\color{olive} }{\color{black} }
\newcommand{\commleaveout}[1]{#1}
}
\DeclareMathOperator{\adj}{adj}
\DeclareMathOperator{\arcosh}{arcosh}

\DeclareMathOperator{\artanh}{artanh}

\DeclareMathOperator{\sgn}{sgn}
\DeclareMathOperator{\tr}{tr}
\DeclareMathOperator{\AC}{AC}
\DeclareMathOperator{\AS}{AS}
\DeclareMathOperator{\AT}{At}
\DeclareMathOperator{\ATT}{AT}
\DeclareMathOperator{\CD}{CD}
\DeclareMathOperator{\Id}{Id}
\DeclareMathOperator{\Log}{Log}
\DeclareMathOperator{\MN}{MN}
\DeclareMathOperator{\NW}{N}
\DeclareMathOperator{\PD}{PD}
\DeclareMathOperator{\SE}{SE}
\DeclareMathOperator{\SEH}{SEH}
\DeclareMathOperator{\SH}{SH}
\DeclareMathOperator{\SL}{SL}

\DeclareMathOperator{\modu}{\,mod\,}
\DeclareMathOperator{\Sin}{\slashed{\mathrm S}in}
\DeclareMathOperator{\Cos}{\slashed{\mathrm C}os}
\DeclareMathOperator{\Tan}{\slashed{\mathrm T}an}
\DeclareMathOperator{\Cot}{\slashed{\mathrm C}ot}

\DeclareMathOperator{\intD}{\mathring D}
\newcommand{\ujnorma}{\slashed N}
\newcommand{\real}{\mathrm{real}}
\newcommand{\complex}{\mathrm{complex}}

\DeclareMathOperator{\Dbar}{D}
\DeclareMathOperator{\ellip}{ell}
\DeclareMathOperator{\hyper}{hyp}
\DeclareMathOperator{\Rea}{Re}
\DeclareMathOperator{\Ima}{Im}

\DeclareMathOperator{\Rexp}{exp_{R}}
\DeclareMathOperator{\Lexp}{exp_{L}}
\DeclareMathOperator{\ad}{ad}
\DeclareMathOperator{\BCH}{BCH}
\DeclareMathOperator{\CRext}{CR^{ext}}
\DeclareMathOperator{\CR}{CR}

\DeclareMathOperator{\DW}{DW}
\DeclareMathOperator{\GL}{GL}
\DeclareMathOperator{\spec}{sp}
\DeclareMathOperator{\dett}{\mathbf{det}}
\DeclareMathOperator{\reC}{\slashed{\mathbf C}}
\DeclareMathOperator{\reD}{\slashed{\mathbf D}}
\DeclareMathOperator{\reG}{\slashed{\mathbf G}}
\DeclareMathOperator{\reP}{\slashed{\mathbf P}}
\DeclareMathOperator{\reL}{\slashed{\mathbf L}}

\DeclareMathOperator{\reW}{\slashed{\mathbf W}}
\DeclareMathOperator{\reX}{\slashed{\mathbf X}}
\newcommand{\qand}{\qquad\text{and}\qquad}
\theoremstyle{definition}
\newtheorem{point}{}[section]
\newtheorem{defin}[point]{Definition}

\newtheorem{remark}[point]{Remark}
\newtheorem{example}[point]{Example}
\newtheorem{exaprop}[point]{Example (Proposition)}

\theoremstyle{plain}

\newtheorem{lemma}[point]{Lemma}
\newtheorem{cor}[point]{Corollary}
\newtheorem{theorem}[point]{Theorem}

\newcommand{\bem}{\begin{bmatrix}}
\newcommand{\eem}{\end{bmatrix}}

\newcommand{\leaveout}[1]{}
\newcommand{\m}{\mathbf}

\newcommand{\eqed}{
\pushQED{\qed}
\qedhere
\popQED
}

\newcommand{\eqedremark}{
\renewcommand{\qedsymbol}{$\triangle$}
\pushQED{\qed}
\qedhere
\popQED
\renewcommand{\qedsymbol}{$\Box$}
}
\newcommand{\qedexer}{  \renewcommand{\qedsymbol}{$\diamondsuit$} \qed \renewcommand{\qedsymbol}{$\Box$}}
\newcommand{\qedremark}{  \renewcommand{\qedsymbol}{$\triangle$} \qed \renewcommand{\qedsymbol}{$\Box$}}

\newcommand{\proofremark}[1]{
\begin{proof}[Remark] #1
\renewcommand{\qedsymbol}{}
\end{proof}
}
\newcommand{\proofremarkqed}[1]{
\begin{proof}[Remark] #1
\renewcommand{\qedsymbol}{$\triangle$}
\end{proof}
\renewcommand{\qedsymbol}{$\Box$}
}
\makeatletter
\newcounter{savesection}
\newcounter{apdxsection}
\renewcommand\appendix{\par
  \setcounter{savesection}{\value{section}}%
  \setcounter{section}{\value{apdxsection}}%
  \setcounter{subsection}{0}%
  \gdef\thesection{\@Alph\c@section}}
\newcommand\unappendix{\par
  \setcounter{apdxsection}{\value{section}}%
  \setcounter{section}{\value{savesection}}%
  \setcounter{subsection}{0}%
  \gdef\thesection{\@arabic\c@section}}
\makeatother
\newcommand{\marginextend}[1]{ \addtolength{\oddsidemargin}{-#1}  \addtolength{\evensidemargin}{-#1}
  \addtolength{\textwidth}{#1}\addtolength{\textwidth}{#1}}
\newcommand{\updownextend}[1]{ \addtolength{\topmargin}{-#1}  \addtolength{\textheight}{#1}
\addtolength{\textheight}{#1}}
\marginextend{1cm}
\updownextend{0cm}
\allowdisplaybreaks[4]
\title[Convergence estimates for the Magnus expansion IIA. $2\times2$ matrices ]{Convergence estimates for the Magnus expansion IIA.
$2\times2$ matrices with operator norm}
\author{Gyula Lakos}
\address{Department of Geometry, Institute of Mathematics, E\"otv\"os University, P\'azm\'any P\'eter s.~1/C,  Budapest, H--1117, Hungary}
\email{lakos@cs.elte.hu}
\begin{document}
\begin{abstract}
We  review and provide simplified proofs related to the Magnus expansion, and improve convergence estimates.
Observations and improvements concerning the Baker--Campbell--Hausdorff expansion are also made.

In this Part IIA, we investigate the case of $2\times2$ matrices with respect to the operator norm.
We consider norm estimates and minimal presentations in terms of the Magnus and BCH expansions.
Some results are obtained in the complex case, but a more complete picture is obtained in the real case.
\end{abstract}
\maketitle
\section*{Introduction to Part IIA}
The present paper is a direct continuation of Part II \cite{L2}.
This assumes general familiarity with Part I \cite{L1} and a good understanding of Part II \cite{L2}.
In Part II, we obtained norm estimates and inclusion theorems with respect to the conformal range.
In Part IIA, we take a closer look to the more treatable case of $2\times2$ matrices.
Some investigations will be about testing or sharpening our earlier results for $2\times2$ matrices;
other ones will deal with minimal Magnus or Baker--Campbell--Hausdorff presentations.
Studying $2\times 2$ matrices may seem to be a modest objective, but, in reality, the computations are not trivial.
(Basic analytical calculations already appear in Magnus \cite{M};
and more sophisticated techniques appear in Michel \cite{Mi} but in the context of the Frobenius norm.
Our computations are similar in spirit but more complicated, as the operator norm is taken seriously.)
The information exhibited here is relatively more complete in the case of $2\times2$ real matrices,
and much more partial in the case of $2\times2$ complex matrices.

In Section \ref{sec:PartIIRev}, we recall some theorems and examples from Part II.
Our results here with respect to $2\times2$ matrices can be viewed relative to these.
In Section \ref{sec:LogA},  we review some technical tools concerning $2\times2$ matrices.
To a variable extent, all later sections use the information discussed here, hence the content of this section is crucial
(despite its elementary nature); however, its content can be consulted as needed later.
Nevertheless, a cursory reading is advised in any case.

Section \ref{sec:Schur} and \ref{sec:BCH} discuss how Schur's formulae and the BCH formula simplify for $2\times2$ matrices, respectively.
In Section \ref{sec:expmag}, we show that for $2\times2$ complex matrices exponentials
only rarely play the ``role of geodesics'' (i. e. parts of minimal Magnus presentations).
In Section \ref{sec:ExamplesBCH}, examples of BCH expansions from $\SL_2(\mathbb R)$ are given with interest in norm growth.
In Section \ref{sec:critbal}, we demonstrate that (appropriately) balanced BCH expansions of $2\times2$ real matrices with cumulative norm $\pi$ are uniformly bounded.
Next, we turn toward BCH minimal presentations of real $2\times2$ matrices (with norm restrictions).
Sections \ref{sec:BCHmoments} and \ref{sec:BCHmin} study the moments of the Schur maps
and apply them to minimal BCH presentations.
In Section \ref{sec:Bound}, the  (critical) ``BCH unit ball'' $\mathcal B_{\frac\pi2,\frac\pi2}$ is described resulting the ``wedge cap''.
At this point several further questions are left open,
but having that prime interesting case seen, we leave the discussion of BCH presentations here.

Then, we start the investigation of Magnus (minimal) expansions for real $2\times2$ matrices.
We find that, in this setting, Magnus minimal presentations are more natural and approachable than BCH minimal presentations.
In Section \ref{sec:PrincLog}, we prove a  logarithmic monotonicity property for the principal disks of $2\times2$ real matrices.
In Section \ref{sec:ExamplesMagnus} we consider some  examples which
test our earlier norm estimates for the Magnus expansion but also help to understand the $2\times2$ real case better
by introducing the ``canonical Magnus developments''.
Section \ref{sec:MagnusGL2} develops a systematic analysis of the $2\times2$ real case.
The first observation is that in this case the Magnus exponent can be read off from the conformal range / principal disk.
Based on the previous examples, we consider (minimal) normal
presentations for $2\times2$ real matrices given as time-ordered exponentials which are not ordinary exponentials.
The conclusion is that those normal presentations are better suited to the geometric description
of  $\widetilde{\GL^+_2}(\mathbb R)$ than the customary exponentials.
In Section \ref{sec:opt22}, we give information about the asymptotics of the optimal norm estimate for $2\times2$ real matrices.
In Section \ref{sec:count22}, however, it is demonstrated that the Magnus exponent of complex $2\times2$ matrices
cannot be read off from the conformal range as simply as in the real case.
(Thus the previous methods cannot transfer that easily.)

\subsection*{Notation and terminology}
Line $y$ of formula $(X)$ will be denoted as $(X/y)$.
If $a,b$ are points of an Euclidean (or just affine) space,
then $[a,b]_{\mathrm e}$ denotes the segment connecting them.
This notation is also applied for (half-)open segments.
It can also be used in $\mathbb R$ conveniently when $a>b$.
Instead of `by the unicity principle of analytic continuation' we will often say `by analytic continuation'
even if the function is already constructed.
In general, we try to follow the notations established in Parts I and II.

\snewpage
\begin{commentx}
\tableofcontents
\end{commentx}
\scleardoublepage\section{Some results from Part II (Review)}
\plabel{sec:PartIIRev}
Here we recall some points of Part II.
It cannot replace the detailed discussion given in Part II, but it serves reference.
~\\
\subsection{Conformal range}
~\\

Suppose that $\mathfrak H$ is a real  Hilbert space.
The logarithmic distance on $\mathfrak H\setminus\{0\}$ is given by
\begin{equation}
d_{\log}(\m v,\m w)=\sqrt{\left(\log \frac{|\m w|_2}{|\m v|_2} \right)^2+\left((\m v,\m w)\sphericalangle\right)^2}.
\plabel{eq:logdist}\end{equation}
\begin{theorem}\plabel{th:cs}
Suppose that $\m z:[a,b]\rightarrow\mathfrak H\setminus\{0\}$ is continuous. Then
\[d_{\log}(\m z(a),\m z(b) )\leq\int_{t\in[a,b]}\frac{|\mathrm d\m z(t)|_2}{|\m z(t)|_2}.\]
In case of equality, $\m z$ is a (not necessarily strictly) monotone subpath of a
distance segment connecting $\m z(a)$ to $\m z(b)$.
\qed
\end{theorem}
For $\m x,\m y \in\mathfrak H\setminus\{0\}$ let $\sphericalangle(\m x,\m y)$ be denote their angle.
This can already be obtained from the underlying real scalar product
$\langle\m x,\m y\rangle_{\mathrm{real}}=\Rea\,\langle\m x,\m y\rangle$.
For $\m x,\m y\in\mathfrak H$, $\m x\neq 0$,  let
\[\m y:\m x=\frac{
\langle \m y,\m x\rangle_{\real}}{|\m x|_2^2}
+\mathrm i\left|\frac{\m y}{|\m x|_2}-  \frac{
\langle \m y,\m x\rangle_{\real}}{|\m x|_2^2}\frac{\m x}{|\m x|_2}\right|_2.\]

For $A\in \mathcal B(\mathfrak H)$, we define the (extended) conformal range as
\[\CRext(A)=\{ A\m x:\m x, \,\overline{(A\m x:\m x)}\,:\, \m x\in\mathfrak H\setminus\{0\}\};\]
and the restricted conformal range as
\[\CR(A)=\{ A\m x:\m x\,:\, \m x\in\mathfrak H\setminus\{0\}\}.\]
(This is a partial aspect of the Davis-Wielandt shell, see Wielandt \cite{Wie}, Davis \cite{D1}, \cite{D2} and also \cite{L2}.)

\begin{theorem}\plabel{th:CRrange} (Time ordered exponential mapping theorem.)

If $\phi$ is $\mathcal B(\mathfrak H)$-valued ordered measure, then
\begin{equation} \CRext(\Lexp \phi)\subset \exp \Dbar(0,\textstyle{\int \|\phi\|_2}),\plabel{eq:CRrange}\end{equation}
and
\begin{equation}  \spec(\Lexp \phi)\subset \exp \Dbar(0,\textstyle{\int \|\phi\|_2}).\plabel{eq:spect}\end{equation}

In particular, if $\int \|\phi\|_2<\pi$, then $\log \Lexp \phi$ is well-defined, and for its spectral radius
\begin{equation} \mathrm r(\log\Lexp \phi) \leq\textstyle{\int \|\phi\|_2}.
\eqed
\plabel{eq:spect2}\end{equation}
\end{theorem}
\snewpage
\subsection{Convergence theorems}
\begin{theorem}\plabel{th:MMNC} (Mityagin--Moan--Niesen--Casas logarithmic convergence theorem.)

If $\phi$ is a $\mathcal B(\mathfrak H)$-valued ordered measure and  $\int \|\phi\|_2<\pi$,
then the Magnus expansion $\sum_{k=1}^\infty \mu_{k,\mathrm L}(\phi)$ is absolute convergent.
In fact, $\log\Lexp(\phi)=\sum_{k=1}^\infty \mu_{k,\mathrm L}(\phi)$ also holds.

The statement also holds if $\mathcal B(\mathfrak H)$ is replaced by any $C^*$-algebra.

[Cf. Sch\"afer \cite{Sch}, Mityagin \cite{Mity}, Moan, Niesen \cite{MN}, Casas \cite{Ca}.]
\qed
\end{theorem}

We say that the ordered measure $\phi$ is a multiple Baker--Campbell--Hausdorff (mBCH) measure,
if, up to reparametrization,  $\phi$ is of form $A_1\m 1\boldsymbol.\ldots\boldsymbol.A_n\m 1$.
In this case, $\phi$ also allows a mass-normalized version
\begin{equation}
\psi=B_1\m 1_{(0,t_1]}\boldsymbol.\ldots\boldsymbol.B_k\m 1_{[t_{k-1},t_k]}
\plabel{eq:mmBCH}
\end{equation}
where $t_i<t_{i+1}$, $\|B_i\|_2=1$ and thus $t_k=\int\|\phi\|_2$.
It is constructed by replacing $A_i$ with $A_i/\|A_i\|_2$ if $\|A_i\|_2\neq0$, and eliminating the term if $A_i=0$.
As it is obtained by a kind a reparametrization, its Magnus expansion is not affected.

\begin{theorem}\plabel{th:better2}
(Finite dimensional critical BCH convergence theorem.)

Let $\mathfrak H$ be a finite dimensional Hilbert space .
Consider the $\mathcal B(\mathfrak H)$ valued mBCH measure $\phi=A_1\m 1\boldsymbol.\ldots\boldsymbol.A_k\m 1$  with cumulative norm
$\|A_1\|_2+\ldots+\|A_k\|_2=\pi$.
Then, the convergence radius of the Magnus (mBCH) expansion of $\phi$ is greater than $1$.
In particular, finite dimensional mBCH expansions with cumulative norm $\pi$ converge.
\qed
\end{theorem}

\begin{lemma}\plabel{lem:common}
Let $\psi$ be a mass-normalized mBCH-measure as in  \eqref{eq:mmBCH}, $\int\|\psi\|_2>0$.

Consider all the Hilbert subspaces $\mathfrak V$ of $\mathfrak H$ such that

(i) $\mathfrak H=\mathfrak V\oplus\mathfrak V^\bot$ is an invariant orthogonal decomposition for all $B_i$.

(ii) $B_1|_{\mathfrak V}=\ldots=B_k|_{\mathfrak V}$, and these are orthogonal (unitary).

Then there is a single maximal such $\mathfrak V$.
\qed
\end{lemma}
In the context of the previous lemma, if the  maximal such $\mathfrak V$ is $0$, then we call $\psi$ reduced.
In particular, this applies if $\bigcap_{i<j}\ker(B_i-B_j)=0$.

\begin{theorem}\plabel{th:better3}
(Finite dimensional logarithmic critical BCH convergence theorem.)

Let $\mathfrak H$ be a finite dimensional Hilbert space.
Consider the $\mathcal B(\mathfrak H)$ valued mass-normalized mBCH measure $\psi$ as in \eqref{eq:mmBCH}
with cumulative norm $\int\|\psi\|_2=\pi$. We claim:

(a) Unless the component operators $B_i$ have a common eigenvector for $\mathrm i$ or $-\mathrm i$ (complex case), or a common eigenblock
$\bem&-1\\1&\eem$ (real case), then  $\log\Lexp(\psi)=\sum_{k=1}^\infty \mu_{k,\mathrm L}(\psi)$ also holds.

(b) If $\psi$ is reduced, then,
for any $t\in\Dbar(0,1)$, $\log\Lexp(t\cdot \psi)=\sum_{k=1}^\infty t^k\mu_{k,\mathrm L}(\psi)$  holds.
Thus, the $log$-able radius of the Magnus (BCH) expansion is also greater than $1$.
\qed
\end{theorem}
~
\snewpage
\subsection{Growth estimates}
~\\

For $p>0$, let us define $0<\ell(p)<\frac\pi2$ as the solution of the equation
\[\ell(p)+p\sin \ell(p)=\frac\pi2.\]
Then $\ell:(0,\infty)\rightarrow (\pi/2,0)_{\mathrm e}$ is a decreasing diffeomorphism.
In particular,
\[\ell(\pi)=0.386519539\ldots\,\,. \]

\begin{theorem}\plabel{th:biestimate2}
Let $\CRext(A)\subset \exp \Dbar(0,p)$, $0< p<\pi$.
Then
\begin{equation}
\|\log A\|_2\leq J(p)=\int_{t=\ell(p)}^{\pi-\ell(p)}
\underbrace{{\frac {p+\sin \left( p\sin  t  \right)-\cos \left( p\sin  t  \right) p\sin t}{2\sin
 \left( p\sin  t   \right) }}}_{JJ(p,t)}\,\mathrm dt.\eqed
 \plabel{eq:biestimate}
\end{equation}
\end{theorem}
\begin{theorem}\plabel{th:Jasympt}
(a) As $p\searrow0$,
\begin{equation}J(p)=p+\frac16\,{p}^{3}-{\frac {1}{72}}{p}^{5}+  \frac{211}{15120} p^7 +O({p}^{9}).\plabel{eq:Jsan1}\end{equation}

(b) As $p\nearrow\pi$,
\begin{equation}
J(p)=
\pi\sqrt{\frac{\pi+p}{\pi-p}}+J_\pi+O({\pi-p})
,
\plabel{eq:Jsan2}
\end{equation}
where
\begin{equation}
J_\pi=-4\tan\ell(\pi)
+\int_{t=\ell(\pi)}^{\pi-\ell(\pi)} \left(JJ(\pi,t)-\frac{2}{\cos^2 t} \right)\,\mathrm dt\notag
\end{equation}
(the integrand extends to a smooth function of $t$). Numerically, $J_\pi=-3.0222\ldots\,\,.$
\qed
\end{theorem}
The function $J(p)$ is not very particular, it can be improved.
An example we have considered to test the effectivity of \eqref{eq:biestimate}--\eqref{eq:Jsan2} was
\snewpage
\begin{example}\plabel{ex:critical}

(The analytical expansion of the Magnus critical case.)

 On the interval $[0,\pi]$, we consider the measure $\Phi$, such that
\[\Phi(\theta)=
\begin{bmatrix}
-\sin2\theta& \cos2\theta\\\cos2\theta&\sin2\theta
\end{bmatrix}
\,\mathrm d\theta|_{[0,\pi]}.\]
Thus, for $p\in(0,\pi)$,
\[\int\|p/\pi\cdot\Phi\|_2=p.\]
Then,
\begin{align}
\mu_{\mathrm L}(p/\pi\cdot \Phi)&=\log\Lexp (p/\pi\cdot\Phi)
\notag\\
&=\pi\left(\frac1{\sqrt{1-(p/\pi)^2}}-1\right)\begin{bmatrix}& -p/\pi-1\\-p/\pi+1&\end{bmatrix}.
\notag
\end{align}
Consequently,
\begin{align}
\|\mu_{\mathrm L}(p/\pi\cdot \Phi)\|_2
&=\pi\left(\frac1{\sqrt{1-(p/\pi)^2}}-1\right)(1+p/\pi)
\notag\\
&=\pi\sqrt{\frac{\pi+p}{\pi-p}}-\pi-p
\notag\\
&=\sqrt2\pi^{3/2} (\pi-p)^{-1/2}-2\pi-\frac{\sqrt2}4\pi^{1/2} (\pi-p)^{1/2}+O(\pi-p),
\plabel{eq:lower1}
\end{align}
as $p\nearrow\pi$.
\qedexer
\end{example}
(We will see, however, that $\Phi|_{[0,p]}$ is a better candidate to test  \eqref{eq:biestimate}--\eqref{eq:Jsan2}.)
\scleardoublepage
\section{Computational background for $2\times2$ matrices}
\plabel{sec:LogA}
For the sake of reference, here we review some facts and conventions
related to $2\times2$ matrices we have already considered in Part II.
However, we also take the opportunity to augment this review
by further observations of analytical nature.
\\

\subsection{The skew-quaternionic form (review)}
~\\

One can write the $2\times2$ matrix $A$ in skew-quaternionic form
\begin{equation}
A= \tilde a\Id_2 +\tilde b\tilde I+\tilde c\tilde J+\tilde d\tilde K\equiv
\tilde a\begin{bmatrix} 1&\\&1\end{bmatrix}
+\tilde b\begin{bmatrix} &-1\\1&\end{bmatrix}
+\tilde c\begin{bmatrix} 1&\\&-1\end{bmatrix}
+\tilde d\begin{bmatrix} &1\\1&\end{bmatrix}.
\plabel{eq:skqform}
\end{equation}

Then,
\[\frac{\tr A}2=\tilde a, \]
and
\[\det A=\tilde a^2+\tilde b^2-\tilde c^2-\tilde d^2.\]
\begin{commentx}
Furthermore,
\[A-\frac{\tr A}2\Id_2=\tilde b\tilde I+\tilde c\tilde J+\tilde d\tilde K; \]

\[A^*=\tilde a\Id_2 -\tilde b\tilde I+\tilde c\tilde J+\tilde d\tilde K;\]

\[\adj A=\tilde a\Id_2 -\tilde b\tilde I-\tilde c\tilde J-\tilde d\tilde K;\]

\[A^{-1}=\frac{\tilde a\Id_2 -\tilde b\tilde I-\tilde c\tilde J-\tilde d\tilde K}{\tilde a^2+\tilde b^2-\tilde c^2-\tilde d^2}.\]
\end{commentx}
~
\subsection{Spectral type (review)}\plabel{ss:SpectType}
~\\

Let us use the notation
\[D_A=\det\left( A-\frac{\tr A}2\Id_2\right)=(\det A)-\frac{(\tr A)^2}4.\]
It is essentially the discriminant of $A$, as the eigenvalues of $A$ are $\frac12\tr A\pm\sqrt{-D_A}$.

In form \eqref{eq:skqform},
\begin{equation}D_A=\tilde b^2-\tilde c^2-\tilde d^2.\plabel{nt:DA}\end{equation}

In the special case of real $2\times2$ matrices, we use the classification

$\bullet$ elliptic case: two conjugate strictly complex eigenvalues,

$\bullet$ parabolic case: two equal real eigenvalues,

$\bullet$ hyperbolic case: two distinct real eigenvalues.

Then, for real $2\times2$ matrices,  $D_A$ measures `ellipticity/parabolicity/hiperbolicity':
If $D_A>0$, then $A$ is elliptic;
if $D_A=0$, then $A$ is parabolic;
if $D_A<0$, then $A$ is hyperbolic.

In the general complex case, there are two main categories: parabolic $(D_A=0)$ and non-parabolic $(D_A\neq 0)$.
\snewpage

\subsection{Principal and chiral disks (review)}\plabel{ss:PrincDisk}
~\\

For $2\times2$ real matrices we can refine the spectral data as follows:
Assume that $A=\bem a&b\\c&d\eem =\tilde a\Id_2 +\tilde b\tilde I+\tilde c\tilde J+\tilde d\tilde K$.
Its principal disk is
\[\PD(A):=\Dbar\left(\frac{a+d}2+\frac{|c-b|}2\mathrm i,
\sqrt{\left(\frac{a-d}2\right)^2+\left(\frac{b+c}2\right)^2} \right)
=\Dbar\left(\tilde a+|\tilde b|\mathrm i,\sqrt{{\tilde c}^2+{\tilde d}^2}\right).\]

This is refined further by the chiral disk
\[\CD(A):=\Dbar\left(\frac{a+d}2+\frac{c-b}2\mathrm i,
\sqrt{\left(\frac{a-d}2\right)^2+\left(\frac{b+c}2\right)^2} \right)
=\Dbar\left(\tilde a+\tilde b\mathrm i,\sqrt{{\tilde c}^2+{\tilde d}^2}\right).\]
The additional data in the chiral disk is the chirality, which is the sign of the twisted trace,
$\sgn\tr\left( \begin{bmatrix}&1\\-1&\end{bmatrix}A\right)=\sgn (c-b)=\sgn\tilde b.$
This chirality is, in fact, understood with respect to a fixed orientation of $\mathbb R^2$.
It does not change if we conjugate $A$ by a rotation, but it changes sign if
we conjugate $A$ by a reflection.
From the properties of the twisted trace, it is also easy too see that $\log$ respects chirality.

One can read off many data from the disks.
For example, if $\PD(A)=\Dbar((\tilde a,\tilde b),r)$, then $\det A=\tilde a^2+\tilde b^2-r^2$.
This is not surprising in the light of
\begin{lemma}\plabel{lem:PDCD}
$\CD$ makes a bijective correspondence between
possibly degenerated disks in $\mathbb C$ and
the orbits of $\mathrm M_2(\mathbb R)$ with respect to
conjugacy by special orthogonal matrices (i. e. rotations).

$\PD$ makes a bijective correspondence between
possibly degenerated disks with center in $\mathbb C^+$
and the orbits of $\mathrm M_2(\mathbb R)$ with respect to
conjugacy by  orthogonal matrices.
\qed
\end{lemma}
The principal / chiral disk is a point if $A$ has the effect of a complex multiplication (that $A$ is a quasicomplex matrix).
In general, matrices $A$ fall into three categories: elliptic, parabolic, hyperbolic;
such that the principal /chiral disk are disjoint, tangent or secant to the real axis, respectively.

\snewpage
\subsection{Conformal range of $2\times2$ real matrices (review)}\plabel{ss:ConfRange}
~\\

The principal / conformal disks turn out to be closely related to the conformal range:
\begin{lemma}\plabel{lem:2key}
Consider the real matrix
$
%\begin{equation}
A=\tilde a\Id_2 +\tilde b\tilde I+\tilde c\tilde J+\tilde d\tilde K.
%\plabel{eq:Areal}\end{equation}
$
We claim:

(a) For $A$ acting on $\mathbb R^2$,
\[\CRext(A^{\mathbb{R}})=\partial\Dbar\left(\tilde a+\tilde b\mathrm i,\sqrt{{\tilde c}^2+{\tilde d}^2}\right)
\cup \partial\Dbar\left(\tilde a-\tilde b\mathrm i,\sqrt{{\tilde c}^2+{\tilde d}^2}\right).\]

(b) For $A$ acting on $\mathbb C^2$,
\begin{align}\notag
\CRext(A^{\mathbb C})=
& \Dbar\left(\tilde a+\tilde b\mathrm i,\sqrt{{\tilde c}^2+{\tilde d}^2}\right)
\setminus \intD\left(\tilde a-\tilde b\mathrm i,\sqrt{{\tilde c}^2+{\tilde d}^2}\right)\\
&\notag\cup \Dbar\left(\tilde a-\tilde b\mathrm i,\sqrt{{\tilde c}^2+{\tilde d}^2}\right)
\setminus \intD\left(\tilde a+\tilde b\mathrm i,\sqrt{{\tilde c}^2+{\tilde d}^2}\right).
\end{align}
This is $\CRext(A^{\mathbb{R}})$ but with the components of $\mathbb C\setminus \CRext(A^{\mathbb{R}})$ disjoint from
$\mathbb R$ filled in.
\qed
\end{lemma}
Thus $\CR(A^{\mathbb{R}})$ is the boundary of the principal (or chiral) disk factored by conjugation.
In terms of hyperbolic geometry (in the Poincar\'e half-space),
$\CR(A^{\mathbb R})$ may yield points or circles in the elliptic case;
asymptotic points or corresponding horocycles (paracycles) in the parabolic case;
lines or pairs of distance lines (hypercycles) in the hyperbolic case.
(In the normal case, it yields points, asymptotic points or lines;
in the non-normal case, it yields circles, asymptotically closed horocycles,
asymptotically closed pairs of distance lines.)
$\CR(A^{\mathbb{C}})$ is $\CR(A^{\mathbb{R}})$  but the $h$-convex closure.

\subsection{Recognizing log-ability}\plabel{ss:RecLogab}
~\\

For finite matrices $\spec(A)\cap\mathbb R=\CRext(A^{\mathrm{(real)}})\cap\mathbb R$.
Consequently, $A$  is $\log$-able if and only if $\CRext(A^{\mathrm{(real)}})\cap(\infty,0]=\emptyset$.
($\CRext$ can be replaced by $\CR$).
Or, for $2\times2$ real matrices, in terms of the principal disk, $A$  is $\log$-able  if and only if $\PD(A)\cap(\infty,0]=\emptyset$.
($\PD$ can be replaced by $\CD$.)
\snewpage

\subsection{Canonical forms of $2\times2$ matrices in skew-quaternionic form}\plabel{ss:CanonForm}
~\\

In $\mathrm M_2(\mathbb R)$,
the effect of conjugation by a rotation matrix $\Id_2\cos\theta+\tilde I\sin\theta$ is  given by
\begin{equation}
a\Id_2+b\tilde I+c\tilde J+d\tilde K\mapsto a\Id_2+b\tilde I+((\cos2\theta)\Id_2+(\sin2\theta)\tilde I)(c\Id_2+d\tilde I)\tilde K.
\plabel{eq:roteff}
\end{equation}
Thus, through conjugation by a rotation matrix, any real $2\times2$ matrix can be brought into shape
\[a\Id_2+b\tilde I+c\tilde J\]
with
\[\text{$a\in\mathbb R$, $b\in\mathbb R$, $c\geq0$}. \]
If we also allow conjugation by $\tilde J$, then
\[b\geq0\]
can be achieved. ($\sgn b$ is the chirality class of the matrix.)

In $\mathrm M_2(\mathbb C)$, using conjugation by rotation matrices,
the coefficient of $\tilde K$ can be eliminated.
Then using conjugation by diagonal unitary matrices, the
phase of the coefficient of $\tilde I$ can be adjusted.
In that way the form
\begin{equation}
a\Id_2+\mathrm e^{\mathrm i\beta}(s_1\tilde I+s_2\tilde J)
\plabel{eq:canonC}
\end{equation}
with
\[\text{$a\in\mathbb C$, $s_1\geq0$, $s_2\geq0$, $\beta\in\mathbb R$}\]
can be achieved.
Conjugating by $\tilde K$,  $\beta\in[0,\pi)$ or a similar restriction can be assumed.
Using the definite abuse of notation
$\sqrt{\tilde K}=\frac{1+\mathrm i}2\Id_2+\frac{1-\mathrm i}2\tilde K$,
we find $\sqrt{\tilde K} \tilde I \sqrt{\tilde K}^{-1}=-\mathrm i\tilde J$
and  $\sqrt{\tilde K} \tilde J \sqrt{\tilde K}^{-1}=-\mathrm i\tilde I$;
and $\sqrt{\tilde K}$ is unitary.
Thus, by a further unitary conjugation, \[\beta\in[0,\pi/2)\] or a similar restriction can be assumed.
(We could also base a canonical form on $\tilde J$ and $\tilde K$, but the
present form is conveniently close to the real case.)

Here normal matrices are characterized by that $s_1=0$ or $s_2=0$.
\\

\subsection{The visualization of certain subsets of $\mathrm M_2(\mathbb R)$}\plabel{ss:Visual22R}
~\\

In what follows we will consider certain subsets $S$  of $\mathrm M_2(\mathbb R)$.
By the rotational effect \eqref{eq:roteff},
an invariant subset $S$ is best to be visualized through the image of the mapping
\[\Xi^{\mathrm{PH}}:\qquad a\Id_2+b\tilde I+c\tilde J+d\tilde K\mapsto (a,b,\sqrt{c^2+d^2}).\]
($\Xi^{\mathrm{PH}}(S)$ can be considered as a subset of the asymptotically closed Poincar\'e half-space.)

In particular, the image of the $0$ centered sphere  in $\mathrm M_2(\mathbb R)$ of radius $\frac\pi2$ is
\[\Xi^{\mathrm{PH}}\left(\partial\mathcal B_{\frac\pi2} \right)=\left\{(a,b,r)\,:\, \sqrt{a^2+b^2}+r=\frac\pi2\right\},\]
a ``conical hat''.

\snewpage
\subsection{Arithmetic consequences of dimension $2$}\plabel{ss:arit2}
~\\

For $2\times 2$ complex matrices $A$,
the Cayley--Hamilton equation
reads as
\begin{equation}
A^2=(\tr A)A-(\det A)\Id_2.
\plabel{eq:arit1}
\end{equation}

Obvious consequences are as follows: For the trace,
\[\tr A^2=(\tr A)^2-2(\det A);\]
%\begin{commentx}
i. e.,
\[\det A=\frac{(\tr A)^2-\tr A^2}2;\]
%\end{commentx}
in the invertible case,
\[A^{-1}=\frac1{\det A}\left((\tr A)\Id_2  -A\right);\]
and, more generally, for the adjugate
\begin{equation}
\adj A=(\tr A)\Id_2-A.
\plabel{eq:adj2}
\end{equation}

Applying \eqref{eq:adj2} to the identity $(\adj \m v)(\adj A)=\adj( A \m v) $, and expanded, we obtain
\begin{equation}
A\m v+\m vA=(\tr A)\m v+(\tr\m v)A+ \tr (A\m v)\Id_2-(\tr A)(\tr\m  v)\Id_2.
\plabel{eq:arit2}
\end{equation}

Multiplying by $A$, and applying \eqref{eq:arit1}, it is easy to deduce that
\begin{equation}
A\m v A=\det (A)\m v+\tr( A\m v)A-(\det A)(\tr \m v)\Id_2.
\plabel{eq:arit3}
\end{equation}
This leads to
\begin{lemma}\plabel{lem:prefive}
(o) Any non-commutative (real) polynomial of $A,\m v$
can be written as linear combinations of
\begin{equation}
\Id_2, A,\m v,[A,\m v]
\plabel{eq:linn}\end{equation}
with coefficients which are (real) polynomials of
\begin{equation}
\text{$\tr A$, $\tr \m v$, $\det A$, $\det \m v$, $\tr (A\m v)$.}
\plabel{eq:poll}
\end{equation}

(a) Any scalar expression built up from $\Id_2,A,\m v$ and $\tr$, $\det$, $\adj$ using algebra operations
can by written as a polynomial of \eqref{eq:poll}.

(b)  Any matrix expression built up from $\Id_2,A,\m v$ and $\tr$, $\det$, $\adj$ using algebra operations
can by written as linear combination of \eqref{eq:linn} with coefficients which are polynomials of \eqref{eq:poll}.
\begin{proof}
(o) Using \eqref{eq:arit1} and \eqref{eq:arit3}, and also their version with the role of $A$ and $\m v$ interchanged, repeatedly,
we arrive to linear combinations $\Id_2,A,\m v, A\m v,\m vA$ (with appropriate coefficients).
Due to, \eqref{eq:arit2}, however, the last two terms can be traded to $[A,\m v]$.
(a--b) are best to be proven by simultaneous induction on the complexity of the algebraic expressions.
In view of (o), the only nontrivial step is when we take determinant of a matrix expression.
This, however, can be resolved by using $\det X=\dfrac{(\tr X)^2-\tr X^2}{2}$.
\end{proof}
\end{lemma}
If we examine the proof of Lemma \ref{lem:prefive}, we see that
it also gives an algorithm for the reduction to standard form.
\begin{example}
\plabel{ex:five}
If we apply Lemma \ref{lem:prefive} with $A$ and $\mathbf v=A^*$,
then \eqref{eq:poll} yields
\begin{equation*}
\tr A,\, \overline{\tr A},\, \det A,\, \overline{\det A},\, \tr (A^*A).%\plabel{eq:polc}
\end{equation*}
(Although $X$ and its conjugate ${\overline X}$ contains the same information, they are arithmetically different objects.)
By taking complex linear combinations, we have
\begin{equation*}
\Rea\tr A,\, \Ima\tr A,\, \Rea\det A,\, \Ima\det A,\, \tr (A^*A).%\plabel{eq:polr}
\end{equation*}
which are real quantities.
These are the ``five data'' associated to $A$ (cf. \cite{L.ell}).
\qedexer
\end{example}

In general, if we allow adjoints in our expressions, then their complexity increases.
We only note the identity
\begin{multline}
\frac{\tr(A^*A)\tr(A^*\m v)}2-(\det A^*)\left(  (\tr A)(\tr\m v)-\tr (A\m v)\right)=\\=
\tr\left(\left(A^*A-\frac{A^*A}{2}\Id_2 \right)\left(A^*\m v-\frac{A^*\m v}{2}\Id_2 \right)\right),
\plabel{eq:NormAv}
\end{multline}
which can be checked by direct computation.
\begin{commentx}
Another identity is
\begin{multline}
|(\tr A)(\tr\m v)-\tr (A\m v)|^2
=\\=
{\tr(A^*A)\tr(\m v^*\m v)+|\tr(A^*\m v)|^2-\tr(A^*A\m v^*\m v)-\tr(A A^*\m v\m v^*)}
.\plabel{eq:specAv}
\end{multline}

\end{commentx}
(Matrix expressions with ``more than two variables'' can also be dealt systematically, but the picture is more complicated.)
\begin{commentx}
However, the following special case is easy to prove:

\begin{lemma}\plabel{lem:refive}
Any complex scalar expression built up from $\Id_2,A,\m v$ and $\tr$, $\det$, $\adj$, $*$ using algebra operations
which is real linear in $\mathbf v$ can by written as  linear combinations of
\begin{gather}
\tr(\m v), \tr(A\m v), \tr(A^*\m v), \tr([A,A^*]\m v)
\plabel{eq:linnQ}\\
\tr(\m v^*), \tr(A\m v^*), \tr(A^*\m v^*), \tr([A,A^*]\m v^*)\notag
\end{gather}
with coefficients which are polynomials of
\begin{equation}
\text{$\tr A$, $\tr A^*$, $\det A$, $\det A^*$, $\tr (A^*A)$.}
\plabel{eq:pollQ}
\end{equation}
That is, as a polynomial \eqref{eq:linnQ} and \eqref{eq:pollQ} but with pure real linear degree $1$ in $\m v$.
\begin{proof}
This can be proven  similarly to (or from) Lemma \ref{lem:refive}.
\end{proof}
\end{lemma}
\end{commentx}

Rearranging \eqref{eq:arit1} into a full square, we obtain
\begin{equation}
\left(A-\frac{\tr A}2\Id_2\right)^2=(-D_A) \Id_2.
\plabel{eq:reduc}
\end{equation}

It is useful to introduce the notation
\[T_{A,\m v}=\frac12\tr\left(\left(A-\frac{\tr A}2\Id_2\right)\left(\m v-\frac{\tr \m v}2\Id_2\right)\right).\]
%\snewpage

By \eqref{eq:reduc},  $D_A=-T_{A,A}$, thus $T_{A,\m v}$ can be thought as a the polarization $D_A$ (times $-1$).
Then one can prove
\begin{equation}
[A,[A,[A,\m v]]]=-4 D_A[A,\m v]
\plabel{eq:tricommut}
;
\end{equation}
\begin{equation}
[A,[\m v,[A,\m v]]]=[\m v,[ A,[A,\m v]]]=4T_{A,\m v}[A,\m v]
\plabel{eq:tricommut2}
;
\end{equation}
\begin{equation}
[\m v,[\m v,[A,\m v]]]=-4 D_{\m v}[A,\m v]
\plabel{eq:tricommut3}.
\end{equation}

Indeed, in \eqref{eq:tricommut}, $A$ can be replaced by $A-\frac{\tr A}2\Id_2$,
in which form, \eqref{eq:reduc} implies \eqref{eq:tricommut}.
Then \eqref{eq:tricommut3} follows by interchange of variables,
and \eqref{eq:tricommut2} follows by polarization in the first two commutator variables.

The commutator analogue of Lemma \ref{lem:prefive} is
\begin{lemma}\plabel{lem:prefivecomm}
Any (real) commutator polynomial of $A,\m v$ with no uncommmutatored terms
can be written as linear combinations of
\begin{equation}
[A,\m v] , [A,[A,\m v]] ,  [\m v,[\m v,A]]
\plabel{eq:linncomm}\end{equation}
with coefficients which are (real) polynomials of
\begin{equation}
\text{$D_A$, $T_{A,\m v} $, $D_{\m v}$.}
\plabel{eq:pollcomm}
\end{equation}
\begin{proof}
Commutator expressions can be reduced to linear combinations of iterated left commutators, whose ends of length four are
either trivial or can be reduced by \eqref{eq:tricommut} or \eqref{eq:tricommut2} (essentially).
After the reduction only commutators of length $2$ or $3$ are left with appropriate coefficients.
Conversely,  by \eqref{eq:tricommut} or \eqref{eq:tricommut2} or \eqref{eq:tricommut3},
the multipliers $D_A$, $T_{A,\m v} $, $D_{\m v}$
can be absorbed to commutators.
\end{proof}
\end{lemma}
\begin{lemma}
\plabel{lem:indep} We find:

(a) In Lemma \ref{lem:prefive}, the entries \eqref{eq:linn} and \eqref{eq:poll} are formally algebraically independent.

(b) In Lemma \ref{lem:prefivecomm}, the entries \eqref{eq:linncomm} and \eqref{eq:pollcomm} are formally algebraically independent.
\begin{proof}
(a) Even the particular case of Example \ref{ex:five} demonstrates that the entries of \eqref{eq:linn} and \eqref{eq:poll}
are algebraically independent over the reals.
This also implies independence in general case over the reals.
As all the algebraic rules are inherently real, this also implies independence over the complex numbers.

(b) It is easy to see that in the statement of Lemma \ref{lem:prefive}, $\det A$, $\tr (A\m v)$, $\det \m v$ can be
replaced by $D_A$, $T_{A,\m v} $, $D_{\m v}$, respectively.
Regarding the base terms, the identities
\begin{equation}
[A,[A,\m v]]
=
-4D_A\left(\m v-\frac{\tr\m v}2\Id_2\right)
-4T_{A,\m v}\left(A-\frac{\tr A}2\Id_2\right)
\plabel{eq:biag1}
\end{equation}
and
\begin{equation}
[\m v,[\m v,\m A]]
=
-4T_{A,\m v}\left(\m v-\frac{\tr\m v}2\Id_2\right)
-4D_{\m v}\left(A-\frac{\tr A}2\Id_2\right)
\plabel{eq:biag2}
\end{equation}
show the independence of the new base terms as
$\left|\begin{matrix}D_A&T_{A,\m v}\\T_{A,\m v}&D_{\m v}\end{matrix}\right|\not\equiv0$.
\end{proof}
\end{lemma}

As a consequence, we have obtained a normal form for (formal) commutator expressions
in the case of $2\times 2$ matrices.
As the absorption rules (\ref{eq:tricommut}--\ref{eq:tricommut3}) are quite simple, this normal form is quite practical.
%\snewpage

We also mention the identities
\begin{equation}
T_{A,[A,\m v]}=0;\plabel{eq:tid1}
\end{equation}
\begin{equation}
T_{A,[A,[A,\m v]]}=0;\plabel{eq:tid2}
\end{equation}
\begin{equation}
T_{A,[\m v,[\m v,A]]}= 4(D_AD_{\m v} -T_{A,\m v}^2);\plabel{eq:tid3}
\end{equation}
which are easy to check by direct computation.
%\snewpage

\subsection{Self-adjointness, conform-unitarity, normality}\plabel{ss:SACUN}
~\\

Note that if $B$ is a $2\times2$ complex self-adjoint matrix, then it is of real hyperbolic or parabolic type.
Thus
\[-D_B\geq0 \qquad\text{if}\qquad B=B^*, \]
with equality if and only if $B$ is a (real) scalar matrix.

In particular, for any $2\times2$ complex matrix $A$,
\[-D_{A^*A}
=\left(\frac{\tr (A^*A)}2\right)^2-|\det A|^2
=\frac12\tr\left(A^*A-\frac{\tr(A^*A)}2\Id_2\right)^2
\quad\geq\quad0.\]
Equality holds if and only $A$ is normal and its eigenvalues have equal absolute values,
i. e. $A$ is conformal-unitary (i. e. it is a unitary matrix times a positive scalar)  or $A$ is the $0$ matrix.
(This follows by considering the unitary-conjugated triangular form of the matrices.)
Note that the $0$ matrix is excluded from being conform-unitary.
In the real case, we can speak about conform-orthogonal matrices.
Then, $2\times2$ conform-orthogonal matrices are either conform-rotations or conform-reflections.

For $2\times2$ complex matrices, in terms of the Frobenius norm,
\[\|[A,A^*] \|_{\mathrm{Frob}}^2=\tr\left([A,A^*]^*[A,A^*]\right)=\tr \,[A,A^*]^2=-2\det [A,A^*]=-2D_{ [A,A^*]}.\]
Thus, $\tr \,[A,A^*]^2=0$, $\det [A,A^*]=0$, $D_{[A,A^*]}=0$
are all equivalent to the normality of $A$.
Actually, one can see, in terms of the operator norm,
$\|[A,A^*] \|_{\mathrm{Frob}}^2=2\|[A,A^*]\|_2^2$.

A somewhat strange quantity is $-D_{\frac{\overline{(\tr A)}A+(\tr A)A^*}2}$.
One can check that it vanishes if and only if $\tr A=0$ or $A$ is conform-unitary.

\begin{lemma}
\plabel{lem:Dexp}
If $A$ is a complex $2\times2$ matrix, then
\[-D_{A^*A}=\left(\frac{\tr (A^*A)}2\right)^2-|\det A|^2;\]
\[-D_{\frac{\overline{(\tr A)}A+(\tr A)A^*}2}=  \dfrac{\tr(A^*A)|\tr A|^2-\det (A^*)\tr(A)^2-\det (A)\tr(A^*)^2}4;\]
\[-D_{\frac12[A,A^*]}=\left(\frac{\tr A^*A}2-\frac{|\tr A|^2}4\right)^2-\left|(\det A)-\frac{(\tr A)^2}4\right|^2.\]
\begin{proof}
Straightforward computation.
\end{proof}
\end{lemma}
\begin{commentx}
\begin{remark}
\plabel{lem:Rexp}
In the terminology of \cite{L.ell}, $-D_{\frac12[A,A^*]}=(U_A)^2-|D_A|^2$.
The quantities above are ``reduced'' in the sense that they depend only on the
conformal range of $A$, cf. \cite{L.ell}.
\qedremark
\end{remark}
\end{commentx}

In the case of real matrices the computations are typically much simpler,
especially in skew-quaternionic form.
For example,
\begin{lemma}
\plabel{lem:DexpSQ}
If $A=\tilde a\Id_2+\tilde b\tilde I+\tilde c\tilde J+\tilde d\tilde K$ is a real $2\times2$ matrix, then
\[-D_{A^*A}=4(\tilde a^2+\tilde b^2)(\tilde c^2+\tilde d^2),\]
\[-D_{\frac{\overline{(\tr A)}A+(\tr A)A^*}2}= 4\tilde a^2(\tilde c^2+\tilde d^2), \]
\[-D_{\frac12[A,A^*]}=4\tilde b^2(\tilde c^2+\tilde d^2).\]
\begin{proof}
Simple computation.
\end{proof}
\end{lemma}

%\snewpage
\subsection{Exponentials (review, alternative)}\plabel{ss:ExpRev}
~\\

Recall that $\Cos$ and $\Sin$  are entire functions on the complex plane such that
\[\Cos(z)=\sum_{n=0}^\infty(-1)^n\frac{z^{n}}{(2n)!}
%\]
\text{\qquad
and
\qquad}
%\[
\Sin(z)=\sum_{n=0}^\infty(-1)^n\frac{z^{n}}{(2n+1)!}.\]
For $x\in\mathbb R$,
\[\Cos(x)=\begin{cases}
\cos\sqrt{x}&\text{if }x>0 \\
1&\text{if }x=0\\
\cosh\sqrt{-x}&\text{if }x<0,
\end{cases}
%\]
\text{\qquad
and
\qquad}
%\[
\Sin(x)=\begin{cases}
\dfrac{\sin\sqrt{x}}{\sqrt{x}}&\text{if }x>0 \\
1&\text{if }x=0\\
\dfrac{\sinh\sqrt {-x}}{\sqrt {-x}}&\text{if }x<0.
\end{cases}\]
\begin{lemma}\plabel{lem:expQuasi}
Let $A$ be a complex $2\times2$ matrix. Then
\[\exp A=\left(\exp\frac{\tr A}2\right)\cdot
\left(
\Cos \left(D_A\right)\Id_2
+\Sin \left(D_A\right)\left(A-\frac{\tr A}2\Id_2\right)
\right).\]
\begin{proof}
$A=\frac{\tr A}2\Id_2+\left(A-\frac{\tr A}2\Id_2\right)$ gives a decomposition to
commuting operators for which, by multiplicativity, the exponential can be computed separately.
In the case of the first summand this is trivial.
In the case of the second summand, the identity \eqref{eq:reduc} and the comparison of the power series implies the statement.
\end{proof}
\end{lemma}
~

\subsection{The differential calculus of $\Cos$ and $\Sin$ (review)}\plabel{ss:SinCosCalc}
~\\

First of all, it is useful to notice that
\[z=\frac{1-\Cos(z)^2}{\Sin(z)^2}\]
(as entire analytic functions).
Then one can easily see that
\[\Cos'(z)=-\frac12\Sin(z)\]
and
\begin{equation*}
%\begin{align*}
\Sin'(z)
%&
=\frac{\Cos(z)-\Sin(z)}{2z}\\
%&
=\frac12  \frac{\Sin(z)^2\cdot(\Cos(z)-\Sin(z))}{1-\Cos(z)^2}
%\end{align*}
\end{equation*}
(as entire analytic functions).
In particular, differentiation will not lead out of the rational field
generated by $\Cos$ and $\Sin$.
~\\
\snewpage

\subsection{A notable differential equation (review)}
\begin{lemma}\plabel{lem:fullprec}
The solution of the ordinary differential equation
\[\frac{\mathrm dA(\theta)}{\mathrm dt}A(\theta)^{-1}=
a\bem&-1\\1&\eem+
b\begin{bmatrix}
-\sin 2c\theta& \cos2c\theta\\\cos2c\theta&\sin2c\theta
\end{bmatrix}
\equiv a\tilde I+ \exp(c\theta\tilde I) b\tilde K \exp(-c\theta\tilde I),
\]
with initial data
\[A(0)=\begin{bmatrix}
1& \\&1
\end{bmatrix} \equiv\Id_2,\]
is given by
\[A(\theta)=F(a\theta,b\theta,c\theta);\]
where
\[F(s,p,w)=  \exp(w\tilde I) \exp (-(w-s)\tilde I+p\tilde K) .\eqed\]

\end{lemma}
In particular,
\[\Lexp\left(\left(a\bem&-1\\1&\eem+
b\begin{bmatrix}
-\sin 2c\theta& \cos2c\theta\\\cos2c\theta&\sin2c\theta
\end{bmatrix}\right)d\theta|_{[0,\pi]}\right)=F(a\pi,b\pi,c\pi).\]

We will also use the special notation
\[W(p,w)=F(0,p,w)\equiv  \exp(w\tilde I) \exp (-w\tilde I+p\tilde K) .\]
~

\subsection{The calculus of $\Cot$}\plabel{ss:CotCalc}
~\\

In what follows, we also use the meromorphic function
\begin{commentx}
s
\[\Tan(z)=\frac{\Sin(z)}{\Cos(z)},\]
and
\end{commentx}
\[\Cot(z)=\frac{\Cos(z)}{\Sin(z)}.\]

\begin{commentx}
$\Tan$ has poles at $z=(k-\frac12)^2\pi^2$, where $k$ is a positive integer, and
\end{commentx}
$\Cot$  has poles at $z=k^2\pi^2$, where $k$ is a positive integer.
\begin{commentx}
We will mainly use the latter function, as $\Cot$ defines an analytic function on $(-\infty,\pi^2)$;
while $\Tan$ has a pole at $\pi^2/4$.
\end{commentx}
It is easy to see that
\[\Cot'(z)=-\frac12\left(1-\frac{\Cot(z)}z +\frac{\Cot(z)^2}z\right).\]

Consequently, $\mathbb R\boldsymbol (z,\Cot(z)\boldsymbol)$ is a differential field (induced from, say, $(-\infty,\pi^2)$),
which is subset of the differential field  $\mathbb R\boldsymbol(\Sin(z),\Cos(z)\boldsymbol)$.
It is not much smaller: As
\[\Sin(z)^2=\frac1{\Cot(z)^2+z},\]
it is not hard to show that the degree of the field extension is
\[\mathbb R\boldsymbol(\Sin(z),\Cos(z)\boldsymbol):\mathbb R\boldsymbol (z,\Cot(z)\boldsymbol)=2.\]
~

\subsection{A collection of auxiliary functions}\plabel{ss:AuxFunct}
~\\
\begin{lemma}
\plabel{lem:CTH}

The expressions

(a)
\[\reC(z)=-2\frac{\Sin'(z)}{\Sin(z)}=\frac{-\Cot(z)+1}{z}\]
\begin{commentx}
\[=\frac{\Sin(z)\cdot(\Sin(z)-\Cos(z))}{1-\Cos(z)^2}
={\frac {-\cos \left( \sqrt {z} \right) \sqrt {z}+\sin \left( \sqrt {z}\right) }{\sin \left( \sqrt {z} \right) z}}\]
\end{commentx}

(b)
\[\reD(z)=2z\reC'(z)+\reC(z)=z\reC(z)^2-2\reC(z)+1=\frac{\Cot(z)^2+z-1}{z}\]
\begin{commentx}
\[=\frac{1-\Sin(z)^2}{1-\Cos(z)^2}
=\frac {  z-\left(\sin \sqrt {z} \right)^2}{{   z\left(\sin \sqrt {z} \right)^2}}\]
\end{commentx}

(c)
\[\reW(z)=2 \reC'(z)=\frac{\Cot(z)^2+\Cot(z)+z-2}{z^2}\]
\begin{commentx}
\[=\frac{\Sin(z)^2(\Cos(z)\Sin(z)-2\Sin(z)^2+1 )}{(1-\Cos(z)^2)^2}\]
\end{commentx}

(d)
\[\reP(z)=\reC(z)^2 - 2 \reC'(z)=\frac{-3\Cot(z)-z+3}{z^2}\]
\begin{commentx}
\[=\frac{\Sin(z)^2(\Cos(z)^2-3\Cos(z)\Sin(z)+3\Sin(z)^2-1 )}{(1-\Cos(z)^2)^2}\]
\end{commentx}

(e)
\[\reG(z)=\reD'(z)=
2z\reC''(z)+3\reC'(z)=
z\reC(z)^3-2\reC(z)^2-2\reC'(z)+\reC(z)\]
\[=\frac{-\Cot(z)^3-z\Cot(z)+1}{z^2}\]
\begin{commentx}
\[=\frac{\Sin(z)\cdot(\Sin(z)^3-\Cos(z))}{(1-\Cos(z)^2)^2} \]
\end{commentx}

(f)
\[\reL(z)= \left(\frac{\reD}{\reC^2}\right)'(z) \cdot\reC(z)^3
%=\reC(z)\reG(z)-\reW(z)\reD(z)
=\reC(z)^3+2\reC'(z)\reC(z)-2\reC(z)\]
\[=\frac{-2\Cot(z)^3-z\Cot(z)^2+3\Cot(z)^2-2z\Cot(z)-z^2+3z-1}{z^3}\]
\begin{commentx}
\[=\frac{\Sin(z)^2\cdot\left( (\Sin(z)-\Cos(z))^2-(1-\Sin(z)^2)^2  \right)}{(1-\Cos(z)^2)^3 }\]
\end{commentx}

(g)
\[\reX(z)=\frac43\reC''(z)={\frac {-2\,{\Cot(z)}^{3}-3\,{\Cot(z)}^{2
}-2\,z\,{\Cot(z)}-3\,{\Cot(z)}-3\,z+8}{3{z}^{3}}}\]
define meromorphic functions $ \reC, \reD, \reW, \reP, \reG, \reL, \reX$
 on the complex plane such that

(i) They are holomorphic on $\mathbb C$ except at $z=(k\pi)^2$ where $k$ is a positive integer.
In those latter points they have poles of order
\[
\text{
$\reC\mapsto 1$, $\reD\mapsto 2$, $\reW\mapsto2 $, $\reP\mapsto 1$, $\reG\mapsto 3$, $\reL\mapsto 3$, $\reX\mapsto 3$
}.
\]

(ii) They are strictly monotone increasing on the interval $(-\infty,\pi^2)$ with range $(0,\infty)$.

(iii) In particular,  they are strictly positive on the interval $(-\infty,\pi^2)$.

\begin{proof}
These facts can be derived using various methods of (complex) analysis.
\end{proof}
\end{lemma}
The following table contains easily recoverable information:
\[
\begin{array}{|c|c|c|c|}
\hline
&&&\\
\text{function}& \text{asymptotics as $x\searrow-\infty$} & \text{value at $x=0$}&\text{ asymptotics at $z\sim k^2\pi^2$ }\\
&  & &\text{ ($k\in\mathbb N\setminus\{0\}$)}\\
&&&\\
\reC&|x|^{-1/2}&\dfrac13& -2\cdot\dfrac1{z-k^2\pi^2}\\
&&&\\
\reD&|x|^{-1}&\dfrac13& 4\pi^2\cdot\dfrac{k^2}{(z-k^2\pi^2)^2}\\
&&&\\
\reW&|x|^{-3/2}&\dfrac{2}{45}& 4\cdot\dfrac{1}{(z-k^2\pi^2)^2}\\
&&&\\
\reP&|x|^{-1}&\dfrac1{15}& -\dfrac{6}{\pi^2}\cdot\dfrac1{k^2(z-k^2\pi^2)}\\
&&&\\
\reG&|x|^{-2}&\dfrac1{15}& -8\pi^2\cdot\dfrac{k^2}{(z-k^2\pi^2)^3}\\
&&&\\
\reL&|x|^{-3}&\dfrac1{135}& -16\cdot\dfrac{1}{(z-k^2\pi^2)^3}\\
&&&\\
\reX&|x|^{-5/2}&\dfrac{16}{2835 }& -\dfrac{16}3\cdot\dfrac{1}{(z-k^2\pi^2)^3}\\
&&&\\
\hline
\end{array}
\]
~

Furthermore, it is also easy to check that the identities
\[  \reC(z)^2  = \reW(z) + \reP(z)  \]
\[  \reC(z) \reD(z) = \reW(z) + \reG(z)  \]
\[  \reC(z) \reG(z) = \reD(z) \reW(z) + \reL (z)  \]
\[  \reP(z) \reD(z) = \reC(z) \reW(z) + \reL(z)   \]
\[  \reG(z) \reP(z) = \reC(z) \reL(z) + \reW(z)^2   \]
\[  \reC(z)^2  \reD(z) = \reC(z) \reW(z) + \reD(z) \reW(z) + \reL(z)   \]
\[ \reC(z)\reP(z)+2\reC(z)\reW(z) = \reL(z)+\reW(z)\]
hold.

\begin{commentx}
Cf.: The rule of differentiation
\[\reC'(z)=\frac{z\reC(z)^2-3\reC(z)+1}{2z}\]
holds. In particular, the rules of simplification
\begin{equation}
z\reC'(z)=\frac{z\reC(z)^2-3\reC(z)+1}2
\plabel{eq:simp1}
\end{equation}
and
\begin{equation}
z\reC''(z)=\frac{z\reC(z)^3-2\reC(z)^2-5\reC'(z)+\reC(z)}2
\plabel{eq:simp2}
\end{equation}
hold.
\end{commentx}

\begin{lemma}
\plabel{lem:DTH}
(a) The function $\left.x\mapsto\frac1{\sqrt{\reD(x)}}\right|_{x\in(-\infty,\pi^2)}$
extends to an analytic function $\mathcal E$ on (a neighbourhood of) $\mathbb R$.
Special values are $\mathcal E(0)=\sqrt3$, $\mathcal E(\pi^2)=0$.

(b) The function $\left.x\mapsto\frac{\reC(x)}{\sqrt{\reD(x)}}\right|_{x\in(-\infty,\pi^2)}$
extends to an analytic function $\mathcal F$ on (a neighbourhood of) $\mathbb R$.
Special values are  $\mathcal F(0)=\frac1{\sqrt3}$, $\mathcal F(\pi^2)=\frac1\pi$.

(c) As analytic functions,
\[
\mathcal E'(z)=-\frac12 \reG(z)\mathcal E(z)^3,
\qand
\mathcal F'(z)=-\frac12 \reL(z)\mathcal E(z)^3.
\]
In particular, $\mathcal E|_{(-\infty,\pi^2]}$ and $\mathcal F|_{(-\infty,\pi^2]}$ are strictly decreasing.
\begin{proof}
(a) It is sufficient to find a compatible extension to $(0,+\infty)$.
One can check that, on the real line, $\reD(x)$ has poles with asymptotics
$\sim\frac{4\pi^2k^2}{(x-k^2\pi^2)^2}$ at $x=k^2\pi^2$ for positive integers $k$;
and otherwise $\reD(x)$ is positive on the real line.
This implies that
\[\mathcal E(x)=\frac{\sgn\sin\sqrt x }{\sqrt{\reD(x)}}\qquad\text{for $x>0$} \]
\begin{commentx}
\[\mathcal E(x)={\frac {\sin \left( \sqrt {x} \right) }{\sqrt {x}}  \sqrt {{\frac {{x}^{2}}{
 \left( \cos \left( \sqrt {x} \right)  \right) ^{2}+x-1}}}}
\]
\end{commentx}
defines an appropriate extension.
The special values are straightforward.

(b) As $\reC$ has only simple poles with asymptotics $\sim\frac{-2}{x-k^2\pi^2}$,
these are cancelled out once multiplied by $\mathcal E(x)$.
The special values are straightforward.

(c) This is straightforward.
\end{proof}
\end{lemma}
~
\snewpage
\subsection{The function $\AC$ (review, alternative)}\plabel{ss:ACRev}
~\\

If $z\in\mathbb C\setminus(-\infty,-1]$, and $t\in [0,1]$, then it is easy to see that
$1+2t(1-t)(z-1)$
is invertible.
(For a fixed $z$, the possible values yield the segment
$\left[1,\frac{1+z}2\right]_{\mathrm e}$.)
Then one can define, for $z\in\mathbb C\setminus(-\infty,-1]$,
\begin{equation}
\AC(z)=\int_{t=0}^1\frac{1}{1+2t(1-t)(z-1)}\,\mathrm dt.
\plabel{eq:ACdef}
\end{equation}

This gives an analytic function in $z\in\mathbb C\setminus(-\infty,-1]$.
\begin{lemma}
\plabel{lem:AC}
(a) In terms of the real domain,
\[\AC(x)=\begin{cases}
\dfrac{\arccos x}{\sqrt{1-x^2}}&\text{if }-1< x<1\\[3mm]
1&\text{if }x=1\\[1mm]
\dfrac{\arcosh x}{\sqrt{x^2-1}}\qquad&\text{if } 1<x.\\
\end{cases}\]

%\begin{commentx}
Note, this can be rewritten as $\AC(x)=\dfrac1{\Sin(\Cos^{-1}(x))}$ for $x\in(-1,+\infty)$.
%\end{commentx}

(b) As an analytic function,
\begin{equation}
\AC'(z)=\frac{z\AC(z)-1}{1-z^2}.
\plabel{eq:ACfun}
\end{equation}

(c) $\AC$ vanishes nowhere on $\mathbb C\setminus(-\infty,-1]$.
$\AC(z)=\pm1$ holds only for $z=1$ with $\AC(1)=1$.

(d) $\AC(z)=1-\frac13(z-1)+O((z-1)^2)$ at $z\sim 1$.

(e) $\AC$ is strictly monotone decreasing on $(-1,+\infty)$ with range $(+\infty,0)_{\mathrm e}$.

\begin{proof}
(b) For $z\neq\pm1$, differentiating under the integral sign, we find
\begin{align*}
\AC'(z)=&\int_{t=0}^1\frac{-2t(1-t)}{(1+2t(1-t)(z-1)  )^2}\,\mathrm dt\\
=&\int_{t=0}^1\frac{ z\dfrac{1}{ 1+2t(1-t)(z-1)}-1}{1-z^2}\,\mathrm dt
+\int_{t=0}^1\frac{ 1-6t(1-t)-4t^2(1-t)^2(z-1)  }{(1+z)(1+2t(1-t)(z-1) )^2}\,\mathrm dt
\\
=&\frac{\displaystyle{ z\int_{t=0}^1\dfrac{1}{ 1+2t(1-t)(z-1)}\,\mathrm dt-1}}{1-z^2}
+\left[\frac{ t(1-t)(1-2t)}{(1+z)( 1+2t(1-t)(z-1))}\right]_{t=0}^1
\\
=&\frac{z\AC(z)-1}{1-z^2}+0.
\end{align*}

(a)
The special value $\AC(1)=1$ is trivial.
Restricted to $x\in(-1,1)$, using \eqref{eq:ACfun}, we find
$\dfrac{\mathrm d}{\mathrm dx}\left(\AC(x)\sqrt{1-x^2}\right)=-\dfrac1{\sqrt{1-x^2}}$.
Considering the primitive functions, we find $\AC(x)=\frac{(\arccos x)+c}{\sqrt{1-x^2}}$
for $x\in(-1,1)$, with an appropriate $c\in\mathbb C$.
As $\AC$ limits to a finite value for $x\nearrow1$, only the case $c=0$ is possible.
Restriction to $(1,+\infty)$ is similar.

(e) is immediate from \eqref{eq:ACdef}.

(c) From \eqref{eq:ACdef}, $\sgn\Ima\AC(z)=-\sgn\Ima z$ is immediate.
Then it is sufficient to restrict to $(-1,+\infty)$.
Then the statement follows from (d).

(d) Elementary analysis.
\end{proof}
\end{lemma}

From Lemma \ref{lem:AC}(a), it is easy to the see that, for $x\in(-\infty,\pi^2)$,
\begin{equation}
\AC(\Cos(x))=\frac1{\Sin(x)}
\plabel{eq:logexpkey}
\end{equation}
holds.
By analytic continuation, it  also holds in an open neighbourhood of $(-\infty,\pi^2)$.

\begin{commentx}
So does
\[\Cos^{-1}(x)=(1-x^2)\AC(x)^2,\]
but this, in fact, extends to $z\in\mathbb C\setminus(-\infty,-1]$,
as long as we consider the principal branch of $\Cos^{-1}$.

The related identities
\[z=\Cos((1-z^2)\AC(z)^2)\]
and
\[1=\AC(z)\Sin((1-z^2)\AC(z)^2)\]
also hold for $z\in\mathbb C\setminus(-\infty,-1]$.
\\
\end{commentx}

\begin{commentx}
We remark that a closely related function, defined
for $z\in\mathbb C\setminus(-\infty,-1]$, is given by
\[\ATT(z)=\frac{\AC\left(\dfrac{1-z}{1+z}\right)}{1+z}.\]
In terms of the real domain, it has the property that
\[\ATT(x)=\begin{cases}
\dfrac{\artanh\sqrt{- x}}{\sqrt{ -x}}&\text{if }-1<x<0,\\
1&\text{if }x=0,\\
\dfrac{\arctan\sqrt x}{\sqrt x}&\text{if }x>0.
\end{cases}\]
One can check that it satisfies the simple differentiation rule
\[\frac{\mathrm d}{\mathrm dz}\ATT(z)=\frac1{2z(1+z)}-\frac{\ATT(z)}{2z},\]
as an analytic function.
\end{commentx}

\snewpage
\subsection{Logarithm (review, alternative)}\plabel{ss:LogRev}
~\\

Recall that, in a Banach algebra,
$A$ is $\log$-able if and only if the spectrum of $A$ is disjoint from $(-\infty,0]$.
In that case, the logarithm of $A$ is defined as
\begin{equation}
\log A=\int_{\lambda=0}^1 \frac{A-1}{\lambda+(1-\lambda)A}\,\mathrm d\lambda
=\int_{t=0}^1 \frac{A-1}{(1-t)+tA}\,\mathrm dt.
\plabel{eq:logdef}
\end{equation}

\begin{lemma} \plabel{lem:logreal}
Let $A$ be a  real $2\times 2$ matrix.

(a) Then  $A$ is a $\log$-able if and only if
 $\det A>0$ and $\dfrac{\tr A}{2\sqrt{\det A}}>-1$.

(b) In the $\log$-able case
\begin{equation}\log A= (\log\sqrt{\det A})\Id_2+\frac{\AC\left(\dfrac{\tr A}{2\sqrt{\det A}}\right)}{\sqrt{\det A}}
\left(A- \frac{\tr A}2\Id_2\right).\plabel{eq:log2}\end{equation}
\end{lemma}

Suppose that $A$ is $n\times n$ complex matrix which is $\log$-able.
Let
\[\sqrt{\dett A}:= \det\sqrt A =\exp\frac{\tr \log A}{2}=\exp\frac12\int_{t=0}^1 \tr\frac{\mathrm d((1-t)\Id+tA )}{(1-t)\Id+tA}.\]
One can see that this is the product of the square root of eigenvalues (with multiplicity).
\begin{lemma}\plabel{lem:logcomplexA}
Suppose that $A$ is $2\times2$ complex matrix which is $\log$-able.

(a) Then $\sqrt{\dett A}\in\mathbb C\setminus (\infty,0]$, $\dfrac{\tr A}{2\sqrt{\dett A}}\in\mathbb C\setminus (-\infty,-1]$.

(b) Furthermore,
the extended form of \eqref{eq:log2} holds:
\begin{equation}\log A= (\log\sqrt{\dett A})\Id_2+\frac{\AC\left(\dfrac{\tr A}{2\sqrt{\dett A}}\right)}{\sqrt{\dett A}}
\left(A- \frac{\tr A}2\Id_2\right).\plabel{eq:log2c}
\eqed\end{equation}
\begin{proof}[Proofs]
Lemma \ref{lem:logreal}(a) and Lemma \ref{lem:logcomplexA}(a) can be obtained from the
examination of the eigenvalues in a relatively straightforward manner.
For the rest, it is sufficient to prove \ref{lem:logcomplexA}(b).
The proof is almost ``tautological'':
Using the differentiation rule \eqref{eq:ACfun},
by a long but straightforward symbolic computation,
one obtains, for $t\in[0,1]$,
\[ \frac{A-\Id_2}{(1-t)\Id_2+tA}=\]\[=\frac{\mathrm d}{\mathrm dt}\left(
\left(
\log\sqrt{\dett ( (1-t)\Id_2+tA)}\right)\Id_2
+\frac{\AC\left(\dfrac{\tr ( (1-t)\Id_2+tA)}{2\sqrt{\dett ( (1-t)\Id_2+tA)}}\right)}{\sqrt{\dett ( (1-t)\Id_2+tA)}}t
\left(A- \frac{\tr A}2\Id_2\right)
\right)\]
(Indeed, the symbolic computation is perfectly valid for, say, $t\sim1$, $A\sim\exp\tilde J$, where
$\dfrac{\tr ( (1-t)\Id_2+tA)}{2\sqrt{\dett ( (1-t)\Id_2+tA)}}\sim\cosh1\neq\pm1$;
and even $\sqrt{\dett\cdot}$ can be replaced by $\sqrt{\det\cdot}$.
Then, by analyticity in $t$ and analyticity in $A$, the general identity valid.)
Integrated, \eqref{eq:log2c} is obtained.
\end{proof}
\end{lemma}
~\snewpage

\subsection{$\AC$ near $-1$, some asymptotics}\plabel{ss:asymptAC}
~\\

From Lemma \ref{lem:AC}(a) (and analytic continuation), it easy to see that
\[\AC(y-1)=\frac{\pi}{\sqrt{y}\sqrt{2-y}}-\AC(1-y)\]
is valid for $y\in\mathbb C\setminus ( (-\infty,0] \cup [2,\infty  ) )$.
This is informative regarding what kind of analytic reparametrizations are useful for $\AC$ near $-1$.
More in terms of real analysis, one has
\begin{lemma}
\plabel{lem:asymptAC}
(a) For $y\in (0,+\infty)$, the function
\[y\mapsto \frac{\AC(y-1)}{\dfrac{\pi}{\sqrt 2\sqrt{y}}}\]
is monotone decreasing with range $(1,0)_{\mathrm e}$.

(b) In particular, for $0<y_1<y_2$,
\[\frac{\AC(y_2-1)}{\AC(y_1-1)}<\frac{\sqrt{y_1}}{\sqrt{y_2}}.\]
\begin{proof}
Part (a) can addressed by elementary analyis.
Part (a) implies part (b).
%Part (b) is an immediate consequence of (a).
\end{proof}
\end{lemma}
~

\snewpage
\subsection{Norms (review)}\plabel{ss:NormRev}

\begin{lemma} \plabel{lem:normcomputeC}
Let $A
%=\begin{bmatrix}a&b\\c&d\end{bmatrix}=\tilde a\Id_2+\tilde b\tilde I+\tilde c\tilde J+\tilde d\tilde K
$ be a real or complex matrix. Then, for the norm,
\begin{align}
\left\|A \right\|_2
&=\sqrt{\frac{\tr(A^*A)}2+\sqrt{-D_{A^*A}}}
\plabel{eq:normcomputeC}
\\\notag
&=\sqrt{\frac{\tr(A^*A)}2+\sqrt{\frac{(\tr(A^*A))^2}4-|\det A|^2}}
\\\notag
&=\frac{\sqrt{\tr(A^*A)+2|\det A|}+\sqrt{\tr(A^*A)-2|\det A|} }2
%\\
%&=\frac{\sqrt{|a|^2+|b|^2+|c|^2+|d|^2+2|ad-bc| }+\sqrt{|a|^2+|b|^2+|c|^2+|d|^2-2|ad-bc|}}2\\
%&=\sqrt{\frac{|\tilde a|^2+|\tilde b|^2+|\tilde c|^2+|\tilde d|^2+ |\tilde a^2+\tilde b^2-\tilde c^2-\tilde d^2| }{2}}\\
%&\qquad +\sqrt{\frac{|\tilde a|^2+|\tilde b|^2+|\tilde c|^2+|\tilde d|^2- |\tilde a^2+\tilde b^2-\tilde c^2-\tilde d^2| }{2}}
;
\end{align}
and, for the co-norm,
\begin{align}
\|A\|_2^-
%&
=\left\|A^{-1} \right\|_2^{-1}
%\plabel{eq:cpnormcomputeC}
%\\\notag
&=\sqrt{\frac{\tr(A^*A)}2-\sqrt{-D_{A^*A}}}
\plabel{eq:cpnormcomputeC}
\\\notag
&=\sqrt{\frac{\tr(A^*A)}2-\sqrt{\frac{(\tr(A^*A))^2}4-|\det A|^2}}
\\\notag
&=\frac{\sqrt{\tr(A^*A)+2|\det A|}-\sqrt{\tr(A^*A)-2|\det A|} }2
%\\
%&=\frac{\sqrt{|a|^2+|b|^2+|c|^2+|d|^2+2|ad-bc| }-\sqrt{|a|^2+|b|^2+|c|^2+|d|^2-2|ad-bc|}}2\\
%&=\sqrt{\frac{|\tilde a|^2+|\tilde b|^2+|\tilde c|^2+|\tilde d|^2+ |\tilde a^2+\tilde b^2-\tilde c^2-\tilde d^2| }{2}}\\
%&\qquad -\sqrt{\frac{|\tilde a|^2+|\tilde b|^2+|\tilde c|^2+|\tilde d|^2- |\tilde a^2+\tilde b^2-\tilde c^2-\tilde d^2| }{2}}
.
\end{align}
In particular,
\begin{equation}
\|A\|_2\cdot\|A\|_2^-=|\det A|
.
\plabel{eq:compmult}
\end{equation}

In the case of real matrices, the results are the same for the Hilbert spaces  $\mathbb R^2$ and $\mathbb C^2$.
\qed
\end{lemma}

For  $2\times2$ matrices, we define the signed co-norm as
\begin{equation}
\lfloor A\rfloor_2=
\begin{cases}
0&\text{if }A=0\\\\
\dfrac{\det A}{\|A\|_2}&\text{if }A\neq 0.\\
\end{cases}
\plabel{eq:signedconorm}
\end{equation}

Then,
\begin{equation}
\left| \left\lfloor A\right\rfloor_2 \right|=\|A\|_2^-,
\end{equation}
and
\begin{equation}
\|A\|_2\cdot \left\lfloor A\right\rfloor_2 =\det A.\plabel{eq:realmult}
\end{equation}
However, we will essentially consider the signed co-norm only for real matrices.
\snewpage
\begin{lemma} \plabel{lem:normcompute}
Let $A=\begin{bmatrix}a&b\\c&d\end{bmatrix}
=\tilde a\Id_2+\tilde b\tilde I+\tilde c\tilde J+\tilde d\tilde K$ be a real matrix. Then
\begin{align}
\left\|A \right\|_2
&=\frac{\sqrt{\tr(A^*A)+2\det A}+\sqrt{\tr(A^*A)-2\det A} }2
\plabel{eq:2norm}\\
&=\frac{\sqrt{(a+d)^2+(c-b)^2}+\sqrt{(a-d)^2+(b+c)^2}}2\notag\\
&=\sqrt{\tilde a^2+\tilde b^2}+\sqrt{\tilde c^2+\tilde d^2}.
\plabel{eq:norm22skew}
\end{align}
On the other hand,
\begin{align}
\|A\|_2^-=\left\|A^{-1} \right\|_2^{-1}
&=\left|\frac{\sqrt{\tr(A^*A)+2\det A}-\sqrt{\tr(A^*A)-2\det A} }2\right|
\plabel{eq:2conorm1}\\
&=\left|\frac{\sqrt{(a+d)^2+(c-b)^2}-\sqrt{(a-d)^2+(b+c)^2}}2\right|\notag\\
&=\left|\sqrt{\tilde a^2+\tilde b^2}-\sqrt{\tilde c^2+\tilde d^2}\right|
.
\notag
\end{align}
It is true that
\begin{align}
\sgn \det A
&=\sgn\frac{\sqrt{\tr(A^*A)+2\det A}-\sqrt{\tr(A^*A)-2\det A} }2
\plabel{eq:2conorm2}\\
&= \sgn \frac{\sqrt{(a+d)^2+(c-b)^2}-\sqrt{(a-d)^2+(b+c)^2}}2\notag\\
&=\sgn\left(\sqrt{\tilde a^2+\tilde b^2}-\sqrt{\tilde c^2+\tilde d^2}\right).
\notag
\end{align}
Furthermore,
\begin{align}
\left\lfloor A\right\rfloor_2=\sgn(\det A )\left\|A \right\|_2^-
&=\frac{\sqrt{\tr(A^*A)+2\det A}-\sqrt{\tr(A^*A)-2\det A} }2
\plabel{eq:2conorm}\\
&=\frac{\sqrt{(a+d)^2+(c-b)^2}-\sqrt{(a-d)^2+(b+c)^2}}2\notag\\
&=\sqrt{\tilde a^2+\tilde b^2}-\sqrt{\tilde c^2+\tilde d^2}
.\eqed
\notag
\end{align}

\end{lemma}
~
\snewpage

\subsection{Directional derivatives}\plabel{ss:DirDer}
~\\

Whenever $f:M\mapsto f(M)$ is a function on an open domain of matrices,
we define the derivative of $f$ at $A$ in direction $\m v$, along smooth curves as
\[\mathrm D_{\m v\text{ at }M=A}\left( f(M)\right)=\lim_{t\rightarrow 0} f(\gamma(t)),\]
whenever it gives the same value for all $\gamma$ such that $\gamma$ is smooth,
$\gamma(0)=A$, $\gamma'(0)=\m v$.
This is a sufficiently flexible notion to deal with some mildly singular $f$.
If $f$ is smooth, then the directional derivatives agrees to the usual multidimensional differential.
E.~g.~
\[\mathrm D_{\m v\text{ at }M=A}\left( M^{-1}\right)=-A^{-1}\m vA^{-1},\]
or
\begin{equation}
\mathrm D_{\m v\text{ at }M=A}\left( \tr M\right)=\tr \m v,
\plabel{eq:basder1}
\end{equation}
\begin{equation}
\mathrm D_{\m v\text{ at }M=A}\left( \det M\right)=\tr \left(\m v\adj A \right),
\plabel{eq:basder22}
\end{equation}
\begin{equation}
\mathrm D_{\m v\text{ at }M=A}\left( \tr(M^*M)\right)=2\Rea\tr (A^*\m v).
\plabel{eq:basder3}
\end{equation}

Note, for $2\times2$ matrices, \eqref{eq:basder22} reads as
\begin{equation}
\mathrm D_{\m v\text{ at }M=A}\left( \det M\right)=(\tr A)(\tr\m v)-\tr (A\m v).
\plabel{eq:basder2}
\end{equation}
Moreover, also for $2\times2$ matrices,
\begin{equation}
\mathrm D_{\m v\text{ at }M=A}\left( D_M\right)=-2 T_{A,\m v}.
\plabel{eq:deet1}
\end{equation}

Furthermore,
\begin{equation}
\mathrm D_{\m v\text{ at }M=A}\left( T_{M,B}\right)=T_{\m v,B}.
\plabel{eq:deet2}
\end{equation}
~

\snewpage
\subsection{Smoothness of norm}\plabel{ss:SmoothNorm}
~\\

Now we examine the smoothness properties of the norm of complex $2\times2$ matrices.
For a complex $2\times2$ matrix $\m v$, we set
\[S(\m v):=\frac12\tr\frac{\m v+\m v^*}2+\sqrt{-D^{\phantom{M}}_{\frac{\m v+\m v^*}2}}.\]
(This is the higher eigenvalue of $\frac{\m v+\m v^*}2$.)
It is easy to see that for
\[\m v=v_a\Id+v_b\tilde I+v_c\tilde J+v_d\tilde K\]
(complex coefficients), it yields
\[S(\m v)=(\Rea v_a)+\sqrt{ (\Ima v_b)^2+(\Rea v_c)^2+(\Rea v_d)^2  }.\]

\begin{lemma}
\plabel{lem:smooth2}
(o) On the space of complex $2\times2$ matrices,
the norm operation is smooth except at matrices $A$ such that the ``norm discriminant''
\[-D_{A^*A}=0.\]

(a) The directional derivative of the function $A\mapsto\|M\|_2$ at $M=0$, along smooth curves, is,
 just \[\mathrm D_{\m v\text{ at }M=0}\left( \|M\|_2\right)=\|\m v\|_2.\]

(b) The directional derivative of the function $M\mapsto\|M\|_2$ at $M=\Id$, along smooth curves, is
\begin{align*}
\mathrm D_{\m v\text{ at }M=\Id}\left( \|M\|_2\right)
&=S(\m v).
\end{align*}

(c) The directional derivative of the function $M\mapsto\|M\|_2$ at any conform-unitary $A$, along smooth curves, is
\[
\mathrm D_{\m v\text{ at }M=A}\left( \|M\|_2\right)=\|A\|_2\cdot S(A^{-1}\m v)=\|A\|_2\cdot S(\m v A^{-1})
=\frac{S(A^*\m v)}{\|A\|_2}=\frac{S(\m v A^*)}{\|A\|_2}
.\]

(d) The directional derivative of the function $M\mapsto\|M\|_2$ at any not conform-unitary $A\neq0$, is
\[\mathrm D_{\m v\text{ at }M=A}\left( \|M\|_2\right)=\]
\[=\frac1{\|A\|_2}\Rea\left(\frac{\tr(A^*\mathbf v)}2 +\frac{\dfrac12
\tr\left(
\left(A^*A-\dfrac{\tr(A^*A)}2\Id_2\right)\left(A^*\mathbf v-\dfrac{\tr(A^*\mathbf v)}2\Id_2\right)
\right)
}{\sqrt{-D_{A^*A}}}\right)\]
\[=\frac{ \|A\|_2\Rea\left(\tr(A^*\m v) \right) -\|A\|_2^{-1}\Rea\left(\det(A^*)\Bigl(\tr(A)\tr(\m v) -\tr(A\m v)\Bigr)\right)
}{2\sqrt{-D_{A^*A}}}.\]
\begin{proof}
(o) The smoothness part follows from the norm expressions
% with encapsulated square roots in Lemma \ref{lem:normcomputeC}.
(\ref{eq:normcomputeC}/1--2).
The non-smoothness part will follow from the explicit expressions for the directional derivatives.
Part (a) is trivial.
Part (b) is best to be done using Taylor expansions of smooth curves.
Part (c) follows from (b) by the conform-unitary displacement argument:
Neighbourhoods of conform-unitary matrices are related to each other
just by multiplication by conform-unitary matrices (left or right, alike).
Such a multiplication, however, just simply scales up the norm.
Part (d) follows from the basic observations (\ref{eq:basder1}--\ref{eq:basder3})
and standard composition rules.
\end{proof}
\end{lemma}
\snewpage
\begin{remark}
\plabel{rem:altsmooth}
If we apply the notation
\[\boldsymbol\{A,B\boldsymbol\}=A^*B+B^*A\]
and
\[\boldsymbol\langle A,B\boldsymbol\rangle=\frac12\tr
\left(
\left(A-\frac{\tr A}2\Id_2\right)\left(B-\frac{\tr B}2\Id_2\right)
\right),\]
then we find
\[\|A\|_2=\sqrt{ \frac12\tr\tfrac12\boldsymbol\{A,A\boldsymbol\}
+\sqrt{ \boldsymbol\langle \tfrac12\boldsymbol\{A,A\boldsymbol\} , \tfrac12\boldsymbol\{A,A\boldsymbol\} \boldsymbol\rangle  }  },\]
and in case Lemma \ref{lem:smooth2} (c),
\[\mathrm D_{\m v\text{ at }M=A}\left( \|M\|_2\right)=
\frac{\dfrac12\tr\tfrac12\boldsymbol\{A,\m v\boldsymbol\}
+
\sqrt{ \boldsymbol\langle \tfrac12\boldsymbol\{A,\m v\boldsymbol\} , \tfrac12\boldsymbol\{A,\m v\boldsymbol\} \boldsymbol\rangle  }
 }{\|A\|_2};\]
and in case Lemma \ref{lem:smooth2} (d),
\[\mathrm D_{\m v\text{ at }M=A}\left( \|M\|_2\right)=
\frac{\dfrac12\tr\tfrac12\boldsymbol\{A,\m v\boldsymbol\}
+
\dfrac{
{ \boldsymbol\langle \tfrac12\boldsymbol\{A,A\boldsymbol\} , \tfrac12\boldsymbol\{A,\m v\boldsymbol\} \boldsymbol\rangle  }
}{
\sqrt{ \boldsymbol\langle \tfrac12\boldsymbol\{A,A\boldsymbol\} , \tfrac12\boldsymbol\{A,A\boldsymbol\} \boldsymbol\rangle  }
}
 }{\|A\|_2};\]
 making the computation more suggestive.
 \qedremark
\end{remark}

For the sake of curiosity, we include the corresponding statement for co-norms.
Here the situation is analogous but slightly more complicated as there is an additional case for non-smoothness.
Indeed, the signed co-norm \eqref{eq:signedconorm} has similar smoothness properties as the norm.
However, when absolute value is taken to obtain the co-norm, then there is an additional source
for non-smoothness: the case when the determinant vanishes.

For a complex $2\times2$ matrix $\m v$, we set
\[S^-(\m v):=\frac12\tr\frac{\m v+\m v^*}2-\sqrt{-D^{\phantom{M}}_{\frac{\m v+\m v^*}2}}.\]
(This is the lower eigenvalue of $\frac{\m v+\m v^*}2$.)
It is easy to see that for
\[\m v=v_a\Id+v_b\tilde I+v_c\tilde J+v_d\tilde K\]
(complex coefficients), it yields
\[S^-(\m v)=(\Rea v_a)-\sqrt{ (\Ima v_b)^2+(\Rea v_c)^2+(\Rea v_d)^2  }.\]

\snewpage
\begin{lemma}
\plabel{lem:cosmooth2}
(o) On the space of complex $2\times2$ matrices,
the co-norm operation is smooth except at matrices $A$ such that
\[\det A=0\]
or
\[-D_{A^*A}=0.\]

(a) The directional derivative of the function $A\mapsto\|M\|_2^-$ at $M=0$, along smooth curves, is,
 just \[\mathrm D_{\m v\text{ at }M=0}\left( \|M\|_2^-\right)=\|\m v\|_2^-.\]

(\'{a}) The directional derivative of the function $M\mapsto\|M\|_2^-$ at $A$ with $\det A=0$, $A\neq0$, along smooth curves, is
\[\mathrm D_{\m v\text{ at }M=A}\left( \|M\|_2^-\right)
=\frac{|(\tr A)(\tr\m v)-\tr (A\m v)|}{\sqrt{\tr (A^*A) }}.\]
(In this case, $\|A\|_2=\sqrt{\tr (A^*A) }$ and $\|A\|_2^-=0$.)

(b) The directional derivative of the function $M\mapsto\|M\|_2^-$ at $M=\Id$, along smooth curves, is
\begin{align*}
\mathrm D_{\m v\text{ at }M=\Id}\left( \|M\|_2^-\right)
&=S^-(\m v).
\end{align*}

(c) The directional derivative of the function $M\mapsto\|M\|_2^-$ at any conform-unitary $A$, along smooth curves, is
\[\mathrm D_{\m v\text{ at }M=A}\left( \|M\|_2^-\right)=\|A\|_2^-\cdot S^-(A^{-1}\m v)=\|A\|_2^-\cdot S^-(\m v A^{-1})
=\frac{S^-(A^*\m v)}{\|A\|_2^-}=\frac{S^-(\m v A^*)}{\|A\|_2^-}
.\]
(In this case, $\|A\|_2=\|A\|_2^-=\sqrt{\tr (A^*A)/2 }$.)

(d) The directional derivative of the function $M\mapsto\|M\|_2^-$ at any not conform-unitary $A\neq0$, is
\[\mathrm D_{\m v\text{ at }M=A}\left( \|M\|_2^-\right)=\]
\[=\frac1{\|A\|_2^-}\Rea\left(\frac{\tr(A^*\mathbf v)}2 -\frac{\dfrac12
\tr\left(
\left(A^*A-\dfrac{\tr(A^*A)}2\Id_2\right)\left(A^*\mathbf v-\dfrac{\tr(A^*\mathbf v)}2\Id_2\right)
\right)
}{\sqrt{-D_{A^*A}}}\right)\]
\[
=\frac{(\|A\|_2^-)^{-1}\Rea\left(\det(A^*)\Bigl(\tr(A)\tr(\m v)-\tr(A\m v) \Bigr)\right) -\|A\|_2^-\Rea\left(\tr(A^*\m v)\right)
}{2\sqrt{-D_{A^*A}}}.\]

\begin{proof}
Similar to the previous statement.
\end{proof}
\end{lemma}

For the signed co-norm,
the corresponding statement to Lemma \ref{lem:smooth2}
is easy to recover from
\[\mathrm D_{\m v\text{ at }M=A}\left( \lfloor M\rfloor_2\right)=
\frac{(\tr A)(\tr\m v)-\tr (A\m v)}{\|A\|_2}
-\frac{\det A}{\|A\|^2_2}\cdot \mathrm D_{\m v\text{ at }M=A}\left( \|M\|_2\right)\]
(for $A\neq0$).
But this counterpart to Lemma \ref{lem:smooth2} is not particularly simple.

However, things are manageable for real $2\times 2$ matrices;
which case is even more simple as the use of $\Rea$ and $\Ima$ can be avoided:

\begin{lemma}
\plabel{lem:realsmooth}
(o) On the space of real $2\times2$ matrices,
the norm and signed co-norm operations are smooth except at matrices $A$ such that
\[-D_{A^*A}=0.\]

(a)
The directional derivative of the function $A\mapsto\|M\|_2$ at $M=0$, along smooth curves, is
\[\mathrm D_{\m v\text{ at }M=0}\left( \|M\|_2\right)=\|\m v\|_2.\]

The directional derivative of the function $A\mapsto\lfloor M\rfloor_2^-$ at $M=0$, along smooth curves, is
\[\mathrm D_{\m v\text{ at }M=0}\left( \lfloor M \rfloor_2\right)=\lfloor\m v\rfloor_2.\]

(b)
The directional derivative of the function $M\mapsto\|M\|_2$ at $M=\Id$, along smooth curves, is
\begin{align*}
\mathrm D_{\m v\text{ at }M=\Id}\left( \|M\|_2\right)
&=S(\m v).
\end{align*}
Here, for $\m v=v_a\Id+v_b\tilde I+v_c\tilde J+v_d\tilde K$ (with real coefficients),
\[S(\m v)=v_a+\sqrt{ ( v_c)^2+( v_d)^2  }.\]

The directional derivative of the function $M\mapsto\lfloor M\rfloor_2$ at $M=\Id$, along smooth curves, is
\begin{align*}
\mathrm D_{\m v\text{ at }M=\Id}\left( \lfloor M\rfloor_2\right)
&=S^-(\m v).
\end{align*}
Here, for $\m v=v_a\Id+v_b\tilde I+v_c\tilde J+v_d\tilde K$ (with real coefficients),
\[S^-(\m v)=v_a-\sqrt{ ( v_c)^2+( v_d)^2  }.\]

(c) The directional derivative of the function $M\mapsto\|M\|_2$ at any conform-unitary $A$, along smooth curves, is
\[\mathrm D_{\m v\text{ at }M=A}\left( \|M\|_2\right)
=\|A\|_2\cdot S(A^{-1}\m v)
=\|A\|_2\cdot S(\m v A^{-1})=\frac{S(A^*\m v)}{\|A\|_2}=\frac{S(\m v A^*)}{\|A\|_2}.\]

The directional derivative of the function $M\mapsto\lfloor M\rfloor_2$ at any conform-unitary $A$, along smooth curves, is
\[\mathrm D_{\m v\text{ at }M=A}\left( \lfloor M\rfloor_2\right)
=\lfloor A\rfloor_2\cdot S^-(A^{-1}\m v)
=\lfloor A\rfloor_2\cdot S^-(\m v A^{-1})=\frac{S^-(A^*\m v)}{\lfloor A\rfloor_2}=\frac{S^-(\m v A^*)}{\lfloor A\rfloor_2}.\]
(In this case, $\|A\|_2=\|A\|_2^-=\sqrt{\tr (A^*A)/2 }$.)

(d)
The directional derivative of the function $M\mapsto\|M\|_2$ at any not conform-unitary $A\neq0$, is
\[\mathrm D_{\m v\text{ at }M=A}\left( \|M\|_2\right)
=\frac{ \|A\|_2\tr(A^*\m v) -\lfloor A\rfloor_2\Bigl( \tr(A)\tr(\m v)-\tr(A\m v) \Bigr)
}{2\sqrt{-D_{A^*A}}}.\]

 The directional derivative of the function $M\mapsto\lfloor M\rfloor_2$ at any not conform-unitary $A\neq0$, is
\[\mathrm D_{\m v\text{ at }M=A}\left( \lfloor M\rfloor_2\right)=\frac{ \| A\|_2 \Bigl(\tr(A)\tr(\m v)-\tr(A\m v) \Bigr)   -\lfloor A\rfloor_2\tr(A^*\m v)
}{2\sqrt{-D_{A^*A}}}.\]
\begin{proof}
It follows from the complex picture.
\end{proof}
\end{lemma}
Next, we concentrate on a special case regarding complex $2\times2$ matrices.
\begin{lemma}
\plabel{lem:DDcomp}
(a) Let $A$ be a complex $2\times2$ matrix. Then
\[D_{A^*A}=D_{\frac{\overline{(\tr A)}A+(\tr A)A^*}2}+D_{\frac12[A,A^*]}.\]

(b) With some abuse of notation, for $\lambda\in[0,1]$, let
\[A(\lambda)=A-\lambda\frac{\tr A}2\Id_2.\]
Then,
\[D_{\frac{\overline{(\tr A(\lambda))}A(\lambda)+(\tr A(\lambda))A(\lambda)^*}2}=
(1-\lambda)^2D_{\frac{\overline{(\tr A)}A+(\tr A)A^*}2},\]
and
\[D_{\frac12[A(\lambda),A(\lambda)^*]}=D_{\frac12[A,A^*]}.\]

(c) The map
\[\lambda\in[0,1]\mapsto \|A(\lambda)\|_2\]
is smooth, with the possible exception of $\lambda=1$.
However, if $A$ is non-normal (or $\tr A=0$), then smoothness also holds at $\lambda=1$.

\begin{proof}
(a), (b) are straightforward computations.
This follows from (a), (b) and the concrete form
$\left\|A \right\|_2=\sqrt{\frac{\tr(A^*A)}2+\sqrt{-D_{A^*A}}}$ of the norm.
Note that the vanishing of $D_{\frac12[A,A^*]}$ is equivalent to normality.
\end{proof}
\end{lemma}

\begin{lemma}
\plabel{lem:Dtrace}
Suppose that $A$ is a complex $2\times2$ matrix.

(a) If $A$ is conform-orthogonal, then
\[\mathrm D_{(\tr A)\Id_2\text{ at }M=A}\left( \|M\|_2\right)=\frac{\dfrac{|\tr A|^2}2}{\|A\|_2}.\]

(b) In the generic case $-D_{A^*A}>0$,
\[\mathrm D_{(\tr A)\Id_2\text{ at }M=A}\left( \|M\|_2\right)=\frac{\dfrac{|\tr A|^2}2+
\dfrac{-D_{\frac{\overline{(\tr A)}A+(\tr A)A^*}2}}{\sqrt{-D_{A^*A}}}
}{\|A\|_2}.\]

(c) Consequently, if $\tr A\neq0$, then
\[\mathrm D_{(\tr A)\Id_2\text{ at }M=A}\left( \|M\|_2\right)>0.\]

%\begin{commentx}
(d) Consequently, if $A\neq0$ is normal, then
\[\mathrm D_{(\tr A)\Id_2\text{ at }M=A}\left( \|M\|_2\right)=\frac{\dfrac{|\tr A|^2}2+
\sqrt{ -D_{\frac{\overline{(\tr A)}A+(\tr A)A^*}2}}
}{\|A\|_2}.\]
%\end{commentx}
\begin{proof}
(a) and (b) are straighforward consequences of Lemma \ref{lem:smooth2}.
(c) is a trivial consequence of (a) and (b).
%\begin{commentx}
So is (d), after considering Lemma \ref{lem:DDcomp}(a).
%\end{commentx}
\end{proof}
\end{lemma}
~
\snewpage
\subsection{Some inequalities for norms}\plabel{ss:InEqNorms}
~\\
\begin{cor}\plabel{lem:detrace}
(a)
If $A$ is a real $2\times2$ matrix, then
\[\left\|A-\frac{\tr A}2\Id_2\right\|_2\leq\|A\|_2\]
and
\[\left\lfloor A-\frac{\tr A}2\Id_2\right\rfloor_2\leq\lfloor A\rfloor_2.\]
The inequalities are strict if $\tr A\neq 0$.

(b) In fact,  if $\tr A\neq 0$, then the maps
\[t\mapsto \left\|A-t\frac{\tr A}2\Id_2\right\|_2\]
and
\[t\mapsto\left\lfloor A-t\frac{\tr A}2\Id_2\right\rfloor_2\]
are strictly decreasing on the interval $t\in[0,1]$.
\begin{proof}
(a) If $A$ is written as in \eqref{eq:skqform}, then
the inequalities trivialize as
$\sqrt{\tilde b^2}+\sqrt{\tilde c^2+\tilde d^2} \leq \sqrt{\tilde a^2+\tilde b^2}+\sqrt{\tilde c^2+\tilde d^2}$
and
$\sqrt{\tilde b^2}-\sqrt{\tilde c^2+\tilde d^2} \leq \sqrt{\tilde a^2+\tilde b^2}-\sqrt{\tilde c^2+\tilde d^2}$.

(b) follows along similar lines.
\end{proof}
\end{cor}
\begin{cor}\plabel{lem:detraceC}
(a)
If $A$ is a complex $2\times2$ matrix, then
\[\left\|A-\frac{\tr A}2\Id_2\right\|_2\leq\|A\|_2.\]
The inequality is strict if $\tr A\neq 0$.

(b) In fact,  if $\tr A\neq 0$, then the map
\[t\mapsto \left\|A-t\frac{\tr A}2\Id_2\right\|_2\]
is strictly decreasing on the interval $t\in[0,1]$.
\begin{proof}
It is sufficient to prove (b), which follows from Lemma \ref{lem:Dtrace}(c).
\end{proof}
\end{cor}
(Similar statement is not true for the unsigned co-norm, not even in the real case.)

\snewpage
\subsection{Moments of linear maps}\plabel{ss:MomentsLin}
~\\

If $\ell:\mathrm M_2(\mathbb C)\rightarrow \mathbb R$ is a (real) linear map, then
it can be represented uniquely by a complex $2\times2$ matrix $A$ such that
\[\ell(\m v)=\frac12\Rea\tr(A^*\m v).\]
This matrix $A$ is the moment associated to $\ell$.

\begin{lemma}
\plabel{lem:MN}
If $-D_{A^*A}\neq0$, then the moment associated to $\m v\mapsto\mathrm D_{\m v\text{ at }M=A}\left( \|M\|_2\right)$ is
\begin{align*}
\MN(A):=
&\frac1{\|A\|_2}\left(A +\frac{ AA^*A-\dfrac{\tr(A^*A)}2A}{\sqrt{-D_{A^*A}}}\right)\\
=&\frac1{\|A\|_2}\left(A +\frac{\det (A)\left(A^*-(\tr A^*)\Id_2\right)  +\dfrac{\tr(A^*A)}2A}{\sqrt{-D_{A^*A}}}\right).
\end{align*}
\begin{proof}
This follows from Lemma \ref{lem:smooth2}(d).
\end{proof}
\end{lemma}
~
\snewpage
\subsection{Possible exponentials from $\mathrm M_2(\mathbb R)$  (review)} \plabel{ss:PosRExp}
~\\
\begin{lemma}
\plabel{lem:PosRExp}
For $A\in\mathrm M_2(\mathbb R)$, let us consider the set
\[\Log^{\mathbb R} A=\{M\in \mathrm M_2(\mathbb R)\,:\, \exp(M)=A \}.\]
The the following hold:

(a) If $A$ is hyperbolic ($D_A<0$) with two positive eigenvalues,
\[\Log^{\mathbb R} A=\{\log A\};\]
where $\log A$ is hyperbolic, with $\spec(\log A)\subset \mathbb R$.

(b) If $A$ is hyperbolic ($D_A<0$) with no two positive eigenvalues, then
\[\Log^{\mathbb R} A=\emptyset.\]

(c) If $A=a\Id_2$ is a positive scalar matrix, then
\[\Log^{\mathbb R} A=\{  (\log a)\Id_2+2k\pi L\,:\, k\in\mathbb Z,\,  L^2=-\Id_2\}.\]
Among its elements, the one with the lowest norm is $\log A= (\log a)\Id_2$.

(d) If $A=a\Id_2$ is a negative scalar matrix, then
\[\Log^{\mathbb R} A=\{  (\log -a)\Id_2+(2k+1)\pi L\,:\, k\in\mathbb Z,\,  L^2=-\Id_2\}.\]
Among its elements, the two elements with the lowest norm are $(\log -a)+\pi\tilde I$ and $(\log -a)-\pi\tilde I$.
($A$ is not log-able).

(e) If $A$ is parabolic ($D_A=0$) but not a positive or negative scalar matrix, then
\[\Log^{\mathbb R} A=\emptyset.\]

(f) If $A$ is elliptic ($D_A>0$), then there are unique $r>0$ and $\phi\in(0,\pi)$ and $I_A\in \mathrm M_2(\mathbb R)$
such that $(I_A)^2=-\Id_2$ and
\[A=r(\cos\phi)\Id_2+r(\sin\phi)I_A.\]
In that case,
\[\Log^{\mathbb R} A=\{  (\log r)\Id_2+ \phi I_A+2k\pi I_A\,:\, k\in\mathbb Z\}.\]
Among its elements, the one with the lowest norm is $\log A= (\log r)\Id_2+ \phi I_A$.
\begin{proof}
$\Log^{\mathbb R} A$ can be recovered by elementary linear algebra.
The comments about the norms follow from the simple norm expression of (\ref{eq:2norm}/3).
(Note that the skew-involutions are all traceless.)
\end{proof}
\end{lemma}

\snewpage
\subsection{Further particularities in the real case}\plabel{ss:RealPart}
~\\

\textbf{(I)}
If $A=\tilde a\Id_2+\tilde b\tilde I+\tilde c\tilde J+\tilde d\tilde K$ is a real $2\times2$ matrix, then
we set
\[\boldsymbol X^+(A):=\tilde a\tilde b\Id_2+(\tilde c^2+\tilde d^2-\tilde a^2)\tilde I-\tilde b\tilde c\tilde J-\tilde b\tilde d\tilde K.\]
This operation is chiral; if $U$ is orthogonal, then
$\boldsymbol X^+(U A U^{-1})=(\det U)\cdot U\boldsymbol X^+( A )U^{-1}$.
One can check that
\[\tr\left(A^*\boldsymbol X^+(A)\right)=0.\]
(Other reasonable choice would be $\boldsymbol X^-(A):=-\boldsymbol X^+(A)$.)

\begin{commentx}
In the setting complex $2\times2$ matrices, the analogous operation is best to
be defined only in the non-selfadjoint case $A\neq A^*$, and only as a two-valued function:
We set
\[\boldsymbol X^\pm(A):=\frac{A^*A^*A+AA^*A^*-AAA-A^*AA^*}{\pm4\sqrt{D_{\frac{A-A^*}2}}}.\]
(If $C$ is nonzero skew-adjoint, then $D_C>0$.)
One can check that
\[\Rea\tr\left(A^*\boldsymbol X^\pm(A)\right)=0\]
(any value chosen).
\end{commentx}

\textbf{(II)}
Assume that $A=a\Id_2+b\tilde I+c\tilde J+d\tilde K$ is a real $2\times2$ matrix such that $\det A>0$.
We
\begin{commentx}
define the $m$-chirality (multiplicative chirality) of $A$ as
\[\mathrm{mchir}\,\,A=\frac{-\tr{(A\tilde I)}}{2\sqrt{\det A}}=\frac{b}{\sqrt{a^2+b^2-c^2-d^2}};\]
and we
\end{commentx}
define the $m$-distortion (multiplicative distortion) of $A$ as
\[\mathrm{mdist}\,\,A=\frac{\sqrt{\tr(A^*A)-2\det A}}{2\sqrt{\det A}}=\frac{\sqrt{c^2+d^2}}{\sqrt{a^2+b^2-c^2-d^2}}.\]
~
\snewpage

\subsection{The flattened hyperboloid model HP and some relatives}\plabel{ss:HypProj}
~\\

(Recall, a review of hyperbolic geometry can be found in Berger \cite{Ber}; for our conventions, see \cite{L2}.)

The flattened hyperboloid model, sometimes also popularized as
the ``Gans model'' (cf. \cite{Gan}), is the ``vertical'' projection of the usual hyperboloid model.
The transcription from the CKB model is
\[\left(x_{\mathrm{HP}},y_{\mathrm{HP}}\right)=
\left(\frac{x_{\mathrm{CKB}}}{ \sqrt{1-(x_{\mathrm{CKB}})^2-(y_{\mathrm{CKB}})^2} },
\frac{y_{\mathrm{CKB}}}{ \sqrt{1-(x_{\mathrm{CKB}})^2-(y_{\mathrm{CKB}})^2} }\right);\]
and the  transcription to the CKB model is
\[\left(x_{\mathrm{CKB}},y_{\mathrm{CKB}}\right)=
\left(\frac{x_{\mathrm{HP}}}{ \sqrt{1+(x_{\mathrm{HP}})^2+(y_{\mathrm{HP}})^2} },
\frac{y_{\mathrm{HP}}}{ \sqrt{1+(x_{\mathrm{HP}})^2+(y_{\mathrm{HP}})^2} } \right).\]
An advantage of HP is that its points are represented by points of $\mathbb R^2$;
its disadvantage is, however, that the asymptotical points of the hyperbolic plane
are represented only by asymptotical points of $\mathbb R^2$.

For technical reasons, we will also consider the $\arctan$-transformed HP model with
\[\left(x_{\mathrm{AHP}},y_{\mathrm{AHP}}\right):=\left(\arctan x_{\mathrm{HP}},y_{\mathrm{HP}}\right).\]
The problem with construction is that it not well-adapted to asymptotic points,
except at $(\pm1,0)^{\mathrm{CKB}}$, where it (or, rather, its inverse) realizes blow-ups.

For this reason, a similar but better construction can be derived from the CKB model,
we let
\begin{align*}
\left(x_{\mathrm{ACKB}},y_{\mathrm{ACKB}}\right):=&\left(\arcsin x_{\mathrm{CKB}},
\frac{y_{\mathrm{CKB}}}{\sqrt{ 1-(x_{\mathrm{CKB}})^2}  }\right)
\\
=&\left(\arctan  \frac{ x_{\mathrm{HP}} }{\sqrt{ 1-(y_{\mathrm{HP}})^2}}  ,
\frac{y_{\mathrm{HP}}}{\sqrt{ 1-(y_{\mathrm{HP}})^2}  }\right).
\end{align*}

It (or, rather, its inverse) realizes a blow-up of unit disk of CKB with to
$[-\frac\pi2,\frac\pi2]\times[-1,1]$.
The points $(\pm1,0)^{\mathrm{CKB}}$ are blown up.

\snewpage
\scleardoublepage\section{Schur's formulae for $2\times2$ matrices}
\plabel{sec:Schur}

\subsection{BCH and Schur's formulae (review)}\plabel{ss:SchurRev}
~\\

Recall, that for formal variables $X$ and $Y$,
\begin{equation}
\BCH(X,Y)=\log((\exp X)(\exp Y)).\plabel{eq:BCHeq}
\end{equation}
Whenever $A$ and $B$ are elements of a Banach algebra $\mathfrak A$, then the convergence of $\BCH(A,B)$ can be considered.
Natural notions of convergence are (a) absolute convergence of terms grouped by joint homogeneity in $A$ and $B$;
(b) absolute convergence of terms grouped by separate homogeneity in $A$ and $B$;
or (c) absolute convergence of terms grouped by non-commutative monomials in $A$ and $B$.
We adopt the first viewpoint as it is equivalent to the convergence of the corresponding Magnus series;
(b) is stricter, and (c) is even stricter and less natural; but even case (c) makes relatively little difference to (a),
as it is discussed in Part I.
Now, if $\BCH(A,B)$ is absolute convergent, then $\exp \BCH(A,B)= (\exp A)(\exp B)$.
Another issue is whether $\BCH(A,B)=\log((\exp A)(\exp B))$ holds.
This latter question is basically about the spectral properties of $(\exp A)(\exp B)$.
As it was discussed in already Part I, if $|A|<\pi$ or just $\spec(A)\subset\{z\,:\,|\Ima z|<\pi\}$, and $t$ is sufficiently small, then
$\BCH(A,tM)=\log((\exp A)(\exp tM))$ and
\begin{equation}
\left.\frac{\mathrm d}{\mathrm dt}\log(\exp(A)\exp(tM))\right|_{t=0}=\beta(-\ad A)M;
\plabel{schur1}
\end{equation}
and similarly, $\BCH(tM,A)=\log((\exp tM)(\exp A))$ and
\begin{equation}
\left.\frac{\mathrm d}{\mathrm dt}\log(\exp(tM)\exp(A))\right|_{t=0}=\beta(\ad A)M;
\plabel{schur2}
\end{equation}
where $\ad(X)Y=[X,Y]=XY-YX $ and
\begin{equation}\beta(x)=\frac x{\mathrm e^x-1}=\sum_{j=0}^\infty \beta_j x^j
=1-{\frac {1}{2}}x+{\frac {1}{12}}{x}^{2}-{\frac {1}{720}}{x}^{4}+\frac1{30240}x^6+\ldots \qquad .\plabel{eq:schur}\notag\end{equation}
\eqref{schur1} and \eqref{schur2} are F. Schur's formulae and they embody partial but practical information about the BCH series.
\\

\snewpage
\subsection{Schur's formulae in the $2\times2$ matrix case}\plabel{ss:Schur22}
~\\

In the $2\times2$ setting, Schur's  formulae can be written in closed form:
\begin{lemma}
\plabel{lem:iven}
If $A,\m v$ are complex matrices, and $\spec(A)\subset\{z\,:\,|\Ima z|<\pi\}$, then
\begin{multline}
\mathrm D_{\m v\text{ at }M=0}\left( \log(\exp(A) \exp(M)) \right)
=\\=
\left.\frac{\mathrm d}{\mathrm dt}\log(\exp(A)\exp(t\m v))\right|_{t=0}
=\m v+\frac12[A,\m v]+\reC(D_A)\cdot\frac14[A,[A,\m v]]\plabel{eq:riven}
\end{multline}
and
\begin{multline}
\mathrm D_{\m v\text{ at }M=0}\left( \log(\exp(M)\exp(A)) \right)=\\=
\left.\frac{\mathrm d}{\mathrm dt}\log(\exp(t\m v)\exp(A))\right|_{t=0}
=\m v-\frac12[A,\m v]+\reC(D_A)\cdot\frac14[A,[A,\m v]]\plabel{eq:liven}
\end{multline}
hold.
\begin{proof}
This is a direct computation based on  the explicit formulas in
Lemma \ref{lem:expQuasi} and Lemma \ref{lem:logcomplexA}, and the rules of derivation.
(The key identity for simplifying back is
\eqref{eq:logexpkey},
which holds only if $x$ is in a certain neighborhood of $(-\infty,\pi^2)$,
but then one can use analytic continuation.)
\end{proof}
\begin{proof}[Alternative proof]
For $A\sim0$, it follows from the traditional form of Schur's formulae
combined with \eqref{eq:tricommut}.
Then the formulas extends by analytic continuation.
\end{proof}
\end{lemma}
\begin{commentx}
\begin{remark}
\plabel{lem:Cthexp}
For $x\sim0$,
\[\reC(x)=\frac13+\frac1{45}\cdot x+\frac4{945}\cdot \frac{x^2}2+\frac2{1575}\cdot \frac{x^3}6+\frac{16}{31185}\cdot \frac{x^4}{24}+O(x^5).
 \eqedremark\]
\end{remark}
\end{commentx}

Simple consequences are
\begin{lemma}
\plabel{lem:LocLog}
If $D_A\notin\{k^2\pi^2\,:\,k\in\mathbb N\setminus\{0\}\}$, then there exists uniquely, a local analytic branch
$\log^{\mathrm{loc}}$ of $\exp^{-1}$ near $\exp A$ such that $\log^{\mathrm{loc}}(\exp A)=A$.
\begin{proof}
The Schur map provides an inverse to the differential of $\exp$, thus the inverse function theorem can be used.
\end{proof}
\end{lemma}

\begin{lemma}
\plabel{lem:MRL}
Assume that $A$ is complex $2\times2$ matrix such that $-D_{A^*A}\neq0$
and
\[\spec(A)\subset\{z\in\mathbb C\,:\,|\Ima z|<\pi\}.\]

(a) Then the moment associated to
$\m v\mapsto\mathrm D_{\m v\text{ at }M=\Id_2}\left( \|\log(\exp(A)\exp(M))\|_2\right)$ is
\[\mathrm{MR}(A):=\MN(A)+ \frac12[A^*,\MN(A)]+\overline{\reC(D_A)}\cdot[A^*,[A^*,\MN(A)]].\]

(b) Then the moment associated to
$\m v\mapsto\mathrm D_{\m v\text{ at }M=\Id_2}\left( \|\log(\exp(M)\exp(A))\|_2\right)$ is
\[\mathrm{ML}(A)=\MN(A)- \frac12[A^*,\MN(A)]+\overline{\reC(D_A)}\cdot[A^*,[A^*,\MN(A)]].\]
\begin{proof}
This follows from Lemma \ref{lem:iven}.
\end{proof}
\end{lemma}
\begin{commentx}
\noindent(The real case will be discussed in greater detail later.)
\end{commentx}

Using Lemma \ref{lem:iven}, we can obtain higher terms in the BCH expansion.
Let $A^{R}(t)=\log(\exp(A)\exp(t\m v))$.
Then
\[\frac{\mathrm d}{\mathrm dt} A^{R}(t) =\m v+\frac12[A^{R}(t),\m v]+\reC(D_{A^{R}(t)})\cdot\frac14[A^{R}(t),[A^{R}(t),\m v]]\]
holds.
Taking further derivatives, we can compute  $\frac{\mathrm d^n}{\mathrm dt^n} A^{R}(t) $.
Due to the rules of derivation, and equations of type (\ref{eq:tid1}--\ref{eq:tid3})
we obtain commutator expressions of $ A^{R}(t) $ , $\m v$,
with coefficients which are polynomials of $\reC^{(n)}(D_{A^{R}(t)} )$ (higher derivatives of $\reC$ as functions)
and $D_{A^{R}(t)}$, $T_{A^{R}(t),\m v}$, $D_{\m v}$.
Specifying $t=0$, we obtain the corresponding higher terms in the BCH expression.
Similar argument applies with $A^{L}(t)=\log(\exp(t\m v)\exp(A))$.
In particular, we obtain
\begin{lemma}\plabel{lem:schur22}
If $A,\m v$ are complex matrices, and $\spec(A)\subset\{z\,:\,|\Ima z|<\pi\}$, then
\begin{multline}
%\begin{align}
\frac{\mathrm d^2}{\mathrm d^2t}\log(
%&
\exp(A)\exp(t\m v))\Bigr|_{t=0}=
%\plabel{schur11}\\\notag=&
\frac{\reC(D_A)+\reD(D_A)}{4}\,[\m v,[\m v,A]]
\plabel{schur11}\\
+\frac{\reC(D_A)}{4}\,\underbrace{[A,[\m v,[\m v,A]]]}_{=-4T_{A,\m v} [A,\m v]}
%\\\notag&
+\frac{2\reP(D_A)+3\reW(D_A)}{16} \, \underbrace{[A,[A,[\m v,[\m v,A]]]]}_{=-4T_{A,\m v} [A,[A,\m v]]}
%\end{align}
\end{multline}
and
\begin{multline}
%\begin{align}
\frac{\mathrm d^2}{\mathrm d^2t}\log(
%&
\exp(t\m v)\exp(A))\Bigr|_{t=0}=
%\plabel{schur22}\\\notag=&
\frac{\reC(D_A)+\reD(D_A)}{4}\,[\m v,[\m v,A]]
\plabel{schur22}\\
-\frac{\reC(D_A)}{4}\,[A,[\m v,[\m v,A]]]
%\\\notag&
+\frac{2\reP(D_A)+3\reW(D_A)}{16} \, [A,[A,[\m v,[\m v,A]]]]
.
%\end{align}
\end{multline}

\begin{commentx}
\begin{align}
\frac{\mathrm d^2}{\mathrm d^2t}\log(&\exp(A)\exp(t\m v))\Bigr|_{t=0}=\plabel{schur11old}
\\\notag=&
\frac{\reC(D_A)+\reD(D_A)}{4}\,[\m v,[\m v,A]]
+\frac{\reC(D_A)}{4}\,[A,[\m v,[\m v,A]]]
\\\notag&+\frac{2\reP(D_A)+3\reW(D_A)}{16} \, [A,[A,[\m v,[\m v,A]]]]
\\\notag=&
\frac{1-\reC(D_A)-\reD(D_A)}{2}\,[\m v,[\m v,A]]
+\frac{\reC(D_A)}{4}\,[A,[\m v,[\m v,A]]]
\\\notag&+\frac{2\reP(D_A)+3\reW(D_A)}{16} \,[ [A,[A,\m v]],[A,\m v] ]
\\\notag=&
\frac16\,[\m v,[\m v,A]]
+\frac{\reC(D_A)}{4}\,[A,[\m v,[\m v,A]]]
\\\notag&+\frac{2\reP(D_A)+3\reW(D_A)}{48}\, \Bigl(2 [A,[A,[\m v,[\m v,A]]]] +[ [A,[A,\m v]],[A,\m v] ] \Bigr).
\end{align}
and
\begin{align}
\frac{\mathrm d^2}{\mathrm d^2t}\log(&\exp(t\m v)\exp(A))\Bigr|_{t=0}=\plabel{schur22old}
\\\notag=&
\frac{\reC(D_A)+\reD(D_A)}{4}\,[\m v,[\m v,A]]
-\frac{\reC(D_A)}{4}\,[A,[\m v,[\m v,A]]]
\\\notag&+\frac{2\reP(D_A)+3\reW(D_A)}{16} \, [A,[A,[\m v,[\m v,A]]]]
\\\notag=&
\frac{1-\reC(D_A)-\reD(D_A)}{2}\,[\m v,[\m v,A]]
-\frac{\reC(D_A)}{4}\,[A,[\m v,[\m v,A]]]
\\\notag&+\frac{2\reP(D_A)+3\reW(D_A)}{16} \,[ [A,[A,\m v]],[A,\m v] ]
\\\notag=&
\frac16\,[\m v,[\m v,A]]
-\frac{\reC(D_A)}{4}\,[A,[\m v,[\m v,A]]]
\\\notag&+\frac{2\reP(D_A)+3\reW(D_A)}{48}\, \Bigl(2 [A,[A,[\m v,[\m v,A]]]] +[ [A,[A,\m v]],[A,\m v] ] \Bigr).
\end{align}
\end{commentx}
hold.
\begin{proof}
Direct computation.
\end{proof}
\end{lemma}
(Although generating functions for the higher terms are known, cf. Goldberg \cite{G}, computing with them is messier.)

\begin{commentx}
\begin{remark}
\plabel{rem:schur22var}
Note that the presentations in Lemma \ref{lem:schur22} depend greatly on the (combinations of) Lie monomials used.
Indeed, the functions $\reD(x),\reW(x),\reP(x)$ are merely convenient polynomials of $x, \reC(x), \reC'(x) ,\reC''(x) $.
As, for example, the identity
\[[ [A,[A,\m v]],[A,\m v] ] - [A,[A,[\m v,[\m v,A]]]] =[\m v,[A,[A,[\m v,A]]]]=-4D_A[\m v,[\m v,A]]\]
holds, rules \eqref{eq:simp1} and \eqref{eq:simp2} enter.
\qedremark
\end{remark}
\snewpage

Some relatively compact presentations for the terms of third order are:
\begin{lemma}
\plabel{rem:schur22third}
\begin{align*}
\left.\frac{\mathrm d^3}{\mathrm d^3t}\log(\exp(A)\exp(t\m v))\right|_{t=0}=\,&
\frac{2\reC(D_A)-\reP(D_A)-3 \reG(D_A)-3\reW(D_A)}{16}
\\&\quad\cdot\Bigl( [[\m v,[A,\m v]],[A,\m v]]  +2[[[A,[A,\m v]],\m v],\m v]  \Bigr)
\\&+\frac{2\reP(D_A)+3\reW(D_A)}{32}
\\&\quad\cdot \Bigl[A, [[\m v,[A,\m v]],[A,\m v]]  +2[[[A,[A,\m v]],\m v],\m v]  \Bigr]
\\&+\frac{6\reW(D_A)+  6\reL(D_A)+  \reC(D_A)\reW(D_A) + 3\reX(D_A)}{64}
\\&\quad\cdot[[\m v, [A,[A,\m v]] ] , [A,[A,\m v]]  ]
\end{align*}
and
\begin{align*}
\left.\frac{\mathrm d^3}{\mathrm d^3t}\log(\exp(A)\exp(t\m v))\right|_{t=0}=\,&
\frac{2\reC(D_A)-\reP(D_A)-3 \reG(D_A)-3\reW(D_A)}{16}
\\&\quad\cdot\Bigl( [[\m v,[A,\m v]],[A,\m v]]  +2[[[A,[A,\m v]],\m v],\m v]  \Bigr)
\\&-\frac{2\reP(D_A)+3\reW(D_A)}{32}
\\&\quad\cdot \Bigl[A, [[\m v,[A,\m v]],[A,\m v]]  +2[[[A,[A,\m v]],\m v],\m v]  \Bigr]
\\&+\frac{6\reW(D_A)+  6\reL(D_A)+  \reC(D_A)\reW(D_A) + 3\reX(D_A)}{64}
\\&\quad\cdot[[\m v, [A,[A,\m v]] ] , [A,[A,\m v]]  ]
\end{align*}
\begin{proof}
Direct computation.
\end{proof}
\end{lemma}
\end{commentx}
Due to the nature of the recursion, one can see that
$\left.\frac{\mathrm d^n}{\mathrm d^nt}\log(\exp(A)\exp(t\m v))\right|_{t=0}$
and
$\left.\frac{\mathrm d^n}{\mathrm d^nt}\log(\exp(t\m v)\exp(A))\right|_{t=0}$,
for $n\geq2$, will be linear combinations of $[A,\m v]$, $[A,[A,\m v]]$,  $[\m v,[\m v,A]]$
with functions of $D_A$, $T_{A,\m v}$, $D_{\m v}$.

Taking this formally, cf. Lemma \ref{lem:prefivecomm}, all this implies
the qualitative statement of the Baker--Campbell--Hausdorff formula,
i. e. that the terms in the BCH expansion, grouped by (bi)degree, are commutator polynomials,
for $2\times2$ matrices.
This was achieved here as the special case of the argument using Schur's formulae.
(Which is a typical argument, see Bonfiglioli, Fulci \cite{BF} for a review of the topic, or \cite{L0} for some additional viewpoints.)
The difference to the general case is that the terms in the expansion (grouped by degree in $\m v$) can be kept in finite
form (meaning the kind of polynomials we have described).
\snewpage
\scleardoublepage\section{The BCH formula for $2\times2$ matrices }
\plabel{sec:BCH}
Let us take another viewpoint on the BCH expansion now.
Recall, in the case of $2\times2$ matrices, as it was demonstrated by Lemma \ref{lem:prefivecomm} and Lemma \ref{lem:indep}, commutator polynomials can be
represented quite simply, and allowing a formal calculations.
It is also a natural question whether those commutator terms of the BCH expansion can be obtained in other efficient manners.
The following formal calculations will address these issues.

Let us use the notation
\begin{equation}\hat A=A-\frac{\tr A}2\Id_2\plabel{eq:Ahatdef}\end{equation}
and
\begin{equation}\hat{\m v}=\m v-\frac{\tr \m v}2\Id_2.\plabel{eq:vhatdef}\end{equation}
Then it is easy to that
\begin{equation}
 \BCH(A,\m v)=\frac{\tr A}2\Id_2+\frac{\tr \m v}2\Id_2+\BCH(\hat A,\hat {\m v}).
\plabel{eq:gagg0}
\end{equation}

Then Lemma \ref{lem:prefive} can be applied to $\BCH(\hat A,\hat {\m v})$.
Due to the traceless of $\hat A$, $\hat {\m v}$, only the
scalar terms $\det \hat A=D_A$, $\det \hat {\m v}=D_{\m v}$, and $\tr(\hat A\hat {\m v})=2T_{A,\m v}$
remain. In the vector terms $[\hat A,\hat {\m v}]=[A,\m v]$.
Thus we arrive to
\begin{align}
\BCH(\hat A,\hat {\m v}) =&\hat g_0(D_A,T_{A,\m v} ,D_{\m v})\, \Id_2\plabel{eq:gaom00}\\\notag
&+\hat g_1(D_A,T_{A,\m v} ,D_{\m v})\,  \hat A\\\notag
&+\hat g_2(D_A,T_{A,\m v} ,D_{\m v})\,\hat{\m v}\\\notag
&+\hat g_3(D_A,T_{A,\m v} ,D_{\m v})\,[A,\m v].
\end{align}
(This is valid for any formal expression of $\hat A, \hat{\m v}$.)

From the (formal) determinant one can deduce that $g_0\equiv 0$.
Considering \eqref{eq:gagg0}, we obtain
\begin{align}
\BCH(\hat A,\hat {\m v}) =&A+\m v+\plabel{eq:gaom01}\\\notag
&+g_1(D_A,T_{A,\m v} ,D_{\m v})\, \hat A  \\\notag
&+g_2(D_A,T_{A,\m v} ,D_{\m v})\,\hat{\m v} \\\notag
&+g_3(D_A,T_{A,\m v} ,D_{\m v})\,[A,\m v] .
\end{align}
where
\[g_1(D_A,T_{A,\m v} ,D_{\m v})=\hat g_1(D_A,T_{A,\m v},D_{\m v})-1,\]
\[g_2(D_A,T_{A,\m v} ,D_{\m v})=\hat g_2(D_A,T_{A,\m v},D_{\m v})-1,\]
\[g_3(D_A,T_{A,\m v} ,D_{\m v})=\hat g_3(D_A,T_{A,\m v},D_{\m v}).\]
\snewpage

This is very close to a commutator expansion. Using \eqref{eq:biag1} and \eqref{eq:biag2}, we obtain
\begin{align}
\BCH(A,\m v)=A+\m v
&+f_1(D_A,T_{A,\m v} ,D_{\m v})\,  [A,\m v]\plabel{eq:gaom2}\\\notag
&+f_2(D_A,T_{A,\m v} ,D_{\m v})\,[A,[A,\m v]]\\\notag
&+f_3(D_A,T_{A,\m v} ,D_{\m v})\,[\m v,[\m v,A]].
\end{align}
where
\[f_1=g_3,\]
and
\begin{equation}
f_2=\frac{T_{A,\m v} g_1-D_{\m v} g_2}{4(D_AD_{\m v} -T_{A,\m v}^2)},
\plabel{eq:rot1}
\end{equation}
\begin{equation}
f_3=\frac{-D_{A} g_1+T_{A,\m v} g_2}{4(D_AD_{\m v} -T_{A,\m v}^2)}.
\plabel{eq:rot2}
\end{equation}

From this, we can see that the qualitative version of the BCH theorem is
a very weak but not entirely trivial statement:
It is equivalent to divisibility in \eqref{eq:rot1} and \eqref{eq:rot2}.
\snewpage
\begin{commenty}
\begin{remark}
\plabel{rem:rep44}
$\BCH(A,\m v)$ or $\BCH(\hat A,\hat {\m v})$ can be expanded in the formal base
\[\Bigl(
\Id_2; A; \m v; [A,\m v]
\Bigr)\]
by multiplication on the left.
But, more conveniently, they can also be expanded in the formal base
\[\widetilde{\mathfrak A}=\Bigl(
\Id_2; \hat A; \hat{\m v}; [A,\m v]
\Bigr)\]
by multiplication on the left.

In that case
\[L_{\widetilde {\mathfrak A}}\hat A=
\begin{bmatrix}
&-D_A&T_{A,\m v}&\\
1&&&-2T_{A,\m v}\\
&&&-2D_A\\
&&\frac12&
\end{bmatrix}
\]
and
\[L_{\widetilde {\mathfrak A}}\hat{\m v}=
\begin{bmatrix}
&T_{A,\m v}&-D_{\m v} &\\
&&&2D_{\m v}\\
1&&&2T_{A,\m v}\\
&-\frac12&&
\end{bmatrix}
.\]
Consequently,
\begin{equation}
L_{\widetilde {\mathfrak A}} \BCH(\hat A,\hat {\m v})=
\BCH\left(
\begin{bmatrix}
&-D_A&T_{A,\m v}&\\
1&&&-2T_{A,\m v}\\
&&&-2D_A\\
&&-\frac12&
\end{bmatrix}
,
\begin{bmatrix}
&T_{A,\m v}&-D_{\m v} &\\
&&&2D_{\m v}\\
1&&&2T_{A,\m v}\\
&-\frac12&&
\end{bmatrix}
\right).
\plabel{eq:lagg0}
\end{equation}
Then the first column of $L_{\widetilde {\mathfrak A}} \BCH(\hat A,\hat {\m v})$
informs about the coefficients of \eqref{eq:gaom00}.

The RHS of
\eqref{eq:lagg0} can be computed effectively by Goldberg's theorem (non-commutator version), cf. Goldberg, \cite{G}
(or by any other suitable method).
At first, the homogeneity relations of the components are not clear, but,
purely formally let us conjugate by the diagonal matrix
\[F=\begin{bmatrix}1&&&\\&\sqrt{D_A}&&\\&&\sqrt{D_{\m v}}&\\&&&2\sqrt{D_A}\sqrt{D_{\m v}}\end{bmatrix}.\]
Then we obtain
\[F(L_{\widetilde {\mathfrak A}}\hat A)F^{-1}=
\begin{bmatrix}
&-\sqrt{D_A}&\frac{T_{A,\m v}}{\sqrt{D_{\m v}}}&\\
\sqrt{D_A}&&&-\frac{T_{A,\m v}}{\sqrt{D_{\m v}}}\\
&&&-\sqrt{D_A}\\
&&\sqrt{D_A}&
\end{bmatrix}
\]
and
\[F(L_{\widetilde {\mathfrak A}}\hat{\m v})F^{-1}=
\begin{bmatrix}
&\frac{T_{A,\m v}}{\sqrt{D_{\m v}}}&-\sqrt{D_{\m v}} &\\
&&&\sqrt{D_{\m v}}\\
\sqrt{D_{\m v}}&&&\frac{T_{A,\m v}}{\sqrt{D_{\m v}}}\\
&-\sqrt{D_{\m v}}&&
\end{bmatrix}
.\]
In that form the general homogeneity of matrices is transparent.

If, in \eqref{eq:lagg0}, the BCH term is expanded in or up to bidegree $(n,m)$,
then the computations in (\ref{eq:gaom00}--\ref{eq:rot2})
are correct in or up to bidegree $(n,m)$.
\qedremark
\end{remark}
%\snewpage
\begin{remark}
\plabel{rem:rep33}
(a) Another method to proceed is as follows.
Firstly, we will describe  $\ad \BCH(A,\mathbf v)$, or, more precisely,
its representation $\ad_{\mathfrak A} \BCH(A,\mathbf v)$ in the formal basis
\[\mathfrak A=\Bigr( \hat A; \hat{\m v};[A,\m v]\Bigl)\]
in order to represent $\ad \BCH(A,\m v)$.
In that case
\[\ad_{\mathfrak A} A=
\begin{bmatrix}
&&-4T_{A,\m v}\\
&&-4D_A\\
0&1&
\end{bmatrix}
\]
and
\[\ad_{\mathfrak A} \m v=
\begin{bmatrix}
&&4D_{\m v}\\
&&4T_{A,\m v}\\
-1&0&
\end{bmatrix}
.
\]

Then
\begin{equation}
\ad_{\mathfrak A} \BCH(A,\mathbf v)=\BCH\left(
\begin{bmatrix}
&&-4T_{A,\m v}\\
&&-4D_A\\
0&1&
\end{bmatrix}
,
\begin{bmatrix}
&&4D_{\m v}\\
&&4T_{A,\m v}\\
-1&0&
\end{bmatrix}
\right).
\plabel{eq:gaoma1}
\end{equation}
This can be computed effectively, again, in arbitrary order.

(In order to see the homogeneity relations, purely formally, let us conjugate by the diagonal matrix
$F=\begin{bmatrix} \sqrt{D_A}&&\\&\sqrt{D_{\m v}}&\\&& 2\sqrt{D_A}\sqrt{D_{\m v}}\end{bmatrix}$.
Then we obtain
\[F(\ad_{\mathfrak A} A)F^{-1}=
\begin{bmatrix}
&&-2\frac{T_{A,\m v}}{\sqrt{D_{\m v}}}\\
&&-2\sqrt{D_A}\\
0&2\sqrt{D_A}&
\end{bmatrix}
\]
and
\[F(\ad_{\mathfrak A} \m v)F^{-1}=
\begin{bmatrix}
&&2\sqrt{D_{\m v}}\\
&&2\frac{T_{A,\m v}}{\sqrt{D_{A}}}\\
-2\sqrt{D_{\m v}}&0&
\end{bmatrix}
.
\]
In that form the homogeneity relations of the matrices are transparent.)

Our objective, is to find $f_1, f_2, f_3$ in \eqref{eq:gaom2}.
Let us consider
\[\mathfrak C=\Bigl(
-[A,\m v];\m v; A\Bigr)
.\]
Then it is easy to check that
\[ \tr\Bigl((\ad_{\mathfrak A} \mathfrak B_i)(\ad_{\mathfrak A} \mathfrak C_i)\Bigr)=\delta_{i,j}\cdot 32(D_A D_{\m v}-T_{A,\m v}^2).\]
Therefore, formally,
\[ f_i(D_A,T_{A,\m v} ,D_{\m v}))= \frac{
\tr\Bigl(\Bigl(\ad_{\mathfrak A} \BCH(A,B)- \ad_{\mathfrak A} A - \ad_{\mathfrak B} B\Bigr)(\ad_{\mathfrak A} \mathfrak C_i)\Bigr)
}{32(D_A D_{\m v}-T_{A,\m v}^2)}.
 \]
Again, the qualitative BCH theorem appears as divisibility.

\snewpage
(b) Alternatively, $\ad \BCH(A,\m v)$ can also be expanded in the formal base
\[\mathfrak B=\Bigl(
[A,\m v];[A,[A,\m v]];[\m v,[\m v,A]]
\Bigr).\]
It is easy to see that
\[\ad_{\mathfrak B} A=
\begin{bmatrix}
&-4D_A&-4T_{A,\m v}\\
1&&\\
0&&
\end{bmatrix}
\]
and
\[\ad_{\mathfrak B} \m v=
\begin{bmatrix}
& 4T_{A,\m v}&4D_{\m v}\\
0&&\\
-1&&
\end{bmatrix}
.
\]
Then
\begin{equation}
\ad_{\mathfrak B} \BCH(A,\mathbf v)=\BCH\left(\begin{bmatrix}
&-4D_A&-4T_{A,\m v}\\
1&&\\
0&&
\end{bmatrix}
,
\begin{bmatrix}
& 4T_{A,\m v}&4D_{\m v}\\
0&&\\
-1&&
\end{bmatrix}\right).
\plabel{eq:gaom1}
\end{equation}
This can be computed effectively in any order.

(In order to see the homogeneity relations, purely formally, let us conjugate by the diagonal matrix
$F=\begin{bmatrix}1&&\\&2\sqrt{D_A}&\\&&2\sqrt{D_{\m v}}\end{bmatrix}$.
Then we obtain
\[F(\ad_{\mathfrak B} A)F^{-1}=
\begin{bmatrix}
&-2\sqrt{D_A}&-2\frac{T_{A,\m v}}{\sqrt{D_{\m v}}}\\
2\sqrt{D_A}&&\\
0&&
\end{bmatrix}
\]
and
\[F(\ad_{\mathfrak B} \m v)F^{-1}=
\begin{bmatrix}
& 2\frac{T_{A,\m v}}{\sqrt{D_{A}}}&2\sqrt{D_{\m v}}\\
0&&\\
-2\sqrt{D_{\m v}}&&
\end{bmatrix}
.
\]
In that form the homogeneity relations of the matrices are transparent.)

Otherwise, the argument is the same, just `$\ad_{\mathfrak A}$' should be replaced by `$\ad_{\mathfrak B}$'.
(Using \eqref{eq:biag1} and \eqref{eq:biag2}, one can rewrite the higher terms of  \eqref{eq:gaom2}
from $\mathfrak B$ to $\mathfrak A$, or the other way around.)
\qedremark
\end{remark}
\end{commenty}
\snewpage
Therefore, in the setting of $2\times2$ matrices, an alternative approach to the BCH theorem
is simply to compute the coefficients $f_2$ and $f_3$.
Due to the simplicity of $2\times2$ matrices, this, indeed, can be done by direct computation:

\begin{lemma}
\plabel{lem:flatBCH2}

Consider the traceless $2\times2$ matrices $\hat A$, $\hat {\m v}$ as before.

(a) Then
\begin{multline*}
\det\left( (1-t)\Id_2+t \exp (\hat A)\exp(\hat{\m v})  \right)=\\=
\underbrace{1+2 t(1-t)\Bigr(\Sin(D_A)\Sin(D_{\m v})T_{A,\m v}+\Cos(D_A)\Cos(D_{\m v})-1\Bigr)
}_{\delta( D_A,T_{A,\m v},D_{\m v},t ):=}
.
\end{multline*}

(b)
Let
\[C(\hat A,\m v):=  \Cos(D_{\m v})\Sin(D_A)\hat A + \Cos(D_A)\Sin(D_{\m v})\hat{\m v}
+\frac12\Sin(D_A)\Sin(D_{\m v})\,[A,\m v]\]

Then, as long as $\exp (\hat A)\exp(\hat{\m v})$ is log-able,
\begin{align}
\BCH(\hat A,\hat{\m v})
&=\AC\bigl( \Cos(D_A)\Cos(D_{\m v})+\Sin(D_A)\Sin(D_{\m v})T_{A,\m v} \bigr)\cdot C(\hat A,\m v)
\plabel{eq:BCH2raw}
\\\notag
&=\left( \int_{t=0}^1\frac{\mathrm dt}{ \delta( D_A,T_{A,\m v},D_{\m v},t )}\right)\cdot C(\hat A,\m v).
\end{align}
\begin{proof}
This follows from our formulas concerning the exponential and the logarithm.
\begin{commentx}

(\ref{eq:BCH2raw}rhs/2)
can be computed directly by \eqref{eq:logdef}. Here
\[\frac{\exp(\hat A)\exp(\hat{\m v})-1}{(1-t)+t\exp(\hat A)\exp(\hat{\m v})}=\frac{\frac12 \frac{\mathrm d}{\mathrm dt}\delta( D_A,T_{A,\m v},D_{\m v},t )
\Id_2+ C(\hat A,\m v)}{\delta( D_A,T_{A,\m v},D_{\m v},t )};\]
from which the result is immediate.
\end{commentx}
\end{proof}
\end{lemma}
We could have formulated the previous statement in the general case,
but it would have been more complicated due to the trace terms.
Taking formally, we note that the determinant term $\delta( D_A,T_{A,\m v},D_{\m v} )$ is a formal perturbation of $1$.
(Substituting $D_A=0$, $T_{A,\m v}=0$, $D_{\m v}=0$ gives $1$).
\snewpage
\begin{theorem}
\plabel{eq:BCHth2}
Let $A,\m v, \hat A, \hat {\m v}$ as before. Let
\begin{multline*}
\eta( D_A,D_{\m v},t )=\\=
1-4t(1-t)\Bigr(1-\bigl(t\Cos(D_A)+(1-t)\Cos(D_{\m v})\bigr)\bigl(t\Cos(D_{\m v})+(1-t)\Cos(D_A)\bigr)\Bigr);
\end{multline*}
\[\varkappa_2( D_A,D_{\m v},t )=\Cos(D_{\m v})+   2t(1-t)(\Cos(D_A)-\Cos(D_{\m v});\]
\begin{equation*}
\xi_2( D_A,D_{\m v},t )=\\=
t(1-t)\varkappa_2( D_A,D_{\m v},t ) \cdot \frac12\Sin(D_A)^2\Sin(D_{\m v});
\end{equation*}
\[\varkappa_3( D_A,D_{\m v},t )= \Cos(D_{A})+   2t(1-t)(\Cos(D_{\m v})-\Cos(D_{A});\]
\begin{equation*}
\xi_3( D_A,D_{\m v},t )=\\=
t(1-t)\varkappa_3( D_A,D_{\m v},t ) \cdot  \frac12\Sin(D_{\m v})^2\Sin(D_{A}).
\end{equation*}

Then
\begin{align}
\BCH(A,\m v)=&A+\m v\plabel{eq:ewa}\\\notag
&+\left( \int_{t=0}^1\frac{\mathrm dt}{ \delta( D_A,T_{A,\m v},D_{\m v},t )}\right)\cdot\frac12\Sin(D_A)\Sin(D_{\m v})\,[A,\m v]
\\\notag
&+\left( \int_{t=0}^1\frac{\xi_2( D_A,D_{\m v},t )}{\eta( D_A,D_{\m v},t )}\cdot
\frac{\mathrm dt}{ \delta( D_A,T_{A,\m v},D_{\m v},t )}\right)\cdot[A,[A,\m v]]
\\\notag
&+\left( \int_{t=0}^1\frac{\xi_3( D_A,D_{\m v},t )}{\eta( D_A,D_{\m v},t )}\cdot
\frac{\mathrm dt}{ \delta( D_A,T_{A,\m v},D_{\m v},t )}\right)\cdot[\m v,[\m v ,A]]
\end{align}
is valid as long as  $A,\m v\sim 0$.

(Note that formally,
$\delta( D_A,T_{A,\m v},D_{\m v},t )$,
$\eta( D_A,D_{\m v},t )$,
$\xi_2( D_A,D_{\m v},t )$,
$\xi_2( D_A,D_{\m v},t )$,
$\Cos(D_A)$, $\Cos(D_{\m v}$, $\Sin(D_A)$, $\Sin(D_{\m v}$,
are all perturbations of $1$ for  $ D_A, T_{A,\m v}, D_{\m v}\sim0$.)
\begin{proof}
As a preparation, let us set
\[\rho(u,x,y)={\frac {\arctan \left( u\tan \left( \frac{\sqrt {x}+\sqrt {y}}2
 \right)  \right) +\arctan \left( u\tan \left(  \frac{\sqrt {x}-\sqrt {y}}2  \right)  \right) }{\sqrt {x}}}
.\]
In this form, this is a formal series in $\sqrt x$ and $\sqrt y$, with coefficients which are polynomials of $u$.
By consideration of parities, we see that it gives an analytic power series convergent for $x,y\sim0$ uniformly in
$u\in[-1,1]$.
(Using the addition rule $\arctan(z)+\arctan (w)=\arctan\left(\frac{w+z}{1-wz}\right)$ for $z,w\sim0$,
and utilizing the functions $\Cos, \Sin,\AC$,
\begin{commentx}
or $\ATT$,
\end{commentx}
 we can express it as a composite analytic function
for $x,y\sim0$ uniformly in $u\in[-1,1]$. In fact, it yields
\begin{multline*}
\rho(u,x,y)=\\=
\frac{2u\Sin(x)\AC\left(\dfrac{\Cos(x)+\Cos(y) +u^2\Cos(x)-u^2\Cos(x)  }{
\sqrt{(\Cos(x)+\Cos(y) +u^2\Cos(x)-u^2\Cos(x))^2  +4u^2(1-\Cos(x)^2)}}\right)}{\sqrt{(\Cos(x)+\Cos(y) +u^2\Cos(x)-u^2\Cos(x))^2  +4u^2(1-\Cos(x)^2)}}.
\end{multline*}
\begin{commentx}
\[\rho(u,x,y)=\frac{2u\Sin(x)\ATT\left(\dfrac{4u^2(1-\Cos(x)^2)}{
(\Cos(x)+\Cos(y) +u^2\Cos(x)-u^2\Cos(x))^2}\right)}{\Cos(x)+\Cos(y) +u^2\Cos(x)-u^2\Cos(x)}.\]
\end{commentx}
But this is not particularly enlightening.)
Note, for $x,y\sim0$,
\[\rho(1,x,y)=1\qand\rho(-1,x,y)=-1.\]

Now, we the start the proof proper.
By Lemma \ref{lem:flatBCH2} and the discussion in this section, we know that
\begin{align}
\BCH&(A,\m v)=A+\m v\plabel{eq:owa}\\\notag
&+\left( \int_{t=0}^1\frac{\mathrm dt}{ \delta( D_A,T_{A,\m v},D_{\m v},t )}\right)\cdot\frac12\Sin(D_A)\Sin(D_{\m v})\,[A,\m v]
\\\notag
&+\left( \int_{t=0}^1
\frac{\Cos(D_{\m v})\Sin(D_A)}{ \delta( D_A,T_{A,\m v},D_{\m v},t )}   \mathrm dt\right)\cdot\hat A
-\underbrace{\left( \int_{t=0}^1 \frac{\mathrm d}{\mathrm dt}\left(\frac{\rho(2t-1,D_A,D_{\m v})}{2}\right)
  \mathrm dt\right)}_{=1}\cdot\hat A
\\\notag
&+\left( \int_{t=0}^1
\frac{\Cos(D_A)\Sin(D_{\m v})}{ \delta( D_A,T_{A,\m v},D_{\m v},t )}\mathrm dt\right)\cdot\hat{\m v}
-\underbrace{\left( \int_{t=0}^1 \frac{\mathrm d}{\mathrm dt}\left(\frac{\rho(2t-1,D_{\m v},D_A)}{2}\right)
  \mathrm dt\right)}_{=1}\cdot\hat {\m v}.
\end{align}
Thus, we have to prove
\begin{equation}
\text{RHS(\ref{eq:ewa}/3)} + \text{RHS(\ref{eq:ewa}/4)}=\text{RHS(\ref{eq:owa}/3)} + \text{RHS(\ref{eq:owa}/4)}.
\plabel{eq:obj}
\end{equation}
As the integrands can be expanded as power series of $D_A,T_{A,\m v},D_{\m v}\sim 0$  (uniformly in $t\in[0,1]$),
it is even sufficient to prove this as for formal power series in $D_A,T_{A,\m v},D_{\m v}$.

Now, it is easy to check that
\[\frac{\mathrm d}{\mathrm dt}\left(\frac{\rho(2t-1,D_A,D_{\m v})}{2}\right)
=\frac{\Sin(D_A)\varkappa_2(D_A,D_{\m v},t)}{\eta(D_A,D_{\m v},t)}\]
and
\[\frac{\mathrm d}{\mathrm dt}\left(\frac{\rho(2t-1,D_{\m v},D_A)}{2}\right)
=\frac{\Sin(D_{\m v})\varkappa_3(D_A,D_{\m v},t)}{\eta(D_A,D_{\m v},t)}.\]
Thus, using  \eqref{eq:biag1} and \eqref{eq:biag2}, one finds that
\[
\left(\frac{\Cos(D_{\m v})\Sin(D_A)}{ \delta( D_A,T_{A,\m v},D_{\m v},t )}-
\frac{\mathrm d}{\mathrm dt}\left(\frac{\rho(2t-1,D_A,D_{\m v})}{2}\right) \right)   \cdot\hat A
\]\[
+\left(\frac{\Cos(D_A)\Sin(D_{\m v})}{ \delta( D_A,T_{A,\m v},D_{\m v},t )}-
\frac{\mathrm d}{\mathrm dt}\left(\frac{\rho(2t-1,D_{\m v},D_A)}{2}\right)\right) \cdot\hat{\m v} \]
\[=  \frac{\xi_2( D_A,D_{\m v},t )}{\eta( D_A,D_{\m v},t )}\cdot
\frac{1 }{ \delta( D_A,T_{A,\m v},D_{\m v},t )}[A,[A,\m v]]
\]\[+
\frac{\xi_3( D_A,D_{\m v},t )}{\eta( D_A,D_{\m v},t )}\cdot
\frac{1}{ \delta( D_A,T_{A,\m v},D_{\m v},t )}[\m v,[\m v ,A]] .\]
Integrated, it yields \eqref{eq:obj}.
\end{proof}
\end{theorem}
Taken formally, the statement and the proof of the previous theorem gives
an alternative demonstration of
the qualitative (commutator expansion) statement of the BCH theorem
(for $2\times 2$ matrices).
This is, however, more for the sake of curiosity;
in practice, it is simpler to use \eqref{eq:BCH2raw}
in order to expand as in \eqref{eq:gaom01}, and obtain \eqref{eq:gaom2} with \eqref{eq:rot1}, \eqref{eq:rot2}.

\snewpage
\scleardoublepage\section{Magnus minimality of mBCH expansions in the $2\times2$ case}
\plabel{sec:expmag}
Recall (from Part I) that in a Banach algebra $\mathfrak A$, for an element $A$ we define its Magnus exponent
as
\[\mathcal M_{\mathfrak A}(A)=\inf\{\smallint|\phi|\,:\, \Lexp(\phi)=A\}\]
(giving $+\infty$ for an empty set).
As ordered measures can be replaced by piecewise constant measures with arbitrarily small increase in the cumulative norm,
the infimum can be taken for mBCH measures. We will say that $\phi$ is Magnus-minimal if
$\mathcal M_{\mathfrak A}( \Lexp(\phi))=\smallint|\psi|$.
Now, the natural generalization of Theorem \ref{th:better3}.(b) would be the following:

(X) ``In the setting of finite dimensional Hilbert spaces, if $\psi$ is a reduced BCH measure, then
$\mathcal M_{\mathcal B(\mathfrak H)}(\Lexp(\psi))<\smallint\|\psi\|_2$ (i. e. $\psi$ is not Magnus minimal).''

The objective of this section is to check (X) for $2\times2$ matrices by direct computation.
(In Section \ref{sec:MagnusGL2}, we will see more informative approach in the real case.)

\begin{theorem}
\plabel{th:NonNormalMin}
(a) Suppose that $A$ is real $2\times2$ matrix, which is not normal.
Then $A\m 1$ is not Magnus-minimal, i. e.
\[\mathcal M_{2\times 2\,\,\real}(A\m 1)<\|A\|.\]
(b) Similar  statement holds in the complex case.
\begin{proof}
We can assume that $A$ is sufficiently small, in particular,
$\spec(A/2)\subset \{z\in\mathbb C\,:\,|\Ima z|<\pi\}$.
For $t\sim0$, $t\in\mathbb R$, let
\[A_1(t)=\log\left( \exp(A/2)  \exp(t[A,A^*])  \right);\]
\[A_2(t)=\log\left( \exp(-t[A,A^*]) \exp(A/2)  \right).\]
Then,
\[\exp A=\exp(A_1(t))\exp(A_2(t));\]
this yields a (m)BCH expansion for $\exp A$.
It is sufficient to prove that
\begin{equation}
\left.\frac{\mathrm d\|A_1(t)\|_2}{\mathrm dt}\right|_{t=0+}+
\left.\frac{\mathrm d\|A_2(t)\|_2}{\mathrm dt}\right|_{t=0+}<0
\plabel{eq:dec1}
\end{equation}
in order to demonstrate non-Magnus minimality.
Using \eqref{eq:riven} and \eqref{eq:tricommut}, we find
\[\underbrace{\mathrm D_{[A,A^*]\,\,\mathrm{at}\,\, M=0}\left(\log\left(\exp(A/2)\exp(M)\right)\right)}_{\m v_1:=}=
\Cot\left(D_{A/2}\right)[A,A^*]+\frac14[A,[A,A^*]];
\]
and, using Lemma \ref{lem:smooth2}(d),
\[\left.\frac{\mathrm d\|A_1(t)\|_2}{\mathrm dt}\right|_{t=0+}=
\mathrm D_{\m v_1 \,\,\mathrm{at}\,\, M=A/2 }\,\|M\|_2=-\frac14\frac{\|A_2\|(-D_{[A,A^*]})}{\sqrt{-D_{A*A}}}<0.\]
(Note that $A/2$ being non-normal implies smoothness for the norm.)
Similarly,
\[\underbrace{\mathrm D_{-[A,A^*]\,\,\mathrm{at}\,\, M=0}\left(\log\left(\exp(M)\exp(A/2)\right)\right)}_{\m v_2:=}=
-\Cot\left(D_{A/2}\right)[A,A^*]+\frac14[A,[A,A^*]];
\]
and,
\[\left.\frac{\mathrm d\|A_2(t)\|_2}{\mathrm dt}\right|_{t=0+}=
\mathrm D_{\m v_2 \,\,\mathrm{at}\,\, M=A/2 }\,\|M\|_2=-\frac14\frac{\|A_2\|(-D_{[A,A^*]})}{\sqrt{-D_{A*A}}}<0.\]
Then \eqref{eq:dec1} holds, and so does the statement.
(Remark: it is useful to follow through this computation for \eqref{eq:canonC} with $s_1,s_2>0$.)
\end{proof}
\end{theorem}
As a consequence, we see that any Magnus minimal mBCH presentation must contain only normal matrices.

\begin{theorem}
\plabel{th:ContractMin}
(a) Suppose that $A$ and $B$  are real $n\times n$ matrices such that $\|A+B\|_2< \|A\|_2+\|B\|_2$.
Then $A\m 1\boldsymbol.B\m 1$ is not Magnus-minimal, i. e.
\[\mathcal M_{n\times n\,\,\real}(\Rexp({A\m 1\boldsymbol.B\m 1)})<\|A\|_2+\|B\|_2.\]
(b) Similar  statement holds in the complex case.
\begin{proof} Appling the BCH formula, we see that
\[\lim_{t\searrow 0} \frac{\|\log(\Rexp(tA,tB))\|_2}{\| tA\|_2+\| tB\|_2}=\frac{\| A+ B\|_2}{\| A\|_2+\| B\|_2}.\]
This implies that for sufficiently small $t$,  $\Rexp({(1-t)A\m 1\boldsymbol.tC\boldsymbol.(1-t)B\m 1)}$,
with $tC=\log(\Rexp(tA,tB))$, will be an alternative presentation with smaller cumulative norm.
\end{proof}
\end{theorem}
\begin{defin}
(a)
We say that the complex normal matrices $A$ and $B$ are aligned if, up to (simultaneous) unitary conjugation,
they are of shape
\[A= r_1{\mathrm e}^{\mathrm i \eta}\begin{bmatrix}1&\\&t_1\end{bmatrix}
\qquad
\text{and}
\qquad
B=r_2{\mathrm e}^{\mathrm i \eta}\begin{bmatrix}1&\\&t_2\end{bmatrix}
,\]
where $r_1,r_2\in(0,+\infty)$, $\eta\in[0,2\pi)$, $t_1,t_2\in \Dbar(0,1)$.

(b) We say that the complex normal matrices  $A$ and $B$ are skew-aligned if, up to simultaneous unitary conjugation,
they are of shape
\[A=r_1\bem \mathrm e^{\mathrm i\eta}\cos t&-\mathrm e^{\mathrm i\theta}\sin t \\\mathrm e^{\mathrm i\eta}\sin t&\mathrm e^{\mathrm i\theta}\cos t\eem
\qquad
\text{and}
\qquad
B=r_2\bem \mathrm e^{\mathrm i\eta}\cos t&-\mathrm e^{\mathrm i\phi}\sin t \\\mathrm e^{\mathrm i\eta}\sin t&\mathrm e^{\mathrm i\phi}\cos t\eem
,\]
where $r_1,r_2\in(0,+\infty)$, $t\in(0,\pi)$ and $\eta,\theta,\phi\in[0,2\pi)$ but $\phi\neq\theta$.
(Note that in this case $A$ and $B$ are conform-unitary.
We also remark that $t\rightsquigarrow \pi-t$, $\eta\rightsquigarrow\eta+\pi$, $\phi\rightsquigarrow\phi+\pi$, $\theta\rightsquigarrow\theta+\pi$
is a symmetry by conjugation with $\bem1&\\&-1\eem$. )
\end{defin}

\begin{lemma}
\plabel{lem:Aligned}
Suppose that $A,B$ are nonzero complex $2\times2$ matrices such that $A$ and $B$ are normal.
Then
\[\|A+B\|_2=\|A\|_2+\|B\|_2\]
holds if and only if
$A$ and $B$ are aligned or skew-aligned.
The aligned and skew-aligned cases are mutually exclusive.
For example, aligned pairs commute and skew-aligned pairs do not commute.
\begin{proof}
If $A$ (or $B$) is not conform unitary, then its norm is realized at an eigenvector.
Due to the additive restriction it must be common eigenvector of $A$ and $B$.
It can be assumed that this common eigenvector is $\bem1\\0\eem$, then due to normality
(the orthogonality of eigenspaces), the matrices are aligned.
If $A$ and $B$ are conform-unitary, then it can be assume that the norm of the sum is taken
$\bem1\\0\eem$ again, and due to the restrictions we have a configuration which is an
skew-alignment but  $t=0$ and $\phi=\theta$ allowed.
The excess cases can be incorporated to the aligned case, what remains is the skew-aligned case.
The commutation statement is easy to check.
\end{proof}
\end{lemma}
\begin{lemma}
\plabel{lem:RealAlign}
Suppose that $A$ and $B$ are nonzero, normal, $2\times2$ real matrices.

(a) $A$ and $B$ are aligned if and only if up to simultaneous conjugation by orthogonal matrices
they are of shape
\[A= r_1{\mathrm e}^{\mathrm i \eta}\begin{bmatrix}1&\\&t_1\end{bmatrix}
\qquad
\text{and}
\qquad
B=r_2{\mathrm e}^{\mathrm i \eta}\begin{bmatrix}1&\\&t_2\end{bmatrix}\]
where $r_1,r_2\in(0,+\infty)$, $\eta\in\{0,\pi\}$, $t_1,t_2\in[-1,1]$
[hyperbolically aligned case], or
up to simultaneous conjugation by orthogonal matrices
they are of shape
\[A= r_1\begin{bmatrix}\cos t&-\sin t\\\sin t&\cos t\end{bmatrix}
\qquad
\text{and}
\qquad
B= r_2\begin{bmatrix}\cos t&-\sin t\\\sin t&\cos t\end{bmatrix}\]
with $r_1,r_2\in(0,+\infty)$, $t\in[0,2\pi)$
[elliptically aligned case].

(The hyperbolic and elliptic cases are not mutually exclusive, but the common case
involves only with scalar matrices.)

(b)  $A$ and $B$ are skew-aligned if and only if up to simultaneous conjugation by orthogonal matrices
they are of shape
\[A=r_1\bem \mathrm e^{\mathrm i\eta}\cos t&-\mathrm e^{\mathrm i\theta}\sin t \\\mathrm e^{\mathrm i\eta}\sin t&\mathrm e^{\mathrm i\theta}\cos t\eem
\qquad
\text{and}
\qquad
B=r_2\bem \mathrm e^{\mathrm i\eta}\cos t&-\mathrm e^{\mathrm i\phi}\sin t \\\mathrm e^{\mathrm i\eta}\sin t&\mathrm e^{\mathrm i\phi}\cos t\eem\]
where $r_1,r_2\in(0,+\infty)$, $t\in(0,\pi)$ and $\eta,\theta,\phi\in\{0,\pi\}$,  $\phi\neq\theta$.
\begin{proof}
This follows from the standard properties of normal (in this case: symmetric and conform-orthogonal) matrices.
\end{proof}
\end{lemma}

\snewpage
\begin{theorem}
\plabel{lem:EllipMin}
(a) Suppose that $A,B$ are real normal $2\times2$ matrices which are  skew-aligned.
Then
$A\m 1\boldsymbol.B\m 1$ is not Magnus-minimal, i. e.
\[\mathcal M_{n\times n\,\,\real}(\Rexp({A\m 1\boldsymbol.B\m 1)})<\|A\|_2+\|B\|_2.\]

(b) Similar  statement holds in the complex case.
\begin{proof}
Assume that $\|A\|_2,\|B\|_2<\pi$ and $C$ is an arbitrary matrix. Then
\[\tilde A(t)=\log(\exp(A)\exp(tC))\]
and
\[\tilde B(t)=\log(\exp(-Ct)\exp(B))\]
make sense and analytic for $t\sim 0$.
Furthermore, $(\exp A)(\exp B)= (\exp \tilde A(t))(\exp \tilde B(t))$.

If
\[\left.\frac{\mathrm d\|\tilde A(t)\|_2}{\mathrm dt}\right|_{t=0+}+\left.\frac{\mathrm d\|\tilde B(t)\|_2}{\mathrm dt}\right|_{t=0+}<0,\]
%(we take right derivatives at $t=0$),
then $A\m 1\boldsymbol. B\m 1$ cannot be Magnus minimal, as it can be replaced $\tilde A(t)\m 1\boldsymbol. \tilde B(t)\m 1$
which is of smaller cumulative norm for small $t>0$.
We are going to use this idea.

By Schur's formulae,
\[ \tilde A'(0)=\left.\frac{\mathrm d}{\mathrm dt}\log(\exp(A)\exp(tC))\right|_{t=0}=\beta(-\ad A)C,\]
\[ \tilde B'(0)=\left.\frac{\mathrm d}{\mathrm dt}\log(\exp(-tC)\exp(B))\right|_{t=0}=-\beta(\ad B)C.\]
(This holds as $\|A\|_2,\|B\|_2<\pi$ was assumed.)

In our situation the matrices $A$ and $B$ are conform-unitary matrices, thus
\[\left.\frac{\mathrm d\|\tilde A(t)\|_2}{\mathrm dt}\right|_{t=0+}=
\|\tilde A(0)\|_2\cdot S( \tilde A'(0)\tilde A(0)^{-1} );\]
cf. Lemma \ref{lem:smooth2}(c); similar formula holds with $\tilde B(t)$.

First, we treat the complex case.
Now, assume that
\[A=r\bem \mathrm e^{\mathrm i\eta}\cos t&-\mathrm e^{\mathrm i\theta}\sin t \\\mathrm e^{\mathrm i\eta}\sin t&\mathrm e^{\mathrm i\theta}\cos t\eem\equiv
r\bem\cos t&-\sin t\\\sin t&\cos t\eem\bem\mathrm e^{\mathrm i\eta}&\\&\mathrm e^{\mathrm i\theta}\eem\]
and
\[B=r\bem \mathrm e^{\mathrm i\eta}\cos t&-\mathrm e^{\mathrm i\phi}\sin t \\\mathrm e^{\mathrm i\eta}\sin t&\mathrm e^{\mathrm i\phi}\cos t\eem\equiv
r\bem\cos t&-\sin t\\\sin t&\cos t\eem\bem\mathrm e^{\mathrm i\eta}&\\&\mathrm e^{\mathrm i\phi}\eem\]
but $A\neq B$.
Note that, in this case, $\phi\not\equiv\theta\modu2\pi$ and
\[1-\cos(\phi-\theta)>0.\]
In our computation $\theta,\phi,\eta$ will be fixed, but $r$ can be taken ``sufficiently small''.

Let
\[C=r\frac1{\frac12(1-\cos(\phi-\theta))}
\bem\cos t&-\sin t\\\sin t&\cos t\eem\bem (\mathfrak a_1+\mathrm i\mathfrak a_2)r\mathrm e^{\mathrm i\eta} &\mathfrak br\\
\mathfrak br&\mathrm e^{\mathrm i\phi}-\mathrm e^{\mathrm i\theta}\eem
,\]
where
\[\mathfrak b=\frac12(2-\mathrm e^{\mathrm i(\theta+\eta)}-\mathrm e^{\mathrm i(\phi+\eta)}    )\sin t
,\]
and
\[\mathfrak a_1=\cos \left( \eta+\theta \right)-\cos \left( \eta+\phi \right)
,\]
and
\begin{align*}\mathfrak a_2=\frac1{\frac12(1-\cos(\phi-\theta))}\Bigl(
 &+\sin \left( \eta+\theta \right)  \left( \cos \left( \eta+
\phi \right) -1 \right)  \left( \cos \left( \eta+\phi \right) -2 \right)\\
&-\sin \left( \eta+\phi \right)  \left( \cos \left( \eta+\theta
 \right) -1 \right)  \left( \cos \left( \eta+\theta \right) -2 \right)\\
 &+\frac12\, \left( \cos \left( \eta+\theta \right) -\cos \left(
\eta+\phi \right)  \right) \sin \left( 2\,\eta+\theta+\phi \right)
\Bigr)
.
\end{align*}

In this case,
\begin{align*}
\left.\frac{\mathrm d\|\tilde A(t)\|_2}{\mathrm dt}\right|_{t=0+}&+\left.\frac{\mathrm d\|\tilde B(t)\|_2}{\mathrm dt}\right|_{t=0+}=\\
&  = \|A\|_2\cdot S(\overbrace{\beta(-\ad A)C \cdot A^{-1}}^{\m v_1:=})+
 \|B\|_2\cdot S(   \overbrace{-\beta(\ad B)C\cdot B^{-1}}^{\m v_2:=})\\
&=-\frac{r^4\sin^2 t}{6}\frac1{\frac12(1-\cos(\phi-\theta))}
(\mathfrak a_1\tilde{\mathfrak a}_1+\mathfrak a_2\tilde{\mathfrak a}_2)+O(r^5)
,
\end{align*}
where   $\mathfrak a_1=\tilde{\mathfrak a}_1$ and $\mathfrak a_2=\tilde{\mathfrak a}_2$.
(In fact, $\mathfrak a_1$ and  $\mathfrak a_2$ were chosen accordingly.)
Regarding the previous computation, we note that
$D_{\frac{\m v_i+\m v_i*}{2}}=1+O(r)$
as $r\searrow0$ (for fixed $\theta,\phi,\eta$).
Thus $\sqrt{D_{\frac{\m v_i+\m v_i*}{2}}}=1+O(r)$ also.
This makes the computation valid for small $r$.

Our next observation is that $(\mathfrak a_1,\mathfrak a_2)\neq(0,0)$.
Indeed, suppose that $\mathfrak a_1=0$.
This means that
\[ \cos \left( \eta+\phi \right) =\cos \left( \eta+\theta \right). \]
As we know
\[\eta+\phi\not\equiv\eta+\theta\modu2\pi,\]
this implies
\[ \sin \left( \eta+\phi \right) =-\sin \left( \eta+\theta \right)\neq 0 \]
and
\[ \cos \left( \eta+\phi \right) =\cos \left( \eta+\theta \right)\neq\pm1. \]
Substituting with these yields
\[\mathfrak a_2=\frac1{\frac12(1-\cos(\phi-\theta))}\Bigl(
2\sin \left( \eta+\theta \right)  \left( \cos \left( \eta+\theta
 \right) -1 \right)  \left( \cos \left( \eta+\theta \right) -2 \right)
\Bigr).
\]
Then $\mathfrak a_2$ is nonzero, as any multiplicative component is nonzero.

By this, we have shown that
\[\left.\frac{\mathrm d\|\tilde A(t)\|_2}{\mathrm dt}\right|_{t=0+}+\left.\frac{\mathrm d\|\tilde B(t)\|_2}{\mathrm dt}\right|_{t=0+}<0\]
if $r$ is sufficiently small (for fixed $\theta,\phi,\eta$).
This already contradicts to Magnus minimality, as we can always restrict to around the
join area between $A$ and $B$.

The real case is not different but we can choose only  $\eta,\theta,\phi\in\{0,\pi\}$.
The resulted matrices are all real (as $\sin(\theta+\eta)$, etc., are all $0$).
\end{proof}
\end{theorem}
\snewpage
\begin{theorem}
\plabel{th:mBCHmin}
(a) Assume $\phi=A_1\m 1_{[0,t_1)}\boldsymbol.\ldots\boldsymbol.A_k\m 1_{[t_{k-1},t_k)}$
is a mass-normalized mBCH measure of real
$2\times2$ matrices, which is Magnus minimal,
i. e. , \[\mathcal M_{n\times n\,\,\real}(\Rexp(\phi))=\smallint\|\phi\|_2.\]
Then there are two possible cases:

(i) Up to conjugation by orthogonal matrices, $\phi$ is of shape
$\phi_1\oplus\phi_2$ such that one of the components is of shape $1\cdot\m 1_{[0,t_k)}$
or $-1\cdot\m 1_{[0,t_k)}$ ($1\times1$ matrix) [hyperbolic case]. Or,

(ii) $\phi$ is constant orthogonal, i. e. it is of shape $U\m 1_{[0,t_k)}$ where
$U$ is orthogonal [elliptic case].

(b) Assume $\phi=A_1\m 1_{[0,t_1)}\boldsymbol.\ldots\boldsymbol.A_k\m 1_{[t_{k-1},t_k)}$
is a mass-normalized mBCH measure of complex
$2\times2$ matrices, which is Magnus minimal,
i. e. , \[\mathcal M_{n\times n\,\,\complex}(\Rexp(\phi))=\smallint\|\phi\|_2.\]
Then, up to conjugation by unitary matrices, $\phi$ is of shape
$\phi_1\oplus\phi_2$ such that one of the components is of shape $u\cdot\m 1_{[0,t_k)}$
where $u\in\partial\Dbar(0,1)$.

As a consequence we see that reduced mass-normalized mBCH measures of $2\times2$ cannot be minimal
(valid in real and complex sense alike).
\begin{proof}
(a) From the previous lemmas we know that the matrices must be normal, one
aligned next to each other.
If one of the matrices is of elliptic type then its neighbour must be equal to it.
Assume that all $A_i$ are of decomposable (parabolic or hyperbolic) type.
Note that scalar matrices are freely movable in the decomposition so we can temporarily assume that
the $A_i$ of hyperbolic type are next to each other.
Then its eigenspace decompositions are the same.
Beyond that we have only scalar matrices, thus decomposability follows.
Minimality implies that at least one of the components is minimal,
thus we have the special shape.
(b) is similar except simpler.
\end{proof}
\proofremarkqed{
Minimality may lead to further restrictions on the orthogonal $U$ or the unit $u$.
However, if $\smallint\|\phi\|_2\leq\pi$, then all the indicated shapes are Magnus minimal.
}
\end{theorem}

\snewpage\scleardoublepage\section{Asymptotics of some BCH expansions from $\SL_2(\mathbb R)$}
\plabel{sec:ExamplesBCH}
For the purposes of this section, we consider some auxiliary functions.
\begin{lemma}\plabel{lem:AST}
(a) The function
\[\AS(x)=\sqrt{\frac{\AC(x)^2-1}{1-x^2}}\]
extends from $(-1,1)\cup(1,+\infty)$ to an analytic function on $\mathbb C\setminus(-\infty,-1]$.
$\AS(0)=\frac{\sqrt3}{3}$; and $\AS$ is nowhere vanishing on $\mathbb C\setminus(-\infty,-1]$.

(b) $\AS$ is monotone decreasing on $(-1,\infty)$ with range $(+\infty,0)_{\mathrm e}$.
\begin{proof}
(a) $\mathbb C\setminus(-\infty,-1]$ is simply connected, and the branchings (or vanishing) for
the square roots can occur at $z=\pm$ or $\AC(z)=\pm$, thus ultimately only for $z=1$.
Thus the statement is sufficient to check at $z=1$ with power series.
\begin{commentx}
In fact, an explicit extension formula is given by
\[\AS(x)=\sqrt{\frac{\AC(z)-1}{1-z}}\cdot\frac{ \sqrt{\AC(z)+1}  }{ \sqrt{1+z}  },\]
as for analytic functions, with the standard square root used.
\end{commentx}

(b) Elementary function calculus.
\end{proof}
\end{lemma}

\begin{lemma}\plabel{lem:AtT}
(a) The function
\[\AT(z)=\frac{\AC(z)-1}{\AS(z)}\]
is analytic on $\mathbb C\setminus(-\infty,-1]$.
$\AT(1)=0$ and $\AT'(0)=-\frac{\sqrt3}{3}$;
$\AT(z)$ vanishes only at $z=0$ for  $z\in \mathbb C\setminus(-\infty,-1]$.

(b) The function $x\mapsto x+\AT(x)$ is monotone increasing on $(-1,+\infty)$ (bijectively);
moreover, it also yields a bijection from $(-1,1]$ to itself.
\begin{proof}
(a) This is a consequence of Lemma \ref{lem:AST}(a).

(b) Elementary function calculus.
\end{proof}
\end{lemma}
\begin{remark}
Consequently, for $z\in\mathbb C\setminus(-\infty,-1]$,
\[
\AC(z)=\frac{1-z^2+\AT(z)^2}{1-z^2-\AT(z)^2}
\]
and
\[
\AS(z)=\frac{2\AT(z)}{1-z^2-\AT(z)^2}
\]
as analytic functions. But, in fact, the denominators vanish only az $z=1$.
\qedremark
\end{remark}
\snewpage
\begin{example}\plabel{ex:loxocomp}
Consider the matrices
\[\tilde J=\begin{bmatrix}1&\\&-1\end{bmatrix},\qquad \tilde I=\begin{bmatrix}0&-1\\1&0\end{bmatrix}.\]
For $\alpha,\beta\in\mathbb C$, let
\[\Upsilon_{\alpha,\beta}=\alpha \tilde J\m 1\boldsymbol.\beta \tilde I\m 1. \]
Then
\[\int\|\Upsilon_{\alpha,\beta}\|_2=|\alpha|+|\beta|.\]

For $|\alpha|+|\beta|<\pi$, we can consider
\begin{align}
\mu_{\mathrm L}(\Upsilon_{\alpha,\beta})&=\log(\exp_{\mathrm L}(\Upsilon_{\alpha,\beta}))
\notag\\
&= \log(\exp(\beta\tilde  I)\exp(\alpha \tilde J))
\notag\\
&=\log \begin{bmatrix}\mathrm e^\alpha\cos\beta&-\mathrm e^{-\alpha}\sin\beta\\\mathrm e^{\alpha}\sin\beta&\mathrm e^{-\alpha}\cos\beta\end{bmatrix}
\notag\\
&=\AC(\cosh\alpha\cos\beta)
\begin{bmatrix}\sinh\alpha\cos\beta&-\mathrm e^{-\alpha}\sin\beta\\\mathrm e^{\alpha}\sin\beta  &-\sinh\alpha\cos\beta \end{bmatrix}.
\notag
\end{align}

If $\alpha,\beta\geq0$, then
\[\|\mu_{\mathrm L}(\Upsilon_{\alpha,\beta})\|_2=\AC(\cosh\alpha\cos\beta)\cdot(\sinh\alpha+\cosh\alpha\sin\beta).\]

(a) Now, for $p\in[0,\pi)$, let
\[\tilde \alpha(p)=p-\pi+\sqrt[3]{\pi^2(\pi-p)},\]
\[\tilde\beta(p)=\pi-\sqrt[3]{\pi^2(\pi-p)}.\]
Then $\tilde \alpha(p),\tilde \beta(p)\geq0$, and
\[\tilde \alpha(p)+\tilde \beta(p)=p.\]
Thus,
\[\int\|\Upsilon_{\tilde\alpha(p),\tilde\beta(p)}\|_2=p. \]
As $p\nearrow\pi$, we see that $\tilde \alpha(p)\searrow0$ (eventually) and $\tilde \beta(p)\nearrow\pi$.
Consequently
\[\lim_{p\rightarrow \pi} \cosh\tilde  \alpha(p)\cos\tilde \beta(p)=-1. \]
In that (elliptic) domain $\AC$ is computed by $\arccos$.
Now, elementary function calculus shows that as $p\nearrow\pi$,
\begin{align}
\|\mu_{\mathrm L}(\Upsilon_{\tilde \alpha(p),\tilde \beta(p)})\|_2&\stackrel\rightarrow=
\frac{\arccos(\cosh\tilde \alpha(p)\cos\tilde \beta(p))}{\sqrt{1-\cosh^2\tilde \alpha(p)\cos^2\tilde \beta(p)}}
(\sinh\tilde \alpha(p)+\cosh\tilde \alpha(p)\sin\tilde \beta(p))
\notag\\
&=\sqrt{\frac{12\pi^{8/3}}{\pi^2+6}} (\pi-p)^{-1/3}+O((\pi-p)^{1/3}).
\notag
\end{align}

We see that in Baker--Campbell--Hausdorff  setting we can produce the asymptotics $O((\pi-p)^{-1/3})$,
although having exponent $-1/3$ instead of $-1/2$ is strange.
\snewpage

(b) It is interesting to see that in the setting of the present example, one cannot do much better:

If we try to optimize $\|\mu_{\mathrm L}(\Upsilon_{\alpha,\beta})\|_2$ for $\alpha+\beta$   $(\alpha,\beta\geq0)$ ,
then, after some computation, it turns out that the best approach is
along a well-defined ridge.
This ridge starts hyperbolic, but turns elliptic.
Its elliptic part is part is parametrized by
$x\in (-1,1]$, and
\[\hat\alpha(x)=\arcosh\left(\frac{\AC(x)+\sqrt{\AC(x)^2-4x(1-x\AS(x))\AS(x)}}{2(1-x\AS(x))}
\right);\]
\[\hat\beta(x)=\arccos\left(\frac{\AC(x)-\sqrt{\AC(x)^2-4x(1-x\AS(x))\AS(x)}}{2\AS(x)}
\right).\]

Then
\[\cosh\hat\alpha(x)\cos\hat\beta(x)=x.\]
Actually, $x=1$ gives a parabolic $\Lexp(\Upsilon_{ \hat\alpha(x),\hat\beta(x)})$, but
for $x\in(-1,1)$ it is elliptic.
Then $\hat\alpha(x),\hat\beta(x)\geq0$.
As $x\searrow-1$, one can see that $\alpha\searrow0$ (eventually) and $\beta\nearrow\pi$;
and, more importantly,
\[\hat\alpha(x)+\hat\beta(x)\nearrow\pi. \]
Now, as $x\searrow-1$,
\[\frac{\arccos x}{\sqrt{1-x^2}}(\sinh\hat \alpha(x)+\cosh\hat \alpha(x)\sin\hat \beta(x))=\pi2^{3/4}(x+1)^{-1/4}+O((x+1)^{1/4}),\]
and
\[\pi-\hat\alpha(x)-\hat\beta(x)=\frac13 2^{3/4}(x+1)^{3/4}+O((x+1)^{5/4}).\]
Hence, using the notation $\hat p(x)=\hat\alpha(x)+\hat\beta(x)$, we find
\[\|\mu_{\mathrm L}(\Upsilon_{ \hat\alpha(x),\hat\beta(x)})\|_2=2\pi3^{-1/3}(\pi- \hat p(x))^{-1/3}+O((\pi- \hat p(x))^{1/3}).\]

This $2\pi3^{-1/3}=4.356\ldots$ is just slightly better than
$\sqrt{\frac{12\pi^{8/3}}{\pi^2+6}}=4.001\ldots$.

(c)
In the previous example let $\alpha=\pi-p, \beta=p$ with $p\in[\pi/2,\pi)$.
Let
\begin{equation}
G(p):=\|\mu_{\mathrm L}(\Upsilon_{ \pi-p,p})\|_2
=\AC(\cosh(\pi-p)\cos p)\cdot(\sinh(\pi-p)+\cosh(\pi-p)\sin p)
.
\plabel{eq:Gdef}
\end{equation}
Then
\[\|\mu_{\mathrm L}(\Upsilon_{ \pi-p,p})\|_2=2\pi\sqrt3(\pi-  p)^{-1}+ \left(\frac\pi2\sqrt3-2\right)  (\pi-  p)  + O((\pi- p)^3).\]

In fact, a special value is
\[G\left(\frac\pi2\right)=\frac\pi2\exp\frac\pi2.\]
Using elementary analysis, one can see that
$G(p)$ is strictly monotone increasing on   $p\in[\pi/2,\pi)$.
\qedexer
\end{example}
\snewpage
%\begin{commentx}

\begin{example}\plabel{ex:ellicomp}
Consider the matrices
\[\tilde P=\begin{bmatrix}0&-1\\&0\end{bmatrix},\qquad \tilde I=\begin{bmatrix}&-1\\1&\end{bmatrix}.\]
For $\alpha,\beta\in\mathbb C$, let
\[\tilde\Upsilon_{\alpha,\beta}=\alpha \tilde P\mathbf 1\boldsymbol.\beta \tilde I\mathbf 1. \]
Then
\[\int\|\tilde\Upsilon_{\alpha,\beta}\|_2=|\alpha|+|\beta|.\]

For $|\alpha|+|\beta|<\pi$,  we can consider
\begin{align}
\mu_{\mathrm L}(\tilde\Upsilon_{\alpha,\beta})&=\log(\exp_{\mathrm L}(\tilde \Upsilon_{\alpha,\beta}))
\notag\\
&= \log(\exp(\beta\tilde  I)\exp(\alpha \tilde P))
\notag\\
&=\log \begin{bmatrix}\cos\beta&-\alpha\cos\beta-\sin\beta\\
\sin\beta&-\alpha\sin\beta+\cos\beta\end{bmatrix}
\notag\\
&=\AC\left(\cos\beta-\frac\alpha2\sin\beta\right)
\begin{bmatrix}\frac\alpha2\sin\beta &\alpha\cos\beta-\sin\beta\\
\sin\beta  &-\frac\alpha2\sin\beta \end{bmatrix}.
\notag
\end{align}

If $\alpha,\beta\geq0$, then
\[\|\mu_{\mathrm L}(\tilde\Upsilon_{\alpha,\beta})\|_2=\AC\left(\cos\beta-\frac\alpha2\sin\beta\right)\cdot
\left(\sin\beta+\frac\alpha2\cos\beta+\frac\alpha2\right).\]

For optimal approach, consider $x\in(-1,1]$, and let
\[\hat\alpha(x)=\frac{2\AT(x)}{\sqrt{1-(x+\AT(x))^2}};\qquad \hat\beta(x)=\arccos\left(x+\AT(x)\right).\]

Then
\[\cos\hat\beta(x)-\frac{\hat\alpha(x)}2\sin\hat\beta(x)=x.\]

As $x\searrow-1$, we have $\alpha\searrow0$ (eventually) and $\beta\nearrow\pi$;
and,  $\hat\alpha(x)+\hat\beta(x)\nearrow\pi.$
Now, as $x\searrow-1$,
\begin{align}
\|\mu_{\mathrm L}(\tilde\Upsilon_{ \hat\alpha(x),\hat\beta(x)})\|_2&=
\frac{\arccos x}{\sqrt{1-x^2}}\left(\sin\hat\beta(x)+\frac{\hat\alpha(x)}2\cos\hat\beta(x)+\frac{\hat\alpha(x)}2\right)
\notag\\
&=2^{1/4}\pi(x+1)^{-1/4}+O((x+1)^{1/4}),
\notag
\end{align}
and
\[\pi-\hat\alpha(x)-\hat\beta(x)=\frac23 2^{1/4}(x+1)^{3/4}+O((x+1)^{5/4}).\]
Hence, using the notation $\hat p(x)=\hat\alpha(x)+\hat\beta(x)$, we find
\[\|\mu_{\mathrm L}(\tilde\Upsilon_{ \hat\alpha(x),\hat\beta(x)})\|_2=\pi(4/3)^{1/3}(\pi- \hat p(x))^{-1/3}+O((\pi- \hat p(x))^{1/3}).\]
This leading coefficient $\pi(4/3)^{1/3}=3.457\ldots$ is worse than the corresponding one in the previous example.
\qedexer
\end{example}
%\end{commentx}

The previous two examples suggest, at least in the regime of $2\times2$ real matrices, two ideas.
Firstly, that for BCH expansions substantially stronger asymptotical estimates apply as compared to general Magnus expansions.
Secondly, that for larger norms in BCH expansions normal matrices are preferred.
These issues will be addressed subsequently, with more machinery.
\snewpage
\scleardoublepage\section{Critical singular behaviour of BCH expansions in the $2\times2$ real case}
\plabel{sec:critbal}
For the purposes of the next statement let
\[U=\{A\in\mathrm M_2(\mathbb R)\:\,\|A\|_2\leq1\}\]
and
\[X= (-\pi,\pi)\times U\times U.\]
According to Theorem \ref{th:better2}, the map
\[\mathcal B:(\alpha,A,B)\mapsto \mu_{\mathrm R}\left(\frac{\pi-\alpha}2A\m 1\boldsymbol.\frac{\pi+\alpha}2B\m 1\right) \]
is defined everywhere. In particular, the subset
\[X_0=
(-\pi,\pi)\times \left\{\tilde I\right\}\times \left\{\tilde I\right\}
\cup
(-\pi,\pi)\times \left\{-\tilde I\right\}\times \left\{-\tilde I\right\}
\]
is included, where $\mathcal B$ takes values $\pi\tilde I$ or $-\pi\tilde I$.
\begin{theorem}
\plabel{th:discont22}
(a) For $(\alpha,A,B)\in X\setminus X_0$,
\[\mathcal B(\alpha,A,B)=\log\left(\exp \left(\frac{\pi-\alpha}2A\right)\exp \left(\frac{\pi+\alpha}2B\right)\right); \]
and $\mathcal B|_{X\setminus X_0}$ is analytic (meaning that it is a restriction from an analytic function
which is defined on an  open subset of $(-\pi,\pi)\times \mathrm M_2(\mathbb R)\times\mathrm M_2(\mathbb R) $ containing
$X\setminus X_0$.)

(b) Let $(\alpha_0,I_0,I_0)\in X_0$.
Then
\begin{align*}
\pi&= \liminf_{(\alpha,A,B)\rightarrow (\alpha_0,I_0,I_0)}  \|  \mathcal B(\alpha,A,B)\|_2<\\
&  <\limsup_{(\alpha,A,B)\rightarrow (\alpha_0,I_0,I_0)}  \|  \mathcal B(\alpha,A,B)\|_2
=\pi\sqrt{\frac{\pi-|\alpha_0|+2\cos\frac{\alpha_0}2}{\pi-|\alpha_0|-2\cos\frac{\alpha_0}2}}.
\end{align*}

(Here $\liminf$ and $\limsup$ can be understood either as for $(\alpha,A,B)\in X\setminus X_0$ or $(\alpha,A,B)\in X$.)
\begin{proof}
(a) This follows from Theorem \ref{th:better3} and the analyticity of $\log$.

(b) In these circumstances, $\tr A,\tr B\rightarrow 0$, and as trace can factorized out, and as
detracing does not increase the norm (Corollary \ref{lem:detrace}), we can restrict to the case $\tr A=\tr B= 0$.
By conjugation invariance we can also assume $I_0= \tilde I$.
Taking \eqref{eq:norm22skew} into consideration,
we can restrict to matrices
\[A=\left(1-\frac{1+t}2\xi\right)\tilde I+\frac{1+t}2\xi\left(
s_1\tilde J+s_2\tilde K
\right) \]
and
\[B=\left(1-\frac{1-t}2\xi\right)\tilde I+\frac{1-t}2\xi\left(
r_1\tilde J+r_2\tilde K
\right), \]
where $0<\xi\sim0$, and  $|t|, \sqrt{s_1^2+s_2^2}, \sqrt{r_1^2+r_2^2}\leq 1$.

Essentially, we have to consider $\xi\searrow0$ while $\alpha\rightarrow \alpha_0$.
(Thus all parameters depend on indices $\lambda\in\Lambda$ which we omit.)
\snewpage
Going through the computations, we find
\begin{multline*}
\left\|\log\left(\exp \left(\frac{\pi-\alpha}2A\right)\exp\left( \frac{\pi+\alpha}2B\right)\right)\right\|_2=(\pi+O(\xi))\cdot
\\
\cdot\frac{\sqrt{ \xi^2(\frac{\pi-t\alpha}2 )^2+O(\xi^3)}+\sqrt{ \xi^2\left(\cos\frac\alpha2\right)^2 \left(
\left(s_1\frac{1+t}2-r_1\frac{1-t}2\right)^2+\left(s_2\frac{1+t}2-r_2\frac{1-t}2\right)^2
\right) +O(\xi^3)}}{\sqrt{ \xi^2\left((\frac{\pi-t\alpha}2)^2- \left(\cos\frac\alpha2\right)^2  \left(
\left(s_1\frac{1+t}2-r_1\frac{1-t}2\right)^2+\left(s_2\frac{1+t}2-r_2\frac{1-t}2\right)^2
\right) \right) +O(\xi^3)}}.
\end{multline*}
Note that due to the convexity of the unit disk,
\[0\leq \underbrace{\left(s_1\frac{1+t}2-r_1\frac{1-t}2\right)^2+\left(s_2\frac{1+t}2-r_2\frac{1-t}2\right)^2}_{E^2:=}\leq1.\]
Furthermore, every value from $[0,1]$ can be realized as $E$ (even if we fix $t$).
Thus
\[\left\|\log\left(\exp\left( \frac{\pi-\alpha}2A\right)\exp \left(\frac{\pi+\alpha}2B\right)\right)\right\|_2=\pi\sqrt{\frac{(\pi-t\alpha)+ \left(\cos\frac\alpha2\right) E}{ (\pi-t\alpha)- \left(\cos\frac\alpha2\right) E}}+O(\xi).\]

From this, the statement follows.
\end{proof}
\end{theorem}
\begin{remark}
Using the notation $p=\frac{\pi+|\alpha|}2=\max\left\{\frac{\pi-\alpha}2, \frac{\pi+\alpha}2\right\}$, we find
\[\pi\sqrt{\frac{\pi-|\alpha|+2\cos\frac{\alpha}2}{\pi-|\alpha|-2\cos\frac{\alpha}2}}
=\pi\sqrt {{\frac {\pi -p+\sin  p  }{\pi -p-\sin  p }}}
={\frac {2\pi \,\sqrt {3}}{\pi -p}}-\,\frac{\pi \,\sqrt {3}}{30} \left( \pi
-p \right)+O((\pi-p)^3)
\]
as $p\nearrow\pi$.
As for crude estimates, one can check that
\[\pi\sqrt {{\frac {\pi -p+\sin  p  }{\pi -p-\sin  p }}}<{\frac {2\pi \,\sqrt {3}}{\pi -p}}<G(p)\]
holds for $p\in[\frac\pi2,\pi)$, cf. \eqref{eq:Gdef}.
\qedremark
\end{remark}
\begin{theorem}
\plabel{th:supbound22}
Let $\alpha_0\in (-\pi,\pi)$,    $p=\frac{\pi+|\alpha_0|}2=\max\left\{\frac{\pi-\alpha_0}2, \frac{\pi+\alpha_0}2\right\}$.
Cf. Example \ref{ex:loxocomp}(c) for the definition of $G(p)$.
Then
\[G(p)\leq\left(\sup_{\substack{(\alpha,A,B)\in X\\ |\alpha|\leq|\alpha_0| }} \|\mathcal B(\alpha,A,B) \|_2\right) <+\infty.\]
\begin{proof}
The lower estimate is obvious from \eqref{eq:Gdef}.
As for the upper estimate,
Let $M$ be so that $M>\pi\sqrt {{\frac {\pi -p+\sin  p  }{\pi -p-\sin  p }}}$.
Then
Theorem \ref{th:discont22} shows that
\[\{ (\alpha,A,B)\in X\,:\, |\alpha|\leq|\alpha_0| , \mathcal B(\alpha,A,B)\geq M  \}\]
is a compact subset of $X$, on which $\mathcal B$ is continuous.
From this compactness, the boundedness of $\|\mathcal B\|_2$ follows.
\end{proof}
\end{theorem}
Thus for $2\times 2$ real matrices, sufficient balancedness implies uniform boundedness for BCH expansions.
The most interesting case is $p=\pi/2$ (i. e. $\alpha_0=0$).
In this case the lower bound is sharp;  using a slightly different formulation,
see Theorem \ref{th:supequal22}.

\snewpage
\scleardoublepage\section{Moments associated to Schur's formulae}
\plabel{sec:BCHmoments}

Let
\[\mathcal S=\{A\in\mathrm M_2(\mathbb R)\,:\,\underbrace{ \spec(A)\subset \{z\in\mathbb C\,:\,|\Ima z|<\pi\}
}_{\equiv D_A<\pi^2\text{ in } \mathrm M_2(\mathbb R)},\,-D_{ A^*A}\neq0\}.\]
We can decompose $\mathcal S$ further; so that
\[\mathcal S^{\mathrm{nn}}=\{A\in\mathcal S\,:\,[A,A^*]\neq0\}\]
is the regular or non-normal (here: non-selfadjoint) interior of $\mathcal S$; and
\[\eth^{\mathrm{par}}\mathcal S=\{A\in\mathcal S\,:\,A=A^*\} \]
is the Schur-parabolic or normal (here: self-adjoint) pseudoboundary.
Then
\[\mathcal S=\mathcal S^{\mathrm{nn}}\,\,\dot\cup\,\,\eth^{\mathrm{par}}\mathcal S.\]

Note that $\mathcal S$ is open in $\mathrm M_2(\mathbb R)$, consequently, for the closure in $\mathrm M_2(\mathbb R)$,
\[\overline{\mathcal S}=\mathcal S\,\,\dot\cup\,\,\partial\mathcal S.\]

We can decompose $\partial \mathcal S\subset\mathrm M_2(\mathbb R)$ so that
\[\partial^0\mathcal S=\{0\}\]
is the zero-boundary;
\[\partial^{\mathrm{hyp}}\mathcal S= \{A\in\mathrm M_2(\mathbb R)\,:\, \spec(A)\subset \{z\in\mathbb C\,:\,|\Ima z|<\pi\},\, \text{ $A$ is a conform-reflexion}\}\]
is the Schur-hyperbolic or conform-reflexional boundary;
\[\partial^{\mathrm{ell}}\mathcal S= \{A\in\mathrm M_2(\mathbb R)\,:\, \spec(A)\subset \{z\in\mathbb C\,:\,|\Ima z|<\pi\},\, \text{ $A$ is a conform-rotation}\}\]
is the Schur-elliptic or conform-rotational boundary;
\[\partial^{\mathrm{dell}}\mathcal S= \{A\in\mathrm M_2(\mathbb R)\,:\, \spec(A)\subset (\mathrm i\pi+\mathbb R)\cup (-\mathrm i\pi+\mathbb R),\, \text{ $A$ is a conform-rotation}\}\]
is the Schur-elliptic degenerate boundary;
\[\partial^{\mathrm{dnn}}\mathcal S= \{A\in\mathrm M_2(\mathbb R)\,:\, \spec(A)\subset (\mathrm i\pi+\mathbb R)\cup (-\mathrm i\pi+\mathbb R),\, \text{ $A$ is not a conform-rotation}\}\]
is the non-normal degenerate boundary.
Then
\[\partial \mathcal S= \partial^0\mathcal S \,\,\dot\cup\,\, \partial^{\mathrm{hyp}}\mathcal S
\,\,\dot\cup\,\, \partial^{\mathrm{ell}}\mathcal S
\,\,\dot\cup\,\, \partial^{\mathrm{dell}}\mathcal S \,\,\dot\cup\,\,\partial^{\mathrm{dnn}}\mathcal S .\]

Some components here can be decomposed further naturally:
$\partial^{\mathrm{ell}}\mathcal S$ can be decomposed to
conform-identity part $\partial^{\mathrm{ell1}}\mathcal S$ and the
generic part $\partial^{\mathrm{ell*}}\mathcal S$;
$\partial^{\mathrm{dell}}\mathcal S$ can be decomposed to
traceless part $\partial^{\mathrm{dell0}}\mathcal S$ and the
generic part $\partial^{\mathrm{dell*}}\mathcal S$; and similarly
$\partial^{\mathrm{dnn}}\mathcal S$ can be decomposed to
traceless part $\partial^{\mathrm{dnn0}}\mathcal S$ and the
generic part $\partial^{\mathrm{dnn*}}\mathcal S$.

We let
\[\mathcal S^{\mathrm{ext}}=\mathrm M_2(\mathbb R)\setminus \overline{\mathcal S},\]
the external set of $\mathcal S$.
This set can be decomposed to the degenerate or normal exterior $\mathcal S^{\mathrm{ellext}}$
and the non-normal exterior $\mathcal S^{\mathrm{nnext}}$.
It is reasonable to discriminate $\widehat X^{\mathrm{ext}}$ further by
whether $D_A\in\{k^2\pi^2\,:\,k\in\mathbb N\setminus\{0\}\}$
holds or not. However, $\mathcal S^{\mathrm{ext}}$ will not be, ultimately, of much interest for us:

Let
\[\mathcal S^{\mathrm{acc}}=\mathcal S^{\mathrm{nn}}\,\,\dot\cup\,\,\eth^{\mathrm{par}}\mathcal S
\,\,\dot\cup\,\, \partial^0\mathcal S
\,\,\dot\cup\,\,\partial^{\mathrm{hyp}}\mathcal S
\,\,\dot\cup\,\, \partial^{\mathrm{ell}}\mathcal S
\,\,\dot\cup\,\, \partial^{\mathrm{dell}}\mathcal S\]
(where $\partial^{\mathrm{dnn}}\mathcal S $ and $\mathcal S^{\mathrm{ext}}$ are not included).
Here the point is that elements of $\partial^{\mathrm{dnn}}\mathcal S$ and $\mathcal S^{\mathrm{ext}}$  are not that
useful to exponentiate, as they can easily be replaced by elements of smaller norm:

\begin{lemma}
\plabel{lem:Sacc}
(a) If $A\in \mathrm M_2(\mathbb R)\setminus {\mathcal S^{\mathrm{acc}}}$, then there exists $B\in{\mathcal S^{\mathrm{acc}}}$ such that $\exp B=\exp A$ and $\|B\|_2<\|A\|_2$.

(b) If $A\in {\mathcal S^{\mathrm{acc}}}\setminus \partial^{\mathrm{dell}}\mathcal S$, there there is only one $B\in {\mathcal S^{\mathrm{acc}}}$ such that $\exp B=\exp A$;
namely $B=A$.

(c)  If $A\in \partial^{\mathrm{dell}}\mathcal S$, there are two possible $B\in {\mathcal S^{\mathrm{acc}}}$ such that $\exp B=\exp A$.
If $A=a\Id_2+\pi\tilde I$ or $A=a\Id_2-\pi\tilde I$, then $B=a\Id_2\pm\pi\tilde I$.
\begin{proof}
This is straightforward from Lemma \ref{lem:PosRExp}.
\end{proof}
\end{lemma}
%\snewpage
Note that $A\in\mathcal S$ if and only if  $\log\exp A=A$ and the norm is smooth at $A$.
(This is the reason for the notation.)
Then
\begin{equation}
\mathrm{MR}_A(\m v)=D_{\m v \text{ at }M=0}\left(M\in \mathrm M_2(\mathbb R)\mapsto   \|\log((\exp A)(\exp M))\|_2\right)
\plabel{eq:MRdef}
\end{equation}
yields a linear operator $\mathrm{MR}_A$.
This can be represented by $\mathrm{MR}(A)\in\mathrm M_2(\mathbb R) $ such that
\[\mathrm{MR}_A(\m v)= \frac12\tr\left(\m v \mathrm{MR}(A)^*\right) .\]
Similarly, we can define
\begin{equation}
\mathrm{ML}_A(\m v)=D_{\m v \text{ at }M=0}\left(M\in \mathrm M_2(\mathbb R)\mapsto   \|\log((\exp M)(\exp A))\|_2\right)
\plabel{eq:MLdef}
\end{equation}
and $\mathrm{ML}(A)\in\mathrm M_2(\mathbb R) $ such that
\[\mathrm{ML}_A(\m v)= \frac12\tr\left(\m v \mathrm{ML}(A)^*\right) .\]
As $\log$, $\exp$, $\|\cdot\|_2$ are locally open at these expressions, we know that $\mathrm{ML}(A),\mathrm{MR}(A)\neq0$
(for $A\in\mathcal S$).
It is easy to see that if $U$ is orthogonal, then
\begin{equation}\mathrm{MR}_{UAU^{-1}}(U\m vU^{-1})=\mathrm{MR}_A(\m v)\plabel{eq:MRinv}\end{equation}
and
\[ \mathrm{MR}(UAU^{-1})=U\mathrm{MR}(A)U^{-1}.\]
Similarly statements holds for $\mathrm{ML}$.
Furthermore,
\begin{equation}\mathrm{MR}_A(\m v)=\mathrm{ML}_{A^*}(\m v^*)\plabel{eq:MRadj}\end{equation}
and
\[\mathrm{MR}(A)=(\mathrm{ML}(A^*))^*.\]

By Lemma \ref{lem:LocLog}, the condition $A\in \mathcal S$ (ie.  $D_{A^*A}\neq0$ \& $D_A<\pi^2$)
can be relaxed to $D_{A^*A}\neq0$ \& $D_A\notin \{k^2\pi^2\,:\, k\in\mathbb N\setminus\{0\}\}$ in the previous discussion.
(Thus it extends to most of $\mathcal S^{\mathrm{ext}}$.)

If $A\in \partial^{0}\mathcal S\cup\partial^{\mathrm{hyp}}\mathcal S\cup\partial^{\mathrm{ell}}\mathcal S$,
then \eqref{eq:MRdef} still makes sense, although it is not linear in $\m v$.
Thus, the (non-linear) forms $\mathrm{MR}_A(\m v)$ and  $\mathrm{ML}_A(\m v)$ are defined,
but we leave $\mathrm{MR}(A)$ and $\mathrm{ML}(A)$ undefined.
Nevertheless, \eqref{eq:MRinv} and \eqref{eq:MRadj} still hold.

If  $A\in \partial^{\mathrm{dell}}\mathcal S$, then we will be content to define
$\mathrm{MR}_A(\m v)=\mathrm{ML}_A(\m v)$ for $\m v\in\mathbb R\Id_2+\mathbb R\tilde I$.
This is done as follows:
We extend $\log$ as $\log^*$ such that for negative scalar matrices $B=-\lambda \Id_2$  it yields
$\log^* B=  (\log\lambda)\Id_2 \pm\pi\tilde I$ (two-valued).
Then $\|\log^* B\|_2$ still makes sense uniquely.
Now, for $\m v\in\mathbb R\Id_2+\mathbb R\tilde I$,  we define
\[\mathrm{MR}_A(\m v)=\mathrm{ML}_A(\m v)=\lim_{t\searrow0}\frac{   \|\log^*(\exp(\m vt)\exp(A))\|_2 -\| A\|_2}{t}\]
($\exp A$ is central).
This argument applies more generally, if $D_A\in \{k^2\pi^2\,:\, k\in\mathbb N\setminus\{0\}\}$,
with respect to an appropriate skew-involution $I_A$; but it will be no interest for us.

Another observation is that $\mathrm{MR}_A(\cdot)$ and $\mathrm{ML}_A(\cdot)$ are never trivial (i. e. identically zero).
Indeed, in each case, the direction $\m v=A$ yields to a differentiable increase or decrease (the latter is for  $A\in \partial^{\mathrm{dell}}\mathcal S$)  in the norm.
We will see concrete expressions for $\mathrm{MR}_A(\m v)$ and $\mathrm{ML}_A(\m v)$ later.

%\snewpage
\begin{lemma}\plabel{lem:momentkif}
Suppose that $A= a\Id_2+ b\tilde I+(r\cos\psi)\tilde J+(r\sin\psi)$
such that
\begin{equation}
a^2+b^2>0,\qquad r>0,
\plabel{eq:dom1}
\end{equation}
and
\begin{equation}
b^2-r^2\notin\{k^2\pi^2\,:\,k\in\mathbb N\setminus\{0\}\}.
\plabel{eq:dom2}
\end{equation}
(Note $A\in\mathcal S$ if and only if \eqref{eq:dom1} and $b^2-r^2<\pi^2$ hold.)

(o) Then
\begin{align*}
\mathrm{MR}(A)=&\hat a\Id_2+\hat b\tilde I+\hat c_{\mathrm R}\tilde J+\hat d_{\mathrm R}\tilde K\\
=&\hat a\Id_2+\hat b\tilde I+(\cos\psi\Id_2+\sin\psi \tilde I)(\breve c\tilde J+\breve d\tilde K),
\end{align*}
and
\begin{align*}
\mathrm{ML}(A)=&\hat a\Id_2+\hat b\tilde I+\hat c_{\mathrm L}\tilde J+\hat d_{\mathrm L}\tilde K\\
=&\hat a\Id_2+\hat b\tilde I+(\cos\psi\Id_2+\sin\psi \tilde I)(\breve c\tilde J-\breve d\tilde K),
\end{align*}
where
\begin{align}\hat a&=\frac{a}{\sqrt{a^2+b^2}}
\plabel{eq:crux}\\
\notag
\hat b&=\frac{b}{\sqrt{a^2+b^2}}\left(1+ \left( \sqrt{a^2+b^2}+r\right) \cdot r\reC(b^2-r^2)\right) \\
\notag
\breve c&=1- \frac{b}{\sqrt{a^2+b^2}}\left( \sqrt{a^2+b^2}+r\right) \cdot b\reC(b^2-r^2)  \\
\notag
\breve d&=-\frac{b}{\sqrt{a^2+b^2}}\left( \sqrt{a^2+b^2}+r\right).
\end{align}

(a) Consequently,
\begin{align}
\hat a^2+\hat b^2- \breve c^2-\breve d^2&=\hat a^2+\hat b^2- \hat c_{\mathrm R}^2-\hat d_{\mathrm R}^2=\hat a^2+\hat b^2- \hat c_{\mathrm L}^2-\hat d_{\mathrm L}^2=\notag\\
&=-\left(\frac{b}{\sqrt{a^2+b^2}}\left(\sqrt{a^2+b^2}+r\right)
\sqrt{\reD(b^2-r^2)}\right)^2 \leq 0,
\plabel{eq:hyprad}
\end{align}
with equality only if $b=0$.
In particular,
\begin{equation}
\breve c^2+\breve d^2= \hat c_{\mathrm R}^2+\hat d_{\mathrm R}^2= \hat c_{\mathrm L}^2+\hat d_{\mathrm L}^2>0.
\plabel{eq:ckbrad}
\end{equation}

(b) If $a>0$ and $b=0$, then $(\hat a,\hat b, \breve c,\breve d)=(1,0,\cos\psi,\sin\psi)$, thus
\[\mathrm{MR}(A)=\mathrm{ML}(A)=\Id_2+\cos\psi\tilde J+\sin\psi\tilde K.\]
If $a<0$ and $b=0$, then $(\hat a,\hat b, \breve c,\breve d)=(-1,0,\cos\psi,\sin\psi)$, thus
\[\mathrm{MR}(A)=\mathrm{ML}(A)=-\Id_2+\cos\psi\tilde J+\sin\psi\tilde K.\]

(c)
Restricted to $A\in\mathcal S$, and $b\neq0$, and to the level set  $\ujnorma=\sqrt{a^2+b^2}+r>0$,  the maps
$\mathrm{MR(A)}$ and $\mathrm{ML(A)}$ are injective.

(d) We can also consider the correspondence induced by conjugation by orthogonal matrices, that is the map
\[\breve{\mathrm M}: (a,b,r)\mapsto \left(\hat a,\hat b,\sqrt{\breve c^2+\breve d^2}\right),\]
where the domain is determined by the restrictions \eqref{eq:dom1} and \eqref{eq:dom2}.

Restricted to the $b^2-r^2<\pi^2$, and $b\neq0$, and to the level set  $\ujnorma=\sqrt{a^2+b^2}+r>0$,
the map $\breve{\mathrm M}$ is injective.

\begin{proof}
(o), (a), and (b) are straightforward computations, cf. Lemma \ref{lem:MRL}.
Only \eqref{eq:ckbrad} requires a particular argument:
By \eqref{eq:hyprad}, we see that $\breve c=\breve d=0$ implies $\hat a= \hat b=0$; however
$\hat a=\breve d=0$ is in contradiction to $a^2+b^2>0$.
Due to the fibration property with respect to conjugation be rotation matrices, (d) will imply (c).
Regarding (d): Assume we have $(\hat a, \hat b,\sqrt{\breve c^2+\breve d^2} )$ and $\ujnorma$ given.
Note that, due to $b^2-r^2<\pi^2$, we know not only $\sgn a=\sgn \hat a$, but  $\sgn b=\sgn\hat b$.
Hence, using $\hat a=\frac{a}{\sqrt{a^2+b^2}}$, the value of $\frac{b}{\sqrt{a^2+b^2}}$ can be recovered.
Using \eqref{eq:hyprad}, we can compute $\sqrt{\reD(b^2-r^2)} $.
As $b^2-r^2<\pi^2$, we can recover $b^2-r^2$.
Then, from $\hat b$, we can recover $r$.
Thus we also know $b^2$. Using $\sgn b$, we can deduce the value of $b$.
As $\frac{a}{\sqrt{a^2+b^2}}:\frac{b}{\sqrt{a^2+b^2}} $  is already known, the value of $a$ can also be recovered.
\end{proof}
\end{lemma}

Motivated by Lemma \ref{lem:EqMin} later, in the context of Lemma \ref{lem:momentkif}, we also define the
normalized expressions
\begin{equation}
\mathrm{MR}^{\mathrm{CKB}}(A)=\frac1{\hat c_{\mathrm R}^2+\hat d_{\mathrm R}^2}\mathrm{MR}(A),\qquad
\mathrm{ML}^{\mathrm{CKB}}(A)=\frac1{\hat c_{\mathrm L}^2+\hat d_{\mathrm L}^2}\mathrm{ML}(A)
\plabel{eq:MRLCKB}
\end{equation}
and
\begin{equation}
\breve{\mathrm M}^{\mathrm{CKB}}(a,b,r)\equiv(\hat a^{\mathrm{CKB}},\hat b^{\mathrm{CKB}} )=\left(\frac{\hat a}{\sqrt{\breve c^2+\breve d^2 }},\frac{\hat b}{\sqrt{\breve c^2+\breve d^2 }} \right).
\plabel{eq:MCKB}
\end{equation}
We are looking for a generalization of Lemma \ref{lem:momentkif}(c)(d) to the normalized expressions.
In that matter, $\breve{\mathrm M}^{\mathrm{CKB}}$ will play crucial role.
The notation refers to the (asymptotically closed) CKB model, as
$\breve{\mathrm M}^{\mathrm{CKB}}$ is supposed to take values in it.
(Indeed, by \eqref{eq:hyprad}, we obtain interior points in the CKB model for $b\neq0$, and asymptotical points for $b=0$.)

For a moment, let us consider the variant
which uses the flattened hyperboloid model given by
\begin{equation}
\breve{\mathrm M}^{\mathrm{HP}}(a,b,r)\equiv
(\hat a^{\mathrm{HP}},\hat b^{\mathrm{HP}} )=\left(\frac{\hat a}{\sqrt{-\hat a^2-\hat b^2+\breve c^2+\breve d^2 }},\frac{\hat b}{\sqrt{-\hat a^2-\hat b^2+\breve c^2+\breve d^2 }} \right).\plabel{eq:HPext0}
\end{equation}
It is not defined for $b=0$ (where the denominators vanish), but this is only a minor annoyance.
Using \eqref{eq:crux}, we can write
\begin{equation}
\breve{\mathrm M}^{\mathrm{HP}}(a,b,r)=
\left(\frac{a}{|b|}\cdot\frac{\dfrac{1}{\sqrt{\reD(b^2-r^2)}}}{\sqrt{a^2+b^2}+r },
(\sgn b)\cdot\left(
\frac{\dfrac{1}{\sqrt{\reD(b^2-r^2)}} }{\sqrt{a^2+b^2}+r}+r \frac{\reC(b^2-r^2)}{\sqrt{\reD(b^2-r^2)}}
\right)\right).
\plabel{eq:HPext}\end{equation}
Using Lemma \ref{lem:DTH}, we can write it as
\begin{equation*}
\breve{\mathrm M}^{\mathrm{HP}}(a,b,r)=
\left(\frac{a}{|b|}\cdot\frac{\mathcal E(b^2-r^2)}{\sqrt{a^2+b^2}+r },
(\sgn b)\cdot\left(
\frac{\mathcal E(b^2-r^2) }{\sqrt{a^2+b^2}+r}+r \mathcal F(b^2-r^2)
\right)\right).
\end{equation*}
This form shows that the expression above is unexpectedly well-defined
for $(a,b,r)\in \mathbb R\times \mathbb R\times [0,+\infty)$ as long as $b\neq0$.
This also applies to $\breve{\mathrm M}^{\mathrm{CKB}}$, as HP can simply written back to CKB.

%\snewpage
Now, we show that $\breve{\mathrm M}^{\mathrm{CKB}}$ extends to the case $b=0$ if we
make a suitable blow-up in $b=0$.
This means that instead of trying the domain
\[\breve{\mathcal M}=\{(a,b,r)\in \mathbb R\times \mathbb R\times [0,+\infty)  \},\]
we will consider the domain
\[\widehat{\mathcal M} =\{(s, r, \theta)\in [0,+\infty)\times [0,+\infty) \times (\mathbb R\modu2\pi)  \}\]
where the canonical correspondence is
given by
\[a=s\cos\theta,\qquad b=s\sin\theta\qquad r=r.\]
(This, somewhat colloquially, describes the map $\widehat{\mathcal M}\rightarrow \breve{\mathcal M}$.)

Thus, we can define $\widehat{\mathrm M}{}^{\mathrm{CKB}}(s,\theta,r)=\breve{\mathrm M}^{\mathrm{CKB}}(s\cos\theta,s\sin\theta,r)$
as long as $s\sin\theta\neq0$.

\begin{lemma}\plabel{lem:momentblow}
 (a)
$\widehat{\mathrm M}^{\mathrm{CKB}}$ extends by the formula
\begin{equation*}
 \widehat{\mathrm M}^{\mathrm{CKB}}(s,r,\theta)=
 \left(\frac{\mathcal A}{\sqrt{\mathcal A^2+\mathcal B^2+\mathcal G^2 }},
 \frac{\mathcal B}{\sqrt{\mathcal A^2+\mathcal B^2+\mathcal G^2  }}\right)
\end{equation*}
with
\[\mathcal A\equiv (\cos\theta)\cdot\mathcal E((s\sin\theta)^2-r^2),\]
and
\[\mathcal B\equiv (\sin\theta)\cdot\underbrace{\left(
\mathcal E((s\sin\theta)^2-r^2)+r(s+r) \mathcal F((s\sin\theta)^2-r^2)
\right)}_{\mathcal B_0\equiv},\]
and
\[\mathcal G\equiv (\sin\theta)\cdot\underbrace{(s+r)}_{\mathcal G_0\equiv},\]
smoothly to the domain $\widehat{\mathcal M}$, i. e. to the domain subject to the conditions
\begin{equation}0\leq s,\qquad 0\leq r,\qquad \theta\in\mathbb R\modu 2\pi.\plabel{eq:BlowDomainCondPre2}\end{equation}

(b) If $\sin\theta= 0$, then $\widehat{\mathrm M}^{\mathrm{CKB}}(s,r,\theta)=(\cos\theta,0)$ (which is $(1,0)$ or $(-1,0)$).

If $s=r=0$, then  $\widehat{\mathrm M}^{\mathrm{CKB}}(0,0,\theta)=(\cos\theta,\sin\theta)$.

However, if $\sin\theta\neq0$ and $s+r>0$, then $\widehat{\mathrm M}^{\mathrm{CKB}}(s,r,\theta)$
is in the open unit disk.

(c) $\widehat{\mathrm M}^{\mathrm{HP}}$ extends by the formula
\begin{equation*}
 \widehat{\mathrm M}^{\mathrm{HP}}(s,r,\theta)=
 \left(\frac{\mathcal A}{|\mathcal G| },\frac{\mathcal B}{|\mathcal G|  }\right)=
 \left(\frac{\mathcal A}{|\mathcal G| },(\sgn\sin\theta)\cdot\frac{\mathcal B_0}{\mathcal G_0  }\right),
\end{equation*}
except it is formally not defined for $s+r=0$ or $\sin\theta=0$.

\begin{proof}
(a) It is sufficient to prove that $\mathcal A^2+\mathcal B^2+\mathcal G^2 $ never vanishes.
If $s=r=0$, then $\mathcal A^2+\mathcal B^2+\mathcal G^2 =\mathcal E(0)^2=3$.
If $s+r>0$, then vanishing requires $\sin\theta=0$; thus $(\cos \theta)^2=1$;
consequently $\mathcal A^2+\mathcal B^2+\mathcal G^2 =\mathcal E(-r^2)^2\geq\mathcal E(0)^2=3$, a contradiction.

(b), (c): Direct computation.
\end{proof}
\end{lemma}

\snewpage
Note that the decomposition $\mathrm M_2(\mathbb R)=\mathcal S^{\mathrm{nn}}\cup \eth^{\mathrm{par}}\mathcal S\cup\ldots$
is invariant for conjugation by rotation matrices, thus it descends to a decomposition
$\breve{\mathcal M}=X^{\mathrm{nn}}\cup \eth^{\mathrm{par}}X\cup\ldots$.
By the blow up map, this induces a decomposition
$\widehat{\mathcal M}=\widehat{X}^{\mathrm{nn}}\cup \eth^{\mathrm{par}}\widehat{X}\cup\ldots$.
We recapitulate the situation for $\widehat{\mathcal M}$. We have

$\bullet$ the subset  $\widehat X^{\mathrm{nn}}$ with $0<r$, $0<s$, $\sin\theta\neq0$, $(s\sin\theta)^2-r^2<\pi$
 as the regular (or non-normal) interior;

$\bullet$ the subset  $\eth^{\mathrm{par}}\widehat X$ with $0<r$, $0<s$, $\sin\theta=0$, $(s\sin\theta)^2-r^2<\pi$
as the Schur-parabolic (or self-adjoint) pseudo-boundary.

$\bullet$ the subset $\partial^{0}\widehat X$  with $r=0$, $s=0$ as the zero boundary;

$\bullet$ the subset $\partial^{\mathrm{hyp}}\widehat X$  with $s=0$, $r>0$
as the Schur-hyperbolic (or conform-reflexion) boundary;

$\bullet$ the subset  $\partial^{\mathrm{ell}}\widehat X$  with $s>0$, $r=0$, $(s\sin\theta)^2-r^2<\pi$
as the Schur-elliptic (or conform-rotational) boundary;

$\bullet$ the subset $\partial^{\mathrm{dell}}\widehat X$ with $r=0$, $(s\sin\theta)^2-r^2=\pi^2$
as the Schur-elliptic degenerate boundary;

$\bullet$ the subset $\partial^{\mathrm{dnn}}\widehat X$ with $r>0$, $(s\sin\theta)^2-r^2=\pi^2$
as the non-normal degenerate boundary;

$\bullet$ the subset $\widehat X^{\mathrm{ext}}$ with $(s\sin\theta)^2-r^2>\pi^2$
as the exterior set.

Thus
\[\widehat {\mathcal M}=
\widehat X^{\mathrm{nn}}
\,\,\dot\cup\,\, \eth^{\mathrm{par}}\widehat X
\,\,\dot\cup\,\, \partial^{0}\widehat X
\,\,\dot\cup\,\, \partial^{\mathrm{hyp}}\widehat X
\,\,\dot\cup\,\, \partial^{\mathrm{ell}}\widehat X
\,\,\dot\cup\,\, \partial^{\mathrm{dell}}\widehat X
\,\,\dot\cup\,\, \partial^{\mathrm{dnn}}\widehat X
\,\,\dot\cup\,\,  \widehat X^{\mathrm{ext}}.
\]

Again, certain components can be decomposed naturally:
$\partial^{\mathrm{hyp}}\widehat X$ can be decomposed to
$\partial^{\mathrm{hyp1}}\widehat X$ with $\sin\theta=0$,
and to  $\partial^{\mathrm{hyp*}}\widehat X$ with $\sin\theta\neq0$.
$\partial^{\mathrm{ell}}\widehat X$ can be decomposed to
$\partial^{\mathrm{ell1}}\widehat X$ with $\sin\theta=0$,
and to  $\partial^{\mathrm{ell*}}\widehat X$ with $\sin\theta\neq0$.
$\partial^{\mathrm{dell}}\widehat X$ can be decomposed to
$\partial^{\mathrm{dell0}}\widehat X$ with $\cos\theta=0$,
and to  $\partial^{\mathrm{dell*}}\widehat X$ with $\cos\theta\neq0$.
$\partial^{\mathrm{dnn}}\widehat X$ can be decomposed to
$\partial^{\mathrm{dnn0}}\widehat X$ with $\cos\theta=0$,
and to  $\partial^{\mathrm{dnn*}}\widehat X$ with $\cos\theta\neq0$.
$\widehat X^{\mathrm{ext}}$ can be decomposed to
the Schur-elliptic part
$\widehat X^{\mathrm{ellext}}$  with $r=0$,
and to the non-normal part
$\widehat X^{\mathrm{nnext}}$  with $r>0$.
(We could further discriminate $\widehat X^{\mathrm{ext}}$ by
whether $(s\sin\theta)^2-r^2\in\{k^2\pi^2\,:\,k\in\mathbb N\setminus\{0\}\}$
holds or not.)

Note that $\partial^{\mathrm{hyp}}\widehat X$ is nontrivially blown up from $\partial^{\mathrm{hyp}}X$,
but otherwise $\widehat{\mathcal M}$ and $\breve{\mathcal M}$ are quite similar.
\snewpage
As long as we avoid $(\pm1,0)^{\mathrm{CKB}}$,
we can pass from CKB to ACKB without trouble.
In the context of Lemma \ref{lem:momentblow}, the critical case is when $\sin\theta=0$.
Technically, however, this requires another (very simple) blow-up in domain.

Let us consider the domain
\begin{multline*}
\widetilde{\mathcal M} =\{(s, r, \theta ,\sigma)\in
[0,+\infty)\times [0,+\infty) \times\\ \times \left(
([0,\pi]\modu2\pi)\times\{+1\}
\cup
([\pi,2\pi]\modu2\pi)\times\{-1\}
 \right)  \}
 \end{multline*}
where the canonical $\widetilde{\mathcal M}\rightarrow\widehat{\mathcal M}$ simply forgets $\sigma$.
Technically, the blow-up simply duplicates (cuts along)
$\partial^{\mathrm{hyp1}}\widehat X\cup \eth^{\mathrm{par}}\widehat X \cup \partial^{\mathrm{ell1}}\widehat X$.
This separates the connected $\widehat{\mathcal M}$ into the two components of $\widetilde{\mathcal M}$.

Now, $\widetilde{\mathrm M}^{\mathrm{ACKB}}(s,r,\theta,\sigma)$ is well-defined as long as $\sin\theta\neq0$.
The analogous statement to Lemma \ref{lem:momentblow} is

\begin{lemma}\plabel{lem:momentblow2}
 (a)
$\widetilde{\mathrm M}^{\mathrm{ACKB}}$ extends by the formula
\begin{equation*}
 \widetilde{\mathrm M}^{\mathrm{ACKB}}(s,r,\theta,\sigma)=
 \left(\arcsin\frac{\mathcal A}{\sqrt{\mathcal A^2+\mathcal B^2+\mathcal G^2 }},\sigma\cdot
 \frac{\mathcal B_0}{\sqrt{(\mathcal B_0)^2+(\mathcal G_0)^2  }}\right)
\end{equation*}

(b)
If $s=r=0$, then  $\widetilde{\mathrm M}^{\mathrm{ACKB}}(0,0,\theta,\sigma)=(\frac\pi2-\theta,\sigma)$.

However, if  $s+r>0$, then $\widetilde{\mathrm M}^{\mathrm{ACKB}}(s,r,\theta,\sigma)$
is in $[-\frac\pi2,\frac\pi2]\times (-1,1)$.

(c) $\widetilde{\mathrm M}^{\mathrm{AHP}}$ extends by the formula
\begin{equation*}
 \widetilde{\mathrm M}^{\mathrm{AHP}}(s,r,\theta,\sigma)=
 \left(\arctan\frac{\mathcal A}{|\mathcal G| },\sigma\cdot\frac{\mathcal B_0}{\mathcal G_0  }\right),
\end{equation*}
with values in $[-\frac\pi2,\frac\pi2]\times \mathbb R $,
except it is formally not defined for $s+r=0$.

\begin{proof}
 Direct computation.
\end{proof}
\end{lemma}

\snewpage
Let us use the notation $\ujnorma=s+r$ (which refers to the norm).
Ultimately, we will be interested in the case when $\ujnorma>0$ and  $\ujnorma$ is fixed.

For the sake of visualization, in Figure \ref{fig:figBA04}, we show
a cross-section of $\widetilde{\mathcal M}|_{\sigma=+1}$  for $\ujnorma=5\pi/3$.

\begin{figure}[H]
  \centering
  \begin{subfigure}[b]{.7\linewidth}
  \begin{overpic}[scale=.5]{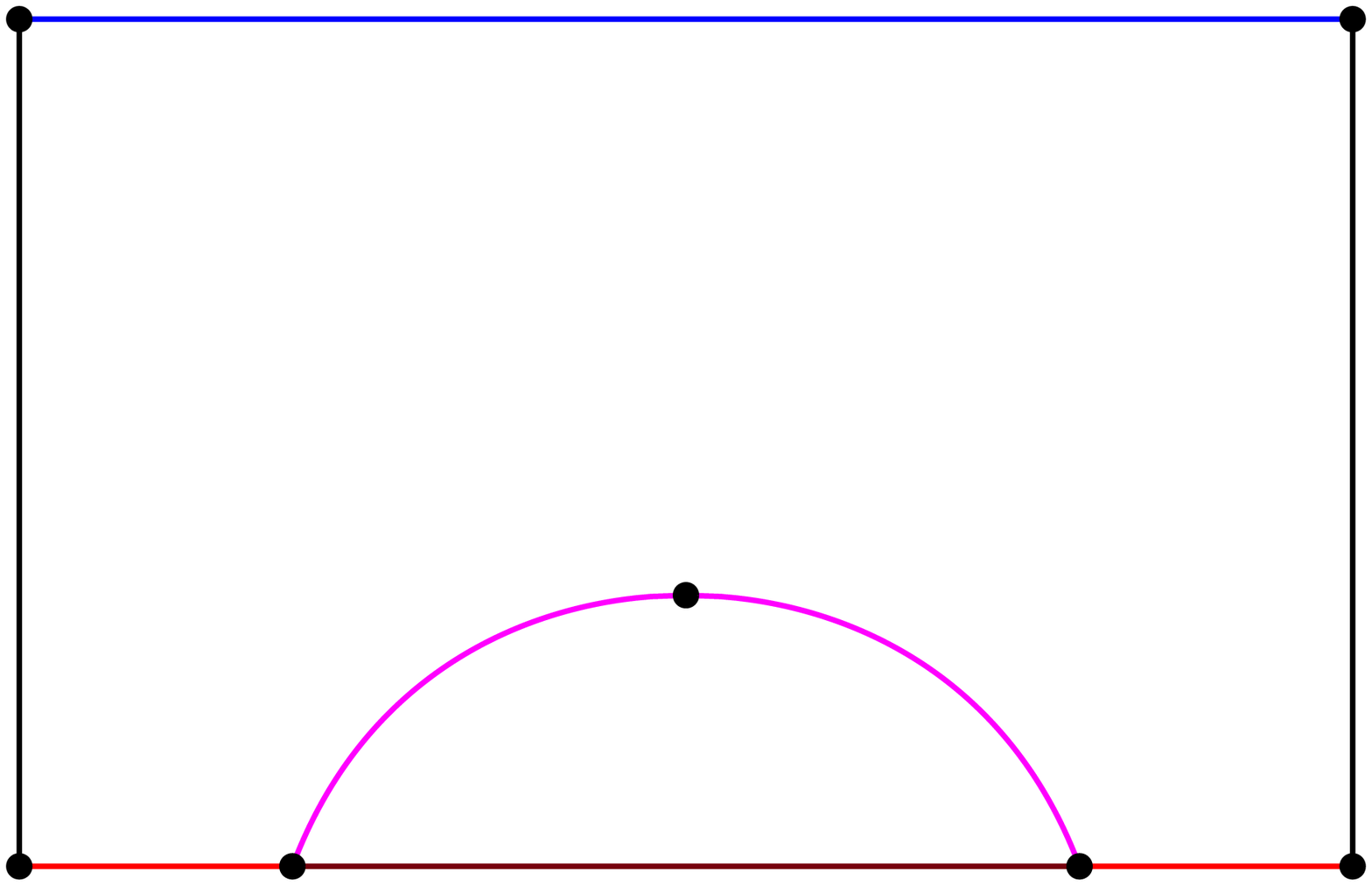}
  \put(100,18){$\partial^{\mathrm{ell1}}\widetilde{X}$}
  \put(100,48){$\eth^{\mathrm{par}}\widetilde{X}$}
  \put(100,78){$\partial^{\mathrm{hyp1}}\widetilde{X}$}
  \put(-12,18){$\partial^{\mathrm{ell1}}\widetilde{X}$}
  \put(-12,48){$\eth^{\mathrm{par}}\widetilde{X}$}
  \put(-12,78){$\partial^{\mathrm{hyp1}}\widetilde{X}$}
  \put(48,13){$\widetilde{X}^{\mathrm{ellext}}$}
  \put(48,27){$\widetilde{X}^{\mathrm{nnext}}$}
  \put(48,42){$\partial^{\mathrm{dnn0}}\widetilde{X}$}
  \put(48,60){$\widetilde{X}^{\mathrm{nn}}$}
  \put(48,82){$\partial^{\mathrm{hyp*}} \widetilde{X}$}
  \put(69,34){$\partial^{\mathrm{dnn*}}\widetilde{X}$}
  \put(19,34){$\partial^{\mathrm{dnn*}}\widetilde{X}$}
  \put(69,13){$\partial^{\mathrm{dell*}}\widetilde{X}$}
  \put(84,13){$\partial^{\mathrm{ell*}}\widetilde{X}$}
  \put(22,13){$\partial^{\mathrm{dell*}}\widetilde{X}$}
  \put(7,13){$\partial^{\mathrm{ell*}}\widetilde{X}$}
  \end{overpic}
    \caption*{Fig. \ref{fig:figBA04}. Cross-section of $\widetilde{\mathcal M}|_{\sigma=+1}$  for $\ujnorma=5\pi/3$ }
  \end{subfigure}
\phantomcaption
\plabel{fig:figBA04}
\end{figure}

In general, the situation is similar for $\ujnorma>\pi$.
For $\ujnorma=\pi$, the degenerate boundary shrinks only to  $\partial^{\mathrm{dell0}}\widetilde{X}$;
and there is no degenerate boundary for $\pi<\ujnorma$.
Now, $\widetilde{\mathcal M}$ contains two such bands,
and $\widehat{\mathcal M}$ is two such band glued together
(all meant for a fixed $\ujnorma>0$).
In this setting $\partial^0 \widetilde{\mathcal M}$ contains two segments,
and $\partial^0 \widehat{\mathcal M}$ is a circle.

Another sort of blow-up which affects the variables $s,r$ is given by passing to the
coordinates $\ujnorma, t$ such that $t\in[0,1]$, and
\begin{equation}
s=\ujnorma t,\qand r=\ujnorma(1-t).
\plabel{eq:tN}
\end{equation}
In merit, as a blow-up, this affects only the case $s=r=0$ of the zero-boundary
(where it is slightly advantageous).
As such, its use is marginal for us.
However, \eqref{eq:tN} is retained as a practical change of coordinates, which is
particularly useful if $\ujnorma$ is kept fixed.

\snewpage
\begin{lemma}
\plabel{lem:momentspec}

In the context of Lemma \ref{lem:momentblow} and Lemma \ref{lem:momentblow2},
the extended actions on some boundary pieces  (restricted for a fixed $\ujnorma>0$, with $\ujnorma=s+r$) are given as follows:

(a) If $s=0$ (so $r=\ujnorma$), i. e. we consider the Schur-hyperbolic boundary $\partial^{\mathrm{hyp}}\widehat X$, then
\[ (\hat a,\hat b, \breve c,\breve d)=\left(\cos\theta , (\coth \ujnorma)\ujnorma\sin\theta,1 ,-\ujnorma\sin\theta\right) .\]
\[\breve c^2+\breve d^2=1+\left(\ujnorma\sin\theta\right)^2.\]
\[\hat a^2+\hat b^2-\breve c^2-\breve d^2=(\sin\theta)^2\left(\left(\frac{\ujnorma}{\sinh\ujnorma}\right)^2-1\right).\]
Here $\theta$ is arbitrary.

The normalized expressions extend as
\[\widehat{\mathrm M}^{\mathrm{CKB}}( 0,\ujnorma,\theta)=
\left(\frac{\cos\theta}{\sqrt{1 +\left(\ujnorma\sin\theta\right)^2}} , \frac{(\coth \ujnorma)\ujnorma\sin\theta}{\sqrt{1 +\left(\ujnorma\sin\theta\right)^2 }}\right)
\]
and (if $\sin\theta\neq0$,)
\[ \widehat{\mathrm M}^{\mathrm{HP}}( 0,\ujnorma,\theta)=
\left(\frac{|\cot\theta|}{\sqrt{1-\left(\dfrac{\ujnorma}{\sinh\ujnorma}\right)^2}}\sgn\cos\theta , \frac{(\coth \ujnorma)\ujnorma}{\sqrt{1-\left(\dfrac{\ujnorma}{\sinh\ujnorma}\right)^2}}\sgn\sin\theta\right).
\]
These curves (in $\theta$) are  two-sided $h$-hypercycles ($h$-distance lines).

(b) If $r=0$  (so $s=\ujnorma$), i. e. we consider the purely elliptic boundary $\partial^{\mathrm{ell}}\widehat X$, then
\[ (\hat a,\hat b, \breve c,\breve d)=\left(\cos\theta , \sin\theta,(\ujnorma\sin\theta) \cot(\ujnorma\sin\theta) ,-\ujnorma\sin\theta\right) ,\]
\[\breve c^2+\breve d^2=\left(\frac{\ujnorma\sin\theta}{\sin\left(\ujnorma\sin\theta\right)}\right)^2,\]
\[\hat a^2+\hat b^2-\breve c^2-\breve d^2=1-\left(\frac{\ujnorma\sin\theta}{\sin\left(\ujnorma\sin\theta\right)}\right)^2 .\]
Here the domain restriction can be expressed as
\[|\sin\theta|<\frac\pi\ujnorma.\]

The normalized expressions extend as
\[\widehat{\mathrm M}^{\mathrm{CKB}}( \ujnorma,0,\theta)=
\left(\frac{\sin\left(\ujnorma\sin\theta\right)} {\ujnorma}\cot\theta,\frac{\sin\left(\ujnorma\sin\theta\right)} {\ujnorma} \right),
\]
and (if $\sin\theta\neq0$,)
\[ \widehat{\mathrm M}^{\mathrm{HP}}( \ujnorma,0,\theta)=
\left( \hat a^{\mathrm{HP}},\hat b^{\mathrm{HP}}  \right)=
\left(\frac{\cos\theta}{\sqrt{\left(\dfrac{\ujnorma\sin\theta}{\sin\left(\ujnorma\sin\theta\right)}\right)^2-1}},
\frac{\sin\theta}{\sqrt{\left(\dfrac{\ujnorma\sin\theta}{\sin\left(\ujnorma\sin\theta\right)}\right)^2-1}}
\right).
\]

These curves (in $\theta$) in the CKB model are radially contracted images of (possibly two open pieces of) the unit circle.
(They give  quasi Cassini ovals.)

(c) If $(s\sin\theta)^2-(\ujnorma-s)^2=\pi^2$, i. e. we consider
$ \partial^{\mathrm{dell}}\widehat X \cup \partial^{\mathrm{dnn}}\widehat X$
(in particular $\ujnorma\geq\pi$), then
\[|\sin\theta|=\frac{\sqrt{\pi^2+(\ujnorma-s)^2}}{s}.\]

\begin{commentx}
Or, conversely
\[s=\frac{\ujnorma^2+\pi^2}{\ujnorma+\sqrt{\ujnorma^2(\sin\theta)^2-\pi^2(\cos\theta)^2}}
,\qquad\text{and}\qquad
r=\frac {\ujnorma^{2} \left( \sin \theta \right) ^{2}-{\pi }^{2}}{\ujnorma \left( \sin \theta \right) ^{2}+
\sqrt {\ujnorma^{2} \left( \sin \theta  \right) ^{2}-{\pi }
^{2} \left( \cos \theta   \right) ^{2}}}\]
\end{commentx}

The domain restrictions are expressed as
\[|\sin\theta|\geq\frac\pi\ujnorma,\qquad\text{or}\qquad
s\geq \frac{\ujnorma^2+\pi^2}{2\ujnorma},\qquad\text{or}\qquad
r\leq \frac{\ujnorma^2-\pi^2}{2\ujnorma}
.\]

The extended actions are given by
\[%\begin{align*}
\widehat{\mathrm M}^{\mathrm{CKB}}( s,r,\theta)
=\left(0,\frac{\frac{\ujnorma-s}\pi}{\sqrt{1+\left(\frac{\ujnorma-s}\pi\right)^2}}\sgn\sin\theta\right)
=\left(0,\frac{\frac{r}\pi}{\sqrt{1+\left(\frac{r}\pi\right)^2}}\sgn\sin\theta\right),
\]%\end{align*}
\begin{commentx}
\[=\left(0,{\frac {\ujnorma\sqrt {{\ujnorma}^{2} \left( \sin  \theta  \right) ^{2}-{\pi }^{2} \left( \cos \theta \right) ^{2}}-{\pi }^{2}}{
\left( {\ujnorma}^{2}+{\pi }^{2} \right) \left(\sin  \theta\right) }}\right), \]
\end{commentx}
and
\[%\begin{align*}
\widehat{\mathrm M}^{\mathrm{HP}}( s,r,\theta)
=\left(0,\frac{\ujnorma-s}\pi\sgn\sin\theta\right)
=\left(0,\frac{r}\pi\sgn\sin\theta\right).
\]%\end{align*}
\begin{commentx}
\[=\left(0,\frac1\pi{\frac {\ujnorma^{2} \left( \sin \theta \right) ^{2}-{\pi }^{2}}{\ujnorma \left( \sin \theta \right) ^{2}+
\sqrt {\ujnorma^{2} \left( \sin \theta  \right) ^{2}-{\pi }
^{2} \left( \cos \theta   \right) ^{2}}}}\sgn\sin\theta\right).\]
\end{commentx}

In particular, the image of $ \partial^{\mathrm{dell}}\widehat X \cup \partial^{\mathrm{dnn}}\widehat X$ is
\[\text{the segment connecting  $\left(0,\frac{\ujnorma^2-\pi^2}{\ujnorma^2+\pi^2}\right)$ and $\left(0,-\frac{\ujnorma^2-\pi^2}{\ujnorma^2+\pi^2}\right)$ in the CKB model},\]
and
\[\text{the segment connecting  $\left(0,\frac{\ujnorma^2-\pi^2}{2\ujnorma\pi}\right)$ and $\left(0,-\frac{\ujnorma^2-\pi^2}{2\ujnorma\pi}\right)$ in the HP model}.\]

(For $\ujnorma=\pi$, this is just the origin, which comes from $ \partial^{\mathrm{dell0}}\widehat X$.
For $\ujnorma>\pi$,  the endpoints come from $ \partial^{\mathrm{dnn0}}\widehat X$,
the origin comes from $ \partial^{\mathrm{dell*}}\widehat X$,
the intermediate points come from $ \partial^{\mathrm{dnn*}}\widehat X$.)

(d) If $\sin\theta=0$, i. e. we consider
$\partial^{\mathrm{hyp1}}\widehat X\cup \eth^{\mathrm{par}}\widehat X \cup\partial^{\mathrm{ell1}}\widehat X$,
or rather $\partial^{\mathrm{hyp1}}\widetilde X\cup \eth^{\mathrm{par}}\widetilde X \cup\partial^{\mathrm{ell1}}\widetilde X$,
then
\[\widehat{\mathrm M}^{\mathrm{ACKB}}( \ujnorma-r,r,\theta,\sigma)=
\left(\frac\pi2\cdot\sgn\cos\theta,\frac{\frac1\ujnorma\mathcal E(-r^2)+r\mathcal F(-r^2)}{\sqrt{1+\left(\frac1\ujnorma\mathcal E(-r^2)+r\mathcal F(-r^2)\right)^2}}
\cdot\sigma\right),\]
and
\[\widehat{\mathrm M}^{\mathrm{AHP}}( \ujnorma-r,r,\theta,\sigma)=
\left(\frac\pi2\cdot\sgn\cos\theta,\left(\frac1\ujnorma\mathcal E(-r^2)+r\mathcal F(-r^2)\right)
\cdot\sigma\right).\]

The absolute values of the seconds coordinates are strictly monotone increasing in $r$,
\[\text{with range
  $\left[\frac{\sqrt3}{\sqrt{3+\ujnorma^2}},\frac{\ujnorma\coth \ujnorma}{\sqrt{1+\ujnorma^2}}\right] $ in
   the ACKB model},\]
and
\[\text{with range
  $\left[\frac{\sqrt3}{\ujnorma},\frac{(\coth \ujnorma)\ujnorma}{\sqrt{1-\left(\dfrac{\ujnorma}{\sinh\ujnorma}\right)^2}}\right] $ in
   the AHP model}.\]
(The lower limits come from $\partial^{\mathrm{ell1}}\widetilde X$,
the upper limits come from $\partial^{\mathrm{hyp1}}\widetilde X$,
the intermediate values come from $\eth^{\mathrm{par}}\widetilde X$.
\begin{proof}
Direct computation.
\end{proof}
\end{lemma}
For the sake of visualization, in Figure \ref{fig:figBA06},
we show pictures about the range of
$\widehat{\mathrm M}^{\mathrm{CKB}}( s,r,\theta) $
restricted to $\overline X\equiv\widehat X^{\mathrm{nn}}\cup\eth^{\mathrm{par}}\widehat X
\cup \partial^{\mathrm{hyp}}\widetilde X \cup \partial^{\mathrm{ell}}\widetilde X
\cup\partial^{\mathrm{dell}}\widetilde X\cup\partial^{\mathrm{dnn}}\widetilde X
$
and to a fixed
$\ujnorma>0$,
with parameter lines induced by $s$ and $\theta$, for the cases $\ujnorma=\frac\pi2$, $\ujnorma=\frac{5\pi}6$,  $\ujnorma=\pi$,  $\ujnorma=\frac{5\pi}3$; and
the image of the relevant boundary pieces in the cases $\ujnorma=\frac{5\pi}3$, $\ujnorma=5\pi$.
Furthermore, in Figure \ref{fig:figBA05},
we also include a picture
about the corresponding range of
$\widetilde{\mathrm M}^{\mathrm{ACKB}}( s,r,\theta,\sigma) $
in the case $\ujnorma=5\pi/6$.

\begin{figure}[H]
  \centering
  \begin{subfigure}[b]{0.4\linewidth}
    \includegraphics[width=2.5in]{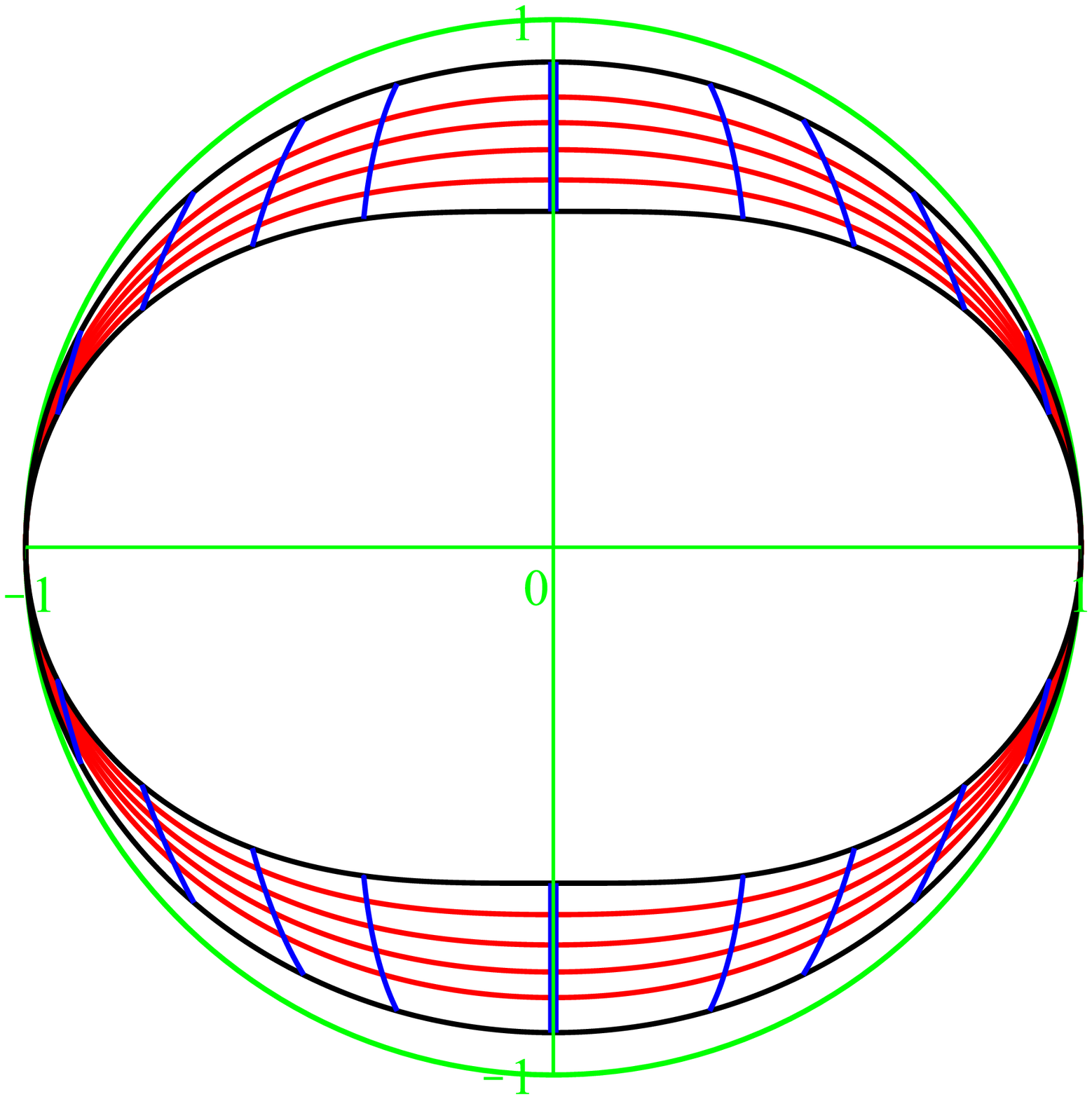}
    \caption*{Fig. \ref{fig:figBA06}(a) $\ujnorma=\pi/2$}
  \end{subfigure}
  \quad
  \begin{subfigure}[b]{0.4\linewidth}
    \includegraphics[width=2.5in]{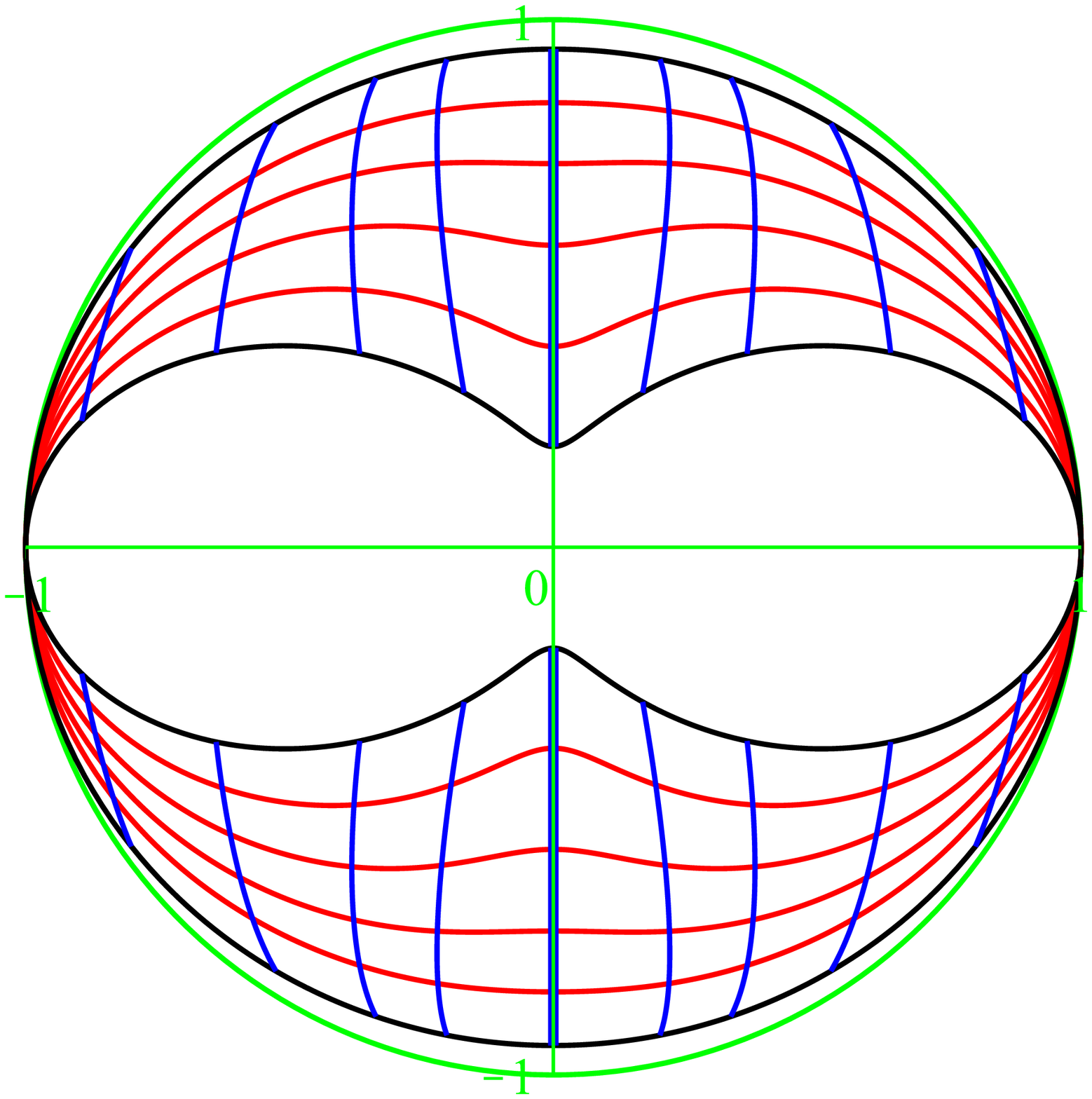}
    \caption*{\ref{fig:figBA06}(b) $\ujnorma=5\pi/6$  }
  \end{subfigure}
\phantomcaption
\plabel{fig:figBA06}
\end{figure}
\begin{figure}[H]
  \ContinuedFloat
  \begin{subfigure}[b]{0.4\linewidth}
    \includegraphics[width=2.5in]{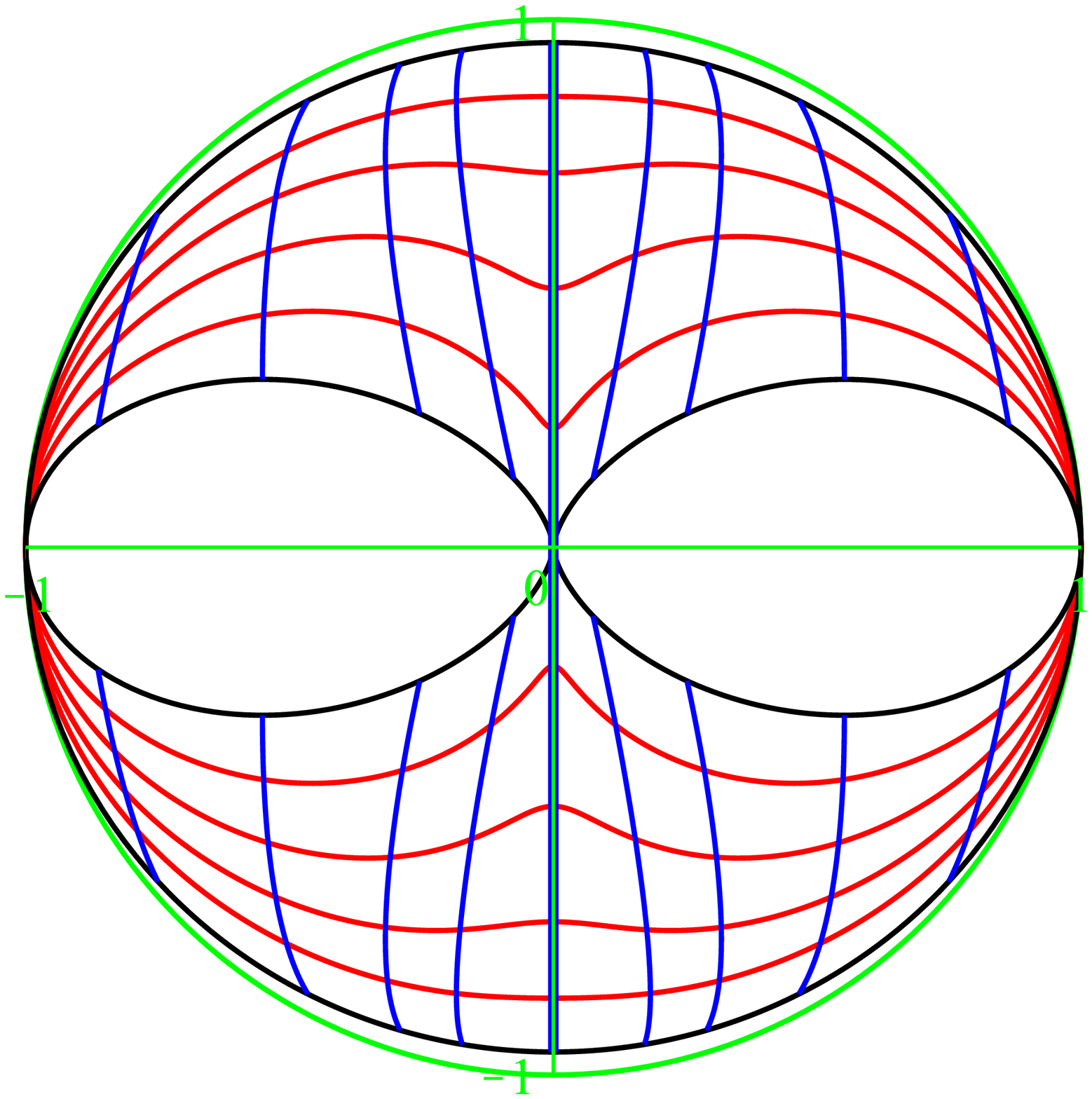}
    \caption*{\ref{fig:figBA06}(c) $\ujnorma=\pi$  }
  \end{subfigure}
  \quad
  \begin{subfigure}[b]{0.4\linewidth}
    \includegraphics[width=2.5in]{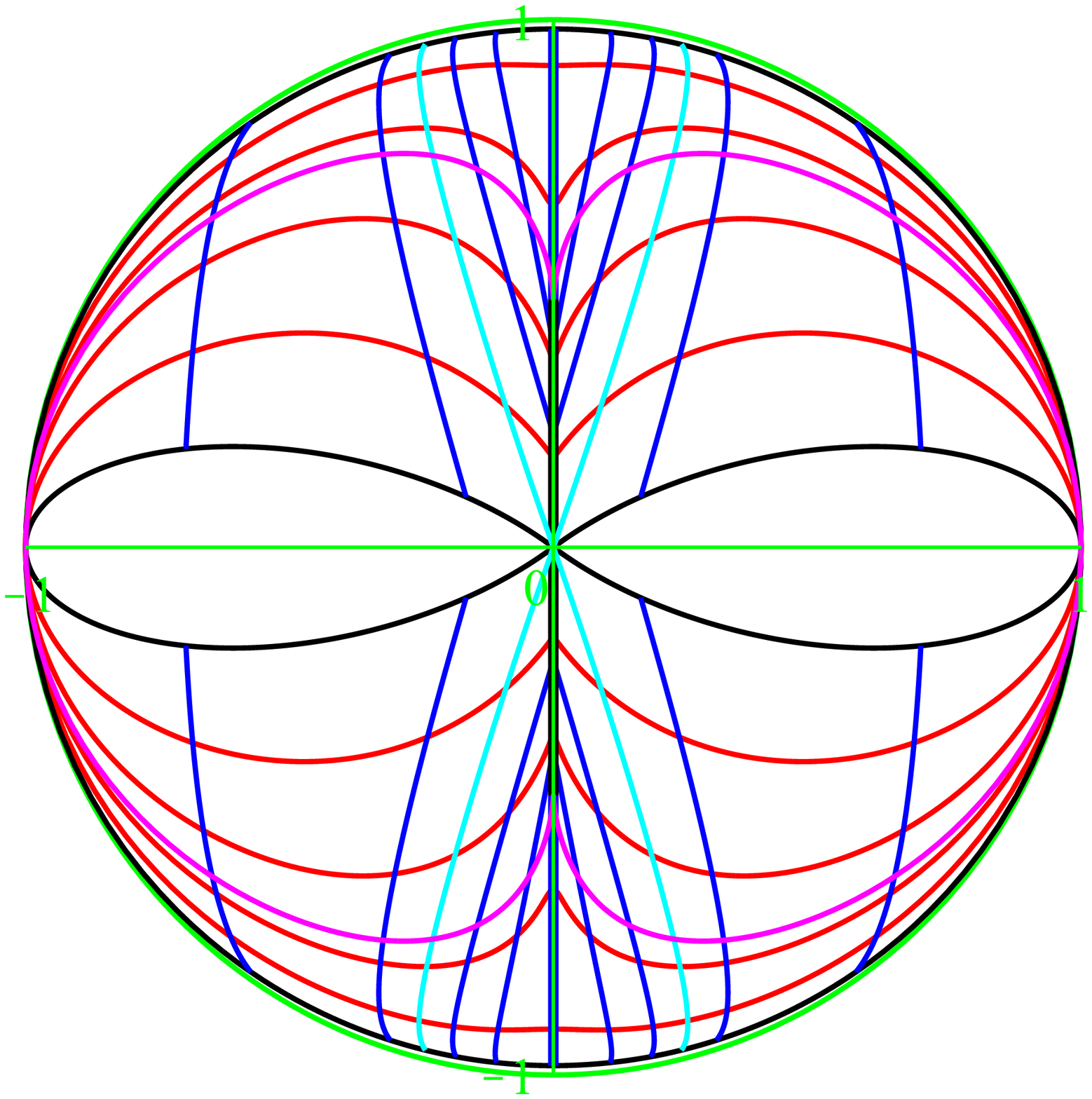}
    \caption*{\ref{fig:figBA06}(d) $\ujnorma=5\pi/3$  }
  \end{subfigure}
\phantomcaption
\end{figure}
~\\[-1cm]
\begin{figure}[H]
  \ContinuedFloat
   \begin{subfigure}[b]{0.4\linewidth}
    \includegraphics[width=2.5in]{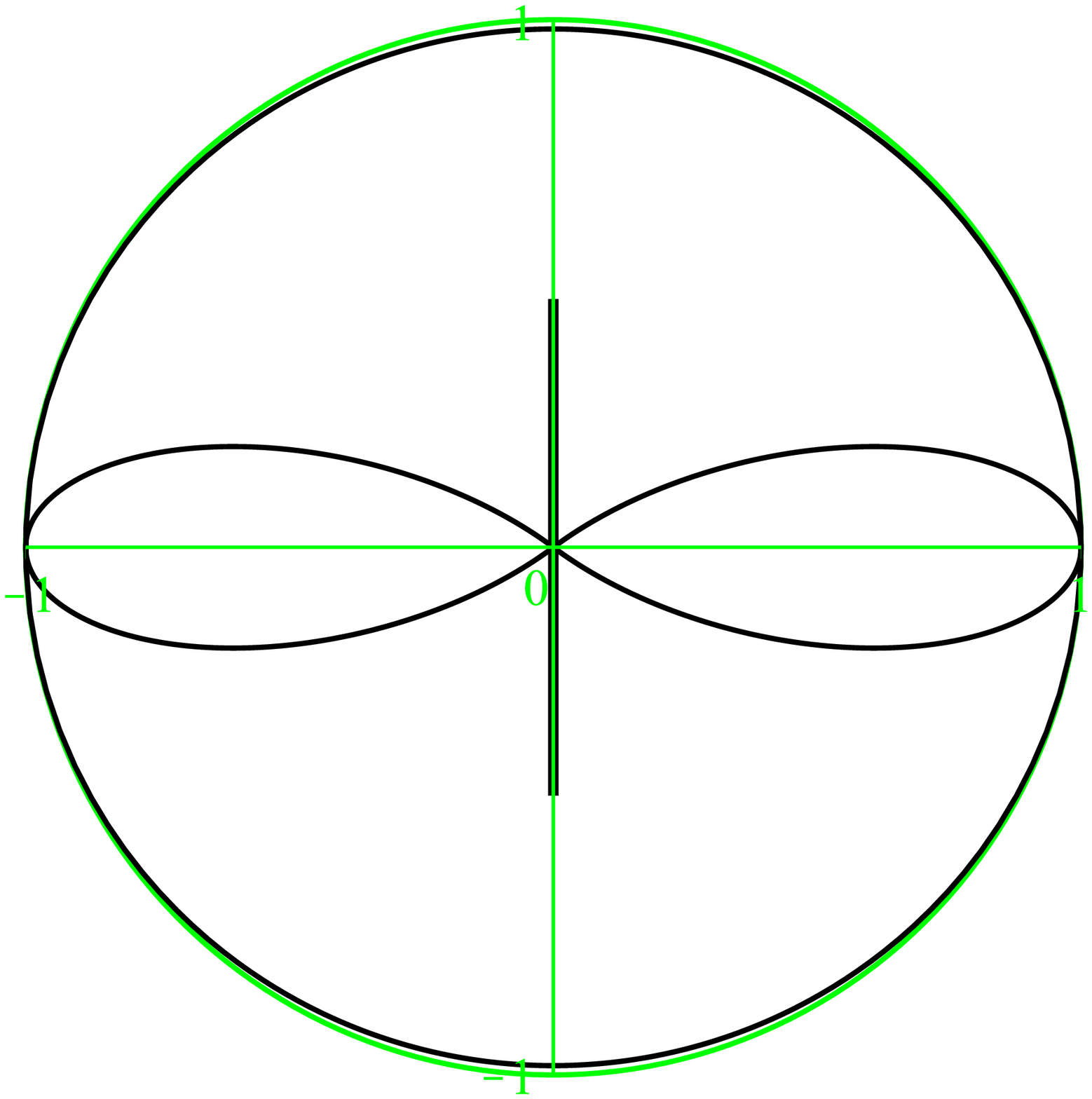}
    \caption*{\ref{fig:figBA06}(e) $\ujnorma=5\pi/3$ (only boundary)}
  \end{subfigure}
  \quad
  \begin{subfigure}[b]{0.4\linewidth}
    \includegraphics[width=2.5in]{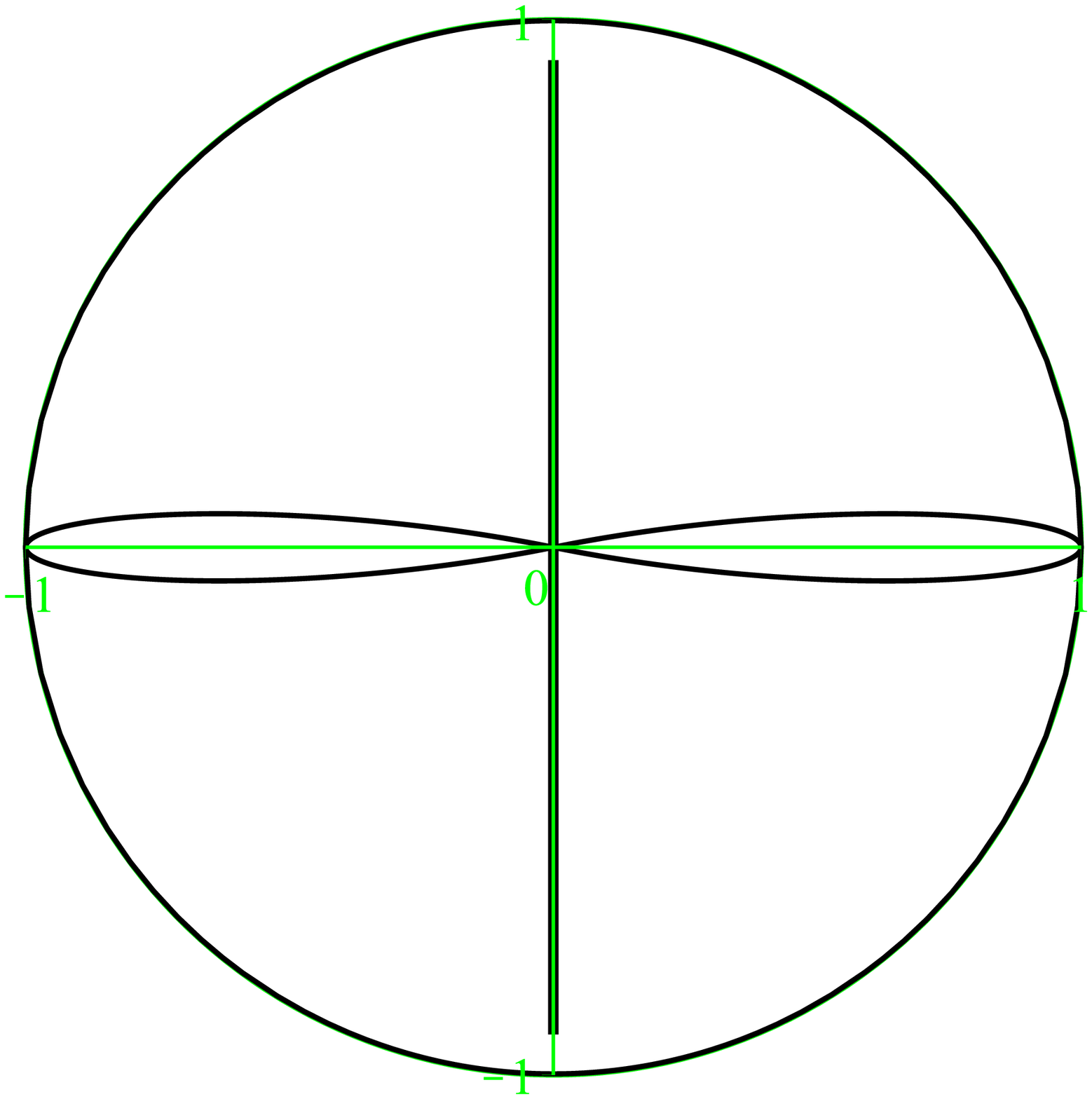}
    \caption*{\ref{fig:figBA06}(f) $\ujnorma=5\pi$ (only boundary)}
  \end{subfigure}
\phantomcaption
\end{figure}
~\\[-1cm]
\begin{figure}[H]
  \centering
  \begin{subfigure}[b]{.7\linewidth}
    \includegraphics[width=4in]{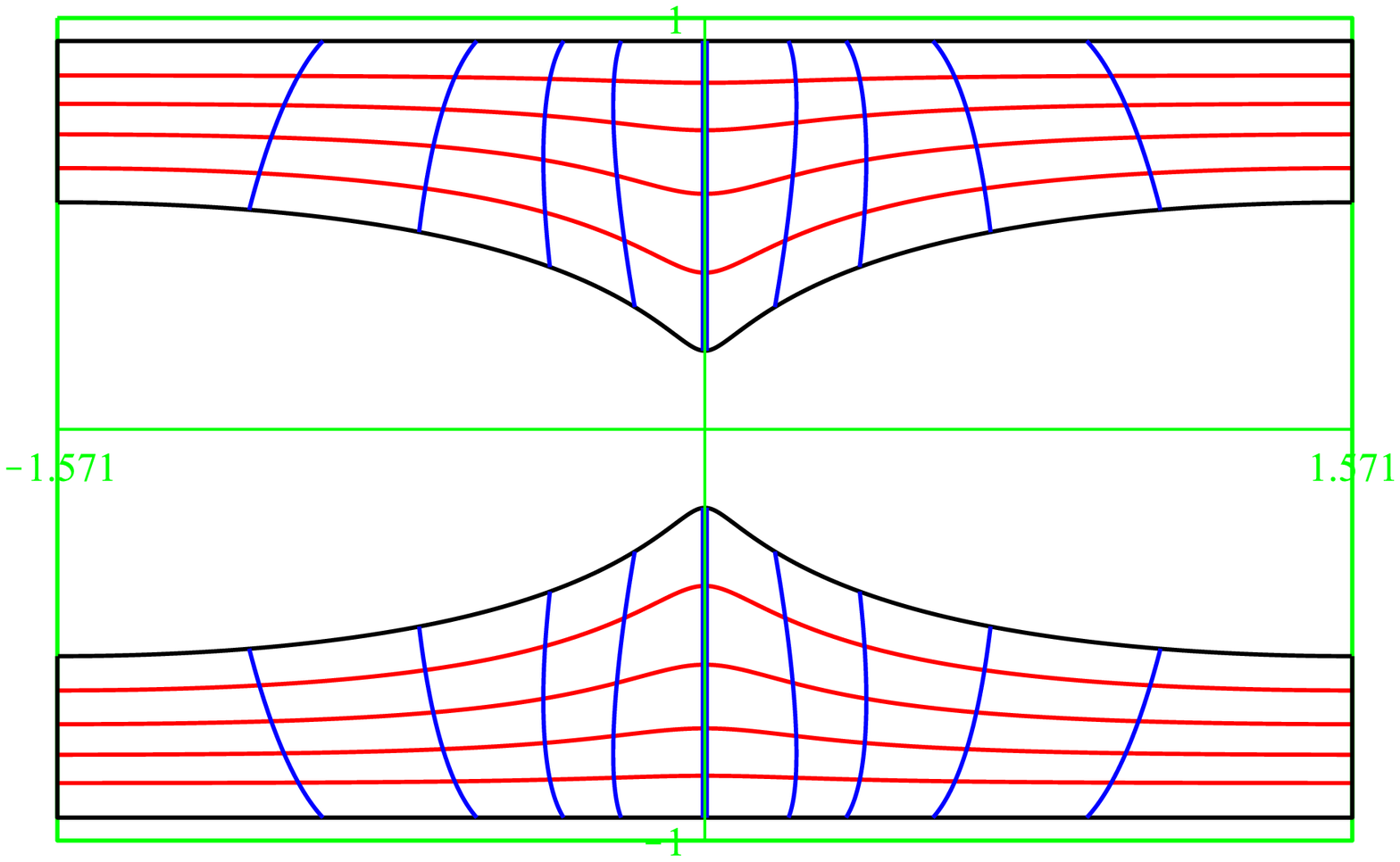}
    \caption*{Fig. \ref{fig:figBA05}. $\ujnorma=5\pi/6$ in ACKB}
  \end{subfigure}
\phantomcaption
\plabel{fig:figBA05}
\end{figure}

We see the following:
Through $\widehat{\mathrm M}^{\mathrm{CKB}}$, the various components of $\overline X$ map as follows:

$\bullet$ The Schur-hyperbolic boundary $\partial^{\mathrm{hyp}}\widehat X$ (i. e.  $s=0$, $t=0$) maps injectively to the outer rim.

$\bullet$ The Schur-elliptic boundary $\partial^{\mathrm{hyp}}\widehat X$  (i. e.  $r=0$, $t=1$) maps injectively to the inner rim
except to the origin.

$\bullet$ The closure of the pseudoboundary
$\partial^{\mathrm{hyp1}}\widehat X\cup \eth^{\mathrm{par}}\widehat X \cup\partial^{\mathrm{ell1}}\widehat X$
(i. e.  $\sin\theta=0$) gives two pinches on the sides.
(This is improved by $\widetilde{\mathrm M}^{\mathrm{ACKB}}$, which is injective on
$\partial^{\mathrm{hyp1}}\widetilde X\cup \eth^{\mathrm{par}}\widetilde X \cup\partial^{\mathrm{ell1}}\widetilde X$.)

$\bullet$ The degenerate boundary maps to the inner slit.
More precisely:

$\quad\circ$ For $\ujnorma<\pi$, there is no degenerate boundary.

$\quad\circ$ For $\ujnorma=\pi$, the degenerate boundary is $\partial^{\mathrm{dell0}}\widehat X$, it maps to the origin.

$\quad\circ$ For $\ujnorma>\pi$, the degenerate boundary is
$\partial^{\mathrm{dell}}\widehat X \cup \partial^{\mathrm{dnn0}}\widehat X \cup\partial^{\mathrm{dnn*}}\widehat X$
(restricted to norm $\ujnorma$).
Then $\partial^{\mathrm{dnn0}}\widehat X $
maps to the upper and lower tips of central slit,
$\partial^{\mathrm{dell*}}\widehat X $ maps to the origin,
$\partial^{\mathrm{dnn*}}\widehat X$ maps to join them.

$\bullet$ The red lines: $\theta$ varies, $s$ (or $r$) is fixed.
The blue lines: $s$ (or $r$) varies, $\theta$ is fixed.
It is quite suggestive that $\widehat X^{\mathrm{nn}}$ maps injectively into two simply connected
domains, but this requires some justification, which is as follows.

In fact, the proof of the next lemma will introduce  some important
computational techniques:

\begin{lemma}\plabel{lem:momentdirinj}
Suppose that $\ujnorma>0$.

(a)
Let $\mathcal S^{\mathrm{nn}}_\ujnorma$ be the subset of  $\mathcal S^{\mathrm{nn}}$ which contains the elements of norm $\ujnorma$.
Then  the operations $\mathrm{MR}^{\mathrm{CKB}}$ and   $\mathrm{ML}^{\mathrm{CKB}}$   (cf. \eqref{eq:MRLCKB}) are injective
on $\mathcal S^{\mathrm{nn}}_\ujnorma$.

(b)
Let $ X^{\mathrm{nn}}_\ujnorma$ be the subset of  $ X^{\mathrm{nn}}$ which is the restriction to  $\ujnorma=r+s$.
Then  the operation $\mathrm{M}^{\mathrm{CKB}}$  (cf. \eqref{eq:MCKB}) is injective
on $X^{\mathrm{nn}}_\ujnorma$.

In fact, its image is the union of two simply connected domains bounded by the image of the boundary pieces $ X^{\mathrm{nn}}$,
laying in the upper and lower half-planes, respectively.

\begin{proof}
Due to the circular fibration property, it is sufficient to prove (b).
By the monotonicity properties of the boundary pieces (cf. Lemma \ref{lem:momentspec})
one can characterize the simply connected domains.
Then one can check that $\widehat{\mathrm{M}}^{\mathrm{CKB}}$ (or $\widehat{\mathrm{M}}^{\mathrm{HP}}$)
start to map into the two simply connected domains near the boundary pieces, which is tedious but doable.
(It requires special arguments only in the case of
 $\partial^{\mathrm{dell0}}\widehat X \cup \partial^{\mathrm{dnn0}}\widehat X $; see comments later.)
Then, for topological reasons, it is sufficient to demonstrate that the Jacobian of $\widehat{\mathrm{M}}^{\mathrm{CKB}}$ (or $\widehat{\mathrm{M}}^{\mathrm{HP}}$)  on
 $X^{\mathrm{nn}}_\ujnorma$ never vanishes.
\snewpage

Then, it is sufficient to consider  the $\mathrm{HP}$ model.
We have to compute the Jacobian of
\[
(\hat a^{\mathrm{HP}},\hat b^{\mathrm{HP}} )=\left(\frac{a}b\cdot\frac{1}{\ujnorma \sqrt{\reD(b^2-r^2)}},
\frac{1}{\ujnorma\sqrt{\reD(b^2-r^2)}}+r\frac{\reC(b^2-r^2)}{\sqrt{\reD(b^2-r^2)}}\right)\]
with
\[a=t\ujnorma\cos\theta,\qquad
b=t\ujnorma\sin\theta,\qquad
r=\ujnorma-t\ujnorma\]
with respect to $t,\theta$ ($\ujnorma$ is fixed) subject to $0<t<1$, $\theta\in(0,\pi)\cup (\pi,2\pi)$.
We find
\begin{align}\frac{\partial(\hat a^{\mathrm{HP}},\hat b^{\mathrm{HP}})}{\partial(t,\theta)}
=\,&-\frac1{\left( \sin \theta \right) ^{2}}\,\frac{\reC}{\reD}\plabel{eq:jacc}\\
&-\,{\frac {1-t \left( \cos \theta \right) ^{2}}{ (\sin\theta)^2}}\, \frac{\reG}{\reD^2}\notag\\
&-\,{  \frac { \left( 1-t \left( \cos\theta  \right) ^{2} \right) }{(\sin\theta)^2 }}\left( 1-t \right)\ujnorma^{2}
\frac{\reL}{\reD^2}\notag\\
&-\,\ujnorma^{2}{t}^{2} \left( \cos \theta \right) ^{2}\frac{\reC\reG}{\reD^2}
\notag
\end{align}
where the arguments of $\reC,\reD,\reG,\reL$ should be $b^2-r^2=(t\ujnorma\sin\theta)^2-(1-t)^2\ujnorma^2$.
Each summand is non-positive;
in fact, strictly negative, with the possible exception of the last one.
Thus the Jacobian is negative.

Moreover, this method with the Jacobian applies to around $\partial^{\mathrm{hyp*}} X$ (i. e. $t=0$, $\sin\theta\neq0$) and
$\partial^{\mathrm{ell*}} X$ (i. e. $t=1$, $\sin\theta\neq0$)
as (\ref{eq:jacc}rhs/1) is always negative.
Moreover,  it also applies to $\partial^{\mathrm{dell*}}\widehat X $ and $\partial^{\mathrm{dnn*}}\widehat X$
where (\ref{eq:jacc}rhs/1)--(\ref{eq:jacc}rhs/3) extend to $0$ zero smoothly but (\ref{eq:jacc}rhs/4)
gives an nonzero contribution.

In fact, if we pass to  $(\hat a^{\mathrm{AHP}},\hat b^{\mathrm{AHP}} )$,
then it becomes regular even near the asymptotical points: It yields
\begin{multline}\frac{\partial(\hat a^{\mathrm{AHP}},\hat b^{\mathrm{AHP}})}{\partial(t,\theta)}
=-\frac{1}{\frac1{\ujnorma^2  \reD}(\cos\theta)^2+(\sin\theta)^2}\Biggl(\,
\,\frac{\reC}{\reD}
+\,(1-t \left( \cos \theta \right) ^{2})\, \frac{\reG}{\reD^2}\notag\\
+\,{  \left( 1-t \left( \cos\theta  \right) ^{2} \right) }\left( 1-t \right)\ujnorma^{2}
\frac{\reL}{\reD^2}
+\,\ujnorma^{2}{t}^{2} \left( \cos \theta \right) ^{2}\left( \sin \theta \right) ^{2}\frac{\reC\reG}{\reD^2}
\Biggr).
\plabel{eq:jacc2}
\end{multline}
This extends smoothly to the case $\sin\theta=0$, i. e. to  $\eth^{\mathrm{par}}\widehat X
\cup \partial^{\mathrm{hyp1}} X \cup \partial^{\mathrm{ell1}} X$.
It is easy to see that, on this set, the Jacobian is strictly negative again.
(Note that in this case $b^2-r^2=-r^2\leq0$, thus $\frac{\reC}{\reD}$ and $\frac{1}{\reD}$ are positive.)

This still leaves the uneasy case of
$\partial^{\mathrm{dell0}}\widehat X \cup \partial^{\mathrm{dnn0}}\widehat X $
regarding boundary behaviour.
However, even those cases can be handled, as failure of indicated behaviour at those point
(1-1 points for each simply connected domain for a fixed $\ujnorma\geq\pi$) would cause irregular behaviour at other points.
\end{proof}

\end{lemma}
\snewpage
The following discussion will be important for us.
For $A\in \mathcal S^{\mathrm{nn}}$, we let
\[\mathrm{NR}(A):=\frac{\boldsymbol X^{+}(\mathrm{MR}(A)) }{\det \mathrm{MR}(A)}
\equiv \frac{\hat b}{\hat a^2+\hat b^2-\breve c^2-\breve d^2}
\left(\hat a\Id_2+\frac{\breve c^2+\breve d^2-\hat a^2}{\hat b}\tilde I-\breve c \tilde J-\breve d K\right).\]

One can immediately see that
\[\mathrm{MR}_A(\mathrm{NR}(A) )=\frac12\tr\left(\mathrm{MR}(A)^* \mathrm{NR}(A)  \right)=0,\]
that is
\[\left. \frac{\mathrm d}{\mathrm dt}\|\log ((\exp A)(\exp( t \mathrm{NR}(A)))\|_2\right|_{t=0}=0.\]
On the other hand,
\begin{lemma}\plabel{lem:cogradient}
\[\left. \frac{\mathrm d^2}{\mathrm dt^2}\|\log ((\exp A)(\exp( t \mathrm{NR}(A)))\|_2\right|_{t=0}<-\frac1\ujnorma<0.\]
\begin{proof} By direct computation one can check that

\begin{align*}\mathrm{LHS}=
&-\frac1\ujnorma\\
&-\frac1\ujnorma(1-t(\cos\theta)^2)\frac{\reG}{\reD^2}\\
&-\ujnorma2(1-t)(1-t(\cos\theta)^2) \frac{\reG\reC}{\reD^2} \\
&-\ujnorma 2t (\sin\theta)^2\frac{\reG}{\reD}\\
&-\ujnorma (1-t)\frac{\reC^2}{\reD}\\
&-\ujnorma^3 (1-t)^2(1-t(\cos\theta)^2)\frac{\reC\reL}{\reD^2}\\
&-\ujnorma^3 t(1-t)(\sin\theta)^2\frac{\reL}{\reD}.
\end{align*}
(The arguments of $\reC,\reD,\reG,\reL$ should be $b^2-r^2\equiv \ujnorma^2t^2(\sin\theta)^2- \ujnorma^2(1-t)^2$.)

Then the statement follows from the observation that every term in the sum is strictly negative on the specified domain.
\end{proof}
\end{lemma}

Here an important observation is that $\mathrm{NR}(A)$ depends solely on $\mathrm{MR}(A)$,
and, by homogeneity, only on $\mathrm{MR}^{\mathrm{CKB}}(A)$.
Similar definitions, statements, and observations can be made regarding $\mathrm{ML}(A)$, and $\mathrm{NL}(A)$.

\begin{commentx}
\begin{remark}\plabel{rem:grad}
A naive gradation (or rather filtration) can be given by
$\reC\sim1$,
$\reD\sim1$,
$\reW\sim2$,
$\reP\sim2$,
$\reG\sim2$,
$\reL\sim3$,
$\ujnorma^2\sim-1$.
This  is induced by behaviour at $z\sim 0$.
\qedremark
\end{remark}
\end{commentx}
\snewpage

\begin{remark}\plabel{rem:momentkif}
Let us address the question that
whether $\mathrm{MR}^{\mathrm{CKB}}(A)$ and $\mathrm{ML}^{\mathrm{CKB}}(A)$
can be extended to \textit{formally} beyond the already discussed cases.

(o) If
$A= a\Id_2+ b\tilde I+(r\cos\psi)\tilde J+(r\sin\psi)$
such that
\begin{equation}
a^2+b^2>0,\qquad r>0,
\end{equation}
then
$\mathrm{MR}^{\mathrm{CKB}}$  extends by the formula
\begin{equation*}
 \mathrm{MR}^{\mathrm{CKB}}(A)=
 \frac{\mathcal A\Id_2+\mathcal B\tilde I+((\cos\psi)\Id_2+(\sin\psi)\tilde I)
 (\mathcal C\tilde J+\mathcal D\tilde K)
 }{\sqrt{\mathcal A^2+\mathcal B^2+\mathcal G^2 }},
\end{equation*}
where
\[\mathcal A\equiv (\cos\theta)\cdot\mathcal E((s\sin\theta)^2-r^2),\]
and
\[\mathcal B\equiv (\sin\theta)\cdot\left(
\mathcal E((s\sin\theta)^2-r^2)+r(s+r) \mathcal F((s\sin\theta)^2-r^2)
\right),\]
and
\[\mathcal C\equiv
\mathcal E((s\sin\theta)^2-r^2)+(\sin\theta)^2\cdot s(s+r) \mathcal F((s\sin\theta)^2-r^2),\]
and
\[\mathcal D\equiv
 (\sin\theta)\cdot(r+s)\mathcal E((s\sin\theta)^2-r^2),\]
and
\[\mathcal G\equiv (\sin\theta)\cdot(s+r).\]
This dispenses with the condition \eqref{eq:dom2}.

\snewpage

(a) If $s>0$, $r=0$, then $\frac{
 (\mathcal C\tilde J+\mathcal D\tilde K)
 }{\sqrt{\mathcal A^2+\mathcal B^2+\mathcal G^2 }}$
can be replaced by an arbitrary $(\cos\xi)\tilde J+(\sin\xi)\tilde K$.
Formally, this leads to
\begin{equation*}
 \mathrm{MR}^{\mathrm{CKB}}(A)=
 \frac{\mathcal A\Id_2+\mathcal B\tilde I
 }{\sqrt{\mathcal A^2+\mathcal B^2+\mathcal G^2 }}+((\cos\psi)\tilde J+(\sin\psi)\tilde K)
\end{equation*}
where
$\psi$ is undecided.
Here
\[\mathcal A\equiv (\cos\theta)\cdot\mathcal E((s\sin\theta)^2),\]
and
\[\mathcal B\equiv (\sin\theta)\cdot\left(
\mathcal E((s\sin\theta)^2)
\right),\]
and
\[\mathcal C\equiv
\mathcal E((s\sin\theta)^2)+(s\sin\theta)^2 \mathcal F((s\sin\theta)^2),\]
and
\[\mathcal D\equiv
 (s\sin\theta)\cdot\mathcal E((s\sin\theta)^2),\]
and
\[\mathcal G\equiv s\sin\theta.\]
(The undecidedness of $\psi$ does not appear $\widehat{\mathrm M}^{\mathrm{CKB}}$.)

(b) If $s=0$, $r>0$, then $\theta$ becomes undecided.
Here
\[\mathcal A\equiv (\cos\theta)\cdot\mathcal E( -r^2),\]
and
\[\mathcal B\equiv (\sin\theta)\cdot\left(
\mathcal E(-r^2)+r^2 \mathcal F(-r^2)
\right),\]
and
\[\mathcal C\equiv
\mathcal E(-r^2),\]
and
\[\mathcal D\equiv
 (\sin\theta)\cdot r\mathcal E(-r^2),\]
and
\[\mathcal G\equiv  r\sin\theta .\]
(The undecidedness of $\theta$ does appear $\widehat{\mathrm M}^{\mathrm{CKB}}$.)

These calculations were formal.
The real geometrical picture is computed in what follows.
``Undecidedness'' will appear as conical degeneracy.
\qedremark
\end{remark}

\begin{lemma}\plabel{lem:momentdeg}
Suppose $A=a\Id_2+b\tilde I$ with $-\pi<b<\pi$  (thus $A\in\partial^{\mathrm{ell}}\mathcal S$)
and $\m v=v_1\Id_2+v_2\tilde I+v_3\tilde J+v_4\tilde K$. Then
\[\mathrm{MR}_A(\m v)=\mathrm{ML}_A(\m v)=\frac{a}{\sqrt{a^2+b^2}}v_1+\frac{b}{\sqrt{a^2+b^2}}v_2
+\frac{b}{\sin b}\sqrt{(v_3)^2+(v_4)^2}.
\]
In particular, if
\[\m v=\hat v_1 \frac{a\Id_2+b\tilde I}{\sqrt{a^2+b^2}}+\hat v_2 \frac{-b\Id_2+a\tilde I}{\sqrt{a^2+b^2}} + \hat v_3\frac{\tilde J}{\frac{b}{\sin b}}+ \hat v_4\frac{\tilde K}{\frac{b}{\sin b}}\]
then
\[\mathrm{MR}_A(\m v)=\mathrm{ML}_A(\m v)=\hat v_1+\sqrt{(\hat v_3)^2+(\hat v_4)^2}. \]
(For $b=0$, $\frac{b}{\sin b}$ should be understood as $1$.)
\begin{proof}
Direct computation.
\end{proof}
\end{lemma}
\snewpage
\begin{lemma}\plabel{lem:momentdeg3}
Suppose $A=c\tilde J+d\tilde K$,  $c^2+d^2>0$,  and $\m v=v_1\Id_2+v_2\tilde I+v_3\tilde J+v_4\tilde K$.
Then
\[\mathrm{MR}_A(\m v)=\frac{cv_3+dv_4}{\sqrt{c^2+d^2}} +\sqrt{v_1^2+\left(\sqrt{c^2+d^2}\frac{\cosh \sqrt{c^2+d^2}}{\sinh \sqrt{c^2+d^2}}v_2-cv_4+dv_3\right)^2}.\]

In particular, if
\begin{align*}
\m v=
&\hat v_1\frac{c\tilde J+d\tilde K}{\sqrt{c^2+d^2}}
+\hat v_2\frac{(\sinh \sqrt{c^2+d^2})\tilde I+(\cosh \sqrt{c^2+d^2})\frac{-d\tilde J+c\tilde K}{\sqrt{c^2+d^2}}}{\dfrac{\sqrt{c^2+d^2}}{\sinh \sqrt{c^2+d^2}}}
\\
&+\hat v_3\Id_2
+\hat v_4 \frac{(\cosh \sqrt{c^2+d^2})\tilde I+(\sinh \sqrt{c^2+d^2})\frac{-d\tilde J+c\tilde K}{\sqrt{c^2+d^2}}}{\dfrac{\sqrt{c^2+d^2}}{\sinh \sqrt{c^2+d^2}}}
,
\end{align*}
then
\[\mathrm{MR}_A(\m v)=\hat v_1+\sqrt{(\hat v_3)^2+(\hat v_4)^2}.\]

Similarly,
\[\mathrm{ML}_A(\m v)=\frac{cv_3+dv_4}{\sqrt{c^2+d^2}} +\sqrt{v_1^2+\left(\sqrt{c^2+d^2}\frac{\cosh \sqrt{c^2+d^2}}{\sinh \sqrt{c^2+d^2}}v_2
+cv_4-dv_3\right)^2}.
\]
In particular, if
\begin{align*}
\m v=
&\hat v_1\frac{c\tilde J+d\tilde K}{\sqrt{c^2+d^2}}+
\hat v_2\frac{(\sinh \sqrt{c^2+d^2})\tilde I+(\cosh \sqrt{c^2+d^2})\frac{d\tilde J-c\tilde K}{\sqrt{c^2+d^2}}}{\dfrac{\sqrt{c^2+d^2}}{\sinh \sqrt{c^2+d^2}}}
\\
&+\hat v_3\Id_2+
\hat v_4 \frac{(\cosh \sqrt{c^2+d^2})\tilde I+(\sinh \sqrt{c^2+d^2})\frac{d\tilde J-c\tilde K}{\sqrt{c^2+d^2}}}{\dfrac{\sqrt{c^2+d^2}}{\sinh \sqrt{c^2+d^2}}}
,
\end{align*}
then
\[\mathrm{ML}_A(\m v)=\hat v_1+\sqrt{(\hat v_3)^2+(\hat v_4)^2}.\]
\begin{proof}
Direct computation.
\end{proof}
\end{lemma}
\begin{lemma}\plabel{lem:seugen}
Suppose that $A_1\in \partial^{\mathrm{dell}}\mathcal S$, and,
in particular, $A_1=a\Id_2+\pi\tilde I$ or $A_1=a\Id_2-\pi\tilde I$.
Let $\m v=v_1\Id_2+v_2\tilde I$.
Then
\[\mathrm{MR}_A(\m v)=\mathrm{ML}_A(\m v)=\frac{av_1-\pi|v_2|}{\sqrt{a^2+\pi^2}}.\]

\begin{proof} Direct computation.
\end{proof}
\end{lemma}
\snewpage
\scleardoublepage\section{BCH minimality}
\plabel{sec:BCHmin}

\begin{defin}
We say that $(A,B)\in\mathrm M_2(\mathbb R)\times \mathrm M_2(\mathbb R)$ is a BCH minimal pair (for $C$),
if for any $(\tilde A,\tilde B)\in\mathrm M_2(\mathbb R)\times \mathrm M_2(\mathbb R)$
such that $(\exp\tilde  A )(\exp \tilde B ) =(\exp A)(\exp B)$ $(=C)$ it holds that
$ \|\tilde A\|_2>\|A\|_2$ or $\|\tilde B\|_2>\|B\|_2$.
\end{defin}
In this section we seek restrictions for BCH minimal pairs.

Any element $C\in\GL^+_2(\mathbb R)$ can easily be perturbed into the product of
two $\log$-able elements (for example, one close to $A$ and one close to $\Id_2$).
By a simple compactness argument we find:
If $C=(\exp A)(\exp B)$, then $C$ allows a minimal pair $(A,B)$ such that $\|\tilde A\|_2\leq \|A\|_2$ and $\|\tilde B\|_2\leq\|B\|_2$.
In minimal pairs, however, we can immediately restrict our attention to some special matrices:

By Lemma \ref{lem:Sacc}, all minimal pairs are from ${\mathcal S^{\mathrm{acc}}}\times{\mathcal S^{\mathrm{acc}}}$.
By the same lemma, we also see that the elements of ${\mathcal S^{\mathrm{acc}}}\times\partial^0\mathcal S $
and  $\partial^0\mathcal S\times{\mathcal S^{\mathrm{acc}}} $ are all minimal pairs.

\begin{defin}
Assume that $(A,B)\in {\mathcal S^{\mathrm{acc}}} \times {\mathcal S^{\mathrm{acc}}}$.
We say that $(A,B)$ is an infinitesimally minimal BCH pair, if we cannot find $\m v\in\mathrm M_2(\mathbb R) $ such that
\begin{equation}
\mathrm{MR}_A(\m v)<0  \qquad\text{and}\qquad \mathrm{ML}_B(-\m v)<0.
\plabel{eq:notBCHmin}
\end{equation}

\end{defin}
\begin{lemma}
\plabel{lem:InfMin}
Assume that $(A,B)\in {\mathcal S^{\mathrm{acc}}} \times {\mathcal S^{\mathrm{acc}}}$,
and $(A,B)$ is not infinitesimally minimal.
Then
one can find $\tilde A,\tilde B \in\mathcal S$ such that
\[\|\tilde A\|_2<\| A\|_2,\qquad\text{and}\qquad \|\tilde B\|_2<\| B\|_2,\]
yet
\[(\exp\tilde  A)(\exp\tilde B)=(\exp A)(\exp B).\]

In particular, infinitesimal minimality is a necessary condition for minimality.
\begin{proof}
If  \eqref{eq:notBCHmin} holds, the let $\tilde A=\log((\exp A)(\exp t\m v))$ and $\tilde B=\log((\exp -t\m v)(\exp B))$
with a sufficiently small $t>0$.
\end{proof}
\end{lemma}
\begin{lemma}
\plabel{lem:EqMin}
Suppose that $A,B \in\mathcal S $ such that  $\mathrm{MR}(A)$ and $\mathrm{ML}(B)$ are not positive multiples of each other.
Then the pair $(A,B)$ is not infinitesimally minimal; in particular,  $(A,B)$ is not minimal.
\begin{proof}
The the nonzero linear functionals $\mathrm{MR}_A(\cdot)$ and $\mathrm{ML}_A(\cdot)$ are not
positive multiples of each other, thus one can find $\m v$ such that \eqref{eq:notBCHmin} holds.
\end{proof}
\end{lemma}

\begin{lemma}\plabel{lem:EAGAN}
If $(A,B)\in\mathcal S\times \mathcal S$ is BCH minimal,
then $(A,B)\in\eth^{\mathrm{par}}\mathcal S\times\eth^{\mathrm{par}}\mathcal S  $.
\begin{proof}
For $A\in \mathcal S^{\mathrm{nn}}$, $\mathrm{MR}^{\mathrm{CKB}}(A)$ and $\mathrm{ML}^{\mathrm{CKB}}(A)$ are
inner points CKB, for $A\in \eth^{\mathrm{par}}\mathcal S$, $\mathrm{MR}^{\mathrm{CKB}}(A)$ and $\mathrm{ML}^{\mathrm{CKB}}(A)$ are
asymptotic points CKB. From this, and the previous lemma,
if $(A,B)$ is infinitesimally minimal, then
$(A,B)\in\eth^{\mathrm{par}}\mathcal S\times\eth^{\mathrm{par}}\mathcal S  $
or $(A,B)\in\mathcal S^{\mathrm{nn}}\times \mathcal S^{\mathrm{nn}}$.
Regarding the second case, by infinitesimal minimality, $\mathrm{MR}^{\mathrm{CKB}}(A)=\mathrm{ML}^{\mathrm{CKB}}(B)$ is required.
Then we can take $\m v=\mathrm{NR}(A)=\mathrm{NL}(B)$, and, by Lemma \ref{lem:cogradient},
$\tilde A=\log((\exp A)(\exp t\m v))$ and $\tilde B=\log((\exp -t\m v)(\exp B))$
will give a counterexample to BCH minimalty with a sufficiently small $t>0$.
\end{proof}
\end{lemma}

\begin{lemma}
\plabel{lem:Ell}
Suppose that $(A,B)\in S^{\mathrm{acc}}\times S^{\mathrm{acc}}$ is an infinitesimally minimal BCH pair. We claim:

(a) If $A\in \partial^{\mathrm{dell}}\mathcal S$, then $B\in \partial^{\mathrm{hyp}}\mathcal S $.

(b) If  $A\in \partial^{\mathrm{ell*}}\mathcal S$, then
$B\in \partial^{\mathrm{ell*}}\mathcal S$ or $B\in \partial^{\mathrm{hyp}}\mathcal S$ or $B\in \mathcal S^{\mathrm{nn}}$.

(c) If  $A \in\partial^{\mathrm{ell0}}\mathcal S$, then
$B\in \partial^{\mathrm{ell0}}\mathcal S$ or $B\in \partial^{\mathrm{hyp}}\mathcal S$ or $B\in \eth^{\mathrm{par}} \mathcal S$.

\begin{proof}
(a) By Lemma \ref{lem:seugen}, for $\m v\sim\pm \tilde I$
the relation holds $\mathrm{MR}_A(\mathbf v)<0$.
(`$\sim$ ' means `approximately'.)
Then the statement is  a  consequence of Lemma \ref{lem:InfMin}, and furthermore, Lemma \ref{lem:momentdeg}.

(b) and (c) follow from the observation that the maps
$\mathrm{MR}_A|_{\mathbb R\Id_2+\mathbb R\tilde I}$ and  $\mathrm{ML}_B|_{\mathbb R\Id_2+\mathbb R\tilde I}$,
if they are linear functionals,
should be proportional to each other (cf. Lemma \ref{lem:InfMin}).
\end{proof}
\end{lemma}
By taking adjoints, the relation described in the previous lemma is symmetric in $A$ and $B$.
This leads to
\begin{theorem}\plabel{th:strat}
For minimal BCH pairs, where all factors are of positive norm, one has
the following incidence possibilities:
\[\begin{array}{c|ccccccc}
&  \mathcal S^{\mathrm{nn}}&\eth^{\mathrm{par}}\mathcal S&\partial^{\mathrm{ell*}}\mathcal S
&\partial^{\mathrm{ell0}}\mathcal S&\partial^{\mathrm{dell}}\mathcal S&\partial^{\mathrm{hyp}}\mathcal S
\\\hline
\mathcal S^{\mathrm{nn}}&\times&\times&\checkmark&\times&\times&\checkmark\\
\eth^{\mathrm{par}}\mathcal S&\times&\checkmark&\times&\checkmark&\times&\checkmark\\
\partial^{\mathrm{ell*}}\mathcal S&\checkmark&\times&\checkmark&\times&\times&\checkmark\\
\partial^{\mathrm{ell0}}\mathcal S&\times&\checkmark&\times&\checkmark&\times&\checkmark\\
\partial^{\mathrm{dell}}\mathcal S&\times&\times&\times&\times&\times&\checkmark\\
\partial^{\mathrm{hyp}}\mathcal S&\checkmark&\checkmark&\checkmark&\checkmark&\checkmark&\checkmark
\end{array}
\]
Here $\checkmark$ means `perhaps possible', and $\times$ means `not possible'.
\begin{proof}
This is a consequence of the previous statements.
\end{proof}
\end{theorem}
Next, we obtain more quantitative restrictions.

\begin{lemma}
\plabel{lem:BAGAN}
(a) For $\ujnorma_1,\ujnorma_2 > 0$ consider the pairs
\[(A,B)=\left(\bem\ujnorma_1&\\&t_1\ujnorma_1\eem  , \bem\ujnorma_2&\\&t_2\ujnorma\eem\right)\]
or
\[(A,B)=\left(-\bem\ujnorma_1&\\&t_1\ujnorma\eem  , -\bem\ujnorma_2&\\&t_2\ujnorma\eem\right)\]
where $t_1,t_2\in(-1,1)$.
Then the pairs $(A,B)$ will produce, up to conjugation by rotation matrices, every (infinitesimally) BCH minimal pair
from $\eth^{\mathrm{par}}\mathcal S \times\eth^{\mathrm{par}}\mathcal S$.
In these cases \[\|\log((\exp A)( \exp B)\|_2=\ujnorma_1+\ujnorma_2.\]

(b) Similar statement holds for $t_1,t_2\in(-1,1]$ with respect to
$(\eth^{\mathrm{par}}\mathcal S \cup \eth^{\mathrm{ell0}}\mathcal S)\times
(\eth^{\mathrm{par}}\mathcal S \cup \eth^{\mathrm{ell0}}\mathcal S)$.
\begin{proof}
(a) By infinitesimal minimality, the matrices should be aligned.
BCH minimality follows from Magnus minimality.
(b) This is a trivial extension of (a).
\end{proof}
\end{lemma}
\snewpage
\begin{lemma}\plabel{lem:ConeDual}
Consider the function given by
\[\alpha(t_1U_1+t_2U_2+t_3U_3+t_4U_4)=t_1+\sqrt{t_3^2+t_4^2}\]
and the linear functional given by
\[\beta( t_1U_1+t_2U_2+t_3U_3+t_4U_4 )=\beta_1t_1+\beta_2t_2+\beta_3t_3+\beta_4t_4.\]

Then we can choose $\m v= t_1U_1+t_2U_2+t_3U_3+t_4U_4$
such that $\alpha(\m v)<0$ and $\beta(-\m v)<0$ unless
\begin{equation}
\beta_2=0\qquad
\text{and}\qquad
\beta_1\geq \sqrt{\beta_3^2+\beta_4^2}
\plabel{eq:coniq}
\end{equation}
holds, in which case finding such a $\m v$ is impossible.
\begin{proof}
If $\beta_2\neq 0$, or $\beta_1<0$, or $0=\beta_1<\sqrt{\beta_3^2+\beta_4^2}$, or $0<\beta_1<\sqrt{\beta_3^2+\beta_4^2}$,
then
\[\m v=-U_1+\frac{\beta_1+1}{\beta_2}U_2\]
or
\[\m v=-U_1,\]
or
\[\m v=-2U_1+\frac{\beta_3}{\sqrt{\beta_3^2+\beta_4^2}}U_3+\frac{\beta_4}{\sqrt{\beta_3^2+\beta_4^2}}U_4\]
or
\[\m v= -\frac{\beta_1+\sqrt{\beta_3^2+\beta_4^2}}{2\beta_1}U_1+\frac{\beta_3}{\sqrt{\beta_3^2+\beta_4^2}}U_3+\frac{\beta_4}{\sqrt{\beta_3^2+\beta_4^2}}U_4\]
respectively, are a good choices.

Conversely, assume that \eqref{eq:coniq} holds and $\alpha(\m v)<0$, i. e.
$t_1<-\sqrt{t_3^2+t_4^2}$, also holds.
Then,
\[\beta_1t_1+\beta_2t_2+\beta_3t_3+\beta_4t_4\leq - \sqrt{\beta_3^2+\beta_4^2}\sqrt{t_3^2+t_4^2}+0+\sqrt{\beta_3^2+\beta_4^2}\sqrt{t_3^2+t_4^2}=0\]
shows that $\beta(-\m v)\geq0$.

(In merit, this is just a simple geometrical statement about  a half-cone.)
\end{proof}
\end{lemma}

\snewpage

\begin{lemma}\plabel{lem:momentdeg2}
Suppose that $A_1=a\Id_2+b\tilde I+c\tilde J\in\mathcal S$, $c>0$, and $B_1=\acute{a}\Id_2+\acute{b}\tilde I\in\partial^{\mathrm{ell}}\mathcal S$,
where $\acute \ujnorma=\sqrt{\acute a^2+\acute b^2}$.
Suppose that $(A_1,B_1)$ is infinitesimally BCH minimal.
Then
\[B_1= \frac{\hat{a}\Id_2+\hat{b}\tilde I}{\sqrt{\hat a^2+\hat b^2}} \,\acute{\ujnorma}\]
and
\begin{equation}
\frac{\hat b\acute{\ujnorma}}{\sin \dfrac{\hat b\acute{\ujnorma}}{\sqrt{\hat a^2+\hat b^2}} }\geq \sqrt{\breve c^2+\breve d^2}.
\plabel{eq:momentdeg2res}
\end{equation}
(For the $\hat b=0$, the LHS is understood as $\sqrt{\hat a^2+\hat b^2} $, but then the condition is vacuous anyway.
Also note that the condition $\left|\dfrac{\hat b\acute{\ujnorma}}{\sqrt{\hat a^2+\hat b^2}} \right|<\pi$
is imposed by  $B_1\in\partial^{\mathrm{ell}}\mathcal S $.)
\begin{proof}
This a consequence of Lemma \ref{lem:ConeDual} and Lemma \ref{lem:momentdeg}.
By (\ref{eq:coniq}/cond1),
\[B_1= \pm\frac{\hat{a}\Id_2+\hat{b}\tilde I}{\sqrt{\hat a^2+\hat b^2}} \,\acute{\ujnorma}.\]
By (\ref{eq:coniq}/cond2), only the sign choice + is valid, and  \eqref{eq:momentdeg2res} holds.
\end{proof}
\end{lemma}

\begin{lemma}\plabel{lem:momentdeg4}
Suppose that $A_1=a\Id_2+b\tilde I+c\tilde J\in\mathcal S$, $c>0$, and $B_1=\acute{c}\tilde J+\acute{d}\tilde K\in\partial^{\mathrm{hyp}}\mathcal S$,
where $\acute\ujnorma=\sqrt{\acute c^2+\acute d^2}$.
If $(A_1,B_1)$ is infinitesimally BCH minimal, then
\[B_1=\left(\sqrt{\breve c^2+\breve d^2- \left(  \frac{\sinh \acute \ujnorma}{\cosh \acute \ujnorma}\hat b\right)^2}\Id_2
 -\frac{\sinh \acute \ujnorma}{\cosh \acute \ujnorma}\hat b\tilde I \right) \frac{\breve c\tilde \Id_2 +\breve d\tilde I}{\breve c^2+\breve d^2}\acute N\tilde J\]
and
\begin{equation}
\breve c^2+\breve d^2\geq \hat a^2+\left(  \frac{\sinh \acute \ujnorma}{\cosh \acute \ujnorma}\hat b\right)^2\frac{\acute\ujnorma^2+1}{\acute\ujnorma^2}.
\plabel{eq:momentdeg4res}
\end{equation}
\begin{proof}
This a consequence of Lemma \ref{lem:ConeDual} and Lemma \ref{lem:momentdeg3}.
By (\ref{eq:coniq}/cond1),
\[B_1=\left(\pm\sqrt{\breve c^2+\breve d^2- \left(  \frac{\sinh \acute \ujnorma}{\cosh \acute \ujnorma}\hat b\right)^2}\Id_2
 -\frac{\sinh \acute \ujnorma}{\cosh \acute \ujnorma}\hat b\tilde I \right) \frac{\breve c\tilde \Id_2 +\breve d\tilde I}{\breve c^2+\breve d^2}\acute N\tilde J.\]
By (\ref{eq:coniq}/cond2), only the sign choice + is valid, and  \eqref{eq:momentdeg4res} holds.
\end{proof}

\end{lemma}

\snewpage
\begin{lemma}\plabel{lem:MinRR}
Suppose that $A_1=c_1\tilde J+d_1\tilde K$ and $A_2=c_2\tilde J+d_2\tilde K$.
Assume that $\slashed N_1=\sqrt{c_1^2+d_1^2}>0$ and $\slashed N_2=\sqrt{c_2^2+d_2^2}>0$.
Let $\psi$ denote the angle between the vectors $(c_1,d_1)$ and $(c_1,d_2)$.
We claim that if $(A_1,A_2)$ is an infinitesimally BCH minimal pair, then $\psi<\pi$,
and
\begin{equation}
\tan\frac\psi2\leq
 {\frac { \left(  \left(\coth \slashed N_1
 \right) +\left(\coth \slashed N_2 \right)
 \right) \slashed N_1\,\slashed N_2}{\slashed N_1\,
\left(\coth \slashed N_1 \right) +\slashed N_2\,\left(\coth \slashed N_2
 \right)  }}.
 \plabel{eq:MinRRres}
\end{equation}
\begin{proof} $\psi<\pi$ is immediate.
If \eqref{eq:MinRRres} does not hold, then with the choice
\[\m v=\frac{A_1}{\slashed N_1}-\frac{A_2}{\slashed N_2}+\frac{\slashed N_1-\slashed N_2}{
\slashed N_1(\coth \slashed N_1)+ \slashed N_2(\coth \slashed N_2)
}\frac{A_1A_2-A_2A_1}{2\slashed N_1\slashed N_2}
,\]
it yields
\[\mathrm{MR}_A(\m v)=\mathrm{ML}_A(-\m v)=
\left( {\frac { \left(  \left(\coth \slashed N_1
 \right) +\left(\coth \slashed N_2 \right)
 \right) \slashed N_1\,\slashed N_2}{\slashed N_1\,
\left(\coth \slashed N_1 \right) +\slashed N_2\,\left(\coth \slashed N_2
 \right)  }}
 -\tan\frac\psi2 \right)\sin\psi <0.\]

(By conjugation we can assume that
\[c_1=\slashed N_1\cos\frac\psi2,\qquad d_1=\slashed N_1\sin\frac\psi2,
\qquad
c_2=\slashed N_2\cos\frac\psi2,\qquad d_2=-\slashed N_2\sin\frac\psi2.\]
Then direct computation yields the statement.)
\end{proof}
\end{lemma}

Assume that $U$ is a subset of $\mathrm M_2(\mathbb R)\times \mathrm M_2(\mathbb R)$ closed under conjugation by rotation matrices.
Then the set $\exp(U)= \{(\exp A)(\exp B)\,:\,(A,B)\in U\}$ is also closed conjugation by rotation matrices.
Thus it can be represented by its image through $\Xi^{\mathrm{PH}}$ accurately.
We say that $\exp(U)$ set has a  factor dimension $p$, if the map $\Xi^{\mathrm{PH}}$
applied to it can be factorized smoothly through a smooth manifold of dimension $p$.

\begin{theorem}
\plabel{stratdim}
Factor dimensions of the BCH minimal subsets of the sets $U_1\times U_2$
restricted to $(A,B)\in U_1\times U_2$, $\|A\|_2=\ujnorma_1$, $\|B\|_2=\ujnorma_2$ are
 \[\begin{array}{c|ccccccc}
&  \mathcal S^{\mathrm{nn}}&\eth^{\mathrm{par}}\mathcal S&\partial^{\mathrm{ell*}}\mathcal S
&\partial^{\mathrm{ell1}}\mathcal S&\partial^{\mathrm{dell}}\mathcal S&\partial^{\mathrm{hyp}}\mathcal S
\\\hline
\mathcal S^{\mathrm{nn}}&&&2&&&2\\
\eth^{\mathrm{par}}\mathcal S&&1&&1&&1\\
\partial^{\mathrm{ell*}}\mathcal S&2&&1&&&1\\
\partial^{\mathrm{ell1}}\mathcal S& &1& &0& &0\\
\partial^{\mathrm{dell}}\mathcal S& & & & &\ &0\\
\partial^{\mathrm{hyp}}\mathcal S&2&1&1&0&0&1
\end{array}
\]
\proofremark{
The numbers given are sort of upper estimates, but one can see that one cannot reduce them generically.
}
\begin{proof}
The previous statements and  simple arguments show this.
E. g. Lemma \ref{lem:momentdeg2} takes care to $\mathcal S^{\mathrm{nn}}\times \partial^{\mathrm{ell*}}\mathcal S$, etc.
\end{proof}
\end{theorem}
Furthermore, note that the image of the restriction of $\eth^{\mathrm{par}}\mathcal S\times \eth^{\mathrm{par}}\mathcal S$
contains the images of the restrictions of  $\eth^{\mathrm{par}}\mathcal S\times \partial^{\mathrm{ell1}}\mathcal S $
and $\partial^{\mathrm{hyp}}\mathcal S \times \eth^{\mathrm{par}}\mathcal S$, which contain
the image of the restriction of
$\partial^{\mathrm{hyp}}\mathcal S\times\partial^{\mathrm{ell1}}\mathcal S$
(and also in reverse order).

\snewpage
Let us define the map
\[\SE_{\ujnorma,\acute{\ujnorma}}:X^{\mathrm{nn}}_{\ujnorma}\rightarrow\GL_2(\mathbb R)\]
by
\[(a,b,r)\mapsto
\exp\left(
a\Id_2+b\tilde I+r\tilde J
\right)
\exp\left(
\frac{\hat{a}\Id_2+\hat{b}\tilde I}{\sqrt{\hat a^2+\hat b^2}} \,\acute{\ujnorma}
\right)
,\]
where $\hat a$ and $\hat b$ are as in Lemma \ref{lem:momentkif} and $\acute\ujnorma\geq0$.

Note that $\log \SE_{\ujnorma,\acute{\ujnorma}}$ is surely well-defined if $0<\ujnorma$ and $\ujnorma+\acute\ujnorma\leq\pi$
(as the BCH presentation is not Magnus-minimal).
\begin{lemma}
\plabel{lem:SElog}
Consider the map $\SE_{\ujnorma,\acute{\ujnorma}}$
with $a=\ujnorma t\cos\theta$, $b=\ujnorma t\sin\theta$, $r=\ujnorma (1-t)$.
(Thus $0<t<1$ and $\sin\theta\neq0$.) We claim:

(a)
\[\mathrm{mdis}\,\SE_{\ujnorma,\acute{\ujnorma}}(a,b,r)=r \Sin(b^2-r^2).  \]

(b)
$\SE_{\ujnorma,\acute{\ujnorma}}(a,b,r)$
is log-able if and only if
\[\Cos(b^2-r^2)\cos\frac{\hat b\acute N}{\sqrt{ \hat a^2+\hat b^2}}
-b\Sin(b^2-r^2)\sin\frac{\hat b\acute N}{\sqrt{ \hat a^2+\hat b^2}} >-1.\]

(c) In the log-able case:
Let $\pi_{23}$ denote the projection to the second and third coordinates.
Then, for the Jacobian,
\[
\frac{\partial\left( \pi_{23}\circ\Xi^{\mathrm{PH}}\circ\log\circ \SE_{\ujnorma,\acute{\ujnorma}} ( \ujnorma t\cos\theta,\ujnorma t\sin\theta,\ujnorma (1-t) )\right)}{\partial(t,\theta)}=
\]
\[(\cos\theta)\cdot\AC\left(\Cos(b^2-r^2)\cos\frac{\hat b\acute N}{\sqrt{ \hat a^2+\hat b^2}}
-b\Sin(b^2-r^2)\sin\frac{\hat b\acute N}{\sqrt{ \hat a^2+\hat b^2}}
\right)\cdot \Sin(b^2-r^2)\cdot\]
\[\]
\[\Biggl(\ujnorma^2 t +\frac{\ujnorma\acute\ujnorma}{(\hat a^2+\hat b^2)^{3/2}}
\Biggl(
1+(1-t)(2-t(\cos\theta)^2)\ujnorma^2 \reC
+t^2(1-t)(\sin\theta)^2(\cos\theta)^2\ujnorma^4\reW
\]\[+
(1-t)(   (1-t)(\sin\theta)^2+(1-t)^2(\cos\theta)^2 +t^2(\sin\theta)^2(\cos\theta)^2 )\ujnorma^4\reC^2
\Biggr)
\Biggr),\]
where the arguments of $\reC$ and $\reW$ should be $b^2-r^2$.

In particular, this is non-vanishing if $\cos\theta\neq0$.
\begin{proof}
(a) This is direct computation.
(b) This is a consequence of Lemma \ref{lem:logreal}(a).

(c) The formula for $\log$ can be applied, then direct computation yields the result.
Regarding non-vanishing, we can see that beyond $\cos\theta$, all multiplicative terms are positive.
\end{proof}
\end{lemma}
\snewpage

Let us define the map
\[\SH_{\ujnorma,\acute{\ujnorma}}:X^{\mathrm{nn}}_{\ujnorma}\rightarrow\GL_2(\mathbb R)\]
by
\begin{multline*}
(a,b,r)\mapsto
\exp\left(
a\Id_2+b\tilde I+r\tilde J
\right)\cdot\\\cdot
\exp\left(
\left(\sqrt{\breve c^2+\breve d^2- \left(  \frac{\sinh \acute \ujnorma}{\cosh \acute \ujnorma}\hat b\right)^2}\Id_2
 -\frac{\sinh \acute \ujnorma}{\cosh \acute \ujnorma}\hat b\tilde I \right) \frac{\breve c\tilde \Id_2 +\breve d\tilde I}{\breve c^2+\breve d^2}\acute N\tilde J
\right),
\end{multline*}
where $\hat a$ and $\hat b$ are as in Lemma \ref{lem:momentkif} and $\acute\ujnorma\geq0$.

Note that $\log \SH_{\ujnorma,\acute{\ujnorma}}$ is surely well-defined if $0<\ujnorma$ and $\ujnorma+\acute\ujnorma\leq\pi$
(as the BCH presentation is not Magnus-minimal).

\snewpage
\begin{lemma}
\plabel{lem:SHlog}
Consider the map $\SH_{\ujnorma,\acute{\ujnorma}}$
with $a=\ujnorma t\cos\theta$, $b=\ujnorma t\sin\theta$, $r=\ujnorma (1-t)$.
(Thus $0<t<1$ and $\sin\theta\neq0$.) We claim:

(a)
\[\mathrm{mdis}\,\SH_{\ujnorma,\acute{\ujnorma}}(a,b,r)=
(\cosh \acute \ujnorma)\Sin(b^2-r^2)\frac{
r  Z +S H
}{\sqrt P}
\]
where
\[S\equiv \frac{\sinh \acute \ujnorma}{\cosh \acute \ujnorma},\]
\[Z\equiv\sqrt{\breve c^2+\breve d^2- \left(  \frac{\sinh \acute \ujnorma}{\cosh \acute \ujnorma}\hat b\right)^2},\]
\[P\equiv\breve c^2+\breve d^2
=1+2(1-t)(\sin\theta)^2\ujnorma^2\reC+(\sin\theta)^2\ujnorma^2\reD+(1-t)^2(\sin\theta)^2\ujnorma^4\reC^2>0,\]
\[H\equiv\frac{\hat b\,\breve c}{\sin\theta} (1-t)+  (\breve c^2+\breve d^2)t
=1+(1-t)(1-t(\cos\theta)^2)\ujnorma^2\reC+t(\sin\theta)^2\ujnorma^2\reD>0,\]
such that the arguments of $\reC$ and $\reD$ should be $b^2-r^2$.

(b)
$\SH_{\ujnorma,\acute{\ujnorma}}(a,b,r)$
is log-able if and only if
\[ (\cosh\acute\ujnorma )\Cos(b^2-r^2)+(\sinh \acute\ujnorma)r\Sin(b^2-r^2)\frac{\hat b\breve d S +  \breve c Z}{P^2} >-1.\]
where $S, Z , P$ are as before.

(c) In the log-able case:
Let $\pi_{23}$ denote the projection to the second and third coordinates.
Then, for the Jacobian,
\[
\frac{\partial\left( \pi_{23}\circ\Xi^{\mathrm{PH}}\circ\log\circ \SH_{\ujnorma,\acute{\ujnorma}} ( \ujnorma t\cos\theta,\ujnorma t\sin\theta,\ujnorma (1-t) )\right)}{\partial(t,\theta)}=
\]
\[(\cos\theta)\cdot\AC\left((\cosh\acute\ujnorma )\Cos(b^2-r^2)+(\sinh \acute\ujnorma)r\Sin(b^2-r^2)\frac{\hat b\breve d S +  \breve c Z}{P^2}
\right)\cdot \Sin(b^2-r^2)\cdot\]
\[\]
\[\frac{\ujnorma^2(\cosh\acute\ujnorma)}{P\sqrt P}\Biggl( rSF +\frac{S^2FH+ t(1-S^2)P^2}Z
\Biggr),\]
where $S,Z,P,H$ are as before, and
\[F\equiv (\cos\theta)^2+t(\sin\theta)^2\ujnorma^2\reC+t^2(\sin\theta)^4\ujnorma^4\reW,\]
and the arguments of $\reC$ and $\reW$ should be $b^2-r^2$.

In particular, this is non-vanishing if $\cos\theta\neq0$.
\begin{proof}
(a) is direct computation.
(b) is a consequence of Lemma \ref{lem:logreal}(a).
(c) The formula for $\log$ can be applied, then direct computation yields the result.
Regarding non-vanishing, we can see that beyond $\cos\theta$, all multiplicative terms are positive.
\end{proof}
\end{lemma}
\snewpage

Let $\SEH_{\ujnorma,\acute{\ujnorma}}$ denote either $\SE_{\ujnorma,\acute{\ujnorma}}$ or $\SH_{\ujnorma,\acute{\ujnorma}}$.
It would be desirable to show that

(SEHC$_1$) ``If $0<\ujnorma,\acute{\ujnorma}$ and $\ujnorma+\acute{\ujnorma}\leq\pi$ then
\[
\|\log\circ\SEH_{\ujnorma,\acute{\ujnorma}}\|_2:X^{\mathrm{nn}}_{\ujnorma}\rightarrow[0,+\infty)
\]
has no local extremum.''

A statement which would provide this can be formulated as follows.
For $(a,b,r)\in\mathbb R^2\times[0,+\infty)$, $\sqrt{a^2+b^2},r>0$
let $\nabla(a,b,r)=\left( \frac{a}{\sqrt{a^2+b^2}}, \frac{b}{\sqrt{a^2+b^2}},1\right)$.
One can recognize it as a sort of the gradient of norm.
Then a stronger statement is

(SEHC$_2$) ``If $0<\ujnorma,\acute{\ujnorma}$ and $\ujnorma+\acute{\ujnorma}\leq\pi$ then
\[
\nabla\left(\Xi^{\mathrm{PH}}\circ\log\circ\SEH_{\ujnorma,\acute{\ujnorma}}( \ujnorma t\cos\theta,\ujnorma t\sin\theta,\ujnorma (1-t) )\right)
\]
and
\begin{multline}\notag
\frac{\partial\left(\Xi^{\mathrm{PH}}\circ\log\circ\SEH_{\ujnorma,\acute{\ujnorma}}( \ujnorma t\cos\theta,\ujnorma t\sin\theta,\ujnorma (1-t) )\right)}{\partial\theta}\times
\\
\times\frac{\partial\left(\Xi^{\mathrm{PH}}\circ\log\circ\SEH_{\ujnorma,\acute{\ujnorma}}( \ujnorma t\cos\theta,\ujnorma t\sin\theta,\ujnorma (1-t) )\right)}{\partial t}
\end{multline}
(the standard vectorial product is meant) are not nonnegatively proportional.''

An even stronger possible statement is as follows.
For $(a,b,r)\in\mathbb R^2\times[0,+\infty)$, let $\pi_{(12)3}$ be defined by
$\pi_{(12)3}(a,b,r)=(\sqrt{a^2+b^2},r)$

(SEHC$_3$) ``If $0<\ujnorma,\acute{\ujnorma}$ and $\ujnorma+\acute{\ujnorma}\leq\pi$, then the Jacobian
\[
\frac{\partial\left( \pi_{(12)3}\circ\Xi^{\mathrm{PH}}\circ\log\circ \SEH_{\ujnorma,\acute{\ujnorma}} ( \ujnorma t\cos\theta,\ujnorma t\sin\theta,\ujnorma (1-t) )\right)}{\partial(t,\theta)}
\]
is nonvanishing.''

This latter statement can be established in particular cases (for $\ujnorma$ and $\acute{\ujnorma}$)
but I do not know a general argument.
\begin{remark}
\plabel{rem:indest}
For $\SE$,

\[0\ll\frac{
\dfrac{\partial\left( \pi_{(12)3}\circ\Xi^{\mathrm{PH}}\circ\log\circ \SE_{\ujnorma,\acute{\ujnorma}} ( \ujnorma t\cos\theta,\ujnorma t\sin\theta,\ujnorma (1-t) )\right)}{\partial(t,\theta)}
}{(\sin\theta)(\cos\theta)(1-t)\AC\left(\Cos(b^2-r^2)\cos\frac{\hat b\acute N}{\sqrt{ \hat a^2+\hat b^2}}
-b\Sin(b^2-r^2)\sin\frac{\hat b\acute N}{\sqrt{ \hat a^2+\hat b^2}}\right)^2}\ll+\infty
\]
seems to be the case  (for $\ujnorma$ and $\acute{\ujnorma}$ fixed).
For $\SH$,
\[0\ll\frac{
\dfrac{\partial\left( \pi_{(12)3}\circ\Xi^{\mathrm{PH}}\circ\log\circ \SE_{\ujnorma,\acute{\ujnorma}} ( \ujnorma t\cos\theta,\ujnorma t\sin\theta,\ujnorma (1-t) )\right)}{\partial(t,\theta)}
}{(\sin\theta)(\cos\theta)(t^2+(\cos\theta)^2)}\cdot\sqrt{t^2+(\sin\theta)^2}\ll+\infty
\]
seems to hold.
These estimates can be established in several special cases.
\qedremark
\end{remark}
\scleardoublepage\section{The balanced critical BCH case}
\plabel{sec:Bound}
Let $\alpha_0\in(-\pi,\pi)$.
Consider
\[\mathcal B_{\frac{\pi-\alpha_0}2,\frac{\pi+\alpha_0}2}=
\left\{\log(\exp(A)\exp(B))\,:\,A,B\in\mathrm M_2(\mathbb R),\|A\|_2\leq \frac{\pi-\alpha_0}2,\|B\|_2\leq \frac{\pi+\alpha_0}2\right\}.  \]

This is not a closed set.
The reason is that $\log$ is not defined at $-\Id$, thus the usual compactness argument does not apply.
This failure in well-definedness affects only two cases, $(A,B)=\left(\frac{\pi-\alpha_0}2\tilde I,\frac{\pi+\alpha_0}2\tilde I \right)$ and
$(A,B)=\left(-\frac{\pi-\alpha_0}2\tilde I,-\frac{\pi+\alpha_0}2\tilde I \right)$.
Yet, the closure $\overline {\mathcal B_{\frac{\pi-\alpha_0}2,\frac{\pi+\alpha_0}2}}$ is larger than $\mathcal B_{\frac{\pi-\alpha_0}2,\frac{\pi+\alpha_0}2}$ by several quasi $\log$-s of $-\Id_2$, as Theorem \ref{th:discont22} shows.
(In this setting, $\log(\exp(A)\exp(B))$ is the same as $\BCH(A,B)$ except the latter one is
well-defined everywhere but still non-continuous at the critical points.
Thus, a quite legitimate version of the set above is its union with $\{\pi \tilde I, -\pi \tilde I \}$.)

\begin{theorem}\plabel{th:exof}
(a) The elements of
\[\overline {\mathcal B_{\frac{\pi-\alpha_0}2,\frac{\pi+\alpha_0}2}}\setminus \mathcal B_{\frac{\pi-\alpha_0}2,\frac{\pi+\alpha_0}2}=
\partial {\mathcal B_{\frac{\pi-\alpha_0}2,\frac{\pi+\alpha_0}2}}\setminus \mathcal B_{\frac{\pi-\alpha_0}2,\frac{\pi+\alpha_0}2}
\] are
exactly the elements $B$ which are are of shape
\[B=b\tilde I+c\tilde J+d\tilde K\]
with
\[b^2-c^2-d^2=\pi^2\]
and
\[\pi\leq |b|+\sqrt{c^2+d^2}\leq\pi\sqrt{\frac{\pi-|\alpha_0|+2\cos\frac{\alpha_0}2}{\pi-|\alpha_0|-2\cos\frac{\alpha_0}2}}.\]

(b) The elements of
\[\partial {\mathcal B_{\frac{\pi-\alpha_0}2,\frac{\pi+\alpha_0}2}}\cap\mathcal B_{\frac{\pi-\alpha_0}2,\frac{\pi+\alpha_0}2}\]
are all of shape $\log((\exp A)(\log B))$ with $\|A\|_2=\frac{\pi-\alpha_0}2$ and $\|B\|_2=\frac{\pi+\alpha_0}2$
such that $(A,B)$ is a BCH minimal pair.

(c) The interior of $\mathcal B_{\frac{\pi-\alpha_0}2,\frac{\pi+\alpha_0}2}$ is connected, containing $0$.
\begin{proof}
(a) The set should be closed to conjugation by orthogonal matrixes, and by continuity of $\exp$, theirs exponentials are $-\Id_2$.
Then only the norms are of question but Theorem \ref{th:discont22} takes care of that.
(b) , (c) These follow from the openness of $\exp$ as long as the domain is restricted to
the spectrum in $\{z\in\mathbb C\,:\,|\Ima z|<\pi\}$.
\end{proof}
\end{theorem}
The set $\mathcal B_{\frac{\pi-\alpha_0}2,\frac{\pi+\alpha_0}2}$ is closed for conjugation by orthogonal rotations,
thus it can be visualized through $\Xi^{\mathrm{PH}}$.
Then $\partial\Xi^{\mathrm{PH}}\left( \mathcal B_{\frac{\pi-\alpha_0}2,\frac{\pi+\alpha_0}2} \right)$, expectedly a
$2$-dimensional object, describes $\left( \mathcal B_{\frac{\pi-\alpha_0}2,\frac{\pi+\alpha_0}2} \right)$.
In the previous section we have devised several restrictions for (infinitesimally) BCH minimal pairs.
In fact, we find that  $\partial\Xi^{\mathrm{PH}}\left( \mathcal B_{\frac{\pi-\alpha_0}2,\frac{\pi+\alpha_0}2} \right)$
must be contained in the union of continuous images of finitely many, at most  $2$ dimensional, manifolds,
which can described explicitly.
Despite this, computation with these objects is tedious.

Therefore we restrict our attention to the case $\mathcal B_{\frac{\pi}2,\frac{\pi}2}$ (i. e. $\alpha_0=0$),
which is the most important for us.
Before giving any argument, by the following Figure \ref{fig:figBA07}, let us show what we will obtain:

\begin{figure}[H]
\centering
   \begin{subfigure}[b]{6.0in}
    \includegraphics[trim=0 200 0 250,clip,width=6.0in]{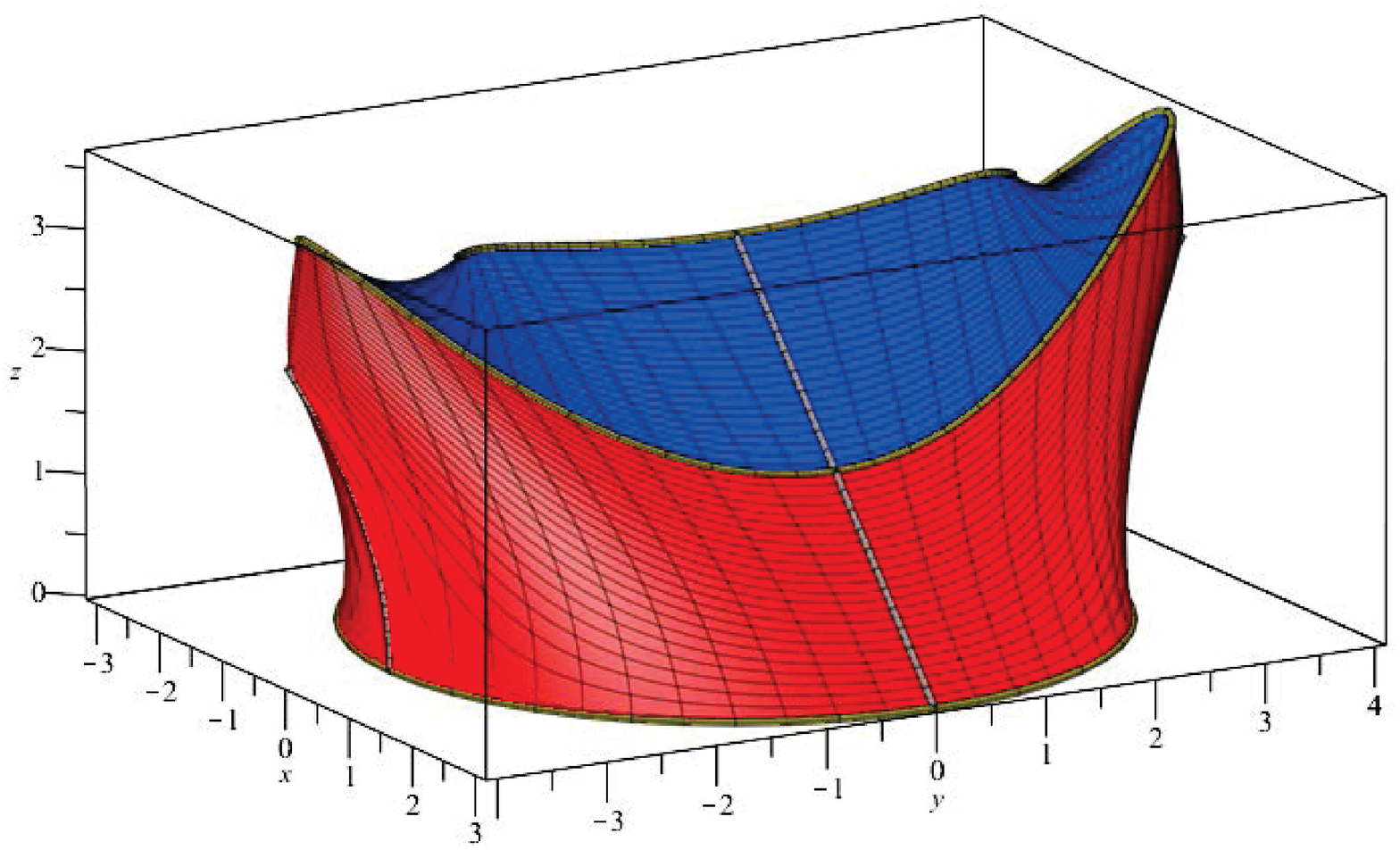}
    \caption*{Fig. \ref{fig:figBA07}(a) $\partial\Xi^{\mathrm{PH}}\mathcal B_{\frac\pi2,\frac\pi2}$}
  \end{subfigure}
\phantomcaption
\plabel{fig:figBA07}
\end{figure}
\begin{figure}[H]
  \ContinuedFloat
   \begin{subfigure}[b]{2.5in}
    \includegraphics[width=2.5in]{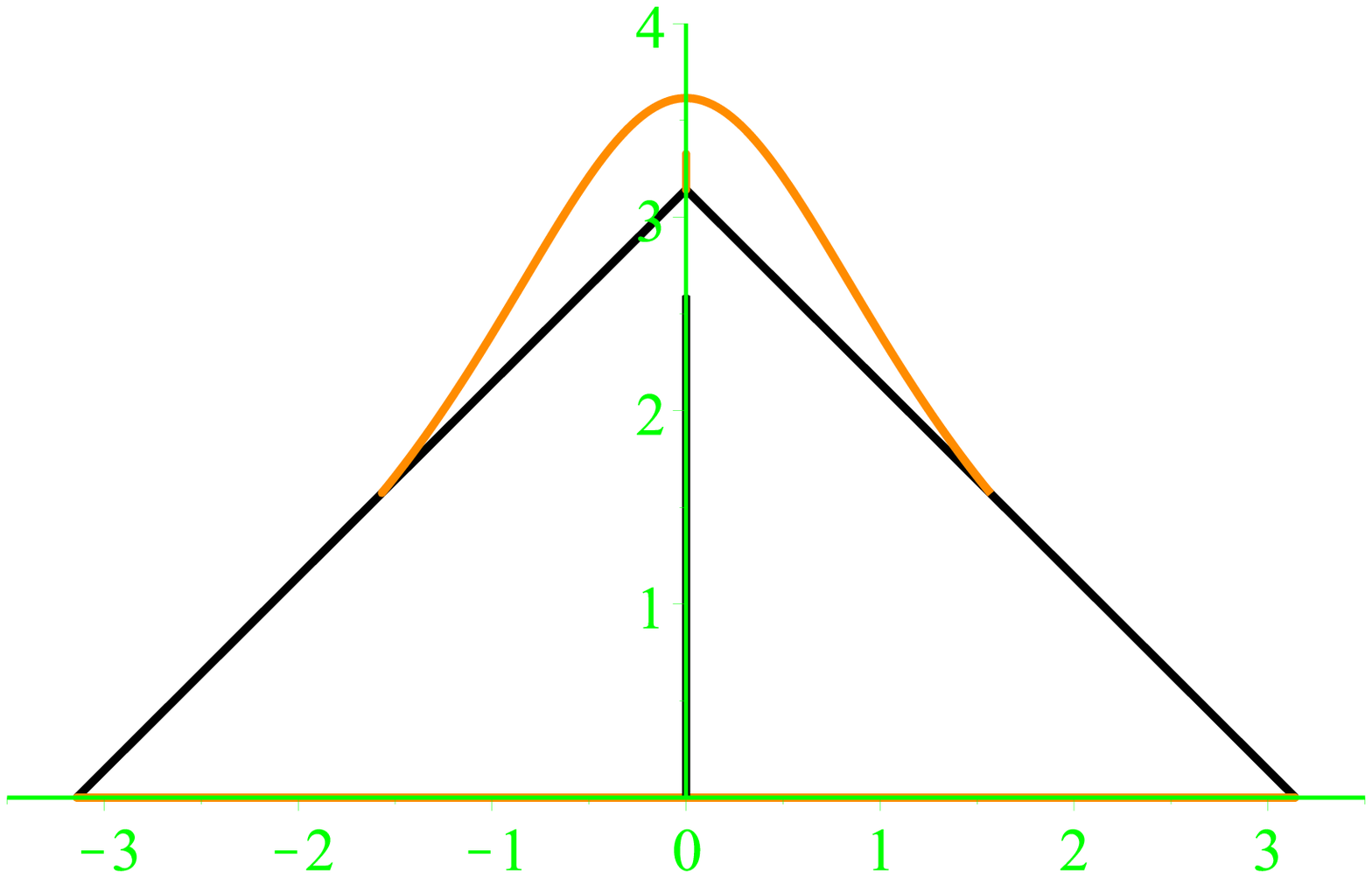}
    \caption*{\ref{fig:figBA07}(b) $xz$-projection}
  \end{subfigure}
   \begin{subfigure}[b]{3.125in}
    \includegraphics[width=3.125in]{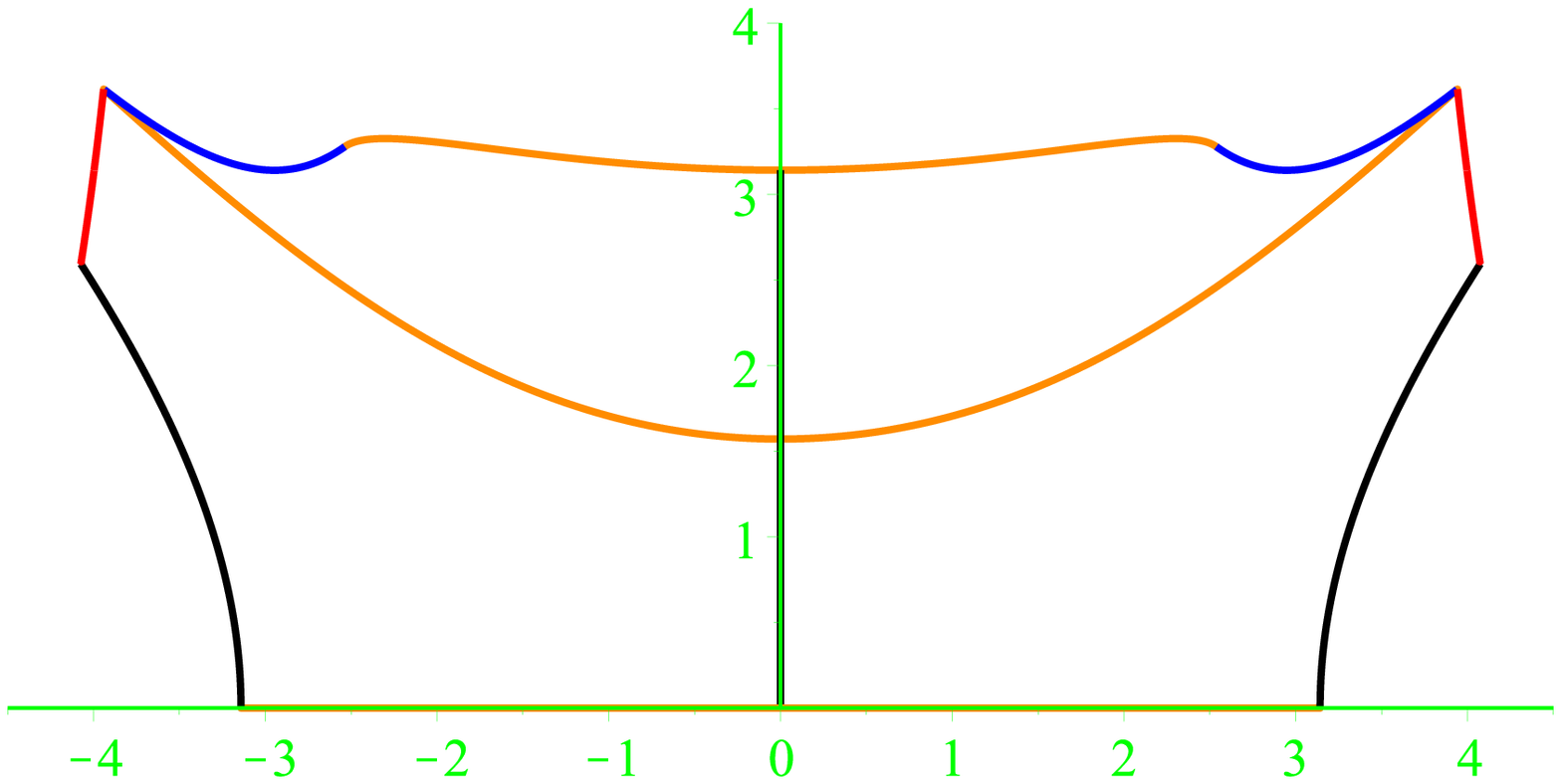}
    \caption*{\ref{fig:figBA07}(c) $yz$-projection}
  \end{subfigure}
\phantomcaption
\end{figure}

$\partial\Xi^{\mathrm{PH}}\mathcal B_{\frac\pi2,\frac\pi2}$ is a ``wedge cap''.
What see is the following (note that $\ujnorma\equiv\pi/2$ here):

$\bullet$ In the middle top, we see the Schur bihyperbolic ridge, the image of the
(infinitesimally) minimal pairs from $\partial^{\mathrm{hyp}}\mathcal S_\ujnorma\times \partial^{\mathrm{hyp}}\mathcal S_\ujnorma$.

$\bullet$ Joining it, we see Schur parabolical segments (in the $xz$-plane), the image of the
(infinitesimally) minimal pairs from $\eth^{\mathrm{par}}\mathcal S_\ujnorma\times \eth^{\mathrm{par}}\mathcal S_\ujnorma$.

$\bullet$ The middle ridge is the  Schur elliptic-hyperbolic ridge, the image of the
(infinitesimally) minimal pairs from
$\partial^{\mathrm{ell}}\mathcal S_\ujnorma\times \partial^{\mathrm{hyp}}\mathcal S_\ujnorma$
or
$\partial^{\mathrm{ hyp}}\mathcal S_\ujnorma\times \partial^{\mathrm{ell}}\mathcal S_\ujnorma$
(transposition invariance shows equality).

$\bullet$ In the front and back see we the closure singularity segments (in the $yz$-plane),
the images coming from limiting to $\pm(\frac\pi2\tilde I, \frac\pi2\tilde I)$

$\bullet$  In the very bottom, we see the  Schur  bielliptic rim, the image of the
(infinitesimally) minimal pairs from
$\partial^{\mathrm{ell1}}\mathcal S_\ujnorma\times \partial^{\mathrm{ell1}}\mathcal S_\ujnorma$
and
$\partial^{\mathrm{ell*}}\mathcal S_\ujnorma\times \partial^{\mathrm{ell*}}\mathcal S_\ujnorma$,
but of $\pm(\frac\pi2\tilde I, \frac\pi2\tilde I)$ which are exceptional
(that is the former case is eliminated).

$\bullet$ The upper, blue area is the  Schur  smooth-hyperbolic area, the image of the
(infinitesimally) minimal pairs from $\mathcal S^{\mathrm{nn}}_\ujnorma\times \partial^{\mathrm{hyp}}\mathcal S_\ujnorma$
or $\partial^{\mathrm{hyp}}\mathcal S_\ujnorma\times\mathcal S^{\mathrm{nn}}_\ujnorma$.
(Again, note transposition invariance.)

$\bullet$ The lower, red area is the  Schur  smooth-elliptic area, the image of the
(infinitesimally) minimal pairs from $\mathcal S^{\mathrm{nn}}_\ujnorma\times \partial^{\mathrm{ell}}\mathcal S_\ujnorma$
or $\partial^{\mathrm{ell}}\mathcal S_\ujnorma\times\mathcal S^{\mathrm{nn}}_\ujnorma$.

\begin{remark}
\plabel{rem:subcritfig}
The situation with $\partial\Xi^{\mathrm{PH}}\mathcal B_{\ujnorma,\ujnorma}$ with $0<\ujnorma<\frac\pi2$ is similar, cf. Figure \ref{fig:figBA8},
\begin{figure}[H]
\centering
   \begin{subfigure}[b]{3.4in}
    \includegraphics[trim=0 100 0 125,clip,width=3.4in]{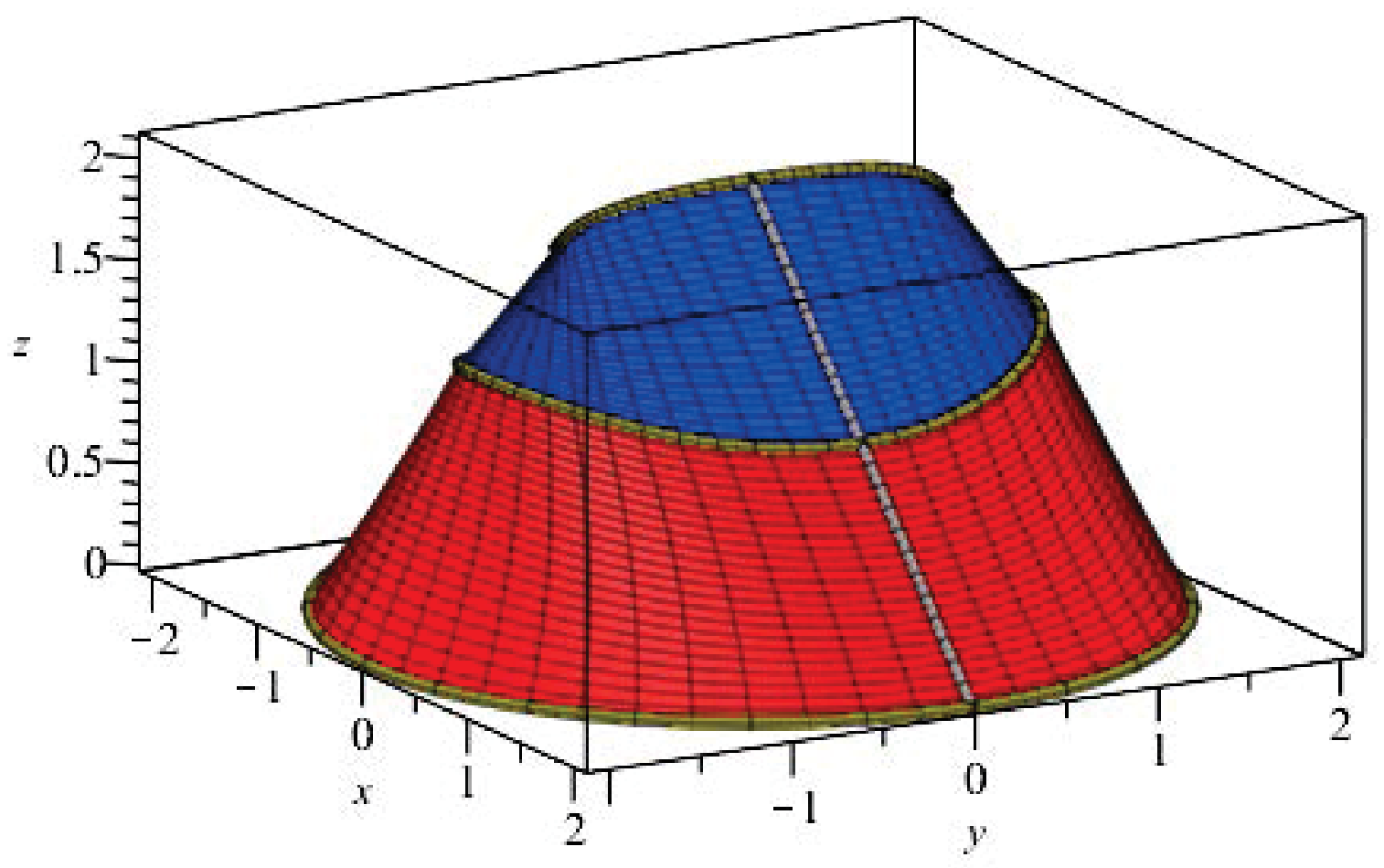}
    \caption*{Fig. \ref{fig:figBA8}(a) $\partial\Xi^{\mathrm{PH}}\mathcal B_{\frac\pi3,\frac\pi3}$}
  \end{subfigure}
\phantomcaption
\plabel{fig:figBA8}
\end{figure}
\begin{figure}[H]
  \ContinuedFloat
   \begin{subfigure}[b]{2.5in}
    \includegraphics[width=2.5in]{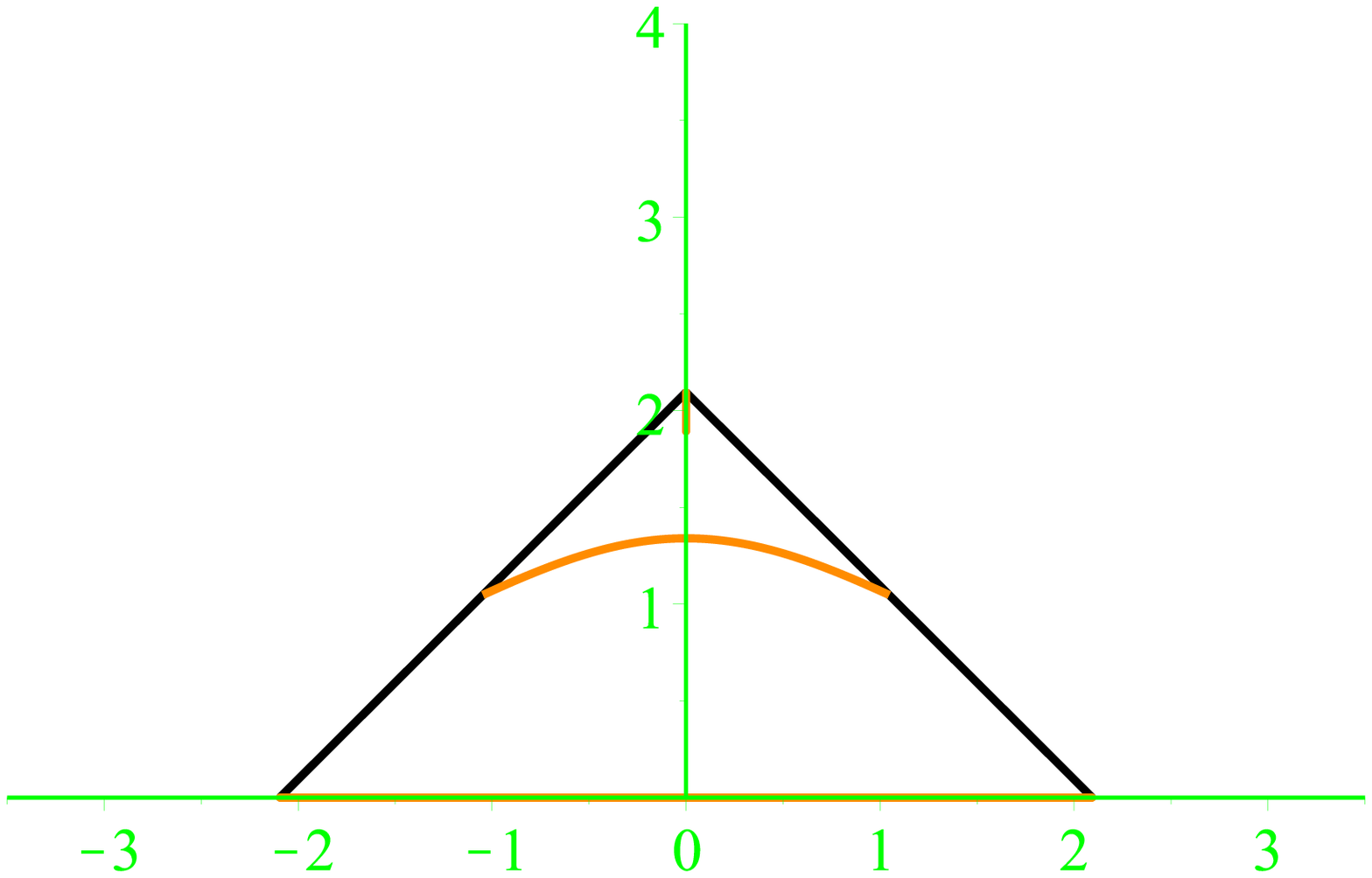}
    \caption*{\ref{fig:figBA8}(b) $xz$-projection}
  \end{subfigure}
   \begin{subfigure}[b]{3.125in}
    \includegraphics[width=3.125in]{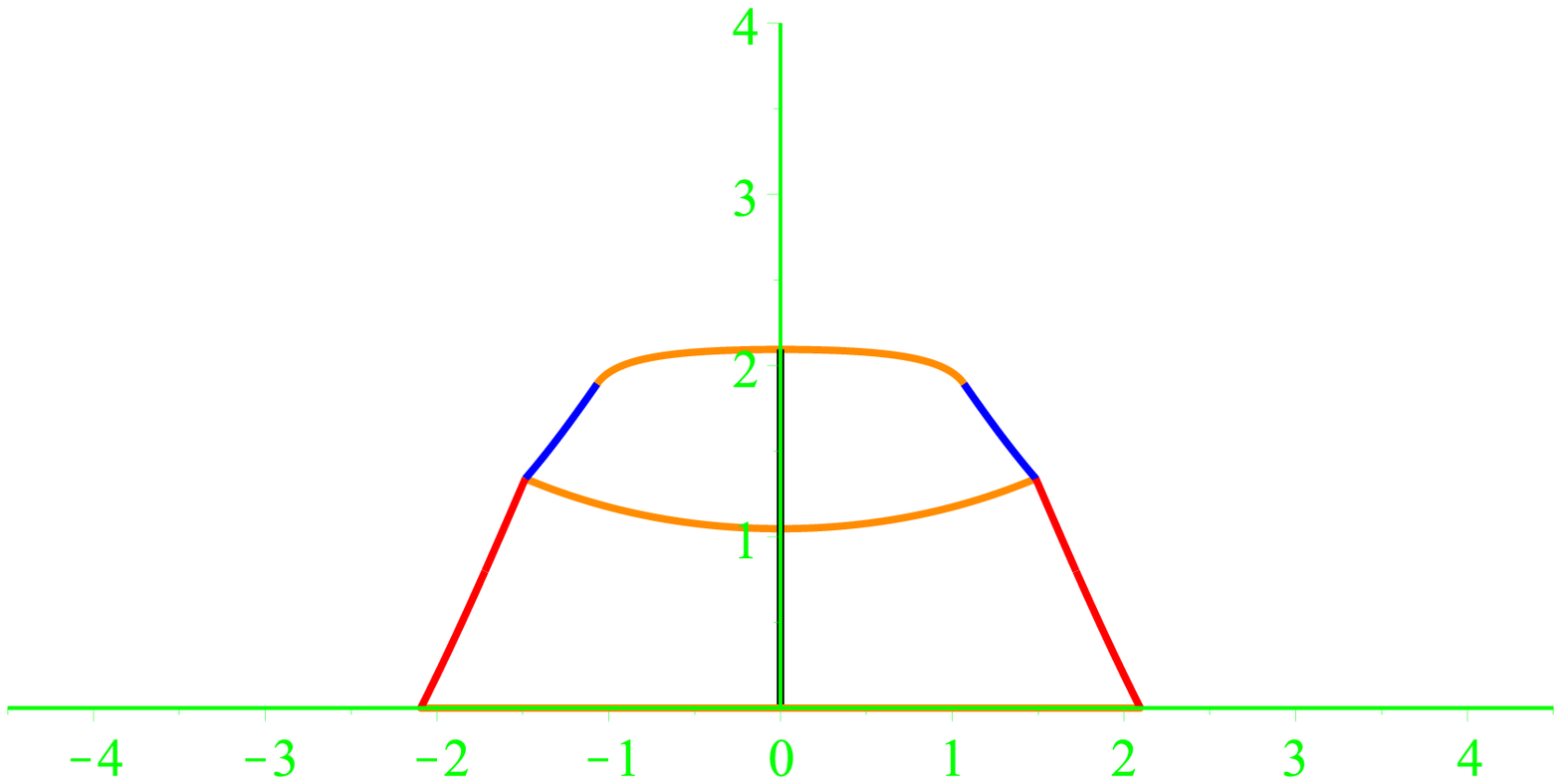}
    \caption*{\ref{fig:figBA8}(c) $yz$-projection}
  \end{subfigure}
\phantomcaption
\end{figure}
\noindent except the ``closure'' singularity does not develop.
\qedremark
\end{remark}

Let us now make the ``statement'' of Figure \ref{fig:figBA07}  more precise.
Let us fix the choice $\ujnorma=\pi/2$.

$\bullet$ We define the Schur bihyperbolic parametrization as
the map
\[\psi\in\left[-2\arctan\ujnorma ,2\arctan\ujnorma \right]\mapsto
\Xi^{\mathrm{PH}}\log\left(
\exp\left(\ujnorma\tilde J\right)
\exp\left(\ujnorma( (\cos\psi)\tilde J +(\sin\psi)\tilde K)\right)
\right)
.\]

$\bullet$ We define the Schur parabolic parametrization as
the map
\[(\sigma,t)\in \{1,-1\}\times(-1,1)\mapsto\Xi^{\mathrm{PH}}\left(
\sigma\cdot\ujnorma\cdot\bem1&\\&t\eem
\right)
.\]

$\bullet$ We define the Schur elliptic*-hyperbolic parametrization as
the map
\[\psi\in(0,\pi)\cup(\pi,2\pi)\mapsto
\Xi^{\mathrm{PH}}\log\left(
\exp\left(\ujnorma( (\cos\psi)\tilde \Id_2 +(\sin\psi)\tilde I)\right)
\exp\left(\ujnorma\tilde J\right)
\right)
.\]

$\bullet$ We define the closure parametrization as
the map
\[ (\sigma,r)\in\{-1,1\}\times\left[0, \frac{2\pi}{\sqrt{\pi^2-4}}\right]\mapsto \left(0,\sigma \cdot \sqrt{\pi^2+r^2},r \right)
.\]

$\bullet$ We define the Schur bielliptic parametrization as
the map
\[\psi\in\left(-\frac\pi2,\frac\pi2\right)\cup\left(\frac\pi2,\frac{3\pi}2\right)\mapsto
\Xi^{\mathrm{PH}} \left(\ujnorma( (\cos\psi)\tilde \Id_2 +(\sin\psi)\tilde I)
\right)
.\]

$\bullet$ We define the Schur smooth-hyperbolic parametrization as
the map
\begin{multline*}
(a,b,r)\in X_\ujnorma\mapsto
\Xi^{\mathrm{PH}}\log\Biggl(
\exp\left(a\Id_2+b\tilde I+r\tilde J\right)\cdot
\\
\cdot
\exp\left(  \left(\sqrt{\breve c^2+\breve d^2- \left(  \frac{\sinh   \ujnorma}{\cosh   \ujnorma}\hat b\right)^2}\Id_2
 -\frac{\sinh   \ujnorma}{\cosh   \ujnorma}\hat b\tilde I \right) \frac{\breve c\tilde \Id_2 +\breve d\tilde I}{\breve c^2+\breve d^2} \ujnorma\tilde J\right)
\Biggr)
.
\end{multline*}
(Here we have used the abbreviations of Lemma \ref{lem:momentkif}.)

$\bullet$ Similarly, we define the Schur smooth-elliptic parametrization as
the map
\[(a,b,r)\in X_\ujnorma\mapsto
\Xi^{\mathrm{PH}}\log\left(
\exp\left(a\Id_2+b\tilde I+r\tilde J\right)
\exp\left(  \frac{\hat{a}\Id_2+\hat{b}\tilde I}{\sqrt{\hat a^2+\hat b^2}}  \ujnorma\right)
\right)
.\]

We will call these maps as the canonical parametrizations in the
$\left(\frac\pi2,\frac\pi2\right)$ case.

\begin{remark}
\plabel{rem:subcrit}
For $0<\ujnorma<\frac\pi2$, the canonical parametrizations can defined similarly in the
$\left(\ujnorma,\ujnorma\right)$ case, except the situation is simpler:
These is no closure parametrization but the Schur bielliptic parametrization
can be defined fully for $[0,2\pi]\modu 2\pi$.
\qedremark
\end{remark}
\begin{theorem}
\plabel{th:critwedge}
(a) Every element of   $\partial\Xi^{\mathrm{PH}}\mathcal B_{\frac\pi2,\frac\pi2}$ occurs the image
of BCH minimal pair with respect norms $\ujnorma_1=\ujnorma_2=\pi/2$ with the exception of the closure singularities.
Conversely, every image of a such a BCH minimal pair or a point of the closure singularity is in
$\partial\Xi^{\mathrm{PH}}\mathcal B_{\frac\pi2,\frac\pi2}$.

(b) The canonical parametrizations taken together map to  $\partial\Xi^{\mathrm{PH}}\mathcal B_{\frac\pi2,\frac\pi2}$
bijectively, and the images fit together topologically as suggested by Figure \ref{fig:figBA07}.
\begin{proof}
First we prove that non BCH minimal points map to the interior.
Indeed as least one component can be replaced by an element of norm less then $\pi/2$.
As the exponential map and the logarithm will be open in these circumstances, perturbing that the entry
yields an open neighborhood in the image.
This proves the first part of (a).
The second part, and, in fact, the rest of the statement, follows by topological reasons if we prove that
the images of the canonical maps fit together topologically as a half-sphere.
For the 1-dimensional canonical parametrizations injectivity and topological incidences are easy to check.
Next one can prove that the 2-dimensional parametrizations limit on their boundary in an
expected manner.
In that blow-ups in ACKB or AHP are instrumental, except at the pairs $\pm(\frac\pi2\tilde I,\frac\pi2\tilde I)$,
continuity (in fact, well-definedness) breaks down.
Yet, we know that the limits exponentiate to $\pm\pi\tilde I$, and
we know the range of the norm by Theorem \ref{th:discont22}.
By this, the closure singularies can be recovered.
Also, the sign relations in (the trace) coordinate $x$ are easy due to the BCH formula.
To the 1-dimensional canonical parametrization we can add the $a=0$ cases of the
smooth-hyperbolic and smooth-elliptic canonical parametrizations.
By the sign relations in $x$, then it is sufficient to prove that
$yz$-projections of the $a\neq0$ parts of smooth-hyperbolic and smooth-elliptic canonical parametrizations
have non-vanishing Jacobians.
Thus the fill out the corresponding regions in Figure \ref{fig:figBA07}(c) as they should.
(That also shows that they fit to the closure singularities properly.)
\end{proof}
\end{theorem}
\begin{theorem}
\plabel{th:supequal22}
\[G\left(\frac\pi2\right)\equiv \frac\pi2\exp\frac\pi2=\left(\sup_{\substack{A,B\in\mathrm M_2(\mathbb R),\\ \|A\|_2=\|B\|_2\leq\pi/2  }} \|\log((\exp A)(\exp B))\|_2\right) .\]
(Here the exceptional cases $A=B=\frac\pi2\tilde I$ and $A=B=-\frac\pi2\tilde I$  do not participate in the supremum.
But, writing $\BCH(A,B)$, they could.)

\begin{proof}[Indication of proof]
Thus, one has to prove
\[\sup\{\|S\|_2\,:\, S\in \mathcal B_{\frac{\pi}2,\frac{\pi}2} \}=\frac\pi2\exp\frac\pi2.\]

According to Theorem \ref{th:exof}, we can replace $\mathcal B_{\frac{\pi}2,\frac{\pi}2} $ by
$\partial\mathcal B_{\frac{\pi}2,\frac{\pi}2}$,
in fact, by $\partial\mathcal B_{\frac{\pi}2,\frac{\pi}2}\cap \mathcal B_{\frac{\pi}2,\frac{\pi}2}$.
By various estimates, we can localize the maxima near to tip of the horns of Figure \ref{fig:figBA07}(a).
By further estimates, one show that the direction of gradients for the Schur smooth-elliptic
parts and  Schur smooth-hyperbolic parts near the tip are inconsistent to the maxima.
Thus it remains to optimize on the Schur elliptic*-hyperbolic part.
There, ultimately, we obtain the tips as expected for the maxima.
\end{proof}
\end{theorem}
\begin{remark}\plabel{rem:tec1}
In the previous proof, the critical part is the reduction to the 1-dimensional boundary pieces.
In a less ad hoc way it is also achieved by showing, say, (SECH${}_3$)
for $\ujnorma=\acute\ujnorma=\frac\pi2$, which is quite doable as a special case.
For 1-dimensional optimizations there are just several methods.
\qedremark
\end{remark}
It seems reasonable to expect that

(SM$_0$) ``If $0<\ujnorma,\acute{\ujnorma}$ and $\ujnorma+\acute{\ujnorma}\leq\pi$, then
the maximum of $\|E\|_2$ for $E\in\mathcal B_{\ujnorma,\acute{\ujnorma}}$ is taken for a traceless pair of matrices.''

Indeed, large norm seems to come from non-commutativity, and thus including tracial parts seems to be quite pointless.
Again, this is quite provable in several special cases, but I do not know a general proof.
It may be possible, however, that a quite simple argument will suffice.
\begin{remark}\plabel{rem:tec2}
The localizations of the maxima, however, are not trivial even in the case $\ujnorma=\acute\ujnorma$.
Indeed, in this setting, there is a critical value
\[C_0=\pi\cdot0.392744\ldots\]
such that
for $0<\ujnorma=\acute\ujnorma<C_0$, the maxima are taken on the Schur bihyperbolic ridge
(cf. Figure \ref{fig:figBA8});
and for $C_0<\ujnorma=\acute\ujnorma\leq\pi$, the maxima are taken on
the Schur elliptic*-hyperbolic ridge; and for $p=C_0$, on both.
\qedremark
\end{remark}

Thus one can expect that

(SM$_+$) ``If $0<\ujnorma,\acute{\ujnorma}$ and $\ujnorma+\acute{\ujnorma}\leq\pi$ but $\ujnorma+\acute{\ujnorma}$ is sufficiently large,  then
the maximum of $\|E\|_2$ for $E\in\mathcal B_{\ujnorma,\acute{\ujnorma}}$ is taken for a traceless Schur elliptic-hyperbolic pair,
the norm of the conform-involution is less or equal than the norm of conform-skew-involution.''

However, this needs further exploration.
\begin{remark}\plabel{rem:tec3}
If this is true, then Example \ref{ex:loxocomp} seems to be sharp case.
\qedremark
\end{remark}
\snewpage
\scleardoublepage\section{Principal disks and logarithm}
\plabel{sec:PrincLog}

\begin{commentx}
For the sake of the next statements using $\tilde a,\tilde b$ instead of $a,b$ would be more appropriate,
but it is probably better to keep the notation simple.
\end{commentx}
\begin{lemma}\plabel{lem:lognorm}
Suppose that $A$ is a $\log$-able real $2\times 2$ matrix with principal disk
\[\PD(A)=\Dbar(a+\mathrm ib,r).\]
In that case,
\begin{equation}
\|\log A\|_2=f_{\mathrm{CA}}(a,b,r)+f_{\mathrm{RD}}(a,b,r),
\plabel{eq:2lognorm}
\end{equation}
and
\begin{equation}
\left\lfloor\log A\right\rfloor_2=f_{\mathrm{CA}}(a,b,r)-f_{\mathrm{RD}}(a,b,r),
\plabel{eq:2logconorm}
\end{equation}
where
\[f_{\mathrm{CA}}(a,b,r)= \sqrt{\Bigl(\log \sqrt{a^2+b^2-r^2} \Bigr)^2+\left(\frac{b\AC\left(\dfrac{a}{\sqrt{a^2+b^2-r^2}}\right)}{\sqrt{a^2+b^2-r^2}}\right)^2}\]
and
\[f_{\mathrm{RD}}(a,b,r)=\frac{r\AC\left(\dfrac{a}{\sqrt{a^2+b^2-r^2}}\right)}{\sqrt{a^2+b^2-r^2}}.\]

In particular, if $\det A=1$, then  $a^2+b^2-r^2=1$, and
$f_{\mathrm{CA}}(a,b,r)=\AC(a)b$, $f_{\mathrm{RD}}(a,b,r)=\AC(a)r$.
\begin{proof}
This is just the combination of
\eqref{eq:log2} and \eqref{eq:2norm}--\eqref{eq:2conorm}, computed explicitly.
\end{proof}
\end{lemma}
%\snewpage
\begin{theorem}\plabel{th:monotone}
Suppose that $A_1,A_2$ are $\log$-able real $2\times 2$ matrices such that
\[\PD(A_1)\subset \PD(A_2).\]
Then
\begin{equation}
\|\log A_1\|_2\leq \|\log A_2\|_2.
\plabel{eq:monotone}
\end{equation}
and
\begin{equation}
\lfloor\log A_1\rfloor_2\geq \lfloor\log A_2\rfloor_2.
\plabel{eq:comonotone}
\end{equation}
\begin{proof}[Remark] The monotonicity of $\|\cdot\|_2$ is strict, except if $\PD(A_1)$
and $ \PD(A_2)$ are centered on the real line
and $\sup\{|\log x|\,:\,x\in\mathbb R\cap \PD(A_1)\}=\sup\{|\log x|\,:\,x\in\mathbb R\cap \PD(A_2)\}$.

The monotonicity of $\lfloor\cdot\rfloor_2$ is strict, except if $\PD(A_1)$
and $ \PD(A_2)$ are centered on the real line
and $\inf\{|\log x|\,:\,x\in\mathbb R\cap \PD(A_1)\}=\inf\{|\log x|\,:\,x\in\mathbb R\cap \PD(A_2)\}$.
\end{proof}
\begin{proof}
Let $f(a,b,r)$ denote the functional expression on the right side of \eqref{eq:2lognorm}.
Then it is a straightforward but long computation to check the identity
\begin{equation}
\left(\frac{\partial f(a,b,r)}{\partial r}\right)^2-
\left(\frac{\partial f(a,b,r)}{\partial a}\right)^2-
\left(\frac{\partial f(a,b,r)}{\partial b}\right)^2
= \left(\frac{f(a,b,r)}{f_{\mathrm{CA}}(a,b,r)} \frac{b\AS\left(\dfrac{a}{\sqrt{a^2+b^2-r^2}}\right)}{{a^2+b^2-r^2}}\right)^2 .
\plabel{eq:conical}
\end{equation}
This is valid, except if $b=0$ and $a=\sqrt{1+r^2}$, the exceptional configurations.
In particular, if $b>0$, then
\[ \left(\frac{\partial f(a,b,r)}{\partial r}\right)^2-
\left(\frac{\partial f(a,b,r)}{\partial a}\right)^2-
\left(\frac{\partial f(a,b,r)}{\partial b}\right)^2
>0 .\]
The  principal disks with $b>0$ form a connected set,
consequently
\begin{equation} \frac{\partial f(a,b,r)}{\partial r}>
\sqrt{\left(\frac{\partial f(a,b,r)}{\partial a}\right)^2+\left(\frac{\partial f(a,b,r)}{\partial b}\right)^2}\plabel{eq:logcone}
\end{equation}
or
\[ \frac{\partial f(a,b,r)}{\partial r}<-
\sqrt{\left(\frac{\partial f(a,b,r)}{\partial a}\right)^2+\left(\frac{\partial f(a,b,r)}{\partial b}\right)^2} \]
should hold globally for $b>0$.
The question is: which one?
It is sufficient to check the sign $\frac{\partial f(a,b,r)}{\partial r}$ at a single place.
Now, it is not hard to check that
\[\frac{\partial f(a,b,r)}{\partial r}\Bigl|_{r=0}=\frac{\AC\left(\frac{a}{\sqrt{a^2+b^2}}\right)}{\sqrt{a^2+b^2}}\]
(except if $a=1,b=0$), which shows that \eqref{eq:logcone} holds.
The meaning of  \eqref{eq:logcone} is that expanding principal disks smoothly with non-real centers
leads to growth in the norm of the logarithm.

Let us return to principal disks $D_i=\PD(A_i)$  in the statement.
If $b_1,b_2>0$, then we can expand the smaller one to the bigger one with non-real centers.
(Indeed, magnify $D_2$ from its lowest point, until the perimeters touch, and then magnify from the touching point.)
This proves the \eqref{eq:monotone} for $b_1,b_2>0$. The general statement follows from the continuity of the norm of the logarithm.
Notice that the norm grows if we can expand through $b>0$.

Regarding \eqref{eq:comonotone}:
Let $f_{\mathrm{co}}(a,b,r)$ denote the functional expression on the right side of \eqref{eq:2logconorm}.
It satisfies the very same equation \eqref{eq:conical} but with $f(a,b,r)$ replaced
by $f_{\mathrm{co}}(a,b,r)$ throughout. However,
\[\frac{\partial f_{\mathrm{co}}(a,b,r)}{\partial r}\Bigl|_{r=0}=
-\frac{\AC\left(\frac{a}{\sqrt{a^2+b^2}}\right)}{\sqrt{a^2+b^2}}.\]
The rest is analogous.
\end{proof}
\end{theorem}
\begin{lemma}\plabel{lem:diskmonotone}
Suppose that  $A_1,A_2$ are real $2\times 2$ matrices. Then
\[\PD(A_1)\subset \PD(A_2)\]
holds if and only if
\[\|A_1+\lambda\Id \|_2\leq \|A_2+\lambda\Id \|_2 \qquad \text{ for all }\lambda\in\mathbb R\]
and
\[\lfloor A_1+\lambda\Id \rfloor_2\geq \lfloor A_2+\lambda\Id \rfloor_2 \qquad
\text{ for all }\lambda\in\mathbb R.\]
\begin{proof}
The norms and co-norms can be read off from the principal disk immediately.
Hence the statement is simple geometry.
\end{proof}
\end{lemma}
\begin{theorem}\plabel{th:diskmonotone}
Suppose that  $A_1,A_2$ are $\log$-able real $2\times 2$ matrices. If
\[\PD(A_1)\subset \PD(A_2),\]
then
\[\PD(\log A_1)\subset \PD(\log A_2).  \]
The monotonicity is strict. Similar statement applies to $\CD$.
\begin{proof}
In this case, the matrices $\mathrm e^{\lambda}A_i$ will also be $\log$-able.
Moreover, $\PD(\mathrm e^{\lambda}A_1)\subset \PD(\mathrm e^{\lambda}A_2)$ holds.
Now, $\log(\mathrm e^{\lambda}A_i)=\log A_i+\lambda\Id$.
By the previous theorem,  $\|\log A_1+\lambda\Id \|_2\leq \|\log A_2+\lambda\Id \|_2$
and $\lfloor A_1+\lambda\Id \rfloor_2\geq \lfloor A_2+\lambda\Id \rfloor_2$ holds for every
$\lambda\in\mathbb R$. According to the previous lemma, this implies the main statement.
The monotonicity is transparent in this case, as both $\log$ and $\exp$ are compatible
with conjugation by orthogonal matrices, hence the orbit correspondence is one-to-one.
$\log$ respects chirality, hence the statement can also be transferred to chiral disks.
\end{proof}
\end{theorem}
\snewpage
\section{Examples: The canonical Magnus developments in $\SL_2(\mathbb R)$}\plabel{sec:ExamplesMagnus}

\begin{example}\plabel{ex:parabolic}
(Magnus parabolic development.)
On the interval $[0,\pi]$,  consider again the measure $\Phi$, such that
\[\Phi(\theta)=
\begin{bmatrix}
-\sin2\theta& \cos2\theta\\\cos2\theta&\sin2\theta
\end{bmatrix}
\,\mathrm d\theta|_{[0,\pi]}.\]
Then, for $p\in[0,\pi)$,
\[\int\|\Phi|_{[0,p]}\|_2=p.\]
Here
\[\Lexp(\Phi|_{[0,p]})=W(p,p)=
\begin{bmatrix}
\cos p&2p\cos p -\sin p\\\sin p&2p\sin p+\cos p
\end{bmatrix}=
(\cos p\Id+\sin p\tilde I)(\Id_2-p\tilde I+p\tilde K).\]
Thus
\[\mu_{\mathrm L}(\Phi|_{[0,p]})=\log \Lexp(\Phi|_{[0,p]})=\AC(\cos p+p\sin p)
\begin{bmatrix}-p\sin p&2p\cos p -\sin p\\\sin p&p\sin p\end{bmatrix}.\]
Consequently,
\[\|\mu_{\mathrm L}(\Phi|_{[0,p]})\|_2=\AC(\cos p+p\sin p)\cdot(\sin p-p\cos p+p).\]
As $p\searrow0$
\begin{equation}
\|\mu_{\mathrm L}(\Phi|_{[0,p]})\|_2=p+{\frac {1}{6}}{p}^{3}-{\frac {1}{72}}{p}^{5}+{\frac {17}{3024}}{p}^
{7}+O \left( {p}^{9} \right).
\plabel{eq:micro}
\end{equation}
As $p\nearrow\pi$,
\begin{equation}
\|\mu_{\mathrm L}(\Phi|_{[0,p]})\|_2=\sqrt2\pi^{3/2} (\pi-p)^{-1/2}-2\pi+\frac{\sqrt{2\pi}(\pi^2-1) }{4} (\pi-p)^{1/2}+O(\pi-p).
\plabel{eq:lower2}
\end{equation}
This is not only better than \eqref{eq:lower1}, but it has the advantage that it can be interpreted
in terms of the solution of a differential equation blowing up.
\qedexer
\end{example}
%\snewpage
\begin{example}\plabel{ex:elliptic}(Magnus elliptic development.)
Let $h\in[0,1]$ be a parameter.
On the interval $[0,\pi]$, consider the measure $\widehat\Phi_h$ such that
\[\widehat\Phi_h(\theta)= (1-h)\begin{bmatrix}& -1\\1&\end{bmatrix}
+h\begin{bmatrix}-\sin2\theta& \cos2\theta\\\cos2\theta&\sin2\theta\end{bmatrix}\,\mathrm d\theta|_{[0,\pi]}. \]
Then, for $p\in[0,\pi)$
\[\int\|\widehat\Phi_{h}|_{[0,p]}\|_2=p.\]
It is easy to see that
\begin{multline*}
\Lexp(\widehat\Phi_{h}|_{[0,p]})=F((1-h)p,hp,p)=\\=\begin{bmatrix}
\cos p&2w\cos p -\sin p\\\sin p&2w\sin p+\cos p
\end{bmatrix}=
(\cos p\Id+\sin p\tilde I)(\Id_2-w\tilde I+w\tilde K).
\end{multline*}
Here $\widehat\Phi_1=\Phi$.
We find that
\[\|\mu_{\mathrm L}(\widehat\Phi_{h}|_{[0,p]})\|=\|\log\Lexp(\widehat\Phi_{h}|_{[0,p]})\|=
\AC(\cos p+hp\sin p)\cdot(\sin p-hp\cos p+hp).
 \]
Thus, if $h\neq 0$, then
\[\lim_{p\nearrow\pi}\|\mu_{\mathrm L}(\widehat\Phi_{h}|_{[0,p]})\|_2=+\infty.\]

It is notable that
\[\CD(\Lexp(\widehat\Phi_{h}|_{[0,p]})) =\Dbar( \mathrm e^{\mathrm ip}-\mathrm i\mathrm e^{\mathrm ip}ph   , ph) ,\]
which is $\CD(\Lexp(\Phi|_{[0,p]}))$ contracted from the boundary point $\mathrm e^{\mathrm ip}$ by factor $h$.
\qedexer
\end{example}
%\snewpage
\begin{example}\plabel{ex:hyperbolic}(Magnus hyperbolic development.)
More generally, let $t$ be a real parameter.
 On the interval $[0,\pi]$ consider the measure $\Phi_{\sin t}$, such that
\[\Phi_{\sin t}(\theta)=
\begin{bmatrix}
-\sin2(\theta\sin t)& \cos2(\theta\sin t)\\\cos2(\theta\sin t)&\sin2(\theta\sin t)
\end{bmatrix}
\,\mathrm d\theta|_{[0,\pi]}.\]
Then, for $p\in[0,\pi)$
\[\int\|\Phi_{\sin t}|_{[0,p]}\|_2=p.\]
$\Phi_1$ is the same as $\Phi$, and $\Phi_{-1}=\tilde K\cdot\Phi_1\cdot\tilde K$. If $t\in(-\pi/2,\pi/2)$, then
\begin{multline}\Lexp(\Phi_{\sin t}|_{[0,p]})=W(p,p\sin t)\\
=(\cos (p\sin t)\Id+\sin (p\sin t)\tilde I)\cdot\left(\cosh(p\cos t)\Id_2+\frac{\sinh(p\cos t)}{\cos t}\Bigl(-\sin t\tilde I+\tilde K\Bigr)\right).
\notag\end{multline}
Consequently,
\begin{multline}
\|\mu_{\mathrm L}(\Phi_{\sin t}|_{[0,p]} )\|_2
=\AC\left({\cosh \left( p\cos t  \right) \cos \left( p\sin t  \right)  +\frac {\sinh \left( p\cos t  \right)
}{\cos t }\sin \left( p\sin t \right) \sin t}\right)
\\
\cdot\left(\left|{\cosh \left( p\cos t  \right) \sin \left( p\sin t  \right) -\frac {\sinh \left( p\cos t  \right)
}{\cos t }\cos \left( p\sin t \right) \sin t}\right|+\frac{\sinh \left( p\cos t \right)}{\cos t }\right).
\plabel{eq:hypmagpre}
\end{multline}

Now, in the special case $p/\pi=\sin t$, we see that
\[\int\|\Phi_{p/\pi}|_{[0,p]}\|_2=p, \]
and
\begin{multline}
\|\mu_{\mathrm L}(\Phi_{p/\pi}|_{[0,p]})\|_2=
\sqrt2\pi^{3/2} (\pi-p)^{-1/2}
-2\pi
+\frac{\sqrt{2\pi}(4\pi^2-3) }{12} (\pi-p)^{1/2}
\\-\frac{4\pi^2-3}3 (\pi-p)^{1}
+\frac{\sqrt2(368\pi^2-840\pi^2-45)}{1440\sqrt\pi} (\pi-p)^{3/2}
+O((\pi-p)^2).
\plabel{eq:lower3}
\end{multline}
This shows that  \eqref{eq:lower2} is not optimal, either.
\qedexer
\end{example}
\begin{remark}
\plabel{rem:posmarois}
For $t\in[-\pi/2,\pi/2]$ and $p\in[0,\pi]$,
\[{\cosh \left( p\cos t  \right) \sin \left( p\sin t  \right) -\frac {\sinh \left( p\cos t  \right)
}{\cos t }\cos \left( p\sin t \right) \sin t}\geq0\]
holds. (Understood as $=0$ for $t=0$. It is also $=0$ for $p=0$.)
\begin{commentx}
The first nontrivial nonzero seems to appear  for $(p,t)=(\boldsymbol z,\pm\pi/2)$.
\end{commentx}
Thus, under these assumptions, the absolute value in \eqref{eq:hypmagpre} is unnecessary.
\qedremark
\end{remark}

\begin{commentx}
In what follows, when we use the terms `Magnus elliptic development' and `Magnus hyperbolic development',
we may allow the case of the Magnus parabolic development.
If we want to exclude it, we say `strictly elliptic'  or `strictly hyperbolic'  development.
\end{commentx}

\scleardoublepage\section{Magnus minimality in the $\GL_2^+(\mathbb R)$ case}\plabel{sec:MagnusGL2}

\begin{theorem}\plabel{th:maximaldisk}
Let $p\in(0,\pi)$.
Consider the family of disks parameterized by $t\in[-\pi/2,\pi/2]$,
such that the centers and radii are
\[\Omega_p(t)=\mathrm e^{\mathrm i p\sin t}\left(\cosh(p\cos t) -\mathrm i\frac{\sinh(p\cos t)\sin t}{\cos t}\right),\]
\[\omega_p(t)=\frac{\sinh(p\cos t)}{\cos t},\]
for $t\neq \pm\pi/2$;
and
\[\Omega_p(\pm\pi/2)=(\cos p+p\sin p )\pm\mathrm i(\sin p-p\cos p),\]
\[\omega_p(\pm\pi/2)=p.\]

(a) The circle $\partial\Dbar(\Omega_p(t),\omega_p(t))$  is tangent to  $\partial\exp\Dbar(0,p)$ at
\[\gamma_p(t)=\mathrm e^{p\cos t+\mathrm ip\sin t}\qquad \text{and}\qquad
\gamma_p(\pi-t \modu2\pi)=\mathrm e^{-p\cos t+\mathrm i p\sin t}.\]
These points are inverses of each other relative to the unit circle.
If the points are equal ($t=\pm\pi/2$), then the disk is the osculating disk at $\gamma_p(t)$.

The disks themselves are orthogonal to the unit circle.
The disks are distinct from each other.
Extending $t\in[-\pi,\pi]$, we have $\Omega_p(t)=\Omega_p(\pi-t \modu2\pi)$,  $\omega_p(t)=\omega_p(\pi-t \modu2\pi)$.

(b)
\[\CD(\Lexp(\Phi_{\sin t}|_{[0,p]})=\CD(W(p,p\sin t))=\Dbar(\Omega_p(t),\omega_p(t)).\]

(c)
The disks $\Dbar(\Omega_p(t),\omega_p(t))$  are the maximal disks in $\exp\Dbar(0,p)$.
The maximal disk $\Dbar(\Omega_p(t),\omega_p(t))$ touches  $\partial\exp\Dbar(0,p)$ only at
$\gamma_p(t)$, $\gamma_p(\pi-t \modu2\pi)$.

\begin{proof}
(a)
The disks are distinct because, the centers are distinct:
For $t\in(-\pi/2,\pi/2)$,
\[\frac{\mathrm d\arg \Omega_p(t)}{\mathrm dt}=\Ima\frac{\mathrm d\log\Omega_p(t)}{\mathrm dt}=
\frac{(p\sin(t)\cosh(p\sin t)-\cosh(p\cos t))\cosh(p\sin t)}{\cosh(p\sin t)^2-\sin^2 t}>0.\]
(Cf. $\int_0^x y\sinh y\mathrm dy=x\cosh x-\sinh x.$)
The rest can easily be checked using the observation
\[\Omega_p(t)=\mathrm e^{p\cos t+\mathrm ip\sin t}-\frac{\sinh(p\cos t)}{\cos t}\mathrm e^{\mathrm i(t+p\sin t)}
=\mathrm e^{-p\cos t+\mathrm ip\sin t}+\frac{\sinh(p\cos t)}{\cos t}\mathrm e^{\mathrm i(-t+p\sin t)}.\]

(b) This is direct computation.

(c) In general, maximal disks touch the boundary curve $\gamma_p$, and any such touching point determines the maximal disk.
(But a maximal disk might belong to different points.)
Due to the double tangent / osculating property the given disks are surely the maximal disks,
once we prove that they are indeed contained in $\exp\Dbar(0,p)$.
However, $\CD(\Lexp(\Phi_{\sin t}|_{[0,p]})=\Dbar(\Omega_p(t),\omega_p(t))$ together
with Theorem \ref{th:CRrange} implies that $\Dbar(\Omega_p(t),\omega_p(t))\subset\exp\Dbar(0,p)$.
The distinctness of the circles implies that they touch the boundary only at the indicated points.
\end{proof}

\snewpage
\begin{proof}[Alternative proof for $\Dbar(\Omega_p(t),\omega_p(t))\subset\exp\Dbar(0,p)$.]
Here we give a purely differential geometric argument.

One can see that the given disks $\Dbar(\Omega_p(t),\omega_p(t))$ are characterized by the following properties:

($\alpha$) If $\gamma_p(t)\neq \gamma_p(\pi-t \modu2\pi)$, then the disk is tangent to $\gamma_p$ at these points.

($\beta$) If $\gamma_p(t)= \gamma_p(\pi-t \modu2\pi)$, i. e. $t=\pm\pi$, then the disk is the osculating disk at $\gamma_p(\pm\pi/2)$.

Now, we prove that  $\Dbar(\Omega_p(t),\omega_p(t))\subset\exp\Dbar(0,p)$.
First, we show that $\Dbar(\Omega_p(0),\omega_p(0))\subset \exp\Dbar(0,p)$.
Indeed, \[\Dbar(\Omega_p(0),\omega_p(0))=\PD \left(\begin{bmatrix}
\mathrm e^p&\\&\mathrm e^{-p}
\end{bmatrix}\right);\] hence, by Theorem \ref{th:CRrange}, the $\log$ of any element of $\Dbar(\Omega_p(0),\omega_p(0))$ is
contained in
\[\PD\left(\log \begin{bmatrix}\mathrm e^p&\\&\mathrm e^{-p}\end{bmatrix}\right)
=\PD\left( \begin{bmatrix}p&\\&{-p}\end{bmatrix}\right)=\Dbar(0,p).\]
Let $L$ be the maximal real number such that $\Dbar(\Omega_p(t),\omega_p(t))\subset\exp\Dbar(0,p) $
for any $t\in[-L,L]$, and $L<\pi/2$. (Due to continuity, there is  a maximum.)
Indirectly, assume that $L<\pi/2$.
Then one of following  should happen:

(i) Besides $\gamma_p(L)$ and $ \gamma_p(\pi-L\modu2\pi)$ there is another pair (due to inversion symmetry)
of distinct points $\gamma_p(\tilde L)$ and $ \gamma_p(\pi-\tilde L\modu2\pi)$, where $\Dbar(\Omega_p(L),\omega_p(L))$ touches the boundary of $\exp\Dbar(0,p)$.

(ii)  $\Dbar(\Omega_p(L),\omega_p(L))$ touches the boundary at $\gamma_p(\pi/2)$ or $\gamma_p(-\pi/2)$.

(iii) $\Dbar(\Omega_p(L),\omega_p(L))$ is osculating at $\gamma_p(L)$ or at $ \gamma_p(\pi-L\modu2\pi)$.

(Symmetry implies that $t=\pm L$ are equally bad.)
Case (i) is impossible, because the given circles are distinct and the characterising properties hold.
Case (ii) is impossible, because, due to  $\omega_p(L)>p$ and the  extremality of $\arg \gamma_p(\pm\pi/2)$,
the situation would imply that   $\Dbar(\Omega_p(L),\omega_p(L))$ strictly contains the
osculating disk at $\gamma_p(\pi/2)$ or $\gamma_p(-\pi/2)$,
which is a contradiction to $\Dbar(\Omega_p(L),\omega_p(L))\subset\exp\Dbar(0,p) $.
Case (iii) is impossible, because for the oriented plane curvature of $\gamma_p$,
\[\varkappa_{\gamma_p}(t)=\frac{1+p\cos t}{p\mathrm e^{p\cos t}}<\frac1{\omega_p(t)}=\frac{\cos t}{\sinh (p\cos t) }\]
holds if $\cos t\neq0$. (In general, $\frac{1+x}{\mathrm e^x}<\frac{x}{\sinh x}$ for $x\neq0$.)
This implies $L=\pi/2$, proving the statement.
\end{proof}
\end{theorem}
In what follows, we will not make much issue out of expressions like $\frac{\sinh px}x$ when $x=0$;
we just assume that they are equal to $p$, in the spirit of continuity.

\snewpage
\begin{theorem}\plabel{th:nonmaximaldisk} Suppose that $p\in(0,\pi)$.
Suppose that $D$ is a disk in  $\exp\Dbar(0,p)$, which touches $\partial\exp\Dbar(0,p)$ at
$\gamma_p(t)=\mathrm e^{p\cos t+\mathrm ip\sin t}$.
Then for an appropriate nonnegative decomposition $p=p_1+p_2$,
\[D=\CD\left(\exp(  p_1(\Id\cos t+\tilde I\sin t) ) \cdot   W(p_2,p_2\sin t) \right).\]
The bigger  the $p_2$ is, the bigger the corresponding disk is.
$p_2=p$  corresponds to the maximal disk, $p_2=0$ corresponds to the point disk.
\begin{proof}
Let $W_{p_1,p_2,t} $ denote the argument of $\CD$.
Then its first component is Magnus exponentiable by norm  $p_1$,
and its second component is Magnus exponentiable by norm  $p_2$.
Thus the principal disk must lie in $\exp\Dbar(0,p)$.
One can compute the center and the radius of the chiral disk (cf. the Remark), and
find that  $\gamma_p(t)$ is on the boundary of the disk.
So, $\CD(W_{p_1,p_2,t})$ must be the maximal $\CD(W_{0,p_1+p_2,t})$ contracted from $\gamma_p(t)$.
One, in  particular, finds that the radius of  $\CD(W_{p_1,p_2,t})$ is
\[\frac{\mathrm e^{p_1+p_2}- \mathrm e^{p_1-p_2}}{2\cos t}=\frac{\mathrm e^{p}}{\cos t}(1-\mathrm e^{-2p_2}).\]
This shows that bigger $p_2$ leads to  bigger disk.
\end{proof}
\proofremarkqed{
\renewcommand{\qedsymbol}{$\triangle$}
It is easy to see that, for $p=p_1+p_2$,
\begin{align*}
&\exp(  p_1(\Id\cos t+\tilde I\sin t) ) \cdot   W(p_2,p_2\sin t)=
\notag\\
&=\mathrm e^{p_1\cos t}\exp((p_1+p_2)\sin t\tilde I)\cdot\left(\cosh(p_2\cos t)\Id_2
+\frac{\sinh(p_2\cos t)}{\cos t}\Bigl(-\sin t\tilde I+\tilde K\Bigr)\right)
\notag\\
&=\Lexp\left(\frac{p_1}{p}\begin{bmatrix}\cos t& -\sin t\\\sin t&\cos t\end{bmatrix}
+\frac{p_2}{p}\begin{bmatrix}-\sin(2\theta\sin t)& \cos(2\theta\sin t)\\
\cos(2\theta\sin t)&\sin(2\theta\sin t)\end{bmatrix}   \mathrm d\theta|_{[0,p]}\right)
\notag\\
&=\Lexp\left(p_1\begin{bmatrix}\cos t& -\sin t\\\sin t&\cos t\end{bmatrix}
+p_2\begin{bmatrix}-\sin(2p\theta\sin t)& \cos(2p\theta\sin t)\\
\cos(2p\theta\sin t)&\sin(2p\theta\sin t)\end{bmatrix}\mathrm d\theta|_{[0,1]}\right).
\notag\qedhere
\end{align*}
}
\end{theorem}
This immediately implies the existence of a certain normal form.
For the sake of compact notation, let
\[\tilde {\mathbb K}:=\{ -\sin\beta \tilde J+\cos\beta \tilde K\,:\,\beta\in[0,2\pi) \},\]
which is  the set of the conjugates of $\tilde K$ by orthogonal matrices.
\snewpage
\begin{theorem}\plabel{th:normalform}
Suppose that $A\in\mathrm M_2(\mathbb R)$ such that $\CD(A)\subset\exp\intD(0,\pi)$.
Assume that $p$ is the smallest real number such that $\CD(A)\subset\exp\Dbar(0,p)$,
and $\CD(A)$ touches $\exp\partial\Dbar(0,p)$ at $\mathrm e^{p(\cos t +\mathrm i\sin t)}$.
Then there is an nonnegative decomposition $p=p_1+p_2$,
and a matrix $\tilde F\in\tilde {\mathbb K}$, such that
\begin{align}
A=\mathrm e^{p_1\cos t}&\exp(p\sin t\tilde I)\cdot\left(\cosh(p_2\cos t)
\Id_2-\frac{\sinh(p_2\cos t)}{\cos t}\sin t\tilde I\right)+\frac{\sinh(p_2\cos t)}{\cos t}\tilde F
\plabel{eq:norarith}\\
&=\Lexp\left( p_1\exp(t\tilde I)+p_2\exp(2p\theta\sin t\tilde I)\cdot\tilde F \,\,\mathrm d\theta|_{[-1/2,1/2]} \right)
\plabel{eq:noruni}\\
&=\Lexp(\exp(t\tilde I) \,\mathrm d\theta|_{[0,p_1]})
\Lexp\left(  \exp((2\theta-p_1-p_2)\sin t \tilde I)\tilde F \, \mathrm d\theta|_{[0,p_2]} \right)
\plabel{eq:norleft}\\
&=\Lexp\left(  \exp((2\theta+p_1-p_2)\sin t \tilde I)\tilde F \, \mathrm d\theta|_{[0,p_2]} \right) \Lexp(  \exp(t\tilde I) \,\mathrm d\theta|_{[0,p_1]}).
\plabel{eq:norright}
\end{align}

The case $p_1=p_2=0$ corresponds to $A=\Id_2$.

The case $p_1>0,p_2=0$ corresponds to point disk case, the expression does not depend on $\tilde F$.

The case $p_1=0,p_2>0$  corresponds to the maximal disk case,
it has degeneracy $t\leftrightarrow \pi - t\modu 2\pi$.

In the general case $p_1,p_2>0$, the presentation is unique in terms of $p_1,p_2,t\modu 2\pi,\tilde F$.
\begin{proof}
This is an immediate consequence of the previous statement and the observation
$(\cos\alpha+\tilde I\sin\alpha)\tilde K(\cos\alpha+\tilde I\sin\alpha)^{-1}=
(\cos2\alpha+\tilde I\sin2\alpha)\tilde K= -\tilde J\sin2\alpha+\tilde K\cos2\alpha.$
\end{proof}
\end{theorem}

In what follows, we use the notation
\[\NW(p_1,p_2,t,\tilde F)\]
to denote the arithmetic expression on the RHS of \eqref{eq:norarith}.
In itself, it just a matrix value, but
the statement above offers three particularly convenient ways (normal forms) to present is as a left-exponential:
\eqref{eq:noruni} is sufficiently nice and compact with norm density $p$ on an interval of  unit length.
\eqref{eq:norleft} and \eqref{eq:norright} are concatenations of intervals of length $p_1$ and $p_2$
with norm density $1$.
One part is essentially a complex exponential, relatively uninteresting;
the other part is the Magnus parabolic or  hyperbolic development of Examples \ref{ex:parabolic}
and \ref{ex:hyperbolic},
but up to conjugation by a special orthogonal matrix, which is the same to say as `up to phase'.

\begin{theorem}\plabel{th:MPexpress}
Suppose that $A\in\mathrm M_2(\mathbb R)$ such that $\CD(A)\subset\exp\intD(0,\pi)$.
Then
\[\mathcal M_{2\times 2\,\,\real}(A)=\inf\{\lambda\in[0,\pi)\,:\, \CD(A)\subset \exp\Dbar(0,\lambda)\}.\]
Or, in other words,
\[\mathcal M_{2\times 2\,\,\real}(A)=\sup\{|\log z|\,:\, z\in\CD(A) \}.   \]
\begin{proof}
Assume that $p$ is the smallest real number such that $\CD(A)\subset\exp\Dbar(0,p)$.
By Theorem \ref{th:CRrange}, $\mathcal M_{2\times 2\,\,\real}(A)$ is at least $p$,
while the left-exponentials of Theorem \ref{th:normalform} does indeed Magnus-exponentiate them
with norm $p$.
\end{proof}
\end{theorem}
\snewpage
Suppose that $A\in\mathrm M_2(\mathbb R)$ such that $\CD(A)\subset\exp\intD(0,\pi)$, $A\neq\Id_2$, $p=\mathcal M_{2\times 2\,\,\real}(A)$.
If $\det A=1$, then $A$ can be of the three kinds:
Magnus elliptic, when $\CD(A)$ touches $\exp\partial\Dbar(0,p)$ at $\mathrm e^{\mathrm ip}$  or
$\mathrm e^{-\mathrm ip}$, but it is not an osculating disk;
Magnus parabolic, when  $\CD(A)$ touches $\exp\partial\Dbar(0,p)$ at $\mathrm e^{\mathrm ip}$  or
$\mathrm e^{-\mathrm ip}$, and  it is an osculating disk;
or Magnus hyperbolic when $\CD(A)$ touches $\exp\partial\Dbar(0,p)$ at two distinct points.
If $\det A\neq 1$ then $\CD(A)$  touches $\exp\partial\Dbar(0,p)$ at a single point,
asymmetrically; we can call these Magnus loxodromic.
We see that Examples \ref{ex:parabolic}, \ref{ex:elliptic}, and \ref{ex:hyperbolic},
cover all the Magnus parabolic, hyperbolic and elliptic cases up to conjugation by an orthogonal matrix.
In general, if $A$ is not Magnus hyperbolic, then it determines a unique Magnus direction
$\cos t+\mathrm i\sin t$ (in the notation Theorem \ref{th:normalform}).
It is the direction of the farthest point of $\{\log z\,:\,z\in \CD(A)\}$ from the origin.
If $A$ is Magnus hyperbolic, then this direction is determined only up to sign in the real part.
\begin{lemma}\plabel{lem:Magnusclassify}
Suppose $A\in\mathrm M_2(\mathbb R)$ such that $\CD(A)\subset\exp\intD(0,\pi)$, $A\neq \Id_2$, $\det A=1$, $\CD(A)=\Dbar((a,b),r)$.
Then $a^2+b^2=r^2+1$ and $a+1>0$.

We claim that $A$ is Magnus hyperbolic or parabolic if and only if
\[2\arctan\frac{r+|b|}{a+1}\leq r.\]
If $A$ is Magnus elliptic or parabolic, then
\[\mathcal M_{2\times 2\,\,\real}(A)=2\arctan\frac{r+|b|}{a+1}.\]
\begin{proof} $\partial \Dbar((a,b),r)$ intersects the unit circle at
\[(\cos\varphi_\pm,\sin\varphi_\pm):= \left(\dfrac{a\pm br}{a^2+b^2},\dfrac{b\mp ar}{a^2+b^2}\right),\]
$\varphi_\pm\in(-\pi,\pi)$.
In particular, $\dfrac{a\pm br}{a^2+b^2}+1>0$; multiplying them, we get $a+1>0$.
Then $\phi_\pm=2\arctan \frac{r\pm b}{a+1}$.
If one them is equal to $r$, then it is a Magnus parabolic case;
if those are smaller than $r$, then it is a Magnus hyperbolic case;
if one of them is bigger than $r$, this it must be a Magnus elliptic case.
(Cf. the size of the chiral disk in Theorem \ref{th:nonmaximaldisk}.)
\end{proof}
\end{lemma}

Recall that we say that the measure $\phi$ is a minimal Magnus presentation for $A$, if
$\Lexp(\phi)=A$ and $\int\|\phi\|_2=\mathcal M_{2\times 2\,\,\real}(A)$.
\begin{lemma}\plabel{lem:MagnusMinEx}
Any element $A\in\GL^+_2(\mathbb R)$ has at least one minimal Magnus presentation.
\begin{proof}
$\GL^+_2(\mathbb R)$ is connected, which implies that any element $A$
has at least one Magnus presentation $\psi$.
If $\int\|\phi\|_2$ is small enough,
then we can divide the supporting interval of $\phi$ into $\lfloor\mathcal M_{2\times 2\,\,\real}(A)/\pi\rfloor$
many subintervals, such that the variation of $\phi$ on any of them is less than $\pi$.
Replace   $\phi$ by a normal form on every such subinterval.
By this we have managed to get a presentation of variation at most    $\int\|\phi\|_2$
by a data from $([0,\pi]\times[0,\pi]\times [0,2\pi]\times \mathbb K) ^{\lfloor\mathcal M_{2\times 2\,\,\real}(A)/\pi\rfloor}$.
Conversely, such a data always gives a presentation, whose $\Lexp$ depends continuously on
the data.
Then the statement follows from a standard compactness argument.
\end{proof}
\end{lemma}
\snewpage
\begin{lemma}\plabel{lem:hyperbolicestimate}
Suppose that $A_\lambda\rightarrow \Id$, such that $A_\lambda$ is Magnus hyperbolic,
but $A_\lambda\neq \Id$ for any $\lambda$.
Suppose that $\CD(A_\lambda)=\Dbar((1+a_\lambda,b_\lambda),r_\lambda)$.

Then, as the sequence converges,
\[\mathcal M_{2\times 2\,\,\real}(A_\lambda)^2=2 a_\lambda+O(\mathrm{itself}^2);\]
or more precisely,
\[\mathcal M_{2\times 2\,\,\real}(A_\lambda)^2=2 a_\lambda-\frac13a_\lambda^2+\frac32\frac{b_\lambda^2}{a_\lambda} +O(\mathrm{itself}^3).\]
\begin{proof}
We can assume that $A_\lambda=W(p_\lambda,p_\lambda\sin  t_\lambda)$.
From the formula of $W(p,p\sin t)$ one can see that $\CD(W(p,p\sin t))$ is an entire function of $x=p\cos t, y=p\sin t$.   One actually finds that the center is
\begin{align}
(1+\hat a(x,y),\hat b(x,y))=&\biggl(1+\frac{x^2+y^2}2+\frac{(x^2-y^2)(x^2+y^2)}{24 }\notag\\
&+\frac{(x^4-10x^2y^2+5y^4)(x^2+y^2)}{720}+O(x,y)^8, \notag\\
&\frac{y(x^2+y^2)}3 + \frac{y(x^2+y^2)(x^2-y^2)}{30} +O(x,y)^7\biggr).
\notag
\end{align}
(One can check that
in the expansion $\hat a(x,y)$, every term is divisible by $(x^2+y^2)$;
in the expansion $\hat b(x,y)$, every term is divisible by $y(x^2+y^2)$.)
Eventually, one finds that
\[p^2=x^2+y^2= 2\hat a(x,y)+O(x,y)^4\]
and
\[p^2=x^2+y^2= 2\hat a(x,y)- \frac13\hat a(x,y)^2+\frac32 \frac{\hat b(x,y)^2}{\hat a(x,y)}  +O(x,y)^6.\qedhere\]
\end{proof}
\end{lemma}
The hyperbolic developments $p\mapsto W(p,p\sin t) $ are uniform motions in the sense
that the increments $ W((p+\varepsilon),(p+\varepsilon)\sin t)  W(p,p\sin t) ^{-1}$
differ from each other by conjugation by orthogonal matrices as $p$ changes.
In fact,  they are locally characterized by the speed $\sin t$,
and a phase, i. e. conjugation by rotations.

\snewpage
\begin{lemma}\plabel{lem:hyperbolicfit}
Assume that $0< p_1,p_2$; $p_1+p_2<\pi$; $t_1,t_2\in[-\pi/2,\pi/2]$; $\varepsilon\in(-\pi/2,\pi/2]$.
On the interval $[-p_1,p_2]$, consider the measure $\phi$ given by
\[\phi(\theta)=\eta(\theta)\,\mathrm d\theta|_{[-p_1,p_2]}, \]
where
\[\eta(\theta)=\begin{cases}
\begin{bmatrix}
-\sin2(\theta\sin t_2)& \cos2(\theta\sin t_2)\\\cos2(\theta\sin t_2)&\sin2(\theta\sin t_2)
\end{bmatrix}
&\text{if }\theta\geq0\\
\begin{bmatrix}
\cos\varepsilon&-\sin\varepsilon\\\sin\varepsilon&\cos\varepsilon
\end{bmatrix}
\begin{bmatrix}
-\sin2(\theta\sin t_1)& \cos2(\theta\sin t_1)\\\cos2(\theta\sin t_1)&\sin2(\theta\sin t_1)
\end{bmatrix}
\begin{bmatrix}
\cos\varepsilon&\sin\varepsilon\\-\sin\varepsilon&\cos\varepsilon
\end{bmatrix}
&\text{if }\theta\leq0.
\end{cases} \]
Then
\[\mathcal M_{2\times 2\,\,\real}(\Lexp(\phi))<p_1+p_2\]
unless $\varepsilon=0$ and $t_1=t_2$.
\begin{proof}
It is sufficient to prove this for a small subinterval around $0$.
So let us take the choice $p_1=p_2=p/2$, $p\searrow0$.
Then
\[\Lexp(\phi|_{[-p/2,p/2]})=W\left(\frac p2,\frac p2\sin t_2\right)
\begin{bmatrix}
\cos\varepsilon&-\sin\varepsilon\\\sin\varepsilon&\cos\varepsilon
\end{bmatrix}
W\left(-\frac p2,-\frac p2 \sin t_1\right)^{-1}
\begin{bmatrix}
\cos\varepsilon&\sin\varepsilon\\-\sin\varepsilon&\cos\varepsilon
\end{bmatrix}.
\]
Let
\[\Dbar((a_p,b_p),r_p)=\CD(\Lexp(\phi|_{[-p/2,p/2]}) ).\]

(i) If $\varepsilon\in(-\pi/2,0)\cup(0,\pi/2)$, then
\[2\arctan\frac{r_p\pm b_p}{a_p+1}-r_p=\mp\frac14\sin(2\varepsilon)p^2+O(p^3).\]
This shows that $\Lexp(\phi|_{[-p/2,p/2]})$ gets Magnus elliptic.
However,
\[\mathcal M_{2\times 2\,\,\real}(\Lexp(\phi|_{[-p/2,p/2]}))=2\arctan\frac{r_p\pm b_r}{a_p+1}=p\cos(\varepsilon) +O(p^2)\]
shows Magnus non-minimality.

(ii)
If $\varepsilon=\pi/2$, $\sin t_1+\sin t_2\neq0$, then
\[2\arctan\frac{r_p\pm b_r}{a_p+1}-r_p=\mp\frac1{12}(\sin t_1+\sin t_2 )p^3+O(p^4).\]
This also shows Magnus ellipticity, and
\[2\arctan\frac{r_p\pm b_r}{a_p+1}=\frac14|\sin t_1+\sin t_2| p^2+O(p^3)\]
shows Magnus non-minimality.

(iii)
If $\varepsilon=\pi/2$, $\sin t_1+\sin t_2=0$, then $\Lexp(\phi|_{[-p/2,p/2]})=\Id_2$.
Hence, full cancellation occurs, this is not Magnus minimal.

(iv) If $\varepsilon=0$, $\sin t_1\neq\sin t_2$, then  $\sin t_1+\sin t_2<2$, and
\[2\arctan\frac{r_p\pm b_p}{a_p+1}-r_p=\frac1{6}(\pm(\sin t_1+\sin t_2 )-2)p^3+O(p^4).\]
This shows that  $\Lexp(\phi|_{[-p/2,p/2]})$ gets Magnus hyperbolic.
Then, assuming Magnus minimality and using the previous lemma, we get a contradiction by
\[\mathcal M_{2\times 2\,\,\real}(\Lexp(\phi|_{[-p/2,p/2]}))^2=p^2- \frac1{48}p^4(\sin t_2-\sin t_1)^2+O(\text{itself}^3)<p^2. \]
This proves the statement.
\end{proof}
\end{lemma}

\begin{lemma}\plabel{lem:parabolicfit}
Assume that $0< p_1,p_2$; $p_1+p_2<\pi$; $t_1\in[-\pi/2,\pi/2)$.
On the interval $[-p_1,p_2]$, let us consider the measure $\phi$ given by
\[\phi(\theta)=\eta(\theta)\,\mathrm d\theta, \]
where
\[\eta(\theta)=\begin{cases}
\tilde I =\begin{bmatrix}
& -1\\1&
\end{bmatrix}
&\text{if }\theta\geq0\\
\begin{bmatrix}
-\sin2(\theta\sin t)& \cos2(\theta\sin t)\\\cos2(\theta\sin t)&\sin2(\theta\sin t)
\end{bmatrix}
&\text{if }\theta\leq0.
\end{cases} \]
Then
\[\mathcal M_{2\times 2\,\,\real}(\Lexp(\phi))<p_1+p_2.\]
\begin{proof}
Again, it is sufficient to show it for a small subinterval around $0$.

(i) Suppose $t\in(-\pi/2,\pi/2)$.
As $p\searrow0$, restrict to the interval
\[\mathcal I_p=\left[-p,\frac{\sinh p\cos t}{\cos t}-p\right].\]
Then
\[ \Lexp(\phi|_{\mathcal I_p})=\exp\left(\tilde I\left(\sin  \frac{\sinh p\cos t}{\cos t }-p\right) \right)
W(-p,-p\sin t )^{-1}.\]
Let
\[\Dbar((a_p,b_p),r_p)=\CD(\Lexp(\phi|_{\mathcal I_p}) ).\]

If we assume Magnus minimality, then
 \[\mathcal M_{2\times 2\,\,\real}(\Lexp(\phi|_{\mathcal I_p}))=\frac{\sinh p\cos t}{\cos t}=r_p.\]
Thus, $\Lexp(\phi|_{\mathcal I_p})$ is Magnus parabolic.
By direct computation, we find
\[2\arctan\frac{r_p+| b_p|}{a_p+1}=
p+\frac1{3}p^3\max(\cos^2t+\sin t-1,-1-\sin t )+O(p^4),\]
in contradiction to
\[\frac{\sinh p\cos t}{\cos t}=p+\frac1{6}p^3(\cos^2t )+O(p^4), \]
which is another way to express $\mathcal M_{2\times 2\,\,\real}(\Lexp(\phi|_{\mathcal I_p}))$ from the density.
(The coefficients of $p^3$ differ for $t\in(-\pi/2,\pi/2)$.)

(ii) Consider now the case $t=-\pi/2$.
\[2\arctan\frac{r_p\pm b_p}{a_p+1}=\pm\frac12p+O(p^2)\]
shows Magnus ellipticity, and
 \[2\arctan\frac{r_p+|b_p|}{a_p+1}=p-\frac1{12}p^3+O(p^4)\]
shows non-minimality.
\end{proof}
\end{lemma}

\snewpage
Now we deal with the unicity of the normal forms as left exponentials.
In the context of Theorem $\ref{th:normalform}$ we call $\ellip(A):=p_1(\cos t+\tilde I\sin t)$
the elliptic component of $A$, and we call $\hyper(A):=p_2$ the hyperbolic length of $A$.

\begin{theorem}\plabel{th:additivity}
Suppose that $A\in\mathrm M_2(\mathbb R)$ such that $\CD(A)\subset\exp\intD(0,\pi)$, and $\phi$
is a minimal Magnus presentation for $A$ supported on $[a,b]$.

Then, restricted to any subinterval $\mathcal I$, the value
$\ellip(\Lexp(\phi|_{\mathcal I}))$ is a multiple of  $\ellip(A)$ by a nonnegative real number.
Furthermore the interval functions
\[\mathcal I\mapsto \mathcal M_{2\times 2\,\,\real}(\Lexp(\phi|_{\mathcal I}))=\smallint \|\phi|_{\mathcal I}\|_2,\]
\[\mathcal I\mapsto \ellip(\Lexp(\phi|_{\mathcal I})),\]
\[\mathcal I\mapsto \hyper(\Lexp(\phi|_{\mathcal I}))\]
are additive. In particular, if $A$ is Magnus hyperbolic or parabolic, then
$\ellip(\Lexp(\phi|_{\mathcal I}))$ is always $0$.
\begin{proof}
Let us divide the supporting interval of $\phi$ into smaller intervals
$\mathcal I_1,\ldots,\mathcal I_s$.
On these intervals let us replace $\phi|_{\mathcal I_k}$  by a left-complex normal form.
Thus we obtain
\[\phi'=\Phi^{(1)}_{\mathcal K_1}\boldsymbol.(\cos t_1+\tilde I\sin t_1)\m 1_{\mathcal J_1}\boldsymbol.
\ldots\boldsymbol.\Phi^{(s)}_{\mathcal K_s}\boldsymbol.(\cos t_s+\tilde I\sin t_s)\m 1_{\mathcal J_s},\]
where $\mathcal J_j$ are $\mathcal K_j$ are some intervals, and $\Phi^{(j)}_{\mathcal K_j}$ are hyperbolic developments (up to conjugation).
(They can be parabolic but for the sake simplicity let us call them hyperbolic.)
Further, rearrange this as
\[\phi''=
\Phi^{\prime(1)}_{\mathcal K_1}\boldsymbol.
\ldots\boldsymbol.\Phi^{\prime(s)}_{\mathcal K_s}
\boldsymbol.(\cos t_1+\tilde I\sin t_1)\m 1_{\mathcal J_1}\boldsymbol.
\ldots\boldsymbol.(\cos t_s+\tilde I\sin t_s)\m 1_{\mathcal J_s}
,\]
where the hyperbolic developments suffer some special orthogonal conjugation but they remain
hyperbolic developments. Now, the elliptic parts
\[\ellip(\Lexp(\phi|_{\mathcal I_j}) )= |\mathcal J_j|(\cos t_j+\tilde I\sin t_j)\]
must be  nonnegatively proportional to each other, otherwise
cancelation would occur when
the elliptic parts are contracted, in contradiction to the minimality of the presentation.
By this, we have proved that in a minimal presentation elliptic parts of
disjoint intervals are nonnegatively proportional to each other.

Suppose that in a division $|\mathcal J_j|\cos t_j\neq 0$ occurs.
Contract the elliptic parts in $\phi''$ but immediately divide them into two equal parts:
\[\phi'''=
\Phi^{\prime(1)}_{\mathcal K_1}\boldsymbol.
\ldots\boldsymbol.\Phi^{\prime(s)}_{\mathcal K_s}
\boldsymbol.(\cos t_j+\tilde I\sin t_j)\m 1_{\mathcal J}\boldsymbol.(\cos t_j+\tilde I\sin t_j)\m 1_{\mathcal J}\boldsymbol
.\]
Now replace everything but the last term by a normal form
\[\phi''''=
\Phi^{\prime(0)}_{\mathcal K_0}
\boldsymbol.(\cos t_0+\tilde I\sin t_0)\m 1_{\mathcal J_0}\boldsymbol.(\cos t_j+\tilde I\sin t_j)\m 1_{\mathcal J}\boldsymbol
.\]
Taking the determinant of the various left-exponential term we find
\[\mathrm e^{ |\mathcal J_0|\cos t_0+|\mathcal J|\cos t_j }= \mathrm e^{ 2|\mathcal J|\cos t_j }.\]
Thus $|\mathcal J_0|\cos t_0\neq 0$, hence, by minimality $t_j=t_0\modu 2\pi$, moreover  $|\mathcal J_0|=|\mathcal J|$.
However, the $\phi'''$ constitutes a normal form (prolonged in the elliptic part), which in this form is unique,
thus, eventually
\begin{equation}\ellip(\Lexp(\phi))=\sum_{j=1}^s \ellip(\Lexp(\phi|_{\mathcal I_j}))
\plabel{eq:additivity}\end{equation}
must hold.

Suppose now that $\sin t_k=1$ or $\sin t_k=-1$ occurs with $|\mathcal J_k|\neq0$.
Consider $\phi''$.
By Magnus minimality and Lemma \ref{lem:hyperbolicfit}, the hyperbolic development must fit
into single hyperbolic development $\Psi_{\mathcal K}$ (without phase or speed change).
Furthermore, by Lemma \ref{lem:parabolicfit}, $\Psi_{\mathcal K}$ must be parabolic fitting
properly to the elliptic parts. Thus $\phi''$, in fact, yields a
normal form $\Psi_{\mathcal K}\boldsymbol.  (\sin t_k)\m 1_{\mathcal J}$.
Then  \eqref{eq:additivity} holds.

The third possibility in $\phi''$ is that all the intervals $\mathcal J_j$ are of zero length.
Then the hyperbolic developments fit into a single development $\Psi_{\mathcal K}$,
but \eqref{eq:additivity} also holds.

Thus \eqref{eq:additivity} is proven.
It implies nonnegative proportionality relative to the total $\ellip(\Lexp(\phi))$.
Now, subintervals of minimal presentations also yield  minimal presentations, therefore
additivity holds in full generality.
Regarding the interval functions, the additivity of $\mathcal M_{2\times 2\,\,\real}$ is trivial, the additivity of $\ellip$
is just demonstrated, and $\hyper$ is just the $\mathcal M_{2\times 2\,\,\real}$ minus the absolute value (norm) of $\ellip$.
\end{proof}
\end{theorem}
\begin{remark}\plabel{rem:Semilocal}
Suppose that $\phi:\mathcal I\rightarrow \mathcal B(\mathfrak H)$ is a measure.
Assume that $\mathcal I_1\subset \mathcal I$ is a subinterval such that
$\smallint \|\phi|_{\mathcal I_1}\|_2<\pi$.
Let us replace $\phi|_{\mathcal I_1}$ by a Magnus minimal presentation
of $\Lexp(\phi|_{\mathcal I_1})$, in order to obtain an other measure $\phi_1$.
Then we call $\phi_1$  a semilocal contraction of $\phi$.

We call $\phi$  semilocally Magnus minimal, if finitely many application of semilocal contractions
does not decrease $\smallint \|\phi\|_2$.
(In this case, the semilocal contractions will not really be contractions, as they are reversible.)
We call $\phi$  locally Magnus minimal, if any application of a semilocal contraction
does not decrease $\smallint \|\phi\|_2$.
It is easy to see that
\[\text{(Magnus minimal) $\Rightarrow$(semilocally Magnus minimal) $\Rightarrow$(locally Magnus minimal)}.\]
The arrows do not hold in the other directions.
For example, $\tilde I\m 1_{[0,2\pi]}$ is semilocally minimal but not Magnus minimal.
Also, $(-\m 1_{[0,1]})\boldsymbol. \Psi_0\boldsymbol.\m 1_{[0,1]}$
is locally Magnus minimal but not semilocally Magnus minimal:
Using semilocal contractions we can move $(-\m 1_{[0,1]})$ and $\m 1_{[0,1]}$
beside each other, and then there is a proper cancellation.

The proper local generalization of Magnus minimality is semilocal Magnus minimality.
If $\phi$ is locally Magnus minimal, the we can define $\ellip(\phi)$ and $\hyper(\phi)$
by taking a finite division of $\{\mathcal I_j\}$ of $\mathcal I$ to intervals
of  variation less than $\pi$, and simply adding  $\ellip(\phi_j)$ and $\hyper(\phi_j)$.
What semilocality is needed for is to show that $\ellip(\phi_{\mathcal I})$
is nonnegatively proportional to $\ellip(\phi)$, and to a proper definition of the Magnus
direction of $\phi$.

Having that,  semilocally  Magnus minimal presentations up to
semilocal contractions behave like Magnus minimal presentations.
They can also be classified as Magnus elliptic, parabolic, hyperbolic, or loxodromic.
(But they are not elements of $\GL_2^+(\mathbb R)$ anymore but presentations.)
In fact, semilocally Magnus minimal presentations up to semilocal contractions
have a very geometrical interpretation, cf. Remark \ref{rem:cover}.
(Interpreted as elements of $\widetilde{\GL_2^+(\mathbb R)}$.)
\qedremark
\end{remark}

\snewpage
As Theorem \ref{th:normalform} suggests, hyperbolic and parabolic developments are rather rigid,
while in other cases there is some wiggling of elliptic parts.
\begin{theorem}\plabel{th:Magnusminimal}
Suppose that $A\neq\Id_2$, $p=\mathcal M_{2\times 2\,\,\real}(A)<\pi$,
 and $\phi$ is a minimal presentation to $A$ supported on the interval $[a,b]$.

(a) Suppose that $A$ is Magnus hyperbolic or parabolic.
Then there are unique elements $t\in[-\pi/2,\pi/2]$ and $\tilde F\in\tilde{\mathbb K}$ such that
\[\Lexp(\phi|_{[a,x]})=\NW\left(0,\int\|\phi|_{[a,x]}\|_2, t, \tilde F\right).\]
Thus, minimal presentations for Magnus hyperbolic and parabolic matrices are unique,
up to reparametrization of the measure.

(b) Suppose that $\CD(A)$ is point disk.
Then there is a unique element $t\in[0,2\pi)$ such that
\[\Lexp(\phi|_{[a,x]})=\exp\left((\Id_2\cos t+\tilde I\sin t)  \int\|\phi|_{[a,x]}\|_2\right).\]
Thus, minimal presentations for quasicomplex matrices are unique,
up to reparametrization of the measure.

(c) Suppose that $A$ is not of the cases above.
Then there are unique elements $t\in[0,2\pi)$, $p_1$, $p_2>0$, $\tilde F\in\tilde{\mathbb K}$
and surjective monotone increasing function $\varpi_i:[a,b]\rightarrow [0,p_i]$
such that
\[\varpi_1(x)+\varpi_2(x)=x-a\]
and
\[\Lexp(\phi|_{[a,x]})=\NW\left(\varpi_1\left(\int\|\phi|_{[a,x]}\|_2\right) ,
\varpi_2\left(\int\|\phi|_{[a,x]}\|_2\right), t, \tilde F \right).\]
Thus, minimal presentations in the general case
are unique, up to displacement of elliptic parts.
\begin{proof}
Divide $[a,b]$ to $[a,x]$ and $[x,b]$, and replace the minimal presentation by norma forms.
They must fit in accordance to minimality.
\end{proof}
\begin{proof}[Remark]
The statement can easily be generalized to semilocally Magnus minimal presentations.
\renewcommand{\qedsymbol}{}
\end{proof}
\end{theorem}
%\snewpage
Theorem \ref{th:Magnusminimal} says that certain minimal Magnus presentations are essentially unique.
Theorems \ref{th:Magnonexp} and \ref{th:Magcritic} will give some explanation to the fact
that it is not easy to give examples for the Magnus expansion blowing up in the critical case
$\smallint\|\phi\|_2=\pi$.
\begin{theorem}\plabel{th:Magnonexp}
Suppose that $A\neq\Id_2$, $p=\mathcal M_{2\times 2\,\,\real}(A)<\pi$,
 and $\phi$ is a minimal presentation to $A$ supported on the interval $[a,b]$.
If $\phi$ is of shape
\[\Lexp(\phi|_{[a,x]})=\exp\left(S \int\|\phi|_{[a,x]}\|_2\right)\]
with some matrix $S$ (i. e., it is essentially an exponential),
then $S$ is a normal matrix with norm $\|S\|_2=1$.
\begin{proof}
Due to homogeneity, $\ellip(\Phi|_{\mathcal I})$ and $\hyper(\Phi|_{\mathcal I})$ must be
proportional to $\mathcal M_{2\times 2\,\,\real}(\Phi|_{\mathcal I})$.
But it is easy to see that (up to parametrization)
only the homogeneous normal densities \eqref{eq:noruni} have this property,
and they are locally constant only if the Magnus non-elliptic component vanishes with $p_2=0$,
or when they are of special hyperbolic type with $\sin t=0$.
\end{proof}
\end{theorem}
(This redevelops Theorem \ref{th:NonNormalMin} in the real case.)

\begin{theorem}\plabel{th:Magcritic}
Suppose that $\phi$ is a measure,
\[\int\|\phi\|_2=\pi,\]
but $\log \Lexp(\phi)$ does not exist.
Then there are uniquely determined elements $t\in\{-\pi,\pi\}$ and $\tilde F\in\tilde{\mathbb K}$,
a nonnegative decomposition $\pi=p_1+p_2$, with $p_2>0$, and
surjective monotone increasing functions $\varpi_i:[a,b]\rightarrow [0,p_i]$
such that
\[\varpi_1(x)+\varpi_2(x)=x-a\]
and
\[\Lexp(\phi|_{[a,x]})=\NW\left(\varpi_1\left(\int\|\phi|_{[a,x]}\|_2\right) ,
\varpi_2\left(\int\|\phi|_{[a,x]}\|_2\right), t, \tilde F \right).\]

Thus, critical cases with $\log$ blowing up  are the
Magnus elliptic and parabolic (but not quasicomplex) developments up to reparametrization and rearrangement of elliptic parts.
\begin{proof}
The presentation must be Magnus minimal, otherwise the $\log$ would be OK.
Divide $[a,b]$ to $[a,x]$ and $[x,b]$, and replace the minimal presentation by normal parts.
They must fit in accordance to minimality.
It is easy to see that in the Magnus hyperbolic / loxodromic cases
$\CD(\Lexp(\phi|_{[a,x]}))$ has no chance to reach $(-\infty,0]$.
The disks are the largest in the Magnus hyperbolic cases, and the
chiral disks $\CD(W(\pi,\pi\sin t))$ of  Magnus strictly hyperbolic developments
do not reach the negative axis.
So the Magnus elliptic and parabolic cases remain but the quasicomplex one is ruled out.
\end{proof}
\end{theorem}
Thus, even critical cases with  $\int\|\phi\|_2=\pi$ are scarce.

\snewpage
\begin{remark}\plabel{rem:cover}
We started this section by investigating matrices $A$ with $\CD(A)\subset\intD(0,\pi)$.
It is a natural question to ask whether the treatment extends to  matrices $A$
with, say,  $\CD(A)\cap(\infty,0]=\emptyset$.
The answer is affirmative.
However, if we consider this question, then it is advisable to take an even bolder step:

We can extend the statements for $A\in\widetilde{\GL^+_2}(\mathbb R)$, the universal cover of ${\GL^+_2}(\mathbb R)$.
This of course, implies that we have to use the covering exponential
$\widetilde{\exp}:\mathrm M_2(\mathbb R)\rightarrow \widetilde{\GL^+_2}(\mathbb R),$
and $\Lexp$ should also be replaced by $\widetilde{\Lexp}$.
Now,  the chiral disks of elements of $\widetilde{\GL^+_2}(\mathbb R)$ live in $\widetilde{\mathbb C}$,
the universal cover of $\mathbb C\setminus\{0\}$.

Mutatis mutandis,  Theorems \ref{th:maximaldisk}, \ref{th:nonmaximaldisk},
\ref{th:normalform} extend in a straightforward manner.
Remarkably, Theorems \ref{th:cs} and \ref{th:CRrange} have versions in this case
(the group acting on the universal cover of $\mathbb C\setminus\{0\}$),
however we do not really need them that much, because chiral disks can be traced directly
to prove a variant of Theorem  \ref{th:MPexpress}.
Elements of $\widetilde{\GL^+_2}(\mathbb R)$ also have minimal Magnus presentations.
In our previous terminology, they are  semilocally Magnus minimal presentations.
In fact, semilocally Magnus minimal presentations up to semilocal contractions
will correspond to elements of $\widetilde{\GL^+_2}(\mathbb R)$.
Their classification to Magnus hyperbolic, elliptic, parabolic, loxodromic, quasicomplex  elements
extends to $\widetilde{\GL^+_2}(\mathbb R)$.
This picture of $\widetilde{\GL^+_2}(\mathbb R)$ helps to understand ${\GL^+_2}(\mathbb R)$.
Indeed, we see that every element of   ${\GL^+_2}(\mathbb R)$ have countably many
semilocally Magnus minimal presentations up to semilocal contractions,
and among those one or two (conjugates) are minimal.
The Magnus exponent of an element of   ${\GL^+_2}(\mathbb R)$ is the
minimal Magnus exponent of its lifts to $\widetilde{\GL^+_2}(\mathbb R)$.
\qedremark
\end{remark}
\begin{example}
\plabel{ex:RealComplexMin}
Let $\boldsymbol z=4.493\ldots$ be the solution of $\tan z=z$ on the interval $[\pi,2\pi]$.
Consider
\[\boldsymbol Z=\begin{bmatrix}-\sqrt{1+\boldsymbol z^2}-\boldsymbol z&\\ &-\sqrt{1+\boldsymbol z^2}+\boldsymbol z\end{bmatrix}.\]
The determinant of the matrix is $1$, we want to compute its real Magnus exponent.
The optimistic suggestion is
$\mathcal M_{2\times 2\,\,\complex}(\boldsymbol Z)=\sqrt{\pi^2+\log(\boldsymbol z+\sqrt{1+\boldsymbol z^2})^2}=3.839\ldots$.
Indeed, in the complex case, or in the doubled real case, this is realizable from
\[\boldsymbol Z=\exp\begin{bmatrix}\log(\boldsymbol z+\sqrt{1+\boldsymbol z^2})
+\pi\mathrm i&\\ &-\log(\boldsymbol z+\sqrt{1+\boldsymbol z^2})+\pi\mathrm i\end{bmatrix}.\]
However, in the real case, there is ``not enough space'' to do this.
A pessimistic suggestion is   $\pi+|\log(\boldsymbol z+\sqrt{1+\boldsymbol z^2})|=5.349\ldots$.
Indeed, we can change sign by an elliptic exponential, and then continue by a hyperbolic exponential.
This,  we know, cannot be optimal.
In reality, the answer is $\mathcal M_{2\times 2\,\,\real}(\boldsymbol Z)=\boldsymbol z=4.493\ldots$.
In fact, $\boldsymbol Z$ is Magnus parabolic, one can check that $\boldsymbol Z\sim W(\boldsymbol z,\boldsymbol z)$.
This is easy to see from the chiral disk.

In this case there are two Magnus minimal representations, because of the conjugational symmetry.
\qedexer
\end{example}

\scleardoublepage
\section{Optimal asymptotical norm estimate for $\mathrm M_2(\mathbb R)$} \plabel{sec:opt22}
If $\psi$ is an ordered $\mathrm M_2(\mathbb R)$ valued measure with cumulative norm
$p=\int\|\psi\|_2$ such that $0<p<\pi$, then, according to the previous sections,
the maximal possible norm of its Magnus expansion $\mu_{\mathrm L}(\psi)$ is realized via maximal disks,
through canonical Magnus parabolic or hyperbolic developments.
Thus this maximal norm is of shape $\|W(p,p\sin t)\|_2$ where $0\leq \sin t\leq 1$.
Thus we have optimize in $\sin t=s(p)$, although at this point it is not clear that $s(p)$ is unique depending $p$.

Using the implicit function theorem, however, it is easy to see that,
for $p\searrow0$, the optimal ridge is defined with $\sin t=s(p)$, where
\[s(p)=1-\frac16p^2+\frac{47}{360}p^4+O(p^6).\]
Then
\begin{equation}
\|\mu_{\mathrm L}(\Phi_{s(p)}|_{[0,p]})\|_2=p+{\frac {1}{6}}{p}^{3}-{\frac {1}{72}}{p}^{5}+{\frac {31}{3024}}{p}^
{7}+O \left( {p}^{9} \right).
\plabel{eq:micro33}
\end{equation}
We see that this is somewhere in the middle, below the upper estimate \eqref{eq:Jsan1} and above the parabolic case \eqref{eq:micro};
but all deviations are of $O(p^7)$.

Using an appropriate reparametrization and the implicit function theorem, one can see that,
for $p\nearrow\pi$, the optimal ridge is defined with $\sin t=s(p)$, where
\[s(p)=\underbrace{1-\frac1\pi\left({\pi-p}\right)^1}_{p/\pi}+\frac{2^{3/2}}{\pi^{3/2}}\left({\pi-p}\right)^{3/2}
-\frac43\left({\pi-p}\right)^2 +O\left( \left({\pi-p}\right)^{5/2}\right).\]
Then
\begin{multline}
\|\mu_{\mathrm L}(\Phi_{s(p)}|_{[0,p]})\|_2=
\sqrt2\pi^{3/2} (\pi-p)^{-1/2}
-2\pi
+\frac{\sqrt{2\pi}(4\pi^2-3) }{12} (\pi-p)^{1/2}
\\-\frac{4\pi^2-3}3 (\pi-p)^{1}
+\frac{\sqrt2(368\pi^2-120\pi^2-45)}{1440\sqrt\pi} (\pi-p)^{3/2}
+O((\pi-p)^2).
\plabel{eq:lower33}
\end{multline}
This is below the upper estimate \eqref{eq:Jsan2} by $O(1)$, and above the
parabolic case \eqref{eq:lower2} by $O((\pi-p)^{1/2})$,
and above the lower estimate \eqref{eq:lower3} by $O((\pi-p)^{3/2})$.

It is natural to ask whether the general upper estimate can be improved to, say
\[
\pi\sqrt{\frac{\pi+p}{\pi-p}}-2\pi+o(1),
\]
as $p\nearrow\pi$.
\begin{remark}
\plabel{rem:opt}
One can show that $s(p)$ yields a well-defined analytic function for $0<p<\pi$.
This is, however, complicated, as it requires global estimates.
\qedremark
\end{remark}
\snewpage
\section{A counterexample for $2\times2$ complex matrices }
\plabel{sec:count22}
According to Theorem \ref{th:MPexpress},   $A\in\mathrm M_2(\mathbb R)$,  $0\leq p<\pi$, and  $\CR(A)\subset\exp\Dbar(0,p) $ implies that
$A=\exp_{\mathrm L} \phi$ with some appropriate $\phi$ such that $\|\phi\|_2\leq p$.
Here we demonstrate that the corresponding statement is not valid for $2\times2$ complex matrices.

Let $0<p<\pi$. In the CKB model, we will consider the Magnus range
\[S_p=\frac{\mathrm{ CKB}}{\mathrm{PH}}\circ\exp((\mathbb C^+)\cap \Dbar(0,p)).\]
Then it is easy to see that its hyperbolic boundary
\[\partial S_p=\frac{\mathrm{ CKB}}{\mathrm{PH}}\circ\exp((\mathbb C^+)\cap  \partial\Dbar(0,p))\]
is parametrized by the curve
\[t\in(0,\pi)\mapsto s_p(t)=\left(\frac{\cos(p\sin t)}{\cosh(p\cos t)},\tanh (p\cos t)\right).\]
The tangent line at $t$ is given by the equation
\begin{align*}0=&
\underbrace{(\sin \left( t \right) )}_{A_p(t):=}x\\
&+\underbrace{(\sin \left( t \right) \sinh \left( p\cos \left( t \right)  \right)
\cos \left( p\sin \left( t \right)  \right) -\cos \left( t \right)
\cosh \left( p\cos \left( t \right)  \right) \sin \left( p\sin \left(
t \right)  \right)
)}_{B_p(t)):=}y\\
&+\underbrace{(-\sin \left( t \right) \cosh \left( p\cos \left( t \right)  \right)
\cos \left( p\sin \left( t \right)  \right) +\cos \left( t \right)
\sinh \left( p\cos \left( t \right)  \right) \sin \left( p\sin \left(
t \right)  \right)
)}_{C_p(t):=}.
\end{align*}
(In this section, $x,y$ are understood as $x_{\mathrm{CKB}},y_{\mathrm{CKB}}$.)
Note that $A_p(\pi-t)=A_p(t)$, $B_p(\pi-t)=-B_p(t)$, $C_p(\pi-t)=-C_p(t)$.

It is not hard to see that the equation of the ellipse tangent to $s_p$ at the three parameter points  $t<\frac\pi2<\pi-t$ is given by
\begin{multline*}
E_{p,t}(x,y)\equiv
(\cos \left( p\sin \left( t \right)  \right) -\cosh \left( p\cos \left( t \right)  \right) x)^2(A_p(t)\cos(p)+C_p(t))^2
\\-\left( \cos \left( p\sin \left( t \right)  \right) -\cosh \left( p\cos
 \left( t \right)  \right) \cos \left( p \right)
\right)^2\left((A_p(t)x+C_p(t))^2-(B_p(t)y)^2 \right)=0
\end{multline*}
(the coefficient of $y^2$ is nonnegative).

\snewpage

\begin{exaprop}
\plabel{ex:count22}
For the choice
\[\left(p_0,t_0\right)=\left(\frac{14}{15}\pi,\frac{7}{15}\pi\right),\]
the elliptical disk
\[E_{p_0,t_0}(x,y)\leq0 \]
is realized as the conformal range of a $2\times2$ complex matrix $X$ in CKB (i. e. as $\DW_{\mathrm{CKB}}(X)$).
Regarding this $X$ (4 different cases up to unitary conjugation),
\[\CR(X)\subset \exp\Dbar(0,p_0)\]
but
\begin{equation}
\mathcal M_{\text{$2\times2$ complex}}(X)>p_0.
\plabel{eq:cont}
\end{equation}
\begin{proof}
One can check that, for the given choice,
the ellipse $E_{p_0,t_0}(x,y)=0$ lies in the interior of the unit circle,
thus it be realized as the boundary of the conformal range of a (purely complex) $2\times2$ complex matrix $X$.
(Cf. Davis \cite{D1},\cite{D2}, Lins, Spitkovsky, Zhong \cite{LSZ}, or, in greater detail, \cite{L.ell}.)
It is  harder to see but this ellipse lies in $S_{p_0}$, intersecting (tangent to) $\partial S_{p_0}$
at $s_{p_0}(t_0), s_{p_0}(\pi/2), s_{p_0}(\pi-t_0)$. (This requires checking the second derivatives near the critical points,
and for the rest numerical considerations are sufficient.)

Assume that now $X=\Lexp\phi$ with $\|\phi\|_2=p_0$.
Then $\phi$ allows an initial restriction $\phi|_I$ such that $\|\phi_I\|_2=p_0/2$.
Let $X_{1/2}=\Lexp(\phi_I)$.
Then, however, $s_{p_0}(t_0),$ $ s_{p_0}(\pi/2),$ $s_{p_0}(\pi-t_0)\in\DW_{\mathrm{CKB}}(X)$ implies that
that  $s_{p_0/2}(t_0),$ $s_{p_0/2}(\pi/2),$ $s_{p_0/2}(\pi-t_0)\in\DW_{\mathrm{CKB}}(X_{1/2})$.
In this case $\DW_{\mathrm{CKB}}(X_{1/2})$ should be the ellipse with boundary  $E_{p_0/2,t_0}(x,y)=0$.
One can, on the other hand, check that this ellipse is not in the unit circle, which is a contradiction.

In fact, if we allow $X=\Lexp\phi$ with $\|\phi\|_2=p_0+\varepsilon$ (with small $\varepsilon$),
then it still means that the boundary of $\DW_{\mathrm{CKB}}(X_{1/2})$ should pass through
certain points near  $s_{p_0/2}(t_0),$ $ s_{p_0/2}(\pi/2),$ $ s_{p_0/2}(\pi-t_0)$
but  contained in $S_{p_0/2}$ (making the tangents also close at the given points).
This still makes the boundary of $\DW_{\mathrm{CKB}}(X_{1/2})$ close to  $E_{p_0/2,t_0}(x,y)=0$, yielding a contradiction.
(Thus \eqref{eq:cont} could be quantified further.)
\end{proof}
\end{exaprop}

\snewpage

\end{document}